\documentclass[11pt]{book}
\usepackage{amsfonts, amsmath,amssymb}
\usepackage[dvips]{graphicx}

\newcommand{\detail}[1]{\par\noi{\bf [Proof detail\ }{#1}
\hfill{\bf ]}\par\noi\hspace{-4pt}}
\renewcommand{\detail}[1]{}

\newcommand{\file}{$\ti{\ }$/tekst/habil/habil.tex\quad}
\renewcommand{\file}{}

\newcommand{\dis}{\displaystyle}
\newcommand{\txt}{\textstyle}

\newcommand{\vc}{\vspace{\baselineskip}}
\newcommand{\med}{\medskip}
\newcommand{\noi}{\noindent}

\newcommand{\halmos}{\rule{1ex}{1.4ex}}
\def \qed {\nopagebreak{\hspace*{\fill}$\halmos$\medskip}}

\newtheorem{theorem}{Theorem}[chapter]
\newtheorem{proposition}[theorem]{Proposition}
\newtheorem{corollary}[theorem]{Corollary}
\newtheorem{conjecture}[theorem]{Conjecture}
\newtheorem{lemma}[theorem]{Lemma}

\newtheorem{remark}[theorem]{Remark}
\newtheorem{defi}[theorem]{Definition}
\newtheorem{note}[theorem]{Remark}
\newcommand{\bt}{\begin{theorem}}
\newcommand{\et}{\end{theorem}}
\newcommand{\bl}{\begin{lemma}}
\newcommand{\el}{\end{lemma}}
\newcommand{\bp}{\begin{proposition}}
\newcommand{\ep}{\end{proposition}}
\newcommand{\bcor}{\begin{corollary}}
\newcommand{\ecor}{\end{corollary}}
\newcommand{\br}{\begin{remark}\rm}
\newcommand{\er}{\end{remark}}
\newcommand{\bcon}{\begin{conjecture}}
\newcommand{\econ}{\end{conjecture}}
\newcommand{\bdf}{\begin{defi}\rm}
\newcommand{\edf}{\hfill$\Diamond$\end{defi}}
\newcommand{\brm}{\begin{note}\rm}
\newcommand{\erm}{\hfill$\Diamond$\end{note}}

\newcommand{\be}{\begin{equation}}
\newcommand{\ee}{\end{equation}}
\newcommand{\ba}{\begin{array}}
\newcommand{\ea}{\end{array}}
\newcommand{\bc}{\be\begin{array}{r@{\,}c@{\,}l}}
\newcommand{\ec}{\end{array}\ee}

\newcommand{\al}{\alpha}
\newcommand{\bet}{\beta}
\newcommand{\ga}{\gamma}
\newcommand{\Ga}{\Gamma}
\newcommand{\de}{\delta}
\newcommand{\De}{\Delta}
\newcommand{\eps}{\varepsilon}
\newcommand{\la}{\lambda}
\newcommand{\La}{\Lambda}
\newcommand{\sig}{\sigma}

\newcommand{\tet}{\theta}
\newcommand{\oo}{\omega}
\newcommand{\om}{\Omega}

\newcommand{\si}{\ensuremath{\sigma}}

\newcommand{\Ai}{{\cal A}}
\newcommand{\Bi}{{\cal B}}
\newcommand{\Ci}{{\cal C}}
\newcommand{\Di}{{\cal D}}
\newcommand{\Ei}{{\cal E}}
\newcommand{\Fi}{{\cal F}}
\newcommand{\Gi}{{\cal G}}
\newcommand{\Hi}{{\cal H}}

\newcommand{\Li}{{\cal L}}
\newcommand{\Mi}{{\cal M}}
\newcommand{\Ni}{{\cal N}}

\newcommand{\Pc}{{\cal P}}
\newcommand{\Qi}{{\cal Q}}
\newcommand{\Ri}{{\cal R}}
\newcommand{\Si}{{\cal S}}
\newcommand{\Ti}{{\cal T}}
\newcommand{\Ui}{{\cal U}}
\newcommand{\Vi}{{\cal V}}
\newcommand{\Wi}{{\cal W}}
\newcommand{\Xc}{{\cal X}}
\newcommand{\Yi}{{\cal Y}}
\newcommand{\Zi}{{\cal Z}}

\newcommand{\R}{{\mathbb R}}
\newcommand{\N}{{\mathbb N}}
\newcommand{\Z}{{\mathbb Z}}

\newcommand{\Q}{{\mathbb Q}}
\newcommand{\T}{{\mathbb T}}

\newcommand{\li}{\langle}
\newcommand{\re}{\rangle}

\newcommand{\desd}{\ensuremath{\Leftrightarrow}}
\newcommand{\volgt}{\ensuremath{\Rightarrow}}
\newcommand{\up}{\uparrow}
\newcommand{\down}{\downarrow}
\newcommand{\sub}{\subset}
\newcommand{\beh}{\backslash}
\newcommand{\co}{{\rm c}}

\newcommand{\isd}{\stackrel{\scriptscriptstyle\Di}{=}}
\newcommand{\asto}[1]{\underset{{#1}\to\infty}{\longrightarrow}}
\newcommand{\Asto}[1]{\underset{{#1}\to\infty}{\Longrightarrow}}
\newcommand{\bplim}{\mbox{\rm bp-lim}}

\newcommand{\ti}{\tilde}
\newcommand{\dgg}{\dagger}
\newcommand{\ov}{\overline}
\newcommand{\un}{\underline}

\newcommand{\pa}{\partial}
\newcommand{\ffrac}[2]{{\textstyle\frac{{#1}}{{#2}}}}
\newcommand{\dif}[1]{\ffrac{\partial}{\partial{#1}}}
\newcommand{\diff}[1]{\ffrac{\partial^2}{{\partial{#1}}^2}}
\newcommand{\difif}[2]{\ffrac{\partial^2}{\partial{#1}\partial{#2}}}

\newcommand{\Diff}[1]{\frac{\partial^2}{{\partial{#1}}^2}}

\newcommand{\di}{\mathrm{d}}
\newcommand{\half}{{[0,\infty)}}
\newcommand{\expo}{\mbox{\large\it e}}
\newcommand{\ex}[1]{\expo^{\,\textstyle{#1}}}
\newcommand{\tr}{{\rm tr}}

\newcommand{\var}{{\rm Var}}
\newcommand{\cov}{{\rm Cov}}
\newcommand{\Pois}{{\rm Pois}}
\newcommand{\Thin}{{\rm Thin}}

\newcommand{\as}{\mbox{a.s.}}

\newcommand{\Widg}{\Wi_{\rm\scriptscriptstyle DG}}
\newcommand{\x}{\mathbf x}
\newcommand{\y}{\mathbf y}
\newcommand{\vb}{\mathbf v}
\newcommand{\hut}{\ov}

\newcommand{\Ys}{\Yi}
\newcommand{\Yp}{Y}
\newcommand{\gal}{\al}

\newcommand{\Tr}[1]{T_{{#1}}}
\newcommand{\expi}{\pi}
\newcommand{\poi}{\oo}
\newcommand{\Poi}{\om}
\newcommand{\lr}{\la}
\newcommand{\lc}{\ga}
\newcommand{\lir}{\la}
\newcommand{\asdto}[2]{\underset{{#1}\down{#2}}{\longrightarrow}}
\newcommand{\Asdto}[2]{\underset{{#1}\down{#2}}{\Longrightarrow}}

\begin{document}

\makeatletter\@addtoreset{equation}{section}
\makeatother\def\theequation{\thesection.\arabic{equation}} 

\renewcommand{\labelenumi}{{\rm (\roman{enumi})}}

\title{\vspace{-3cm}\Huge Extinction versus unbounded growth\\
{\Large Habilitation Thesis of the University Erlangen-N\"urnberg}}
\author{\Large Jan M. Swart\vspace{6pt}\\
{\' UTIA}\\
{Pod vod\'arenskou v\v e\v z\' i 4}\\
{18208 Praha 8}\\
{Czech Republic}\\
{e-mail: swart@utia.cas.cz}}
\date{{\small\file}\large January 31, 2007\\[20pt]
\parbox{400pt}{{\bf\large  Abstract} Certain Markov processes, or deterministic evolution equations, have the property that they are dual to a stochastic process that exhibits {\em extinction versus unbounded growth}, i.e., the total mass in such a process either becomes zero, or grows without bounds as time tends to infinity. If this is the case, then this phenomenon can often be used to determine the invariant measures, or fixed points, of the process originally under consideration, and to study convergence to equilibrium. This principle, which has been known since early work on multitype branching processes, is here demonstrated on three new examples with applications in the theory of interacting particle systems.}
}
\maketitle



{\setlength{\parskip}{-2pt}\tableofcontents}

\newpage

\chapter{Introduction}

\section{Interacting particle systems}

This habilitation thesis treats three subjects from {\em probability theory}, and more precisely, from the field of {\em interacting particle systems}. The binding element is a common technique used to study these subjects, which gives the title to this thesis, which finds its origin in multitype {\em branching theory}, and which is applied here both to branching processes and to processes which do not have the branching property, but still are in some ways similar to branching processes, although in other aspects of their behavior they are completely different. In this introductory section, we zoom out a bit more than is usual in a research paper, and take a look at the whole area of probability theory, and the fields of interacting particle systems and branching theory in particular, to see how they arose historically and how they are related.

Probability theory established itself as a mathematical discipline relatively late in history. Its origins are often traced back to an exchange of letters between Pascal and Fermat in the mid-17th century \cite{Apo69}, although some mention Cardano, one century earlier. The theory was not put on a firm axiomatic basis until the monograph by Kolmogorov in 1933 \cite{Kol33}, who based it on abstract measure theory, which had been developed in the preceding decades following the work of Lebesgue at the turn of the century. Because of these foundations, some authors claim that probability theory is a subfield of measure theory. Although there are measures all over the place, this is probably as justified as saying that algebra is a subfield of linear algebra.

When one tries to look for reasons why probability theory rose so late (why, for example, did the Greeks show no interest?), one is reminded of Einstein's remark `Gott w\"urfelt nicht' (God doesn't gamble). Even today, many people, including some mathematicians, associate mathematics primarely with beautiful structures that are entirely fixed, like a Penrose tiling, while an infinite random structure of the type that occurs in percolation theory evokes a certain disdain: `Why, that can be anything!'. Actually, it can't.

The reason is that once random structures get large, many events tend to get extremely improbable, until in the limit, for infinite systems, their probability is actually zero. The example that everybody knows are the laws of large numbers, which pertain to sums of independent identically distributed random variables. Closely related to this is the central limit theorem, which describes exactly how much randomness is left in the limit, and what the limit distribution is. Once a colleague asked what I was just working on. After hearing my explanation, his reaction was: so you are trying to prove a sort of central limit theorem? The answer is both yes and no.

Indeed, most of probability theory seems to be occupied with proving that certain things are certain in the limit that the system size, or time, or both tend to infinity, and that other things have a limit law.\footnote{I have to add a caveat here for statisticians, who are sometimes treated as probabilists, and sometimes as a species of their own, who from a practical point of view also have a lively interest in small samples, and, generally speaking, seem to be more interested in doing things and managing things, while the probabilist sensu strictu just sits down and tries to understand.} Yet, the methods needed to prove these limit statements are in general completely different from those used in the case of independent random variables. The independent case being well-understood, probabilists nowadays investigate systems of highly dependent components. And while there is just one way in which things can be independent, there are many ways in which things can depend on each other.

Seen from this point of view, the ``theory of interacting particle systems'' sounds like the natural culmination point of all of probability theory. That is not quite true. In fact, the classical book by Liggett called `Interacting Particle Systems' \cite{Lig85} was translated into Russian as `Markovskije Processy s Lokalnym Vzaimodejstvijem' (Markov Processes with Local Interaction), which captures the subject more precisely. Interacting particle systems are always situated in space, which is often $\Z^d$, sometimes $\R^d$, and sometimes another discrete or continuous structure that is in some way translation invariant. At each point in this space, there is some local Markov process going on, that is inherently random, and interacts with the Markov processes surrounding it. Although this interaction is only local, in the long run information can spread arbitrarily far, and therefore it is the long-time behavior of the process that is usually of interest.

This description of interacting particle systems excludes many other dependent systems, such as random walks in random environment, self-enforced and self-avoiding random walks, cellular automata and other deterministic evolutions, random matrices, and percolation theory, although many of these topics have close links with interacting particle systems. It also excludes, unrighteously, interacting particle systems in quantum probability. And, finally, it excludes other active areas of probabilistic research, such as abstract theory of Markov processes and semigroups, stochastic evolution equations, stochastic analysis, and more.

The origin of the field of interacting particle systems lies in 19-th century physics, when scientists like Bolzmann, Van der Waals, and others started to look for the molecular basis of thermodynamics. Thus, the original motivation was to study particles moving around in $\R^3$ according to the deterministic rules of classical Hamiltonian dynamics, or, later, its quantummechanical counterpart, which in a sense is both deterministic and inherently random. The mathematical problems arising from continuous space and deterministic motion being too difficult, people turned to models on lattices, that moreover have a local source of randomness. This class of models is still extremely rich, and apart from their original physical motivation, it was found that models of this type can be used to model many other interesting phenomena in a variety of applications in, for example, biology, sociology, and random network theory. Of the four classical models from \cite{Lig85}, namely the Ising model, voter model, contact process, and exclusion process, only the first and last have a clear physical motivation.

As a mathematical discipline, the field of interacting particle systems started around 1970. Again, compared to other branches of mathematics, this is very recent. This time, the reasons lie probably not only in a lack of interest (after all, the physical problems had been around for a century by that time) but also in the inherent difficulty of the subject. Certain special results date back further, to the mid 40ies; this includes work on multitype branching processes, percolation, and the famous Onsager solution of the 2-dimensional equilibrium Ising model. Gradually, people had to get used to the fact that interacting particle systems rarely allow for explicit solutions, and that very little can be said about them in general. Rather, even the simplest-looking among them required the development of new tools suited exactly for them, and many naive questions remained open for many years.

The systems of interest (interacting particle systems) and the main questions (limit laws for large system sizes and large times) being defined now, we can focus on some more specific topics. The first topic we would like to mention, which motivates much of the work done in the field, is that of {\em phase transitions}. Originally referring to the phenomenon that certain substances (as a general rule with exceptions: pure chemical substances) can either be in a gaseous, fluid, or solid phase, and change abruptly between these phases as the temperature or pressure pass a certain point, the concept has subsequently been generalized to include more phases (e.g.\ graphite versus diamond) and then to describe the general phenomenon that many-particle systems may drastically change their behavior when certain parameters pass certain tresholds, called {\em critical points}.

Phase transitions are a central topic for a number of reasons. First of all, since finite systems running for a finite time generally depend continuously on their parameters, mathematically ideal phase transitions occur only in the limit that the system size, and time, are sent to infinity, and therefore are the typical sort of phenomenon that justifies the study of large or infinite systems. Second, detailed information about them is often hard to get, since they are out of reach of most expansion techniques that tell us something about very high or low values of our parameters. In other words, phase transitions are difficult, and therefore prestigious. The third and most important reason is probably the belief, supported by nonrigorous theory developed by theoretical physicists, that phase transitions are highly universal. Thus, different interacting particle systems may have the `same' phase transition. Although the exact parameter values where this phase transition takes place may differ from one model to the other, zooming in on these phase transitions, and at the same time zooming out in space (and time, if we are not in equilibrium) should always yield roughly the same picture. This can for example be seen from the {\em critical exponents} of these phase transitions, which describe how certain quantities behave according to a certain power law as the critical point is approached. The classical paper in physics on this topic is \cite{WK74}.

Trying to prove results about {\em critical phenonema} that take place at, or in the immediate vicinity of the critical points, in particular, the calculation of critical exponents, has been a big aim behind much work done on interacting particle systems. Progress has been slow. In a number of cases, expansion techniques, such as the lace expansion, have been used to show that certain systems have `trivial' exponents, that are the same as those for other, noninteracting systems. Recently, important progress has been made on critical exponents for two-dimensional systems having conformally invariant scaling limits. The key object in this work is the Stochastic Loewner Equation \cite{Law05}. Apart from these two cases (the `trivial' critical exponents and those from conformal field theory) there is still little process.

Where, in all of this, is the present habilitation thesis situated? No critical exponents will be calculated in what follows, but we will see critical phenomena, and even some universality. In any case, there will be phase transitions around, and we will prove limit laws as time and system size are sent to infinity. A repeating theme in the proofs will be the exploitation of the simple observation that in certain particle systems, the number of particles either becomes zero, or tends to infinity. As far as I am aware off, this idea was first used in multitype branching theory.

The theory of branching processes started with a paper by Galton and Watson in 1874 \cite{WG74}, who studied the problem of the extinction of noble names. The problem drew new interest with the rise of probability theory in the 30-ies and with the study of nuclear chain reactions, which led to the study of multitype processes. It was only in the mid-70-ies, when people started to consider $\Z^d$ as the space of types, that the first branching processes were studied that might truly be called interacting particle systems. Even as such, they hardly deserve the name, since they consist of particles independently hopping around on a lattice, that moreover independently of each other split into more particles or die. The only way in which dependencies arise, which make the model interesting, is through the fact that certain `families' of particles all descend from one and the same `ancestor'. Basic questions about their ergodic behavior were solved by Kallenberg \cite{Kal77} using his famous `backward tree technique'. We will use this technique in Section~\ref{002sec}. It is moreover closely linked to the work in Chapter~\ref{C:typ} of this thesis. The main technique that unites all chapters, however, is the use of `extinction versus unbounded growth', as will be explained in the next section.

\section{Extinction versus unbounded growth}

Certain Markov processes, or deterministic evolution equations, have the property that they are dual to a stochastic process that exhibits {\em extinction versus unbounded growth}, i.e., the total mass in such a process either becomes zero, or grows without bounds as time tends to infinity. If this is the case, then this phenomenon can often be used to determine the invariant measures, or fixed points, of the process originally under consideration, and to study convergence to equilibrium. In this section, we demonstrate this principle, in the historicaly correct order, first on multitype branching processes, and then on the contact process.

\subsection{Extinction versus unbounded growth in branching theory}\label{S:eubr}

Consider a collection of particles of $n$ different types. Assume that each particle of type $i\in\{1,\ldots,n\}$ gives with {\em birth rate} $b_{ij}$ birth to a particle of type $j\in\{1,\ldots,n\}$, and dies with {\em death rate} $d_i$. We will assume that $b_{ij}>0$ and $d_i>0$ for all $i,j$. Let $Y_t(i)$ denote the number of particles of type $i$ at time $t\geq 0$. Then $Y=(Y_t)_{t\geq 0}$ is a Markov process in $\N^n$, which in the usual terminology is called a {\em continuous-time multitype binary branching process}. We write $P^y$ for  the law of $Y$ started in $Y_0=y$ and denote expectation with respect to $P^y$ by $E^y$. It is well-known that
\be\label{Ydual}
E^y\Big[\prod_{i=1}^n(1-u_0(i))^{Y_t(i)}\Big]=\prod_{i=1}^n(1-u_t(i))^{y(i)}\qquad(t\geq 0),
\ee
whenever $u_t=(u_t(1),\ldots,u_t(n))$ is a $[0,1]^n$-valued solution to the system of differential equations
\be\label{difu}
\dif{t}u_t(i)=\sum_{j=1}^nb_{ij}u_t(j)(1-u_t(i))-d_iu_t(i)\qquad(t\geq 0,\ i\in\{1,\ldots,n\}).
\ee
The map that gives $(1-u_t)$ as a function of $(1-u_0)$ and $t$ is what is classically known as the {\em generating function} of the branching process $Y$ (at time $t$). We prefer to work with $u_t$ (and not $1-u_t$) since this will simplify formulas later on.

Formula (\ref{Ydual}) has a useful interpretation in terms of thinning. By definition, a {\em thinning} of a particle configuration $y\in\N^n$ with a vector $v\in[0,1]^n$ is the random particle configuration obtained from $y$ in the following manner. Independently for each particle, we decide with probability $v(i)$ (depending on the type $i$ of the particle) whether we will keep it; with the remaining probability $1-v(i)$ we throw this particle away. If we denote the thinned collection of particles resulting from this procedure by $\Thin_v(y)$, then the left-hand side of (\ref{Ydual}) is just the probability that the configuration $\Thin_{u_t}(Y_t)$ contains no particles. Since the right-hand side of (\ref{Ydual}) has a similar interpretation, we may rewrite (\ref{Ydual}) as
\be\label{Ythin}
P^y[\Thin_{u_0}(Y_t)=0]=P[\Thin_{u_t}(y)=0]\qquad(t\geq 0).
\ee
The relation (\ref{Ydual}), or its rewrite (\ref{Ythin}), are an example of a {\em duality relation}, where the dual of the Markov process $Y$ is in this case the deterministic process $u$.

Using this duality relation, we can deduce information about $Y$ from $u$, and vice versa. To demonstrate this, we will show how the fact that the process $Y$ exhibits extinction versus unbounded growth gives information about the fixed points of the $n$-dimensional differential equation (\ref{difu}).

It is not hard to see that
\be
\dif{t}E[Y_t(i)]=\sum_{j=1}^nM_{ji}E[Y_t(j)]\qquad(t\geq 0),
\ee
where $M_{ji}=b_{ji}-\de_{ij}d_i$ $(i,j=1,\ldots,n)$. Since by adding a constant multiple of the identity, we can make $M$ into a matrix with strictly positive entries, it follows from the Perron-Frobenius theorem that $M$ has a maximal eigenvalue, say $\la$, that corresponds to a positive right and left eigenvector, which are the only nonnegative eigenvectors. If $\la<0$, we say that the branching process $Y$ is {\em subcritical}, if $\la=0$ we say that it is {\em critical}, and if $\la>0$ we say that it is {\em supercritical}. In the subcritical and critical cases, $Y$ {\em dies out}, i.e.,
\be
P^y\big[\exists t\geq 0\mbox{ s.t.\ }Y_s=0\ \forall s\geq t\big]=1\qquad(y\in\N^n).
\ee
(Note that since there is no spontaneous creation of particles, the zero configuration is a trap for the Markov process $Y$.) On the other hand, in the supercritical case, on which we focus from now on, $Y$ {\em survives} with positive probability, i.e.,
\be\label{survprob}
P^y\big[Y_t\neq 0\ \forall t\geq 0]>0\qquad(y\in\N^n,\ y\neq 0).
\ee
Indeed, the probability in (\ref{survprob}) is given by $1-\prod_{i=1}^n(1-p(i))^{y(i)}$, where
\be\label{pdefi}
p(i):=P^{\de_i}\big[Y_t\neq 0\ \forall t\geq 0]>0\qquad(i=1,\ldots,n),
\ee
and $\de_i$ denotes the particle configuration with just one particle of type $i$.

We claim that $p$ is the only nonzero fixed point of the differential equation (\ref{difu}), and the limit point started from any nonzero initial condition. To prove this, we observe that $Y$ exhibits {\em extinction versus unbounded growth}, in the following sense:
\be\label{Yexgro}
P^y\big[\exists t\geq 0\mbox{ s.t.\ }Y_s=0\ \forall s\geq t\quad\mbox{or}\quad\lim_{t\to\infty}|Y_t|=\infty\big]=1\qquad(y\in\N^n),
\ee
where $|y|:=\sum_{i=1}^ny(i)$ denotes the total number of particles in a particle configuration $y\in\N^n$. Why does (\ref{Yexgro}) hold? We will not give a formal proof here, but just explain the main idea. (For a more formal approach, see Lemma~\ref{L:sufexgro} below.) Since we are assuming that the death rates $d_i$ are all positive, it is not hard to show that
\be
\inf_{|y|\leq K}P^y\big[\exists t\geq 0\mbox{ s.t.\ }Y_s=0\ \forall s\geq t\big]>0\qquad(K\geq 0).
\ee
Indeed, if the process $Y$ is started with no more than $K$ particles, then there is a positive chance that all these particles die before they have a chance to branch, and therefore the probability that the process dies out can be estimated from below uniformly in all particle configurations with no more than $K$ particles. Now imagine that the number of particles $|Y_t|$ is less than $K$ at a (random) sequence of times tending to infinity. Then the process would infinitely often have a (uniformly) positive chance to die out in the next time interval of a certain length, and therefore it would eventually have to die out. Since this is true for any $K$, the only way for the process to escape extinction is to let the number of particles tend to infinity. 

We now show how extinction versus unbounded growth (formula (\ref{Yexgro})) implies that any solution of (\ref{difu}) with $u_0\neq 0$ satisfies
\be\label{utop}
\lim_{t\to\infty}u_t=p,
\ee
where $p$ is defined in (\ref{pdefi}). Note that $P[\Thin_v(\de_i)\neq 0]=v(i)$ $(v\in[0,1]^n)$, and therefore, by (\ref{Ythin}),
\be\label{utex}
u_t(i)=P^{\de_i}[\Thin_{u_0}(Y_t)\neq 0]\qquad(t\geq 0,\ i=1,\ldots,n).
\ee
Since we are assuming that $b_{ij}>0$ for all $i,j$, it is easy to see from (\ref{utex}) that $u_0\neq 0$ implies $u_t(i)>0$ for all $i=1,\ldots,n$ and $t>0$, so by a restart argument we may without loss of generality assume that $u_0(i)>0$ for all $i=1,\ldots,n$. 

Using (\ref{utex}) once more, and using extinction versus unbounded growth (formula (\ref{Yexgro})), we see that for large $t$ there are up to an event with small probability only two situations to be considered. Either $Y_t=0$, in which case $\Thin_{u_0}(Y_t)=0$, or $|Y_t|$ is large, in which case, by the fact that $u_0(i)>0$ for all $i$, we know that $\Thin_{u_0}(Y_t)$ is with large probability nonzero. Therefore, $P^{\de_i}[\Thin_{u_0}(Y_t)\neq 0]\cong P[Y_t\neq 0]$ for large $t$, and taking the limit $t\to\infty$ in (\ref{utex}) we arrive at (\ref{utop}). This proves that $p$ is the only nonzero fixed point of the differential equation (\ref{difu}), and the limit point started from any nonzero initial condition.

In a discrete time setting (but with much more general branching mechanisms), the result (\ref{utop}), including a proof based on extinction versus unbounded growth, can be found in Harris \cite[Theorem~II.7.2]{Har63}, who ascribes it to Everett and Ulam \cite{EU48}.

It is not hard to see that the positivy assumptions on the rates $b_{ij}$ and $d_i$ can be weakened considerably. In fact, it suffices if at least one of the $d_i$ is nonzero, and if the $b_{ij}$ are {\em irreducible}, in the sense that for each $i,j\in\{1,\ldots,n\}$, there exist $k_0,\ldots,k_m$ with $k_0=i$, $k_m=j$, and $b_{k_{l-1},k_l}>0$ for all $l=1,\ldots,m$.

\subsection{Extinction versus unbounded growth in the contact process}\label{S:euco}

The standard, nearest neighbor $d$-dimensional {\em contact process} is a Markov process $\eta=(\eta_t)_{t\geq 0}$ taking values in the space of all subsets of $\Z^d$, with the following description. If $i\in\eta_t$, then we say that the {\em site} $i\in\Z^d$ is {\em infected} at time $t\geq 0$, otherwise such a site is called {\em healthy}. Infected sites become healthy with rate $1$. Healthy sites become infected with infection rate $\lambda$ times the number of neighboring infected sites. Here, we say that $i,j\in\Z^d$ are neighbors if $|i-j|=1$.

It is useful to think about the contact process as a frustated branching process. Think of infected sites as being occupied by a particle. Then each particle tries with rate $\la$ to give birth to a particle at each neighboring site. If, however, that site is already occupied by a particle, the birth fails.

Indeed, it is easy to see that $|\eta_t|$, the total number of infected sites, can be bounded from above by a binary branching process with branching rate $2d\la$ and death rate $1$. In particular, if $\la\leq 1/(2d)$, this branching process is (sub)critical, and hence the contact process {\em dies out}. On the other hand, with considerably more effort, it is possible to show that for suffiently large $\la$, the contact process {\em survives} with positive probability, i.e.,
\be
P^A[\eta_t\neq\emptyset\ \forall t\geq 0]>0\qquad(A\neq\emptyset).
\ee
It is easy to show that two contact processes $\eta,\ti\eta$ with infection rates $\la,\ti\la$ can be coupled such that $\eta_t\leq\ti\eta_t$, so it follows that there exists a {\em critical infection rate} $0<\la_{\rm c}<\infty$ such that the contact process dies out for $\la<\la_{\rm c}$ and survives (with positive probability) for $\la>\la_{\rm c}$. The question whether the contact process survives {\em at} $\la=\la_{\rm c}$ was open for almost 15 years; its solution by Bezuidenhout and Grimmett in \cite{BG90} was a major milestone in the development of the theory of the contact process.

We will not touch this subject here, but rather show how the fact that the contact process exhibits extinction versus unbounded growth, together with self-duality, can be used to prove that if the contact process survives, then it has a unique nontrivial homogeneous invariant law. Here, we say that a probability law on the space of all subsets of $\Z^d$ is {\em nontrivial} if it gives zero probability to the empty set, and (spatially) {\em homogeneous} if it is invariant under translations.

It is well-known that the contact process is {\em self-dual}, in the following sense. Fix an infection rate $\la$, and for $A\sub\Z^d$, let $\eta^A$ denote the contact process with this infection rate started in the initial state $\eta^A_0=A$. Then
\be\label{contdual}
P[\eta^A_t\cap B=\emptyset]=P[A\cap\eta^B_t=\emptyset]\qquad(t\geq 0,\ A,B\sub\Z^d).
\ee
Since the contact process is an attractive spin system, it follows from standard theory that it has an {\em upper invariant law} $\ov\nu$, which is the largest invariant law in the sense of stochastic ordering, and the limit law as $t\to\infty$ of the process started with all sites infected:
\be
\Li(\eta^{\Z^d}_t)\Asto{t}\ov\nu.
\ee
Using the self-duality (\ref{contdual}) we can give a useful characterization of $\ov\nu$. Let $\eta^{\Z^d}_\infty$ be a random variable with law $\Li(\eta^{\Z^d}_\infty)=\ov\nu$. Then
\be\label{durep}
P[\eta^{\Z^d}_\infty\cap A=\emptyset]=\lim_{t\to\infty}P[\Z^d\cap\eta^A_t=\emptyset]=P[\exists t\geq 0\mbox{ s.t.\ }\eta^A_t=\emptyset]
\ee
for all finite $A\sub\Z^d$. Since $\Li(\eta^{\Z^d}_t)$ is homogeneous for each $t\geq 0$, so is $\ov\nu$. Using (\ref{durep}) and survival, it is not hard to show that $\ov\nu$ is nontrivial. We claim that it is the only invariant law with this property and moreover, that
\be\label{lawcon}
\Li(\eta_t)\Asto{t}\ov\nu
\ee
when $\eta$ is a contact process started in any initial law $\Li(\eta_0)=\mu$ that nontrivial and homogeneous. To prove this, we observe that the contact process exhibits extinction versus unbounded growth in the following sense:
\be\label{conexgro}
P\big[\exists t\geq 0\mbox{ s.t.\ }\eta^A_t=\emptyset\quad\mbox{or}\quad\lim_{t\to\infty}|\eta^A_t|=\infty\big]=1\qquad(A\sb\Z^d),
\ee
where $|A|$ denotes the cardinality of a set $A$. The proof is basically the same as in the case of multitype branching (see formula (\ref{Yexgro})). Since it may happen that all infected sites become healthy before any further infection has taken place, it is easy to show that
\be
\inf_{|A|\leq K}P\big[\exists t\geq 0\mbox{ s.t.\ }\eta^A_t=\emptyset\big]>0\qquad(K\geq 0).
\ee
Thus, the probability that the process will die out can be estimated from below uniformly in all configurations with at most $K$ infected sites, and therefore the only way for the process to avoid extinction is to let the number of infected sites tend to infinity.

Now let $\Li(\eta_0)=\mu$ be nontrivial and homogeneous. Then, with a bit of trouble, it is possible to show that for each $t>0$, the law $\Li(\eta_t)$ has the property that
\be\label{partev}
\lim_{K\to\infty}\sup_{|A|\leq K}P[\eta_t\cap A_n=\emptyset]=0.
\ee
Therefore, by a restart argument, we may without loss of generality assume that $\Li(\eta_0)$ has this property. Self-duality (formula (\ref{contdual})) tells us that
\be\label{ett}
P[\eta_t\cap A=\emptyset]=P[\eta_0\cap\eta^A_t=\emptyset]\qquad(t\geq 0),
\ee
where $\eta_0$ and $\eta^A_t$ are independent. If $t$ is large, then in evaluating the right-hand side of (\ref{ett}), by extinction versus unbounded growth (\ref{conexgro}), up to an event with small probability we need to consider only two cases. Either $\eta^A_t=\emptyset$, in which case $\eta_0\cap\eta^A_t=\emptyset$, or $|\eta^A_t|$ is large, in which case $\eta_0\cap\eta^A_t$ is with high probability not empty since $\Li(\eta_0)$ has the property (\ref{partev}). It follows that $P[\eta_0\cap\eta^A_t=\emptyset]\cong P[\eta^A_t=\emptyset]$ for large $t$, and taking the limit $t\to\infty$ in (\ref{ett}), using (\ref{durep}), we see that
\be
\lim_{t\to\infty}P[\eta_t\cap A=\emptyset]=P[\eta^{\Z^d}_\infty\cap A=\emptyset],
\ee
for all finite $A\sub\Z^d$, which proves (\ref{lawcon}).

This argument is due to Harris \cite[Theorem~9.2]{Har76}, who builds on earlier work of Vasil'ev, Vasershtein, Leontovich, and others. It can also be found in Ligget's book \cite[Theorem~VI.4.8]{Lig85}.

\section{Overview of the habilitation thesis}

\subsection{Branching processes in renormalization theory}

Certain problems in the study of a special type of interacting particle system, namely {\em linearly interacting catalytic Wright-Fisher diffusions}, lead one to study a special continuous-mass continuous- type space branching process, namely, the {\em super-Wright-Fisher diffusion}. This is a Markov process $\Ys=(\Ys_t)_{t\geq 0}$, taking values in the space of finite measures on $[0,1]$, whose transition probabilities are uniquely characterized by its Laplace functionals
\be\label{Lapl}
E^\mu\big[\ex{-\li\Ys_t,u_0\re}\big]=\ex{-\li\mu,u_t\re}\qquad(t\geq 0),
\ee
where $\li\mu,f\re:=\int\! f\,\di\mu$ and $u$ is a mild solution of the semilinear Cauchy equation
\be\label{Cau}
\dif{t}u_t(x)=\ffrac{1}{2}x(1-x)\diff{x}u_t(x)+\gal u_t(x)(1-u_t(x))\qquad(t\geq 0),
\ee
with $u_0$ any nonnegative continuous function on $[0,1]$. One should think of (\ref{Lapl}) and (\ref{Cau}) as continuous analogues of (\ref{Ydual}) and (\ref{difu}), respectively, where the finite type space $\{1,\ldots,n\}$ has been replaced by $[0,1]$ and the space $\N^n$ of all $n$-type particle configurations has been replaced by the space $\Mi[0,1]$ of all finite measures on $[0,1]$. We can think of $\Ys_t$ as describing a population, consisting of many particles each of which has a very small mass, such that each particle performs a Wright-Fisher diffusion on $[0,1]$, that is, the Markov process in $[0,1]$ whose generator is (the closure of) the operator $\ffrac{1}{2}x(1-x)\diff{x}$, and in addition, particles branch in such a way that the offspring of a bit of mass $\di m$ at position $x$ during a time interval of length $\di t$ produces offspring with mean $(1+\gal\di t)\di m$ and variance $\gal\di t$.

The way how the super Wright-Fisher diffusion $\Ys$ arises in a renormalization analysis of systems of linearly interacting catalytic Wright-Fisher diffusions will be explained in Chapter~\ref{C:ren} For the moment, we take the process in (\ref{Lapl}) for granted, and ask about fixed point(s) and long-time convergence of solutions $u$ to the Cauchy equation (\ref{Cau}). We would like to play the same game as in Section~\ref{S:eubr} and use extinction versus unbounded growth of $\Ys$ to prove convergence of $u$. Apart from the technical complications arising from continuous type space and continuous mass, we meet a more fundamental problem: our underlying motion, the Wright-Fisher diffusion, is not irreducible, i.e., it is not possible to get with positive probability from any point to any other point in the type space.

Indeed, the Wright-Fisher diffusion $Y$ has two traps: $0$ and $1$, and the process started in any initial state satisfies
\be
P\big[\exists\tau<\infty,\ r\in\{0,1\}\mbox{ s.t.\ }Y_t=r\ \forall t\geq\tau\big]=1,
\ee
i.e., the process gets trapped in finite time. For the measure-valued process $\Ys$, this means that with positive probability, in the long run most of the mass gets concentrated in $0$, or $1$, or both. Whether there is also a positive probability that there remains some mass in $(0,1)$ turns out to depend on the parameter $\gal$. For $\gal>1$, the answer is yes; otherwise it is no. As a result, we have to prove extinction versus unbounded growth on each of the part of the type space $\{0\},\{1\}$, and $(0,1)$, and we find three or four (depending on $\gal$) different nonzero fixed points of (\ref{Cau}), each with their own domain of attraction. 

This analysis carried out in Sections~\ref{S:WFi}--\ref{S:WFiii} of Chapter~\ref{C:ren}. There, a similar analysis is carried out also for a related branching process in discrete time, the description of which is somewhat complicated. An important tool in this analysis is the use of  embedded particle systems, as explained in Section~\ref{exex}. The results in this chapter are joint work with Klaus Fleischmann (WIAS, Berlin). Part of this has been published in \cite{FSsup}.

\subsection{Branching-coalescing particle systems}

Consider a model of binary branching random walks, i.e., a collection
of particles situated on a lattice $\La$, where each particle moves
independently of the others according to a continuous time random walk
that jumps from site $i\in\La$ to site $j$ with rate $a(i,j)$, each
particle splits with a branching rate $b\geq 0$ into two new
particles, created on the position of the old one, and each particle
dies with a death rate $d\geq 0$. Let $X_t(i)$ denotes the number of
particles at time $t\geq 0$ at the site $i\in\La$ and write
$X_t:=(X_t(i))_{i\in\La}$. Then, in analogy with (\ref{Ydual}), one
has
\be\label{gendu}
E^x\Big[\prod_{i=1}^n(1-u_0(i))^{X_t(i)}\Big]
=\prod_{i=1}^n(1-u_t(i))^{x(i)}\qquad(t\geq 0),
\ee
whenever $u_t=(u_t(1),\ldots,u_t(n))$ is a $[0,1]^\La$-valued solution
to the system of differential equations
\be\label{difu2}
\dif{t}u_t(i)=\sum_ja(j,i)(u_t(j)-u_t(i))
+bu_t(i)(1-u_t(i))-du_t(i)
\ee
$(t\geq 0,\ i\in\La)$. For each $f\in[0,1]^\La$, set $U_tf:=u_t$
$(t\geq 0)$ where $u$ solves (\ref{difu}) with initial condition
$u_0=f$; then $(U_t)_{t\geq 0}$ is the {\em generating semigroup} of
the branching process $X=(X_t)_{t\geq 0}$.

What happens if in the branching system $X$ we also allow for {\em
coalescence} of particles, i.e., if we let each pair of particles,
present on the same site, coalesce with rate $2c$ (with $c\geq 0$) to
one particle? In this case, we lose the branching property, i.e., we obtain a truly interacting system of particles. It turns out that although there is now no longer a generating semigroup in the classical sense, if we replace the deterministic evolution in (\ref{difu}) by the system of stochastic differential equations (SDE's)
\bc\label{usde}
\di u_t(i)&=&\dis\sum_ja(j,i)(u_t(j)-u_t(i))\,\di t
+bu_t(i)(1-u_t(i))\,\di t-du_t(i)\,\di t\\
&&\dis+\sqrt{2cu_t(i)(1-u_t(i))}\,\di B_t(i)
\qquad\qquad(t\geq 0,\ i\in\La),
\ec
then formula (\ref{gendu}) generalizes to the case with coalescence in the sense that
\be\label{codu}
E\Big[\prod_{i=1}^n(1-u_0(i))^{X_t(i)}\Big]
=E\Big[\prod_{i=1}^n(1-u_t(i))^{X_0(i)}\Big]\qquad(t\geq 0).
\ee
The duality (\ref{codu}) is due to \cite{Shi81,SU86}. It turns out
that the behavior of branching-coalescing particle systems of the type
we have just described is very similar to that of the contact
process. In fact, the history of this type of models seems to be as
least as old as that of the contact process. In particular, our model
is a special case of Schl\"ogl's first model \cite{Sch72}.

Given the similarity of $X$ with a contact process, and the similarity of the duality (\ref{codu}) with the self-duality of the contact process (\ref{contdual}), one can try to mimick the proof of (\ref{lawcon}) in the present set-up. This was done by Shiga and Uchiyama in \cite{SU86} for solutions $u$ to the system of SDE's (\ref{usde}). More precisely, they used extinction versus unbounded growth for the particle system $X$ to prove that the law of the system of SDE's $u$, started in any nontrivial homogeneous initial law, converges for $t\to\infty$ to the upper invariant law of $u$. 

We note that if the death rate $d$ is positive, then the probability that the process $X$ will get extinct can be estimated from below uniformly in all configurations with at most $K$ particles. Therefore, extinction versus unbounded growth for $X$ follows by the same argument as in Sections~\ref{S:eubr} and \ref{S:euco}. If $d=0$, the process cannot get extinct. In this case, it is not completely trivial to show that the number of particles tends to infinity, which forced the authors of \cite{SU86} to make some additional technical assumptions.

In Chapter~\ref{C:braco}, we turn the duality (\ref{codu}) around, and use extinction versus unbounded growth for the system of SDE's $u$ to prove that the law of the particle system $X$ started in any nontrivial homogeneous initial law, converges for $t\to\infty$ to the upper invariant law of $X$. This also involves some technical difficulties, since we need to show that the continuous system $u$ may hit zero in finite time, and we need to show that $X$ has an upper invariant law, which means that we must show that $X$ can be started with infinitely many particles at every site. 

These problems can be overcome, however, and we end up with results that are stronger than those in \cite{SU86}. Additional tools that we use are a self-duality for the system of SDE's $u$, as well as the fact that the particle system $X$ can be obtained from $u$ by Poissonization. This is joint work with Siva Athreya (Bangalore), and has been published in \cite{AS05}.

\subsection{The contact process seen from a typical site}

In the last chapter of this thesis, we return to the classical contact process, but instead of studying the process started in a nontrivial homogeneous initial law as in Section~\ref{S:euco}, we wish to study the process started in finite initial states. It is known that questions about this sort of initial states are much more difficult than those about homogeneous initial laws. Nevertheless, a lot is known for the standard, nearest neighbor process on $\Z^d$. A central technical tool in this work is a dynamical block technique due to \cite{BG90}, which shows that the contact process, whenever it survives, can be compared with oriented percolation with an arbitrary high parameter. This technique finds its origin in older (although published later) work on unoriented percolation \cite{GM90,BGN91}.

While this technique has been very successful for the symmetric nearest-neighbor contact process on $\Z^d$, and can no doubt be extended to short-range contact processes on the same lattice, it is not obvious if it can be adapted to asymmetric processes, or to other lattices than $\Z^d$. Nevertheless, the study of contact processes on other lattices than $\Z^d$ is interesting both from a theoretical and practical poiint of view. The theoretical motivation comes from analogies with unoriented percolation on general transitive graphs, which has proved to be a fruitful topic (see, e.g., \cite{BLPS99}). For unoriented percolation, it is known that it is important whether the underlying lattice is amenable (such as $\Z^d$) or not (e.g.\ a regular tree). Work on the contact process on regular trees by \cite{Pem92,DS95,Lig96,Sta96} makes one suspect that a similar dichotomy could hold for the contact process.

In Chapter~\ref{C:typ}, we study contact processes on general countable groups $\La$. We use a technique from the theory of branching processes, namely Palm measures, to show that indeed, certain aspects of the behavior of the contact process started in finite initial states depend on a property of the underlying lattice. The property that turns out to be important is whether $\La$ has subexponential growth, which is in fact a bit stronger than amenability.

Somewhat surprisingly, it turns out that in this context, extinction versus unbounded growth can again be of use to us. We will see that the local law of the process as seen from a typical `Palmed' infected site at a typical late time can approximately be described by a monotone, translation invariant, harmonic function of the contact process. It is not hard to see that if $\eta^\La_\infty$ is a random variable with law $\Li(\eta^\La_\infty)=\ov\nu$, the upper invariant law, then
\be
f(A):=P[\eta^\La_\infty\cap A\neq\emptyset]
\ee 
also defines an (a priori different) monotone, translation invariant, harmonic function $f$. The key argument in Chapter~\ref{C:typ} uses extinction versus unbounded growth, plus duality, to show that this is up to a multiplicative constant the only such function. This extends the classical result, outlined in Section~\ref{S:euco}, that $\ov\nu$ is the only nontrivial homogeneous invariant law.

\chapter[Renormalization of catalytic WF-diffusions]{Renormalization of catalytic\\ Wright-Fisher diffusions}\label{C:ren}

\section{Introduction}\label{intro}

\subsection{Linearly interacting diffusions}

Let $D\sub\R^d$ be open and convex, let $\ov D$ denote its closure, and assume that $0\in\ov D$. Let $\La$ be a countably infinite group, with group action denoted by $(\xi,\eta)\mapsto\xi\eta$ and unit element $0$. Let $a:\La\times\La\to\R$ be summable and invariant with respect to left multiplication in the group, i.e.,
\be
\sum_{\eta\in\La}|a(\xi,\eta)|<\infty\quad\mbox{and}\quad a(\xi,\eta)=a(\zeta\xi,\zeta\eta)\quad(\xi,\eta,\zeta\in\La),
\ee
and assume that $a$ is irreducible in the sense that for all $\De\sub\La$ with $\De\neq\emptyset,\La$, there exist $\xi\in\De$ and $\eta\in\La\beh\De$ such that either $a(\xi,\eta)\neq 0$ or $a(\eta,\xi)\neq 0$. We assume moreove that
\be\label{bposcor}
a(\xi,\eta)\geq 0\qquad(\xi\neq\eta).
\ee
Consider a collection $\x=(\x_\xi)_{\xi\in\La}$ of $\ov D$-valued processes, solving the martingale problem for the operator
\be
\Ai f(x):=\sum_{\eta,\xi\in\La}a(\eta,\xi)\sum_{i=1}^dx_{\eta,i}\dif{x_{\xi,i}}f(x)+\sum_{\xi\in\La}\sum_{i,j=1}^dw_{ij}(x_\xi)\difif{x_{\xi,i}}{x_{\xi,j}}f(x),
\ee
where we write $x=(x_\xi)_{\xi\in\La}$ and $x_\xi=(x_{\xi,1},\ldots,x_{\xi,d})$ for a point $x\in\ov D^\La$, and the domain of $\Ai$ consists of all functions on $\ov D^\La$ that depend only on finitely many coordinates through a $\Ci^{(2)}$ function of compact support. It is well-known that $\ov D^\La$-valued (weak) solutions to a system of SDE's of the form
\be\label{linSDE}
\di\x_\xi(t)=\sum_{\eta\in\La}a(\eta,\xi)\x_\eta(t)\di t+\sqrt{2}\sig(\x_\xi(t))\di B_\xi(t)\qquad(t\geq 0,\ \xi\in\La),
\ee
solve the martingale problem for $\Ai$, were $(B_\xi)_{\x\in\La}$ is a system of independent $d'$-dimensional Brownian motions, and the $d\times d'$ matrix-valued function $\sig$ is continuous and satisfies
\be
\sum_{k=1}^{d'}\sig_{ik}(x)\sig_{jk}(x)=w_{ij}(x).
\ee
Conversely (see \cite[Theorem~5.3.3]{EK} for the finite dimensional case), every solution to the martingale problem for $\Ai$ can be represented as a solution to the SDE (\ref{linSDE}), where there is some freedom in the choice of the root $\sig$ of the diffusion matrix $w$.

Equation (\ref{linSDE}) says that $\x$ is a system of linearly interacting $d$-dimensional diffusions. As a result of assumption (\ref{bposcor}), the linear drift causes the components $(\x_\xi)_{\xi\in\La}$ to be positively correlated.

Set
\be\label{ladef}
\la:=a(0,0)-\sum_\xi a(0,\xi).
\ee
For reasons that will become clear in a moment (see formula (\ref{mom})~(i) and the remarks below it), if $\la>0$, we have to assume that $\ov D$ is a cone in order for solutions of (\ref{linSDE}) to exist. Under suitable assumptions on the diffusion matrix $w$, it can then be shown that the system of SDE's (\ref{linSDE}) defines a strong Markov process in a Ligget-Spitzer space $\Ei_\ga(\La)$, defined as
\be
\Ei_\ga(\La):=\big\{x\in\ov D^\La:\sum_{\xi\in\La}\ga_\xi|x_\xi|<\infty\big\},
\ee
where $(\ga_\xi)_{\xi\in\La}$ are strictly positive constants such that $\sum_{\xi\in\La}\ga_\xi<\infty$ and $\sum_{\eta\in\La}a(\eta,\xi)\ga_\eta\leq K\ga_\xi$ $(\xi\in\La)$, for some $K<\infty$. The Markov process $\x$ is uniquely defined by the lattice $\La$, the {\em interaction kernel} $a$, the domain $D$, and the {\em diffusion matrix} $w$.

Basic information about the process $\x$ can be obtained by calculating its mean and covariances. Consider a random walk $R=(R_t)_{t\geq 0}$ on $\La$ that jumps from a point $\xi$ to a point $\eta$ with rate $a(\xi,\eta)$ $(\xi\neq\eta)$. This random walk is called the {\em underlying motion} of $\x$. Set
\be
P_t(\xi,\eta):=P^\xi[R_t=\eta].
\ee
and recall the definition of $\la$ in (\ref{ladef}). Write $\x_\xi(t)=(\x_{\xi,1}(t),\ldots,\x_{\xi,d}(t))$. Then
\be\ba{rr@{\,}c@{\,}l}\label{mom}
{\rm (i)}&\dis E[\x_{\xi,i}(t)]&=&\dis e^{\la t}\sum_{\eta\in\La}P_t(\eta,\xi)E[\x_{\eta,i}(0)],\\[5pt]
{\rm (ii)}&\dis\cov(\x_{\xi,i}(t),\x_{\eta,j}(t))&=&\dis e^{2\la t}\sum_{\zeta,\vartheta}P_t(\zeta,\xi)P_t(\vartheta,\eta)\cov(\x_{\zeta,i}(0),\x_{\vartheta,j}(0))\\
&&&\dis+\int_0^te^{2\la s}\sum_\zeta P_s(\zeta,\xi)P_s(\zeta,\eta)E[w_{ij}(\x_\zeta(t-s))]\di s.
\ec
$(t\geq 0,\ \xi,\eta\in\La,\ 1\leq i,j\leq d)$. Let us start the process $\x$ in an initial law $\Li(\x(0))$ that is {\em homogeneous} in the sense that it is invariant with respect to left multiplication in the group, i.e., $\Li((\x_\xi(0))_{\xi\in\La})=\Li((\x_{\zeta\xi}(0))_{\xi\in\La})$ for each $\zeta\in\La$. Then, as a function of the parameter $\la$, the process $\x$ experiences a {\em phase transition} at $\la=0$. If $\la<0$, then in many examples it can be shown that the process started in any homogeneous initial law converges, as $t\to\infty$, to a unique homogeneous invariant law $\nu$. Letting $t\to\infty$ in (\ref{mom})~(i) we see that $\int\nu(\di x)x_{\xi,i}=0$ for each $\xi\in\La$, $i=1,\ldots,d$. On the other hand, as one may guess from (\ref{mom})~(i), for $\la>0$ the process becomes unstable in the sense that the process started in a nonzero homogeneous initial state does not converge to an invariant law, but grows exponentially.

In the {\em critical} case $\la=0$, the long-time behavior of $\x$ is more subtle. Let us call
\be\label{eff}
\pa_w D:=\{x\in\ov D:w_{ij}(x)=0\ \forall i,j=1,\ldots,d\}
\ee
the {\em effective boundary} of $D$ (associated with $w$). Note that $\pa_w D$ is the set of traps of the process $\x$, in the sense that the process started in a constant initial state $\x_\xi(0)=\tet$ $(\xi\in\La)$ with $\tet\in\pa_w D$ satisfies $\x_\xi(t)=\tet$ $(t\geq 0,\ \xi\in\La)$. Let us say an initial law $\Li(\x(0))$ is {\em nontrivial} if $P[\exists\tet\in\pa_w D\ \mbox{s.t.}\ \x_\xi(0)=\tet\ \forall \xi\in\La]=0$. 

A natural question is whether $\x$ has homogeneous nontrivial invariant laws. In order to guess the answer to this question, we must look at the covariance formula (\ref{mom})~(ii). We observe that
\be\label{Green}
G(\xi,\eta):=\int_0^\infty\sum_\zeta P_t(\zeta,\xi)P_t(\zeta,\eta)\di t=E\Big[\int_0^\infty\!1_{\txt\{R^{\dgg,\xi}_t=\ti R^{\dgg,\eta}_t\}}\di t\Big]
\ee
is the expected time spent together by two independent random walks $R^{\dgg,\xi}$ and $\ti R^{\dgg,\eta}$, started in $R^{\dgg,\xi}_0=\xi$ and $\ti R^{\dgg,\eta}_0=\eta$, and jumping from a point $\xi$ to a point $\eta$ with the reversed jump rates $a^\dgg(\xi,\eta):=a(\eta,\xi)$. If $\La$ is an abelian group, with group action denoted by $(\xi,\eta)\mapsto \xi+\eta$, then the difference $R^{\dgg,\xi}_t-\ti R^{\dgg,\eta}_t$ is itself a random walk, with symmetrized jump rates $a_{\rm s}(\xi,\eta):=a(\xi,\eta)+a(\eta,\xi)$, and $G$ is finite if and only this random walk is recurrent. In particular, this is true for finite range jump kernels on $\Z^n$ if and only if $n\leq 2$.

It follows from (\ref{mom})~(ii) that the process $\x$ cannot have nontrivial homogeneous invariant laws with finite second moments if $G(0,0)=\infty$. Indeed, it has been verified for a number of examples of finite range models on $\Z^n$, that $\x$ has nontrivial homogeneous invariant laws if and only if $n>2$. More precisely, in the transient case $n>2$, the process has a nontrivial homogeneous invariant law with mean $\tet$ for each $\tet\in\ov D\beh\pa_w D$, which is the limit law of the process started in any spatially ergodic initial law with mean $\tet$. This type of behavior is called {\em stable behavior}. On the other hand, in the recurrent case $n\leq 2$, the only homogeneous invariant laws of the process are the delta-measures $\de_\tet$ on constant configurations $\tet\in\pa_w D$. In this case, the law of the process started from a spatially ergodic initial law with mean $\tet\in\ov D\beh\pa_w D$ converges, as time tends to infinity, to a convex combination of these delta measures. This means that there are regions in space of growing size, called {\em clusters}, where the process is approximately constant and equal to some $\tet\in\pa_w D$. This type of behavior is called {\em clustering}. 

A general result on stable behavior for $d=1$ (i.e., for one-dimensional domains $D$) can be found in \cite{Shi92}. A general result on clustering for $d=1$ can be found in \cite{CFG96}. Some (weak) general results in dimensions $d\geq 2$ for bounded domains $D$ can be found in \cite{Swa00}. Below, we list some explicit examples that have been treated in the literature.\med

\noi
{\em The Ornstein-Uhlenbeck process} $\ov D=\R$, $w(x)=\al>0$. This is a Gaussian model that has been studied in \cite{Deu89}. This reference also contains results for the subcritical case $\la<0$.\med

\noi
{\em The super-random walk} $\ov D=\half$, $w(x)=\al x$, with $\al>0$. This is the discrete space analogue of the well-known super-Brownian motion \cite{Daw77,Daw93,Eth00}. Both the super-random walk and the super-Brownian motion are continuous-mass branching processes. For these models, the dichotomy between stable behavior and clustering can be proved with the help of Kallenberg's backward tree technique \cite{Kal77,GW91}.\med

\noi
{\em The stepping stone model} $\ov D=[0,1]$, $w(x)=\al x(1-x)$, with resampling parameter $\al>0$. This model, on rather general lattices, has been treated by Shiga \cite{Shi80,Shi80b}, who also gives results for the subcritical case $\la<0$. The diffusion function $w(x)=x(1-x)$ is called the {\em Wright-Fisher diffusion function} and is motivated by applications in population dynamics. Generalizations to other diffusion functions $w:[0,1]\to\R$ that satisfy $w(0)=w(1)=0$ and $w>0$ on $(0,1)$ can be found in \cite{NS80,CG94}. The multidimensional {\em Wright-Fisher diffusion matrix} $w_{ij}(x):=x_i(\de_{ij}-x_j)$ on $\ov D:=\{x\in\R^d:x_i\geq 0,\ \sum_{i=1}^dx_i\leq 1\}$ can be treated with the help of Donnelly and Kurtz's look-down construction \cite{DK96,GLW05}.\med

\noi
{\em Catalytic branching} $\ov D=\half^2$, $w(x)=\left(\ba{@{}cc@{}}\al x_1 & 0\\ 0 & \bet x_1x_2\ea\right)$, with $\al,\bet>0$. This model has been studied in \cite{Pen04}. A continuous space version of this model, the catalytic super-Brownian motion, has been studied in \cite{DF97a,DF97b,EF98,FK99}. A discrete particle version of this model has been studied in \cite{GKW99}.\med

\noi
{\em Mutually catalytic branching} $\ov D=\half^2$, $w(x)=\left(\ba{@{}cc@{}}\al x_1x_2 & 0\\ 0 & \bet x_1x_2\ea\right)$, with $\al,\bet>0$. This model has been studied in \cite{DP98}. Its continuous-space analogue, the mutually catalytic super-Brownian motion, has recieved a lot of attention \cite{DEFMPX02a,DEFMPX02b,DF02,DFMPX03}.\med

\noi
{\em Catalytic Wright-Fisher diffusions} $\ov D=[0,1]^2$, $w(x)=\left(\ba{@{}cc@{}}\al x_1(1-x_1) & 0\\ 0 & p(x_1)x_2(1-x_2)\ea\right)$, where $\al>0$ and the {\em catalyzing function} $p:[0,1]\to\half$ Lipschitz continuous. This model, with the first component replaced by a voter model (which heuristically corresponds to taking $\al=\infty$) has been studied in \cite{GKW01}. This model will also be the main subject of our present chapter.\med

\noi
In the clustering regime (i.e., the case $\La=\Z^n$ with $n\leq 2$, or more generally the case where the quantity $G(0,0)$ from (\ref{Green}) is infinite), it is an interesting problem to determine the {\em clustering distribution}
\be\label{clusdis}
\lim_{t\to\infty}\Li(\x_0(t))
\ee
of the process started in a constant initial state $\x_\xi(0)=\tet$ $(\xi\in\La)$, for all $\tet\in\ov D$. If this limit exists, then it will be concentrated on the effective boundary $\pa_w D$. In dimension $d=1$, when $\pa_w D$ consists of the finite endpoints of the interval $D$, the clustering distribution is trivial. In particular, if $D=[0,1]$, then as a result of (\ref{mom})~(i), it is $\tet\de_1+(1-\tet)\de_0$.

More generally, for any bounded domain $D$ in dimensions $d\geq 1$, let $H_w$ denote the class of $w$-harmonic functions, i.e., functions $h\in\Ci^{(2)}(\ov D)$ satisfying $\sum_{ij}w_{ij}(x)\difif{x_i}{x_j}h(x)=0$ on $D$. Assume that $H_w$ has the property that
\be\label{invhar}
T^c_{x,t}h(H_w)\sub H_w\quad(t\geq 0,\ c>0,\ x\in\ov D),
\ee
where
\be
T^c_{x,t}h(y):= h(x+(y-x)e^{-ct})\quad(t\geq 0,\ c>0,\ x\in\ov D)
\ee
is the semigroup with generator $\sum_{i=1}^dc(x_i-y_i)\dif{y_i}$, i.e., the generator of a deterministic process with a linear drift with strength $c$ towards $x$. Under this assumption, it has been shown in \cite{Swa00} that (\ref{mom})~(i), in the critical case $\la=0$, can be generalized to
\be\label{hexp}
E[h(\x_{\xi,i}(t))]=\sum_{\eta\in\La}P_t(\eta,\xi)E[h(\x_{\eta,i}(0))]\qquad(t\geq 0,\ h\in H_w),
\ee
and this is enough to determine the clustering distribution uniquely. Indeed, the limit in (\ref{clusdis}) must be the unique $H_w$-harmonic measure on $\pa_wD$ with mean $x$. If (\ref{invhar}) holds then we say that $w$ has {\em invariant harmonics}. Diffusion matrices on higher-dimensional domains do not in general have invariant harmonics; this applies in particular to catalytic Wright-Fisher diffusions if the catalyzing function $p$ satisfies $p(0)=0$ and $p(1)>0$.

To get an idea of what the clustering distribution could be in general, we need to analyze the behavior of $\x$ on large space and time scales. We start with the large space-time behavior of the usual stepping stone model.

\subsection{Large space-time behavior}

The behavior of the stepping stone model on $\Z^n$, with resampling parameter $\al$, on large spatial and temporal scales can be studied with the help of its moment dual, a system of rate $\al$ coalescing random walks. In fact, it is in particular the $\al\to\infty$ limit of these models that has been studied in detail, that is, the voter model and its dual, a system of immediately coalescing random walks. A good reference is \cite{CG86}.

In this section, we will especially be interested in the case $n=2$, which is the critical dimension for random walk to be recurrent. Indeed, a 2-dimensional random walk $(R_t)_{t\geq 0}$ is recurrent, but it is only barely so. This is expressed, for example, in the fact that the quantity
\be
E\big[\int_0^t1_{\txt\{R_s=0\}}\big]
\ee
tends very slowly to infinity as $t\to\infty$. (For a precise definition of critical recurrence, see \cite[formula~(1.15)]{Kle96}.) As a result, on $\Z^2$ we see critical phenomenon associated with the phase transition between recurrence and transience.

Let $\x$ be a finite-range stepping stone model on $\Z^2$, started in a constant configuration $\x_\xi(0)=\tet$ $(\xi\in\Z^2)$, for some $\tet\in[0,1]$. Let
\be
\De^t_s:=[0,t^{\frac{1}{2}e^{-s}}]^2\cap\Z^2\qquad(s,t\geq 0)
\ee
be a block of volume $t^{e^{-s}}$, and let
\be\label{xst}
\x^s(t):=\frac{1}{|\De^t_s|}\sum_{\xi\in\De^t_s}\x_\xi(t)\qquad(s,t\geq 0)
\ee
be the average of $\x(t)$ over $\De^t_s$. By combining \cite[Theorem~5]{CG86} and \cite[Theorem~2]{FG94} as described in \cite[Proposition~3.1]{GKW01}, it follows that
\be\label{difclus}
\Li\big((\x^s(t))_{s\geq 0}\big)\overset{\rm f.d.d.}{\underset{t\to\infty}{\Longrightarrow}}(\y_s)_{s\geq 0},
\ee
where $(\y_s)_{s\geq 0}$ is a Wright-Fisher diffusion, i.e., a solution to $\di\y_s=\sqrt{\y_s(1-\y_s)}\di B_s$, started in $\y_0=\tet$. Here f.d.d.\ denotes convergence in finite dimensional distributions. (The question whether the convergence in f.d.d.\ can be replaced by weak convergence in path space is the subject of ongoing research.) Formula (\ref{difclus}) shows how block averages at late times $t$ change as we zoom in in space. Very large block avarages, over blocks of volume $t$, still show the original intensity $\tet$ that the process $\x$ was starting in. As we zoom in on smaller blocks of volume $t^{e^{-s}}$, with $s\geq 0$, the block averages change in a random way, until after some random time, the Wright-Fisher diffusion $\y_s$ hits $0$ or $1$, (with probabilities $1-\tet$ or $\tet$, respectively), and from that random scale on, the block avarages are constant.

Note that the long-time behavior of the limiting diffusion $\y$ in (\ref{difclus}) gives us the clustering distribution (\ref{clusdis}). It seems likely that similar results hold for other models as well; however, the limiting diffusion in (\ref{difclus}) will not always be the Wright-Fisher diffusion. To find out what the limit could be more generally, it is helpful to replace the lattice $\Z^2$ by the hierarchical group, as explained in the next section.

\subsection{Hierarchically interacting diffusions}\label{hier}

For any $N\geq 2$, the {\em hierarchical group} with freedom $N$ is
the set $\om_N$ of all sequences $\xi=(\xi_1,\xi_2,\ldots)$, with
coordinates $\xi_k$ in the finite set $\{0,\ldots,N-1\}$, which are
different from $0$ only finitely often, equipped with componentwise
addition modulo $N$. Setting
\be
\|\xi\|:=\min\{n\geq 0:\xi_k=0\ \forall k>n\}\qquad(\xi\in\om_N),
\ee
$\|\xi-\eta\|$ is said to be the {\em hierarchical distance}
between two sites $\xi$ and $\eta$ in $\om_N$.

Let $\x^N=(\x^N_\xi)_{\xi\in\om_N}$ be a critical system of linearly interacting diffusions on $\om_N$ with interaction kernel given by
\be\label{aN}
a_N(\xi,\eta):=\sum_{k=\|\xi-\eta\|}^\infty\frac{c_{k-1}}{N^{2k-1}}\quad(\xi\neq\eta),\quad a_N(\xi,\xi):=-\sum_{\eta\neq\xi}a_N(\xi,\eta),
\ee
where $(c_k)_{k\geq 0}$ are positive {\em migration constants} such that the quantity $\sum_\xi a_N(0,\xi)=\sum_k c_k/N^k$ is finite. The random walk associated with $a_N$ is recurrent if and only if
\be\label{rec}
\sum_{k=0}^\infty\frac{1}{d_k}=\infty,\qquad\mbox{where}
\quad d_k:=\sum_{n=0}^\infty\frac{c_{k+n}}{N^n}
\ee
(see \cite{DG93a,Kle96}; a similar problem is treated in \cite{DE68}).

Let $\De_k(\xi):=\{\eta:\|\xi-\eta\|\leq k\}$ denote the {\em $k$-block} around $\xi$ and let
\be\label{kblock}
\x^k_\xi(t):=\frac{1}{|\De_k(\xi)|}\sum_{\eta:\|\xi-\eta\|\leq k}\x_\eta(t)
\qquad\qquad(k\geq 0).
\ee 
denote the {\em $k$-block average} around $\xi$. The sequence $(\x^0_0(t),\x^1_0(t),\ldots)$ of block-averages around the origin is called the {\em interaction chain}. Heuristic arguments suggest that in the {\em local mean field limit} $N\to\infty$, the interaction chain converges to a certain well-defined Markov chain. In order to charcterize this chain, we need a few definitions.
\bdf\label{defi}
{\bf (Renormalization class and transformation)} Let $D\sub\R^d$ be nonempty, convex, and open, and let $\ov D$ be its closure. Let $\Wi$ be a collection of continuous functions $w$ from $\ov D$ into the space $M^d_+$ of symmetric non-negative definite $d\times d$ real matrices, such that $\la w\in\Wi$ for every $\la>0$, $w\in\Wi$. We call $\Wi$ a  {\em prerenormalization class} on $\ov D$ if the following three conditions are satisfied:
\begin{enumerate}
\item For each constant $c>0$, $w\in\Wi$, and $x\in\ov D$, the
martingale problem for the operator $A^{c,w}_x$ is well-posed, where
\be\label{Adef}
A^{c,w}_x f(y):=\sum_{i=1}^d c\,(x_i-y_i)\dif{y_i}f(y)
+\sum_{i,j=1}^dw_{ij}(y)\difif{y_i}{y_j}f(y)\qquad(y\in\ov D),
\ee
and the domain of $A^{c,w}_x$ is the space of real functions on $\ov
D$ that can be extended to a twice continuously differentiable
function on $\R^d$ with compact support.
\item For each $c>0$, $w\in\Wi$, and $x\in\ov D$, the martingale
problem for $A^{c,w}_x$ has a unique stationary solution with
invariant law denoted by $\nu^{c,w}_x$.
\item For each $c>0$, $w\in\Wi$, $x\in\ov D$, and $i,j=1,\ldots,d$,
one has $\dis\int_{\ov D}\nu^{c,w}_x(\di y)|w_{ij}(y)|<\infty$.
\end{enumerate}
If $\Wi$ is a prerenormalization class, then we define for each $c>0$
and $w\in\Wi$ a matrix-valued function $F_cw$ on $\ov D$ by
\be\label{Fc}
F_cw(x):=\int_{\ov D}\nu^{c,w}_x(dy)w(y)\qquad(x\in\ov D).
\ee
We say that $\Wi$ is a {\em renormalization class} on $\ov D$ if in addition:
\begin{enumerate}\addtocounter{enumi}{3}
\item For each $c>0$ and $w\in\Wi$, the function $F_cw$ is an element of $\Wi$.
\end{enumerate}
If $\Wi$ is a renormalization class and $c>0$, then the map
$F_c:\Wi\to\Wi$ defined by (\ref{Fc}) is called the 
{\em renormalization transformation} on $\Wi$ with {\em migration
constant} $c$. In (\ref{Adef}), $w$ is called the {\em diffusion
matrix} and $x$ the {\em attraction point.}
\edf
For any renormalization class $\Wi$ and any sequence of (strictly)
positive migration constants $(c_k)_{k\geq 0}$, we define {\em
iterated renormalization transformations} $F^{(n)}:\Wi\to\Wi$, as
follows:
\be\label{Fn}
F^{(n+1)}w:=F_{c_n}(F^{(n)}w)\quad(n\geq 0)\quad\mbox{with}
\quad F^{(0)}w:=w\qquad(w\in\Wi_{\rm cat}).
\ee
We set $s_0:=0$ and
\be\label{sn}
s_n:=\sum_{k=0}^{n-1}\frac{1}{c_k}\qquad(1\leq n\leq\infty).
\ee
With these definitions, we can formulate the following conjecture about the behavior of the interaction chain in the local mean field limit $N\to\infty$.
\begin{conjecture}\label{Ichain}
Let $\Wi$ be a renormalization class. Fix $w\in\Wi$, $\tet\in D$, and positive numbers $(c_k)_{k\geq 0}$ such that for $N$ large enough, $\sum_k c_k/N^k<\infty$. For all $N$ large enough, let $\x^N$ be a solution to (\ref{linSDE}) on $\La=\om_N$ with $a=a_N$ from (\ref{aN}), and assume that $t_N$ are constants such that, for some $n\geq 1$, $\lim_{N\to\infty}N^{-n}t_N=T\in\half$. Then
\be\label{toI}
\Big(\x^{N,n}_0(t_N),\ldots,\x^{N,0}_0(t_N)\Big)
\underset{N\to\infty}{\dis\Longrightarrow}(I^{w}_{-n},\ldots,I^{w}_0),
\ee
where $(I^{w}_{-n},\ldots,I^{w}_0)$ is a Markov chain with transition laws
\be\label{transspec}
P[I^{w}_{-k}\in\di y|I^{w}_{-k-1}=x]=\nu^{c_k,F^{(k)}w}_x(\di y)
\qquad(x\in\ov D,\ 0\leq k\leq n-1)
\ee
and initial state
\be\label{tinit}
I^{w}_{-n}=\y_T,\quad\mbox{where}\quad\di\y_t=c_n(\tet-\y_t)\di t
+\sqrt{2}\sig^{(n)}(\y_t)\di B_t,\quad\y_0=\tet,
\ee
and $\sig^{(n)}$ is a root of the diffusion matrix $F^{(n)}w$.
\end{conjecture}
Rigorous versions of conjecture~\ref{Ichain} have been proved for
renormalization classes on $\ov D=[0,1]$ and $\ov D=\half$ in
\cite{DG93a,DG93b}. See \cite{DG96,DGV95} for similar results. Note
that the Markov chain $I^{w}=(I^{w}_{-n},\ldots,I^{w}_0)$ is a sort of
analogue of the block averages $(\x^s(t))_{s\geq 0}$ defined in
(\ref{xst}). As we will see below, for appropriate choices of the
constants $(c_k)_{k\geq 0}$, the discrete chain $I^{w}$ can be
approximated by a diffusion, in the spirit of (\ref{difclus}). In
order to see this, we need a few facts about renormalization classes.
To keep things as simple as possible, we specialize to renormalization
classes on bounded domains, although much of what we will say, with
some modifications here and there, can be generalized to unbounded
domains.

\subsection{Renormalization classes}\label{gensec}

In this section, we describe some elementary properties that hold
generally for (pre-) renormalization classes on bounded domains. The
proofs of Lemmas~\ref{contlem}--\ref{clusK} can be found in
Section~\ref{gensub} below.

Fix a prerenormalization class $\Wi$ on a set $\ov D$ where
$D\sub\R^d$ is open, bounded, and convex. Then $\Wi$ is a subset of
the cone $\Ci(\ov D,M^d_+)$ of continuous $M^d_+$-valued functions on
$\ov D$. We equip $\Ci(\ov D,M^d_+)$ with the topology of uniform
convergence. We let $\Mi_1(\ov D)$ denote the space of probability
measures on $\ov D$, equipped with the topology of weak
convergence. Our first lemma says that the equilibrium measures
$\nu^{c,w}_x$ and the renormalized diffusion matrices $F_cw(x)$ are
continuous in their parameters.
\bl{\bf(Continuity in parameters)}\label{contlem}
\begin{itemize}
\item[{\bf (a)}] The map $(x,c,w)\mapsto\nu^{c,w}_x$ from
$\ov D\times(0,\infty)\times\Wi$ into $\Mi_1(\ov D)$ is continuous.
\item[{\bf (b)}] The map $(x,c,w)\mapsto F_cw(x)$ from
$\ov D\times(0,\infty)\times\Wi$ into $M^d_+$ is continuous.
\end{itemize}
\el
In particular, $x\mapsto\nu^{c,w}_x$ is a continuous probability
kernel on $\ov D$, and $F_cw\in\Ci(\ov D,M^d_+)$ for all $c>0$ and
$w\in\Wi$. Recall from Definition~\ref{defi} that $\la w\in\Wi$ for
all $w\in\Wi$ and $\la>0$. The reason why we have included this
assumption is that it is convenient to have the next scaling lemma
around, which is a consequence of time scaling.
\bl{\bf(Scaling property of renormalization transformations)}\label{schaal}
One has
\be\left.\ba{rr@{\,}c@{\,}l}\label{schafo}
{\rm (i)}&\nu^{\la c,\la w}_x&=&\nu^{c,w}_x\\
{\rm (ii)}&F_{\la c}(\la w)&=&\la F_cw\\
\ea\qquad\right\}\quad(\la,c>0,\ w\in\Wi,\ x\in\ov D).
\ee
\el
The following simple lemma will play a crucial role in what follows.
\bl{\bf (Mean and covariance matrix)}\label{numom}
For all $x\in\ov D$ and $i,j=1,\ldots,d$, the mean and covariances of
$\nu^{c,w}_x$ are given by
\be\ba{rr@{\;}c@{\;}l}\label{moments2}
{\rm (i)}&\dis\int_{\ov D}\nu^{c,w}_x(\di y)(y_i-x_i)&=&0,\\
{\rm (ii)}&\dis\int_{\ov D}\nu^{c,w}_x(\di y)(y_i-x_i)(y_j-x_j)
&=&\frac{1}{c}F_cw_{ij}(x).
\ec 
\el
Recall the definition of the effective boundary associated with a
diffusion matrix $w$ in (\ref{eff}). The next lemma says that the effective
boundary is invariant under renormalization.
\bl[Invariance of effective boundary]\label{effinv}
One has $\pa_{F_cw}D=\pa_wD$ for all $w\in\Wi$, $c>0$.
\el
{F}rom now on, let $\Wi$ be a renormalization class, i.e., $\Wi$
satisfies also condition~(iv) from Definition~\ref{defi}. Fix a
sequence of (positive) migration constants $(c_k)_{k\geq 0}$. By
definition, the {\em iterated probability kernels} $K^{w,(n)}$
associated with a diffusion matrix $w\in\Wi$ (and the constants
$(c_k)_{k\geq 0}$) are the probability kernels on $\ov D$ defined
inductively by
\be\label{Kdef}
K^{w,(n+1)}_x(\di z):=
\int_{\ov D}\nu^{c_n,F^{(n)}w}_x(\di y)K^{w,(n)}_y(\di z)
\quad(n\geq 0)\quad\mbox{with}\quad K^{w,(0)}_x(\di y):=\de_x(\di y),
\ee
with $F^{(n)}$ as in (\ref{Fn}). Note that $K^{w,(n)}$ is the transition probability from time $-n$ to time $0$ of the interaction chain in the local mean-field limit (see Conjecture~\ref{Ichain}):
\be\label{Kdef2}
K^{w,(n)}_x(\di y):=P[I^{w}_0\in\di y|I^{w}_{-n}=x]
\qquad\quad(x\in\ov D,\ n\geq 0).
\ee
Note moreover that
\be\label{KF}
F^{(n)}w(x)=\int_{\ov D} K^{w,(n)}_x(\di y)w(y)\qquad(x\in\ov D,\ n\geq 0).
\ee
The next lemma follows by iteration from Lemmas~\ref{contlem} and
\ref{numom}. It their essence, this lemma
and Lemma~\ref{clusK} below go back to \cite{BCGH95}.
\bl[Basic properties of iterated kernels]\label{basic}
For each $w\in\Wi$, the $K^{w,(n)}$ are continuous probability kernels
on $\ov D$. Moreover, for all $x\in\ov D$, $i,j=1,\ldots,d$, and $n\geq 0$,
the mean and covariance matrix of $K^{w,(n)}_x$ are given by
\be\ba{rr@{\;}c@{\;}l}\label{moments}
{\rm (i)}&\dis\int_{\ov D} K^{w,(n)}_x(\di y)(y_i-x_i)&=&0,\\
{\rm (ii)}&\dis\int_{\ov D} K^{w,(n)}_x(\di y)(y_i-x_i)(y_j-x_j)
&=&s_nF^{(n)}w_{ij}(x).
\ec
\el
We equip the space $\Ci(\ov D,\Mi_1(\ov D))$ of continuous probability
kernels on $\ov D$ with the topology of uniform convergence (since
$\Mi_1(\ov D)$ is compact, there is a unique uniform structure on
$\Mi_1(\ov D)$ generating the topology). For `nice' renormalization
classes, it seems reasonable to conjecture that the kernels
$K^{w,(n)}$ converge as $n\to\infty$ to some limit $K^{w,\ast}$ in
$\Ci(\ov D,\Mi_1(\ov D))$. If this happens, then formula
(\ref{moments})~(ii) tells us that the rescaled renormalized diffusion
matrices $s_nF^{(n)}w$ converge uniformly on $\ov D$ to the covariance
matrix of $K^{w,\ast}$.

We will mainly be interested in the case that
$\lim_{n\to\infty}s_n=\infty$. Indeed, if the iterated kernels
converge to a limit $K^{w,\ast}$, then this condition
guarantees that this limit is concentrated on the effective boundary:
\bl[Concentration on the effective boundary]\label{clusK}
If $s_n\asto{n}\infty$, then for any $f\in\Ci(\ov D)$ such that $f=0$ on
$\pa_w D$:
\be\label{concen}
\lim_{n\to\infty}\sup_{x\in\ov D}
\Big|\int_{\ov D}K^{w,(n)}_x(\di y)f(y)\Big|=0.
\ee
\el
Note that $s_n\to\infty$ if and only if $\sum_k1/c_k=\infty$. We can
think of this condition as the $N\to\infty$ limit of the condition
$\sum_k 1/d_k=\infty$ in (\ref{rec}). Thus, the condition
$s_n\to\infty$ guarantees that the corresponding system of linearly
interacting diffusions on the hierarchical group with migration
constants $(c_k)_{k\geq 0}$ {\em clusters in the local mean field
limit}.

Most of the discussion in this section carries over to renormalization
classes on unbounded $D$, but in this case, the second moments of the
iterated kernels $K^{w,(n)}$ may diverge as $n\to\infty$. As a result,
because of formula (\ref{moments})~(ii), the $s_n$ may no longer be
the right scaling factors to find a nontrivial limit of the
renormalized diffusion matrices; see, for example, \cite{BCGH97}.

\subsection{Rescaled transformations}

We return to renormalization classes on bounded domains, and focus our
attention on the clustering regime $s_n\to\infty$. Since we expect
$s_nF^{(n)}w$ to converge to a limit (namely, the covariance matrix of
$K^{w,\ast}$), we will use Lemma~\ref{schaal} to convert the rescaled
iterates $s_nF^{(n)}$ into (usual, not rescaled) iterates of another
transformation. For this purpose, it will be convenient to modify the
definition of our scaling constants $s_n$ a little bit. Fix some
$\bet>0$ and put
\be\label{hutsn}
\hut s_n:=\bet+s_n\qquad(n\geq 0).
\ee
Define {\em rescaled renormalization transformations} $\hut F_\ga:\Wi\to\Wi$ by
\be\label{hatFc}
\hut F_\ga w:=(1+\ga)F_{1/\ga}w\qquad(\ga>0,\ w\in\Wi).
\ee
Using (\ref{schafo})~(ii), one easily deduces that
\be\label{sFw}
\hut s_nF^{(n)}w=\hut F_{\ga_{n-1}}\circ\cdots\circ\hut F_{\ga_0}(\bet w)
\qquad(w\in\Wi,\ n\geq 1),
\ee
where
\be\label{gan}
\ga_n:=\frac{1}{\hut s_nc_n}\qquad(n\geq 0).
\ee

We can reformulate the condition $s_n\to\infty$ from Lemma~\ref{clusK}
in terms of the constants $(\ga_n)_{n\geq 0}$. Indeed, it is not hard
to check\footnote{To see this, let $\ov s_\infty\in(0,\infty]$ denote
the limit of the $\ov s_n$ and note that on the one
hand,\label{footnote} $\sum_n1/(\hut s_nc_n)\geq\sum_n\log(1+1/(\hut
s_nc_n))=\log(\prod_n\hut s_{n+1}/\hut s_n)=\log(\hut s_\infty/\hut
s_1)$, while on the other hand $\sum_n 1/(\hut
s_nc_n)\leq\prod_n(1+1/(\hut s_nc_n))=\prod_n\hut s_{n+1}/\hut
s_n=\hut s_\infty/\hut s_1$.} that the following three conditions are
equivalent:
\be\label{sumga}
{\rm (i)}\quad s_n\asto{n}\infty,\qquad{\rm (ii)}
\quad\hut s_n\asto{n}\infty,\qquad{\rm (iii)}\quad\sum_n\ga_n=\infty.
\ee
In view of (\ref{sFw}), it is natural to assume that the $\ga_n$
converge to a limit $\ga^\ast\in[0,\infty]$. Since $\hut s_{n+1}/\hut
s_n=1+\ga_n$, it is not hard to see that the following conditions are
equivalent:
\be\label{gamma}
{\rm (i)}\quad\frac{s_{n+1}}{s_n}\asto{n}1+\ga^\ast,
\qquad{\rm (ii)}\quad\frac{\hut s_{n+1}}{\hut s_n}\asto{n}1+\ga^\ast,
\qquad{\rm (iii)}\quad\ga_n\asto{n}\ga^\ast.
\ee
If $0<\ga^\ast<\infty$, then, in the light of (\ref{sFw}), we expect
$\hut s_nF^{(n)}w$ to converge to a fixed point of the transformation
$\hut F_{\ga^\ast}$. If $\ga^\ast=0$, the situation is more complex.
In this case, we expect the orbit $\hut s_nF^{(n)}w\mapsto\hut
s_{n+1}F^{(n+1)}w\mapsto\cdots$, for large $n$, to approximate a
continuous flow, the generator of which is
\be\label{gener}
\lim_{\ga\to 0}\ga^{-1}\Big(\hut F_\ga w-w\Big)(x)
=\ffrac{1}{2}\sum_{i,j=1}^dw_{ij}(x)\difif{x_i}{x_j}w(x)+w(x)\qquad(x\in\ov D).
\ee
To see that the right-hand side of this equation equals the left-hand
side if $w$ is twice continuously differentiable, one needs a Taylor
expansion of $w$ together with the moment formulas (\ref{moments2})
for $\nu^{1/\ga,w}_x$. Under condition condition~(\ref{sumga})~(iii),
we expect this continuous flow to reach equilibrium.

In the light if these considerations, we are led to at the following
general conjecture.
\begin{conjecture}\label{con}{\bf (Limits of rescaled renormalized
diffusion matrices)}
Assume that $s_n\to\infty$ and $s_{n+1}/s_n\to 1+\ga^\ast$ for some
$\ga^\ast\in[0,\infty]$. Then, for any $w\in\Wi$,
\be\label{concon}
s_nF^{(n)}w\asto{n}w^\ast,
\ee
where $w^\ast$ satisfies
\be\ba{rr@{\,}c@{\,}ll}\label{afp}
{\rm (i)}&\dis\hut F_{\ga^\ast}w^\ast&=&\dis w^\ast
\qquad&\mbox{if}\ \ 0<\ga^\ast<\infty,\\[5pt]
{\rm (ii)}&\dis\ffrac{1}{2}\sum_{i,j=1}^d
w^\ast_{ij}(x)\difif{x_i}{x_j}w^\ast(x)+w^\ast(x)&=&\dis 0
\qquad(x\in\ov D)\quad&\mbox{if}\ \ \ga^\ast=0,\\[5pt]
{\rm (iii)}&\dis\lim_{\ga\to\infty}\ov F_\ga w^\ast&=&w^\ast
\qquad&\mbox{if}\ \ \ga^\ast=\infty.
\ec
\end{conjecture}
We call (\ref{afp})~(ii), which is in some sense the $\ga^\ast\to 0$
limit of the fixed point equation (\ref{afp})~(i), the {\em asymptotic
fixed point equation}. A version of formula (\ref{afp})~(ii) occurred
in \cite[formula~(1.3.5)]{Swa99} (a minus sign is missing there).

In particular, one may hope that for a given effective boundary, the
equations in (\ref{afp}) have a unique solution. Our main result
(Theorem~\ref{main} below) confirms this conjecture for a
renormalization class of catalytic Wright-Fisher diffusions and for
$\ga^\ast<\infty$. In Section~\ref{numer} below, we discuss numerical
evidence that supports Conjecture~\ref{con} in the case $\ga^\ast=0$
for other renormalization classes on compacta as well.

\subsection{Diffusive clustering}

Assuming that the rescaled renormalized diffusion matrices
$s_nF^{(n)}w$ converge to a limit $w^\ast$, we can make a guess about
the limit of the iterated probability kernels $K^{w,(n)}$.
\begin{conjecture}{\bf(Limits of iterated probability kernels)}\label{Kcon}
Assume that $s_nF^{(n)}w\to w^\ast$ as $n\to\infty$. Then, for any $w\in\Wi$,
\be
K^{w,(n)}\asto{n}K^\ast,
\ee
where $K^\ast$ has the following description:
\renewcommand{\labelenumi}{\rm(\roman{enumi})}
\begin{enumerate}
\item If $0<\ga^\ast<\infty$, then
\be
K^\ast_x=\lim_{n\to\infty}P^x[I^{\ga^\ast}_n\in\,\cdot\,],
\ee
where $(I^{\ga^\ast}_n)_{n\geq 0}$ is the Markov chain with
transition law $P[I^{\ga^\ast}_{n+1}\in\cdot\,|I^{\ga^\ast}_n=x]
=\nu^{1/\ga^\ast,w^\ast}$. 
\item If $\ga^\ast=0$, then
\be\label{Kcrit}
K^\ast_x=\lim_{t\to\infty}P^x[I^0_t\in\,\cdot\,],
\ee
where $(I^0_s)_{s\geq 0}$ is the diffusion process with generator
$\sum_{i,j=1}^dw^\ast_{ij}(y)\difif{y_i}{y_j}$.
\item If $\ga^\ast=\infty$, then
\be
K^\ast_x=\lim_{\ga\to\infty}\nu^{1/\ga,w^\ast}_x.
\ee
\end{enumerate}
\end{conjecture}
If $\ga^\ast<\infty$, this conjecture is motivated by the observation that in this case, the Markov chain $(I^w_{-n},\ldots,I^w_0)$ from Conjecture~\ref{Ichain} is approximately time homogeneous for $n\to\infty$. The case $\ga^\ast=0$ is of particular interest. In this case $I^w_{-n},I^w_{-n+1},\ldots$ converges, in the right scaling, to the diffusion $(I^0_s)_{s\geq 0}$ with diffusion matrix $w^\ast$. This is a sort of analogon of the diffusive clustering result (\ref{difclus}). Based on this analogy, we can make one more conjecture.
\begin{conjecture}{\bf (Clustering distribution on $\Z^2$)}\label{Z2con}
Let $D\sub\R^d$ be open, bounded, and convex, and let $\Wi$ be a
renormalization class on $\ov D$. Assume that the asymptotic fixed
point equation (\ref{afp})~(ii) has a unique solution $w^\ast$ in
$\Wi$. Let $\sig$ be a continuous root of a diffusion matrix
$w\in\Wi$. Let $\x=(\x_\xi)_{\xi\in\Z^2}$ be a $\ov D^{\Z^2}$-valued
process, solving the system of SDE's
\be
\di\x_\xi(t)=\sum_{\eta:\,|\eta-\xi|=1}
\!\big(\x_\eta(t)-\x_\xi(t)\big)\,\di t+\sig(\x_\xi(t))\di B_\xi(t),
\ee
with initial condition $\x_\xi(0)=\tet\in\ov D$ $(\xi\in\Z^2)$. Then
\be
\Li(\x_0(t))\Asto{t}P[I^0_\infty\,|\,I^0_0=\tet]\qquad(\xi\in\Z^2),
\ee
where $(I^0_s)_{s\geq 0}$ is the diffusion with generator
$\sum_{ij}w^\ast_{ij}(y)\difif{y_i}{y_j}$.
\end{conjecture}

\subsection{Numerical solutions to the asymptotic fixed point equation}
\label{numer}

Let $t\mapsto w(t,\,\cdot\,)$ be a solution to the continuous flow
with the generator in (\ref{gener}), i.e., $w$ is an $M^d_+$-valued
solution to the nonlinear partial differential equation
\be\label{flow}
\dif{t}w(t,x)=\ffrac{1}{2}\sum_{i,j=1}^dw_{ij}(t,x)\difif{x_i}{x_j}w(t,x)
+w(t,x)\qquad(t\geq 0,\ x\in\ov D).
\ee
Solutions to (\ref{flow}) are quite easy to simulate on a computer. We
have simulated solutions for all kind of diffusion matrices (including
nondiagonal ones) on the unit square $[0,1]^2$, with the effective
boundaries 1--6 depicted in Figure~\ref{fixfig}. For all initial
diffusion matrices $w(0,\,\cdot\,)$ we tried, the solution converged
as $t\to\infty$ to a fixed point $w^\ast$. In all cases except case~6,
the fixed point was unique. The fixed points are listed in
Figure~\ref{fixfig}. The functions $p^\ast_{0,1,0}$ and $q^\ast$ from
Figure~\ref{fixfig} are plotted in Figure~\ref{graph}.

\begin{figure}
\begin{center}
\setlength{\unitlength}{1cm}
\begin{tabular}{|c|c|c|}
\hline
case & effective boundary & fixed points $w^\ast$ of (\ref{flow})\\
\hline
1& \begin{picture}(1,1.1)(.2,.2)
   \put(-.2,-.2){\framebox(1,1)[c]{}}
   \put(-.2,-.2){\circle*{.17}}
   \put(-.2,.8){\circle*{.17}}
   \put(.8,-.2){\circle*{.17}}
   \put(.8,.8){\circle*{.17}}  
   \end{picture}
& $\left(\ba{@{}cc@{}}x_1(1-x_1)&0\\0&x_2(1-x_2)\ea\right)$\\
2& \begin{picture}(1,1.1)(.2,.2)
   \put(-.2,-.2){\framebox(1,1)[c]{}}
   \put(-.2,-.2){\circle*{.17}}
   \put(-.2,.8){\circle*{.17}}
   \put(.8,-.2){\circle*{.17}}
   \put(.8,.8){\circle*{.17}}  
   \linethickness{.17cm}
   \put(-.2,-.2){\line(0,1){1}}
   \end{picture}
 & $\left(\ba{@{}cc@{}}x_1(1-x_1)&0\\0&p^\ast_{0,1,0}(x_1)x_2(1-x_2)\ea\right)$\\
3& \begin{picture}(1,1.1)(.2,.2)
   \put(-.2,-.2){\framebox(1,1)[c]{}}
   \put(-.2,-.2){\circle*{.17}}
   \put(-.2,.8){\circle*{.17}}
   \put(.8,-.2){\circle*{.17}}
   \put(.8,.8){\circle*{.17}}  
   \linethickness{.17cm}
   \put(-.2,-.2){\line(0,1){1}}
   \put(-.2,-.2){\line(1,0){1}}
   \end{picture}
 & $\left(\ba{@{}cc@{}}q^\ast(x_1,x_2)&0\\0&q^\ast(x_2,x_1)\ea\right)$\\
4& \begin{picture}(1,1.1)(.2,.2)
   \put(-.2,-.2){\framebox(1,1)[c]{}}
   \put(-.2,-.2){\circle*{.17}}
   \put(-.2,.8){\circle*{.17}}
   \put(.8,-.2){\circle*{.17}}
   \put(.8,.8){\circle*{.17}}  
   \linethickness{.17cm}
   \put(-.2,-.2){\line(0,1){1}}
   \put(.8,-.2){\line(0,1){1}}
   \end{picture}
 & $\left(\ba{@{}cc@{}}x_1(1-x_1)&0\\0&0\ea\right)$\\
5& \begin{picture}(1,1.1)(.2,.2)
   \put(-.2,-.2){\framebox(1,1)[c]{}}
   \put(-.2,-.2){\circle*{.17}}
   \put(-.2,.8){\circle*{.17}}
   \put(.8,-.2){\circle*{.17}}
   \put(.8,.8){\circle*{.17}}  
   \linethickness{.17cm}
   \put(-.2,-.2){\line(0,1){1}}
   \put(-.2,-.2){\line(1,0){1}}
   \put(.8,-.2){\line(0,1){1}}
   \end{picture}
 & $\left(\ba{@{}cc@{}}x_1(1-x_1)1_{\{x_2>0\}}&0\\0&0\ea\right)$\\
6& \begin{picture}(1,1.1)(.2,.2)
   \put(-.2,-.2){\framebox(1,1)[c]{}}
   \put(-.2,-.2){\circle*{.17}}
   \put(-.2,.8){\circle*{.17}}
   \put(.8,-.2){\circle*{.17}}
   \put(.8,.8){\circle*{.17}}  
   \linethickness{.17cm}
   \put(-.2,-.2){\line(0,1){1}}
   \put(-.2,-.2){\line(1,0){1}}
   \put(.8,-.2){\line(0,1){1}}
   \put(-.2,.8){\line(1,0){1}}
   \end{picture}
 & $g^\ast(x_1,x_2)\left(\ba{@{}cc@{}} m_{11}&m_{12}\\m_{21}&m_{22}\ea\right)$\\[14pt]
\hline
\end{tabular}
\caption[Fixed points of the flow (\ref{flow}).]
{Fixed points of the flow (\ref{flow}).}\label{fixfig}
\end{center}
\end{figure}

The fixed points for the effective boundaries in cases 1,2, and 4 will
be described in Theorem~\ref{main} below. In particular,
$p^\ast_{0,1,0}$ is the function from Theorem~\ref{main}~(c). The
simulations suggest that the domain of attraction of these fixed
points (within the class of ``all'' diffusion matrices on $[0,1]^2$)
is actually a lot larger than the classes for which we are able to prove
convergence in Theorem~\ref{main}.

The function $q^\ast$ from case~3 satisfies $q^\ast(x_1,1)=x_1(1-x_1)$
and is zero on the other parts of the boundary. In contrast to what
one might perhaps guess in view of case~2, $q^\ast$ is {\em not} of
the form $q^\ast(x_1,x_2)=f(x_2)x_1(1-x_1)$ for some function $f$.

Case~5 is somewhat degenerate since in this case the fixed point is
not continuous.

The only case where the fixed point is not unique is case~6. Here, $m$
can be any positive definite matrix, while $g^\ast$, depending on $m$,
is the unique solution on $(0,1)^2$ of the equation
$1+\frac{1}{2}\sum_{i,j=1}^2m_{ij}\difif{x_i}{x_i}g^\ast(x)=0$, with
zero boundary conditions. Some diffusion matrices that are in the
domain of attraction of these fixed points are described in
Theorem~\ref{isoTheorem} below. The simulations indicate that the true
domain of attraction is much larger than what can be proved (and
includes nonisotropic matrices).

\begin{figure}
\begin{center}
\includegraphics[width=4cm,height=4cm]{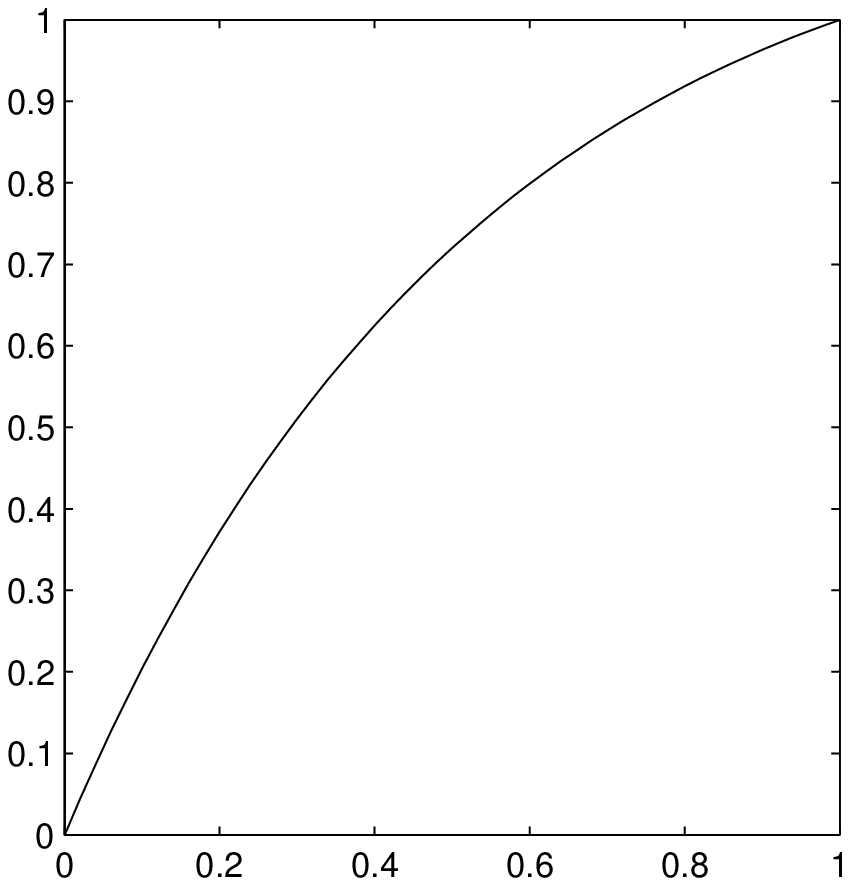}
\hspace{1cm}
\includegraphics[width=5cm,height=4cm]{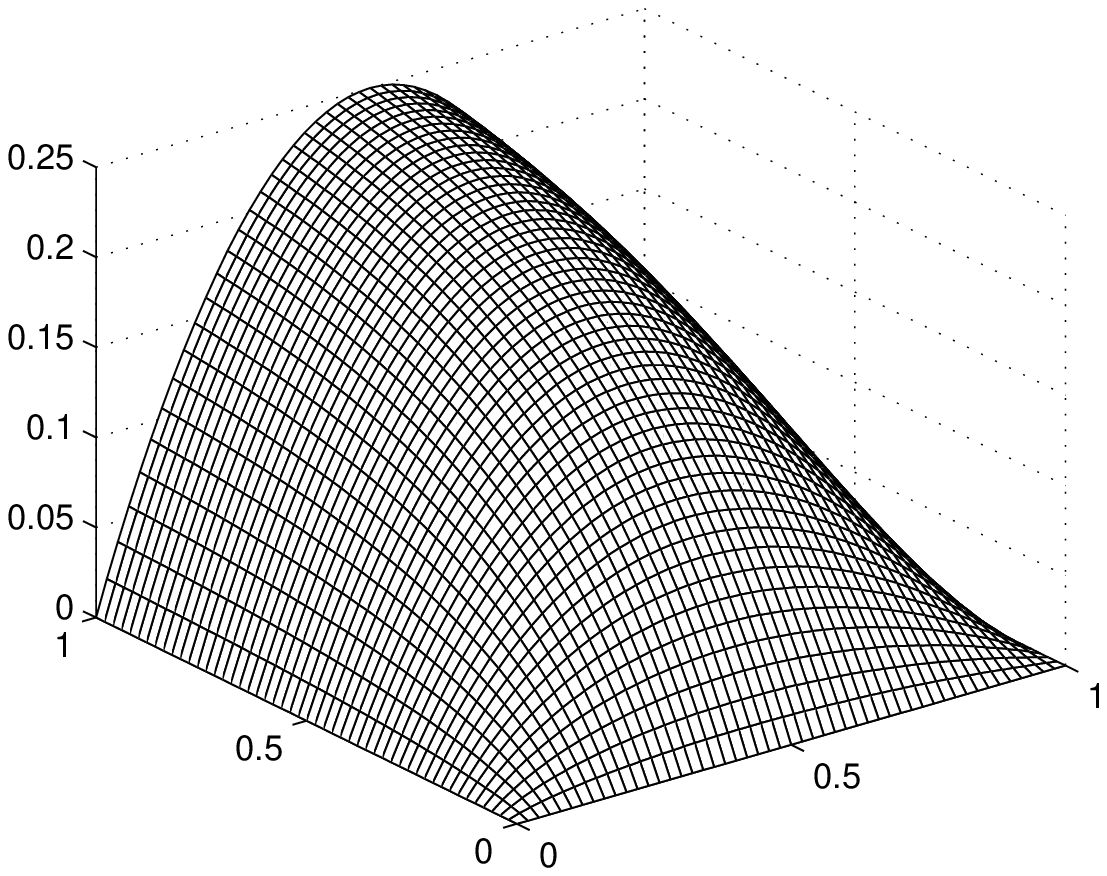}
\caption[The functions $p^\ast_{0,1,0}$ and $q^\ast$ from cases~2 and 3
of Figure~\ref{fixfig}.]{The functions $p^\ast_{0,1,0}$ and $q^\ast$ from
cases~2 and 3 of Figure~\ref{fixfig}.}\label{graph}
\end{center}
\end{figure}

\subsection{Known results}\label{prevsec}

In this section we discuss some results that have been derived previously
for renormalization classes on compact sets.
\bt\label{WFTheorem}
{\bf\cite{BCGH95,DGV95} (Universality class of Wright-Fisher models)}
Let $D:=\{x\in\R^d:x_i>0\ \forall i,\ \sum_{i=1}^dx_i<1\}$, and let
$\{e_0,\ldots,e_d\}$, with $e_0:=(0,\ldots,0)$ and
$e_1:=(1,0,\ldots,0),\ldots,\ e_d:=(0,\ldots,0,1)$ be the extremal
points of $\ov D$. Let $w^\ast_{ij}(x):=x_i(\de_{ij}-x_j)$ $(x\in\ov
D$, $i,j=1,\ldots,d)$ denote the standard Wright-Fisher diffusion
matrix, and assume that $\Wi$ is a renormalization class on $\ov D$
such that $w^\ast\in\Wi$ and $\pa_w\ov D=\{e_0,\ldots,e_d\}$ for all
$w\in\Wi$. Let $(c_k)_{k\geq 0}$ be migration constants such that
$s_n\to\infty$ as $n\to\infty$. Then, for all $w\in\Wi$, uniformly on
$\ov D$,
\be\label{toWF}
s_n F^{(n)}w\asto{n} w^\ast.
\ee
\et
The convergence in (\ref{toWF}) is a consequence of Lemmas~\ref{basic}
and \ref{clusK}: The first moment formula (\ref{moments})~(i) and
(\ref{concen}) show that $K^{w,(n)}_x$ converges to the unique
distribution on $\{e_0,\ldots,e_d\}$ with mean $x$, and by the second
moment formula (\ref{moments})~(ii) this implies the convergence of
$s_nF^{(n)}w$.

In order for the iterates in (\ref{toWF}) to be well-defined,
Theorem~\ref{WFTheorem} {\em assumes} that a renormalization class
$\Wi$ of diffusion matrices $w$ on $\ov D$ with effective boundary
$\{e_0,\ldots,e_d\}$ is given. The problem of finding a nontrivial
example of such a renormalization class is open in dimensions greater
than one. In the one-dimensional case, however, the following result
is known.
\bl{\bf \cite{DG93b} (Renormalization class on the unit interval)}
The set
\be\label{WWF}
\Widg:=\{w\in\Ci[0,1]:w=0\mbox{ on }\{0,1\},\ w>0\mbox{ on }(0,1),
\ w\mbox{ Lipschitz}\}
\ee
is a renormalization class on $[0,1]$.
\el
About renormalization of isotropic diffusions, the following result
is known. Below, $\pa D:=\ov D\beh D$ denotes the topological boundary of $D$.
\bt\label{isoTheorem}
{\bf\cite{HS98} (Universality class of isotropic models)}
Let $D\sub\R^d$ be open, bounded, and convex and let $m\in M^d_+$ be
fixed and (strictly) positive definite. Set
$w^\ast_{ij}(x):=m_{ij}g^\ast(x)$, where $g^\ast$ is the unique
solution of $1+\ffrac{1}{2}\sum_{ij}m_{ij}\difif{x_i}{x_j}g^\ast(x)=0$
for $x\in D$ and $g^\ast(x)=0$ for $x\in\pa D$. Assume that $\Wi$ is a
renormalization class on $\ov D$ such that $w^\ast\in\Wi$ and such
that each $w\in\Wi$ is of the form
\be
w_{ij}(x)=m_{ij}g(x)\qquad(x\in\ov D,\ i,j=1,\ldots,d),
\ee
for some $g\in\Ci(\ov D)$ satisfying $g>0$ on $D$ and $g=0$ on $\pa
D$. Let $(c_k)_{k\geq 0}$ be migration constants such that
$s_n\to\infty$ as $n\to\infty$. Then, for all $w\in\Wi$, uniformly on
$\ov D$,
\be\label{togast}
s_n F^{(n)}w\asto{n}w^\ast.
\ee
\et
The proof of Theorem~\ref{isoTheorem} follows the same lines as the
proof of Theorem~\ref{WFTheorem}, with the difference that in this
case one needs to generalize the first moment formula
(\ref{moments})~(i) in the sense that $\int_{\ov D} K^{w,(n)}_x(\di
y)h(y)=h(x)$ for any $m$-harmonic function $h$, i.e., $h\in\Ci^{(2)}(\ov D)$
satisfying $\sum_{ij}m_{ij}\difif{x_i}{x_j}h(x)=0$ for $x\in D$. The
kernel $K^{w,(n)}_x$ now converges to the $m$-harmonic measure on $\pa
D$ with mean $x$, and this implies (\ref{togast}).

Again, in dimensions $d\geq 2$, the problem of finding a `reasonable'
class $\Wi$ satisfying the assumptions of Theorem~\ref{isoTheorem} is
so far unresolved. The problem with verifying conditions (i)--(iv)
from Definition~\ref{defi} in an explicit set-up is that (i) and (ii)
usually require some smoothness of $w$, while (iv) requires that one
can prove the same smoothness for $F_cw$, which is difficult.

The proofs of Theorems~\ref{WFTheorem} and \ref{isoTheorem} are both
based on invariant harmonics (see (\ref{invhar})). Since diffusion
matrices of catalytic Wright-Fisher diffusions do not in general have
invariant harmonics, in order to prove our main result
(Theorem~\ref{main} below), we will need quite different techniques.

Closely related to this is the fact that in the renormalization
classes from Theorems~\ref{WFTheorem} and \ref{isoTheorem}, the unique
attraction point $w^\ast$ does not depend on the parameter $\ga^\ast$
from (\ref{gamma}). As a result, it turns out that the class $\{\la
w^\ast:\la>0\}$ is a fixed shape. Here, for any prerenormalization
class $\Wi$, a {\em fixed shape} is a subclass $\hat\Wi\sub\Wi$ of the
form $\hat\Wi=\{\la w:\la>0\}$ with $0\neq w\in\Wi$, such that
$F_c(\hat\Wi)\sub\hat\Wi$ for all $c>0$. The next lemma, which will be
proved in Section~\ref{gensub} below, describes how fixed shapes for
renormalization classes on compact sets typically arise.
\bl{\bf(Fixed shapes)}\label{L:shape} Assume that for each
$0<\ga^\ast<\infty$, there is a $0\neq w^\ast=w^\ast_{\ga^\ast}\in\Wi$
such that $s_nF^{(n)}w\asto{n}w^\ast_{\ga^\ast}$ whenever $w\in\Wi$,
$s_n\to\infty$, and $s_{n+1}/s_n\to 1+\ga^\ast$. Then:\med

\noi
{\bf (a)} $w^\ast_{\ga^\ast}$ is the unique solution in $\Wi$ of
equation (\ref{afp})~(i).\med

\noi
{\bf (b)} If $w^\ast=w^\ast_{\ga^\ast}$ does not depend on $\ga^\ast$, then
\be\label{fix}
F_c(\la w^\ast)=(\ffrac{1}{\la}+\ffrac{1}{c})^{-1}w^\ast\qquad(\la,c>0).
\ee
Moreover, $\{\la w^\ast:\la>0\}$ is the unique fixed shape in $\Wi$.\med

\noi
{\bf (c)} If the $w^\ast_{\ga^\ast}$ for different values of $\ga^\ast$
are not constant multiples of each other, then $\Wi$ contains no fixed shapes.
\el
In our main result (Theorem~\ref{main} below), we will describe a renormalization class which we believe contains no fixed shape.

\section{Catalytic Wright-Fisher diffusions}

\subsection{Main result}

Motivated by the previous sections, we will now take the abstract definition of a renormalization class as our starting point, and study iterated renormalization transformations on one such class. Earlier work of this sort has been done in \cite{BCGH95,BCGH97,HS98,Sch98,CDG04}. The subject of our study will be the following renormalization class on $[0,1]^2$.
\bdf\label{Wcatdef}
{\bf (Renormalization class of catalytic Wright-Fisher diffusions)}
We set $\Wi_{\rm cat}:=\{w^{\al,p}:\al>0,\ p\in\Hi\}$, where
\be\label{walp}
w^{\al,p}(x):=\left(\ba{@{}cc@{}}\al x_1(1-x_1)&0\\
0&p(x_1)x_2(1-x_2)\ea\right)\qquad(x=(x_1,x_2)\in[0,1]^2),
\ee
and
\be\label{Hdef}
\Hi:=\{p:p\mbox{ a real function on }[0,1],\ p\geq 0,
\ p\mbox{ Lipschitz continuous}\}.
\ee
Moreover, we put
\be\label{Hlr}
\Hi_{l,r}:=\{p\in\Hi:\; 1_{\{p(0)>0\}}=l,\ 1_{\{p(1)>0\}}=r\}\qquad(l,r=0,1),
\ee
and set $\Wi^{l,r}_{\rm cat}:=\{w^{\al,p}:\al>0,
\ p\in\Hi_{l,r}\}\quad(l,r=0,1)$.
\edf
Solutions $\y=(\y^1,\y^2)$ to the martingale problem for $A^{c,w^{\al,p}}_x$ (recall (\ref{Adef})) can be represented as solutions to the SDE
\be\ba{rr@{\,}c@{\,}l}\label{catsde}
{\rm (i)}&\dis\di\y^1_t&=&\dis c\,(x_1-\y^1_t)\di t
+\sqrt{2\al\y^1_t(1-\y^1_t)}\di B^1_t,\\[5pt]
{\rm (ii)}&\dis\di\y^2_t&=&\dis c\,(x_2-\y^2_t)\di t
+\sqrt{2p(\y^1_t)\y^2_t(1-\y^2_t)}\di B^2_t.
\ec
We call $\y^1$ the Wright-Fisher {\em catalyst} with {\em resampling
rate} $\al$ and $\y^2$ the Wright-Fisher {\em reactant} with {\em
catalyzing function} $p$.

Here is our main result:
\bt{\bf (Main result)}\label{main}\newline
{\bf (a)} The set $\Wi_{\rm cat}$ is a renormalization class on
$[0,1]^2$ and $F_c(\Wi^{l,r}_{\rm cat})\sub\Wi^{l,r}_{\rm cat}$
$(c>0,\ l,r=0,1)$.\vc

\noi
{\bf (b)} Fix (positive) migration constants $(c_k)_{k\geq 0}$ such that
\be\label{sga}
{\rm (i)}\quad s_n\asto{n}\infty\qquad\mbox{and}\qquad{\rm (ii)}
\quad\frac{s_{n+1}}{s_n}\asto{n}1+\ga^\ast
\ee
for some $\ga^\ast\geq 0$. If $w\in\Wi^{l,r}_{\rm cat}$ $(l,r=0,1)$,
then uniformly on $[0,1]^2$,
\be\label{wast}
s_nF^{(n)}w\asto{n}w^\ast,
\ee
where the limit $w^\ast$ is the unique solution in $\Wi^{l,r}_{\rm cat}$
to the equation
\be\ba{rr@{\,}c@{\,}ll}\label{wiga}
{\rm (i)}&\dis(1+\ga^\ast)F_{1/\ga^\ast}w^\ast&=&\dis w^\ast\qquad
&\mbox{if}\ \ \ga^\ast>0,\\[5pt]
{\rm (ii)}&
\dis\ffrac{1}{2}\sum_{i,j=1}^2w^\ast_{ij}(x)\difif{x_i}{x_j}w^\ast(x)+w^\ast(x)
&=&\dis 0\qquad(x\in[0,1]^2)\quad&\mbox{if}\ \ \ga^\ast=0.
\ec
{\bf (c)} The matrix $w^\ast$ is of the form $w^\ast=w^{1,p^\ast}$,
where $p^\ast=p^\ast_{l,r,\ga^\ast}\in\Hi_{l,r}$
depends on $l,r,$ and~$\ga^\ast$. One has
\be
p^\ast_{0,0,\ga^\ast}\equiv 0\quad\mbox{and}
\quad p^\ast_{1,1,\ga^\ast}\equiv 1\qquad\mbox{ for all }\ga^\ast\geq 0.
\ee
For each $\ga^\ast\geq 0$, the function $p^{\ast}_{0,1,\ga^\ast}$ is
concave, nondecreasing, and satisfies $p^\ast_{0,1,\ga^\ast}(0)=0$,
$p^\ast_{0,1,\ga^\ast}(1)=1$. By symmetry, analoguous statements hold
for $p^\ast_{1,0,\ga^\ast}$.
\et
Conditions (\ref{sga})~(i) and (ii) are satisfied, for example, for
$c_k=(1+\ga^\ast)^{-k}$. Note that the functions
$p^\ast_{0,0,\ga^\ast}$ and $p^\ast_{1,1,\ga^\ast}$ are independent of
$\ga^\ast\geq 0$. We believe that on the other hand,
$p^\ast_{0,1,\ga^\ast}$ is not constant as a function of $\ga^\ast$,
but we have not proved this.\footnote{In support of this, if $\Ui_\ga$ $(\ga>0)$ are transformations such that $\ov F_\ga ^{1,p}=w^{1,\Ui_\ga p}$ (see (\ref{hutform}) below), then a heuristic calculation for $p=p^\ast_{0,1,0}$ yields $\Ui_\ga p(x)=p(x)+\ga^2 x(1-x)\big\{\ffrac{1}{2}p''(x)-\ffrac{4}{3}(p'(x))^2-\ffrac{4}{3}xp'''(x)\big\}+O(\ga^3)$, which implies that $p^\ast_{0,1,0}\neq p^\ast_{0,1,\ga^\ast}$ for $\ga^\ast$ small enough.}
If this is confirmed, then by Lemma~\ref{L:shape}, it follows that
$\Wi^{0,1}_{\rm cat}$, unlike all renormalization classes studied
previously, contains no fixed shapes.

\detail{If $p=p^\ast_{0,1,0}$, then a heuristic calculaton yields
\be\label{gaex}
\Ui_\ga p(x)=p(x)+\ga^2 x(1-x)\big\{\ffrac{1}{2}p''(x)-\ffrac{4}{3}(p'(x))^2-\ffrac{4}{3}xp'''(x)\big\}+O(\ga^3)
\ee
as $\ga\to0$. Note that $\ffrac{1}{2}p''(x)-\ffrac{4}{3}(p'(x))^2<0$ for $x\in(0,1)$. By (\ref{ddp2}) and (\ref{ddp4}) below, $p''(0)=-2p'(0)$ and $p''(1)=-2p'(1)$. Hence, by the concavity of $p$, one has $p''(0)<p''(1)$ which shows that $p'''(x)>0$ for some $x\in(0,1)$, and therefore the right-hand side of (\ref{gaex}) is different from $p(x)$ for some $x\in(0,1)$ and $\ga$ small enough. It follows that $p$ is not a fixed point under $\Ui_\ga$.}

The function $p^\ast_{0,1,0}$ is the unique nonnegative
solution to the equation
\be\label{pdi}
\ffrac{1}{2}x(1-x)\diff{x}p(x)+p(x)(1-p(x))=0\qquad(x\in[0,1])
\ee
with boundary conditions $p(0)=0$ and $p(1)>0$. This function occurred
before in the work of Greven, Klenke, and Wakolbinger
\cite[formulas~(1.10)--(1.11)]{GKW01}, who studied linearly
interacting catalytic Wright-Fisher diffusions catalyzed by a voter
model. They believe their results to hold for a Wright-Fisher catalyst
too, i.e., for a model of the form
\bc\label{Z2}
\di\x^1_\xi(t)&=&\dis\sum_{\eta:\,|\eta-\xi|=1}\!
\big(\x^1_\eta(t)-\x^1_\xi(t)\big)\,\di t
+\sqrt{2\al\x^1_\xi(t)(1-\x^1_\xi(t))}\,\di B^1_\xi(t),\\[5pt]
\di\x^2_\xi(t)&=&\dis\sum_{\eta:\,|\eta-\xi|=1}\!
\big(\x^2_\eta(t)-\x^2_\xi(t)\big)\,\di t
+\sqrt{2p(\x^1_\xi(t))\x^2_\xi(t)(1-\x^2_\xi(t))}\,\di B^2_\xi(t),
\ec
where $\al>0$ is a constant, $p$ is a nonnegative function on $[0,1]$
satisfying $p(0)=0$ and $p(1)>0$, but they could not prove this due to
certain technical difficulties that a $[0,1]$-valued catalyst would
create, compared to the simpler $\{0,1\}$-valued voter model. They
determined the clustering distribution of their model on $\Z^2$, which
turns out to coincide with the prediction made based on
renormalization theory in Conjecture~\ref{Z2con}, with
$w^\ast=w^{1,p^\ast_{0,1,0}}$ as in our Theorem~\ref{main}.

The work in \cite{GKW01} not only provides the main motivation for the
present chapter, but also inspired some of our techniques for proving
Theorem~\ref{main}. This concerns in particular the proof of
Proposition~\ref{identprop} below, which makes the connection between
renormalization transformations and a branching process. We hope that
conversely, our techniques may shed some light on the problems left
open by \cite{GKW01}, in particular, the question whether their
results stay true if the voter model catalyst is replaced by a
Wright-Fisher catalyst. It seems plausible that their results may not
hold for the model in (\ref{Z2}) if the catalyzing function $p$ grows
too fast at $0$. On the other hand, our proofs suggest that $p$ with a
finite slope at $0$ should be OK. (In particular, while deriving
formula (\ref{intext3}) below, we use that $p$ can be bounded from above by
$r_+h_{0,1}$ for some $r_+>0$, which requires that $p$ has a finite
slope at $0$.)

\subsection{Open problems}\label{probsec}

The general program of studying renormalization classes in the sense
of Definition~\ref{defi} contains a wealth of open problems. In our
proofs, we make heavy use of the single-way nature of the catalyzation
in (\ref{catsde}), in particular, the fact that $\y^1$ is an
autonomous process which allows one to condition on $\y^1$ and
consider $\y^2$ as a process in a random environment created by
$\y^1$.  As soon as one leaves the single-way catalytic regime one
runs into several difficulties, both technically (it is hard to prove
that a given class of matrices is a renormalization class in the sense
of Definition~\ref{defi}) and conceptually (it is not clear when
solutions to the asymptotic fixed shape equation (\ref{afp})~(ii) are
unique). Therefore, it seems at present hard to verify the complete
picture for renormalization classes on the unit square that arises
from the numerical simulations described in Section~\ref{numer} and
Figures~\ref{fixfig} and~\ref{graph}, unless one or more essential new
ideas are added.

In this context, the study of the nonlinear partial differential
equation (\ref{flow}) and its fixed points seems to be a challenging
problem. This may be a hard problem from an analytic point of view,
since the equation is degenerate and not in divergence form. For the
renormalization class $\Wi_{\rm cat}$, the quasilinear equation
(\ref{flow}) reduces to the semilinear equation (\ref{cau}), which is
analytically easier to treat and moreover has a probabilistic
interpretation in terms of a superprocess. We do not
know whether solutions to equation (\ref{flow}) can in general be
represented in terms of a stochastic process of some sort.

Even for the renormalization class $\Wi_{\rm cat}$, several
interesting problems are left open. One of the most urgent ones is to
prove that the functions $p^\ast_{0,1,\ga^\ast}$ are not constant in
$\ga^\ast$, and therefore, by Lemma~\ref{L:shape}~(c), $\Wi^{0,1}_{\rm
cat}$ contains no fixed shapes. Moreover, we have not investigated the
iterated renormalization transformations in the regime
$\ga^\ast=\infty$.  Also, we believe that the convergence in
(\ref{Ucon})~(ii) does not hold if the condition that $p$ is Lipschitz
is dropped, in particular, if $p(0)=0$ and $p$ has an infinite slope
at $0$. For $p\in\Hi_{0,0}$, it seems
plausible that a properly rescaled version of the iterates
$\Ui^{(n)}p$, with $\Ui_\ga$ as in (\ref{Uiga}) below, converges to a
universal limit, but we have not investigated this either. Finally, we
have not investigated the convergence of the iterated kernels
$K^{w,(n)}$ from (\ref{Kdef}) (in particular, we have not verified
Conjecture~\ref{Kcon}) for the renormalization class $\Wi_{\rm cat}$.

Our methods, combined with those in \cite{BCGH95}, can probably be
extended to study the action of iterated renormalization
transformations on diffusion matrices of the following more general
form (compared to (\ref{walp})):
\be\label{walp2}
w(x)=\left(\ba{@{}cc@{}}g(x_1)&0\\0&p(x_1)x_2(1-x_2)\ea\right)
\qquad(x=\in[0,1]^2),
\ee
where $g:[0,1]\to\R$ is Lipschitz, $g(0)=g(1)=0$, $g>0$ on $(0,1)$,
and $p\in\Hi$ as before. This would, however, require a lot of extra
technical work and probably not generate much new insight. The
numerical simulations mentioned in Section~\ref{numer} suggest that
many diffusion matrices of an even more general form than
(\ref{walp2}) also converge under renormalization to the limit points
$w^\ast$ from Theorem~\ref{main}, but we don't know how to prove this.

In the next sections, we will show that for the renormalization class
$\Wi_{\rm cat}$, the rescaled renormalization transformations $\ov
F_\ga$ from (\ref{hatFc}) can be expressed in terms of the log-Laplace
operators of a discrete time branching process on $[0,1]$. This will
allow us to use techniques from the theory of spatial branching
processes to verify Conjecture~\ref{con} for the renormalization class
$\Wi_{\rm cat}$ in the case $\ga^\ast<\infty$.

\subsection{Poisson-cluster branching processes}

We first need some concepts and facts from branching theory. Finite
measure-valued branching processes (on $\R$) in discrete time have
been introduced by Ji\v rina \cite{Jir64}. We need to consider only a
special class. 

Let $E$ be a separable, locally compact, and metrizable space. We let
$\Ci(E)$ and $B(E)$ denote the spaces of all continuous, and bounded
Borel measurable, real functions on $E$, respectively. We put
$\Ci_+(E):=\{f\in\Ci(E):f\geq 0\}$ and define $B_+(E)$ analogously. We
let $\Mi(E)$ denote the space of all finite measures on $E$, equipped
with the topology of weak convergence. The subspace of probability
measures is denoted by $\Mi_1(E)$. For $\mu\in\Mi(E)$ and $f\in B(E)$
we use the notation $\li\mu,f\re:=\int_E f\,\di\mu$ and
$|\mu|:=\mu(E)$.

We call a continuous map $\Qi$ from $E$ into $\Mi_1(\Mi(E))$ a
{\em continuous cluster mechanism}. By definition, an $\Mi(E)$-valued
random variable $\Xc$ is a {\em Poisson cluster measure} on $E$ with
locally finite {\em intensity measure} $\mu$ and continuous cluster
mechanism $\Qi$, if its log-Laplace transform satisfies
\be\label{randclust}
-\log E\big[\ex{-\li\Xc,f\re}\big]=\int_E\!\mu(\di x)\Big(1-\int_{\Mi(E)}\!\!\Qi(x,\di\chi)\ex{-\li\chi,f\re}\Big)\quad(f\in B_+(E)).
\ee
For given $\mu$ and $\Qi$, such a Poisson cluster measure exists, and
is unique in distribution, provided that the right-hand side of
(\ref{randclust}) is finite for $f=1$. It may be constructed as
$\Xc=\sum_i\chi_{x_i}$, where $\sum_i\de_{x_i}$ is a (possibly
infinite) Poisson point measure with intensity $\mu$, and given
$x_1,x_2,\ldots$, the $\chi_{x_1},\chi_{x_2},\ldots$ are independent
random variables with laws
$\Qi(x_1,\,\cdot\,),\Qi(x_2,\,\cdot\,),\ldots$, respectively.

Now fix a finite sequence of functions $q_k\in\Ci_+(E)$ and continuous
cluster mechanisms $\Qi_k$ ($k=1,\ldots,n$), define
\be\label{Vk}
\Ui_kf(x):=q_k(x)\Big(1-\int_{\Mi(E)}\!\!\Qi_k(x,\di\chi)\ex{-\li\chi,f\re}\Big)\qquad(x\in E,\ f\in B_+(E),\ k=1,\ldots,n),
\ee
and assume that
\be\label{fincon}
\sup_{x\in E}\Ui_k1(x)<\infty\qquad(k=1,\ldots,n).
\ee
Then $\Ui_k$ maps $B_+(E)$ into $B_+(E)$ for each $k$, and for each
$\Mi(E)$-valued initial state $\Xc_0$, there exists a
(time-inhomogeneous) Markov chain $(\Xc_0,\ldots,\Xc_n)$ in $\Mi(E)$,
such that $\Xc_k$, given $\Xc_{k-1}$, is a Poisson cluster measure
with intensity $q_k\Xc_{k-1}$ and cluster mechanism $\Qi_k$.
It is not hard to see that the process started in $\mu$ satisfies
\be\label{Vi}
E^\mu\big[\ex{-\li\Xc_n,f\re}\big]
=\ex{-\li\mu,\Ui_1\circ\cdots\circ\Ui_n f\re}
\qquad(\mu\in\Mi(E),\ f\in B_+(E)).
\ee
We call $\Xc=(\Xc_0,\ldots,\Xc_n)$ the {\em Poisson-cluster branching
process} on $E$ with {\em weight functions} $q_1,\ldots,q_n$ and
cluster mechanisms $\Qi_1,\ldots,\Qi_n$. The operator $\Ui_k$ is
called the {\em log-Laplace operator} of the transition law from
$\Xc_{k-1}$ to $\Xc_k$. Note that we can write (\ref{Vi}) in the
suggestive form
\be\label{Vi2}
P^\mu\big[\Pois(f\Xc_n)=0\big]
=P\big[\Pois\big((\Ui_1\circ\cdots\circ\Ui_n f)\mu\big)=0\big].
\ee
Here, if $\mu$ is an $\Mi(E)$-valued random variable, then
$\Pois(\mu)$ denotes an $\Ni(E)$-valued random variable such that
conditioned on $\mu$, $\Pois(\mu)$ is a Poisson point measure with
intensity $\mu$.

\subsection{The renormalization branching process}\label{RBP}

We will now construct a Poisson-cluster branching process on $[0,1]$
of a special kind, and show that the rescaled renormalization
transformations on $\Wi_{\rm cat}$ can be expressed in terms of the
log-Laplace operators of this branching process.

By Lemma~\ref{ergo} below, for each $\ga>0$ and $x\in[0,1]$, the SDE
\be\label{Yclx}
\di\y(t)=\ffrac{1}{\ga}\,(x-\y(t))\di t+\sqrt{2\y(t)(1-\y(t))}\di B(t),
\ee
has a unique (in law) stationary solution. We denote this solution by
$(\y^\ga_x(t))_{t\in\R}$. Let $\tau_\ga$ be an independent exponentially
distributed random variable with mean $\ga$, and set
\be\label{Zidef}
\Zi^\ga_x:=\int_0^{\tau_\ga}\de_{\y^\ga_x(-t/2)}\di t\qquad(\ga>0,\ x\in[0,1]).
\ee
Define constants $q_\ga$ and continuous (by Corollary~\ref{Ugacont} below)
cluster mechanisms $\Qi_\ga$ by
\be\label{Qqdef}
q_\ga:=\ffrac{1}{\ga}+1\qquad\mbox{and}\qquad\Qi_\ga(x,\,\cdot\,)
:=\Li(\Zi^\ga_x)\qquad(\ga>0,\ x\in[0,1]),
\ee
and let $\Ui_\ga$ denote the log-Laplace operator with (constant)
weight function $q_\ga$ and cluster mechanism $\Qi_\ga$, i.e.,
\be\label{Uiga}
\Ui_\ga f(x):=q_\ga\Big(1-\int_{\Mi([0,1])}\!\!\Qi_\ga(x,\di\chi)
\ex{-\li\chi,f\re}\Big)\qquad(x\in[0,1],\ f\in B_+[0,1],\ \ga>0).
\ee
We now establish the connection between renormalization transformations
on $\Wi_{\rm cat}$ and log-Laplace operators.
\bp{\bf (Identification of the renormalization transformation)}\label{identprop}
Let $\hut F_\ga$ be the rescaled renormalization transformation on
$\Wi_{\rm cat}$ defined in (\ref{hatFc}). Then
\be\label{hutform}
\hut F_\ga w^{\txt 1,p}=w^{\txt 1,\Ui_\ga p}\qquad(p\in\Hi,\ \ga>0).
\ee
\ep
Fix a diffusion matrix $w^{\al,p}\in\Wi_{\rm cat}$ and migration constants
$(c_k)_{k\geq 0}$. Define constants $\hut s_n$ and $\ga_n$ as in (\ref{hutsn})
and (\ref{gan}), respectively, where $\bet:=1/\al$. Then
Proposition~\ref{identprop} and formula (\ref{sFw}) show that
\be\label{renbra}
\hut s_nF^{(n)}w^{\al,p}=w^{\txt 1,\Ui_{\ga_{n-1}}\circ\cdots
\circ\Ui_{\ga_0}(\frac{p}{\al})}.
\ee
Here $\Ui_{\ga_{n-1}},\ldots,\Ui_{\ga_0}$ are the log-Laplace
operators of the Poisson-cluster branching process
$\Xc=(\Xc_{-n},\ldots,\Xc_0)$ with weight functions
$q_{\ga_{n-1}},\ldots,q_{\ga_0}$ and cluster mechanisms
$\Qi_{\ga_{n-1}},\ldots,\Qi_{\ga_0}$. We call $\Xc$ (started at some
time $-n$ in an initial law $\Li(\Xc_{-n})$) the {\em renormalization
  branching process}. By formulas (\ref{Vi}) and (\ref{renbra}), the
study of the limiting behavior of rescaled iterated renormalization
transformations on $\Wi_{\rm cat}$ reduces to the study of the
renormalization branching process $\Xc$ in the limit $n\to\infty$.

\subsection{Convergence to a time-homogeneous process}

Let $\Xc=(\Xc_{-n},\ldots,\Xc_0)$ be the renormalization branching
process introduced in the last section. If the constants
$(\ga_k)_{k\geq 0}$ satisfy $\sum_n\ga_n=\infty$ and
$\ga_n\to\ga^\ast$ for some $\ga^\ast\in\half$, then $\Xc$ is almost
time-homo\-geneous for large $n$. More precisely, we will prove the
following convergence result.
\bt{\bf(Convergence to a time-homogenous branching process)}\label{supercon}
Assume that $\Li(\Xc_{-n})\Asto{n}\mu$ for some probability law $\mu$
on $\Mi([0,1])$.\med

\noi
{\bf (a)} If $0<\ga^\ast<\infty$, then
\be\label{disco}
\Li(\Xc_{-n},\Xc_{-n+1},\ldots)\Asto{n}
\Li(\Yi^{\ga^\ast}_0,\Yi^{\ga^\ast}_1,\ldots),
\ee
where $\Yi^{\ga^\ast}$ is the time-homogenous branching process with
log-Laplace operator $\Ui_{\ga^\ast}$ in each step and initial law
$\Li(\Yi^{\ga^\ast}_0)=\mu$.\med

\noi
{\bf (b)} If $\ga^\ast=0$, then
\be\label{coco}
\Li\Big(\big(\Xc_{-k_n(t)}\big)_{t\geq 0}\Big)\Asto{n}
\Li\Big(\big(\Yi^0_t\big)_{t\geq 0}\Big),
\ee
where $\Rightarrow$ denotes weak convergence of laws on path space,
$k_n(t):=\min\{k:0\leq k\leq n$, $\sum_{l=k}^{n-1}\ga_l\leq t\}$, and
$\Yi^0$ is the superprocess on $[0,1]$ with underlying motion
generator $\ffrac{1}{2}x(1-x)\diff{x}$ and activity and growth
parameter both identically $1$, started in the initial law
$\Li(\Yi^0_0)=\mu$.
\et
We call the superprocess $\Yi^0$ from part~(b) the {\em
super-Wright-Fisher diffusion}. It is the time-homo\-ge\-neous Markov
process in $\Mi[0,1]$ with continuous sample paths, whose Laplace
functionals are given by
\be
E^\mu\big[\ex{-\li\Yi^0_t,f\re}\big]=\ex{-\li\mu,\Ui^0_tf\re}
\qquad(\mu\in\Mi[0,1],\ f\in B_+[0,1],\ t\geq 0),
\ee
where $\Ui^0_tf=u_t$ is the unique mild solution of the semilinear
Cauchy equation
\be
\left\{\ba{r@{\,}c@{\,}l}\label{cau}
\dif{t}u_t(x)&=&\ffrac{1}{2}x(1-x)\diff{x}u_t(x)+u_t(x)(1-u_t(x))
\quad(t\geq 0,\ x\in[0,1]),\\
u_0&=&f.\ea\right.
\ee
For a further study of the renormalization branching process $\Xc$
and its limiting processes $\Yi^{\ga^\ast}$ ($\ga^\ast\geq 0$) we
will use the technique of embedded particle systems, which we explain
in the next section.

\subsection{Weighted and Poissonized branching processes}

In this section, we explain how from a Poisson-cluster branching
process it is possible to construct other branching processes by
weighting and Poissonization. We first need to introduce spatial
branching particle systems in some generality.

Let $E$ again be separable, locally compact, and metrizable.  We set
$\Ci_{[0,1]}(E):=\{f\in\Ci(E):0\leq f\leq 1\}$ and define
$B_{[0,1]}(E)$ analogously. We write $\Ni(E)$ for the space of finite
counting measures, i.e., measures of the form
$\nu=\sum_{i=1}^m\de_{x_i}$ with $x_1,\ldots,x_m\in E$ ($m\geq 0$). We
interpret $\nu$ as a collection of particles, situated at positions
$x_1,\ldots,x_m$. For
$\nu\in\Ni(E)$ and $f\in B_{[0,1]}(E)$, we adopt the notation
\be
f^{\txt 0}:=1\quad\mbox{and}\quad f^{\txt\nu}:=\prod_{i=1}^mf(x_i)
\quad\mbox{when}\quad\nu=\sum_{i=1}^m\de_{x_i}\quad(m\geq 1).
\ee
We call a continuous map $x\mapsto Q(x,\,\cdot\,)$ from $E$ into
$\Mi_1(\Ni(E))$ a {\em continuous offspring mechanism.}

Fix continuous offspring mechanisms $Q_k$ ($1\leq k\leq n$), and
let $(X_0,\ldots,X_n)$ be a Markov chain in $\Ni(E)$ such that,
given that $X_{k-1}=\sum_{i=1}^m\de_{x_i}$, the next step of the
chain $X_k$ is a sum of independent random variables with laws
$Q_k(x_i,\,\cdot\,)$ ($i=1,\ldots,m$). Then
\be\label{V}
E^\nu\big[(1-f)^{\txt X_n}\big]=(1-U_1\circ\cdots\circ U_nf)^{\txt\nu}
\qquad(\nu\in\Ni(E),\ f\in B_{[0,1]}(E)),
\ee
where $U_k:B_{[0,1]}(E)\to B_{[0,1]}(E)$ is defined as
\be\label{VV}
U_kf(x):=1-\int_{\Ni(E)}\!\!Q^k(x,\di\nu)(1-f)^{\txt\nu}
\qquad(1\leq k\leq n,\ x\in E,\ f\in B_{[0,1]}(E)).
\ee
We call $U_k$ the {\em generating operator} of the transition law from
$X_{k-1}$ to $X_k$, and we call $X=(X_0,\ldots,X_n)$ the {\em
branching particle system} on $E$ with generating operators
$U_1,\ldots,U_n$. It is often useful to write (\ref{V}) in the
suggestive form
\be\label{V2}
P^\nu\big[\Thin_f(X_n)=0\big]=P\big[\Thin_{U_1\circ\cdots\circ U_n f}(\nu)
=0\big]\qquad(\nu\in\Ni(E),\ f\in B_{[0,1]}(E)).
\ee
Here, if $\nu$ is an $\Ni(E)$-valued random variable and $f\in
B_{[0,1]}(E)$, then $\Thin_f(\nu)$ denotes an $\Ni(E)$-valued random
variable such that conditioned on $\nu$, $\Thin_f(\nu)$ is obtained
from $\nu$ by independently throwing away particles from $\nu$, where
a particle at $x$ is kept with probability $f(x)$. One has the
elementary relations
\be\label{elrel}
\Thin_f(\Thin_g(\nu))\isd\Thin_{fg}(\nu)\quad\mbox{and}
\quad\Thin_f(\Pois(\mu))\isd\Pois(f\mu),
\ee
where $\isd$ denotes equality in distribution. 

We are now ready to describe weighted and Poissonized branching
processes. Let $\Xc=(\Xc_0,\ldots,\Xc_n)$ be a Poisson-cluster
branching process on $E$, with continuous weight functions
$q_1,\ldots,q_n$, continuous cluster mechanisms $\Qi_1,\ldots,\Qi_n$,
and log-Laplace operators $\Ui_1,\ldots,\Ui_n$ given by (\ref{Vk}) and
satisfying (\ref{fincon}). Let $\Zi^k_x$ denote an $\Mi(E)$-valued
random variable with law $\Qi_k(x,\,\cdot\,)$. Let $h\in\Ci_+(E)$ be
bounded, $h\neq 0$, and put $E^h:=\{x\in E:h(x)>0\}$. For $f\in
B_+(E^h)$, define $hf\in B_+(E)$ by $hf(x):=h(x)f(x)$ if $x\in E^h$
and $hf(x):=0$ otherwise.
\bp{\hspace{-1.6pt}\bf(Weighting of Poisson-cluster branching processes)\hspace{-1pt}}\label{weightprop}
Assume that there exists a constant $K<\infty$ such that $\Ui_kh\leq
Kh$ for all $k=1,\ldots,n$. Then there exists a Poisson-cluster
branching process $\Xc^h=(\Xc^h_0,\ldots,\Xc^h_n)$ on $E^h$ with
weight functions $(q_1^h,\ldots,q^h_n)$ given by $q^h_k:=q_k/h$,
continuous cluster mechanisms $\Qi^h_1,\ldots,\Qi^h_n$ given by
\be\label{Qih}
\Qi^h_k(x,\,\cdot\,):=\Li(h\Zi^k_x)\qquad(x\in E^h),
\ee
and log-Laplace operators $\Ui^h_1,\ldots,\Ui^h_n$ satisfying
\be\label{htrafo1}
h\,\Ui^h_kf:=\Ui_k(hf)\qquad(f\in B_+(E^h)).
\ee
The processes $\Xc$ and $X^h$ are related by
\be\label{weighting}
\Li(\Xc^h_0)=\Li(h\Xc_0)\quad\mbox{implies}\quad\Li(\Xc^h_k)
=\Li(h\Xc_k)\qquad(0\leq k\leq n).
\ee
\ep
\bp{\bf(Poissonization of Poisson-cluster branching processes)}\label{Poisprop}
Assume that $\Ui_kh\leq h$ for all $k=1,\ldots,n$. Then there exists
a branching particle system $X^h=(X^h_0,\ldots,X^h_n)$ on $E^h$ with
continuous offspring mechanisms $Q^h_1,\ldots,Q^h_n$ given by
\be\label{Zxdef}
Q^h_k(x,\,\cdot\,):=\frac{q_k(x)}{h(x)}P\big[\Pois(h\Zi^k_x)\in\cdot\,\big]
+\Big(1-\frac{q_k(x)}{h(x)}\Big)\de_0(\,\cdot\,)
\qquad(x\in E^h),
\ee
and generating operators $U^h_1,\ldots,U^h_n$ satisfying
\be\label{htrafo}
hU^h_kf:=\Ui_k(hf)\qquad(f\in B_{[0,1]}(E^h)).
\ee
The processes $\Xc$ and  $X^h$ are related by
\be\label{Poissonization}
\Li(X^h_0)=\Li(\Pois(h\Xc_0))\quad\mbox{implies}\quad\Li(X^h_k)
=\Li(\Pois(h\Xc_k))\qquad(0\leq k\leq n).
\ee
\ep
Here, the right-hand side of (\ref{Zxdef}) is always a probability
measure, despite that it may happen that $q_k(x)/h(x)>1$. The
(straightforward) proofs of Propositions~\ref{weightprop} and
\ref{Poisprop} can be found in Section~\ref{Poissec} below. If
(\ref{weighting}) holds then we say that $\Xc^h$ is obtained from
$\Xc$ by {\em weighting} with density $h$. If (\ref{Poissonization})
holds then we say that $X^h$ is obtained from $\Xc$ by {\em
Poissonization} with density $h$. Proposition~\ref{Poisprop} says
that a Poisson-cluster branching process $\Xc$ contains, in a way,
certain `embedded' branching particle systems $X^h$. Poissonization
relations for superprocesses and embedded particle systems have
enjoyed considerable attention, see \cite{FStrim} and references
therein.

A function $h\in B_+(E)$ such that $\Ui_kh\leq h$ is called {\em
$\Ui_k$-super\-harmonic}. If the reverse inequality holds we say
that $h$ is {\em $\Ui_k$-sub\-harmonic}. If $\Ui_kh=h$ then $h$ is
called {\em $\Ui_k$-harmonic}.

\subsection{Extinction versus unbounded growth for embedded particle
systems}\label{exex}

In this section we explain how embedded particle systems can be used
to prove Theorem~\ref{main}. Throughout this section $(\ga_k)_{k\geq 0}$
are positive constants such that $\sum_n\ga_n=\infty$ and
$\ga_n\to\ga^\ast$ for some $\ga^\ast\in\half$, and
$\Xc=(\Xc_{-n},\ldots,\Xc_0)$ is the renormalization branching process
on $[0,1]$ defined in Section~\ref{RBP}. We write
\be
\Ui^{(n)}:=\Ui_{\ga_{n-1}}\circ\cdots\circ\Ui_{\ga_0}.
\ee
In view of formula (\ref{renbra}), in order to prove Theorem~\ref{main},
we need the following result.
\bp{\bf (Limits of iterated log-Laplace operators)}\label{Uit}
Uniformly on $[0,1]$,
\be\ba{rr@{\,}c@{\,}ll}\label{Ucon}
{\rm (i)}&\dis\lim_{n\to\infty}\Ui^{(n)}p&=&1&
\dis\qquad(p\in\Hi_{1,1}),\\[5pt]
{\rm (ii)}&\dis\lim_{n\to\infty}\Ui^{(n)}p&=&0&
\dis\qquad(p\in\Hi_{0,0}),\\[5pt]
{\rm (iii)}&\dis\lim_{n\to\infty}\Ui^{(n)}p&=&p^{\ast}_{0,1,\ga^\ast}&
\dis\qquad(p\in\Hi_{0,1}),
\ec
where $p^{\ast}_{0,1,\ga^\ast}:[0,1]\to[0,1]$ is a function depending
on $\ga^\ast$ but not on $p\in\Hi_{0,1}$.
\ep
In our proof of Proposition~\ref{Uit}, we will use embedded particle
systems $X^h=(X^h_{-n},\ldots,X^h_0)$ obtained from $\Xc$ by Poissonization
with certain $h$ taken from the classes $\Hi_{1,1}$, $\Hi_{0,0}$, and
$\Hi_{0,1}$. Below, $P^{-n,\de_x}$ denotes the law of the process started
at time $-n$ with one particle at $x$.
\bl{\bf(Embedded particle system with $h_{1,1}$)}\label{11lem}
The constant function $h_{1,1}(x):=1$ is $\Ui_\ga$-harmonic for each
$\ga>0$. The corresponding embedded particle system $X^{h_{1,1}}$ on
$[0,1]$ satisfies
\be\label{expl}
P^{-n,\de_x}\big[|X^{h_{1,1}}_0|\in\cdot\,\big]\Asto{n}\de_\infty
\ee
uniformly\footnote{Since $\Mi_1[0,\infty]$ is compact in the topology
of weak convergence, there is a unique uniform structure compatible
with the topology, and therefore we can unambiguously talk about uniform
convergence of $\Mi_1[0,\infty]$-valued functions (in this case,
$x\mapsto P^{-n,\de_x}\big[|X^{h_{1,1}}_0|\in\cdot\,\big]$).} for
all $x\in[0,1]$.
\el
In (\ref{expl}) and similar formulas below, $\Rightarrow$ denotes
weak convergence of probability measures on $[0,\infty]$. Thus,
(\ref{expl}) says that for processes started with one particle
on the position $x$ at times $-n$, the number of particles at
time zero converges to infinity as $n\to\infty$.
\bl{\bf(Embedded particle system with $h_{0,0}$)}\label{00lem}
The function $h_{0,0}(x):=x(1-x)$ $(x\in[0,1])$ is
$\Ui_\ga$-super\-harmonic for each $\ga>0$. The corresponding
embedded particle system $X^{h_{0,0}}$ on $(0,1)$ is critical
and satisfies
\be\label{ext}
P^{-n,\de_x}\big[|X^{h_{0,0}}_0|\in\cdot\,\big]\Asto{n}\de_0
\ee
locally uniformly for all $x\in(0,1)$.
\el
Here, we say that a branching particle system $X$ is {\em critical} if each
particle produces on average one offspring (in each time step and
independent of its position). Formula (\ref{ext}) says that the
embedded particle system $X^{h_{0,0}}$ gets extinct during the time
interval $\{-n,\ldots,0\}$ with probability tending to one as
$n\to\infty$. We can summarize Lemmas~\ref{11lem} and \ref{00lem} by
saying that the embedded particle system associated with $h_{1,1}$
grows unboundedly while the embedded particle system associated with
$h_{0,0}$ becomes extinct as $n\to\infty$.

We will also consider an embedded particle system $X^{h_{0,1}}$ for a
certain $h_{0,1}$ taken from $\Hi_{0,1}$. It turns out that this
system either gets extinct or grows unboundedly, each with a positive
probability. In order to determine these probabilities, we need to
consider embedded particle systems for the time-homogeneous processes
$\Yi^{\ga^\ast}$ ($\ga^\ast\in\half$) from (\ref{disco}) and
(\ref{coco}). If $h\in\Hi_{0,1}$ is $\Ui_{\ga^\ast}$-superharmonic for
some $\ga^\ast>0$, then Poissonizing the process $\Yi^{\ga^\ast}$ with
$h$ yields a branching particle system on $(0,1]$ which we denote by
$Y^{\ga^\ast,h}=(Y^{\ga^\ast,h}_0,Y^{\ga^\ast,h}_1,\ldots)$.
Likewise, if $h\in\Hi_{0,1}$ is twice continuously differentiable and
satisfies
\be\label{contsup}
\ffrac{1}{2}x(1-x)\diff{x}h(x)-h(x)(1-h(x))\leq 0,
\ee
then Poissonizing the super-Wright-Fisher diffusion $\Yi^0$ with $h$
yields a continuous-time branching particle system on $(0,1]$, which
we denote by $Y^{0,h}=(Y^{0,h}_t)_{t\geq 0}$. For example, for $m\geq 4$,
the function $h(x):=1-(1-x)^m$ satisfies (\ref{contsup}).
\bl{\bf(Embedded particle system with $h_{0,1}$)}\label{01lem}
The function $h_{0,1}(x):=1-(1-x)^7$ is $\Ui_\ga$-super\-harmonic for
each $\ga>0$. The corresponding embedded particle system $X^{h_{0,1}}$
on $(0,1]$ satisfies
\be\label{expext}
P^{-n,\de_x}\big[|X^{h_{0,1}}_0|\in\cdot\,\big]\Asto{n}
\rho_{\ga^\ast}(x)\de_\infty+(1-\rho_{\ga^\ast}(x))\de_0,
\ee
locally uniformly for all $x\in(0,1]$, where
\be\label{rhodef}
\rho_{\ga^\ast}(x):=\left\{\ba{ll}\dis P^{\de_x}[Y^{\ga^\ast,h_{0,1}}_k\neq 0
\ \forall k\geq 0]\quad&(0<\ga^\ast<\infty),\\[5pt]
\dis P^{\de_x}[Y^{0,h_{0,1}}_t\neq 0\ \forall t\geq 0]\quad
&(\ga^\ast=0).\ea\right.
\ee
\el
We now explain how Lemmas~\ref{11lem}--\ref{01lem} imply
Proposition~\ref{Uit}. In doing so, it will be more convenient to work
with weighted branching processes than with Poissonized branching
processes. A little argument (which can be found in Lemma~\ref{exgro}
below) shows that Lemmas~\ref{11lem}--\ref{01lem} are equivalent to
the next proposition.
\bp{\bf(Extinction versus unbounded growth)}\label{P:exgr}
Let $h_{1,1}$, $h_{0,0}$, and $h_{0,1}$ be as in
Lemmas~\ref{11lem}--\ref{01lem}. For $\ga^\ast\in\half$,
put $p^\ast_{1,1,\ga^\ast}(x):=1$, $p^\ast_{0,0,\ga^\ast}(x):=0$
$(x\in[0,1])$, and
\be\label{pdef}
p^\ast_{0,1,\ga^\ast}(0):=0\quad\mbox{and}
\quad p^\ast_{0,1,\ga^\ast}(x):=h_{0,1}(x)\rho_{\ga^\ast}(x)\qquad(x\in(0,1]),
\ee
with $\rho_{\ga^\ast}$ as in (\ref{rhodef}). Then, for $(l,r)=(1,1),(0,0)$,
and $(0,1)$,
\be\label{intexp}
P^{-n,\de_x}\big[\li\Xc_0,h_{l,r}\re\in\cdot\,\big]\Asto{n}
\ex{-p^\ast_{l,r,\ga^\ast}(x)}\de_0
+\big(1-\ex{-p^\ast_{l,r,\ga^\ast}(x)}\big)\de_\infty,
\ee
uniformly for all $x\in[0,1]$.
\ep
Formula (\ref{intexp}) says that the weighted branching process
$\Xc^{h_{l,r}}$ exhibits a form of extinction versus unbounded
growth. More precisely, for large $n$ the total mass of
$h_{l,r}\Xc_0$ is close to $0$ or $\infty$ with high probability.\vc

\noi
{\bf Proof of Proposition~\ref{Uit}} By (\ref{Vi}),
\be\label{Vi3}
\Ui^{(n)}p(x)=-\log E^{-n,\de_x}\big[\ex{-\li\Xc_0,p\re}\big]
\qquad(p\in B_+[0,1],\ x\in[0,1]).
\ee
We first prove formula (\ref{Ucon})~(ii). For $(l,r)=(0,0)$,
formula (\ref{intexp}) says that
\be\label{intext}
P^{-n,\de_x}[\li\Xc_0,h_{0,0}\re\in\cdot\,]\Asto{n}\de_0
\ee
uniformly for all $x\in[0,1]$. If $p\in\Hi_{0,0}$, then we can
find $r>0$ such that $p\leq rh_{0,0}$. Therefore, (\ref{intext})
implies that for any $p\in\Hi_{0,0}$,
\be\label{intext2}
P^{-n,\de_x}[\li\Xc_0,p\re\in\cdot\,]\Asto{n}\de_0.
\ee
By (\ref{Vi3}) it follows that
\be\label{innie}
\Ui^{(n)}p(x)=-\log E^{-n,\de_x}\big[\ex{-\li\Xc_0,p\re}\big]\asto{n}0,
\ee
where the limits in (\ref{intext2}) and (\ref{innie}) are uniform in
$x\in[0,1]$. This proves formula (\ref{Ucon})~(ii). To prove formula
(\ref{Ucon})~(iii), note that for any $p\in\Hi_{0,1}$ we can choose
$0<r_-<r_+$ such that $r_-h_{0,1}\leq p+h_{0,0}\leq r_+h_{0,1}$.
Therefore, (\ref{intexp}) implies that
\be\label{intext3}
P^{-n,\de_x}[\li\Xc_0,p\re+\li\Xc_0,h_{0,0}\re\in\cdot\,]\Asto{n}
\ex{-p^\ast_{0,1,\ga^\ast}(x)}\de_0
+\big(1-\ex{-p^\ast_{0,1,\ga^\ast}(x)}\big)\de_\infty.
\ee
Using moreover (\ref{intext}), we see that
\be\label{intext4}
P^{-n,\de_x}[\li\Xc_0,p\re\in\cdot\,]\Asto{n}
\ex{-p^\ast_{0,1,\ga^\ast}(x)}\de_0
+\big(1-\ex{-p^\ast_{0,1,\ga^\ast}(x)}\big)\de_\infty.
\ee
By (\ref{Vi3}), it follows that
\be
\Ui^{(n)}p(x)=-\log E^{-n,\de_x}\big[\ex{-\li\Xc_0,p\re}\big]\asto{n}
p^\ast_{0,1,\ga^\ast}(x)
\ee
where all limits are uniform in $x\in[0,1]$. This proves
(\ref{Ucon})~(iii). The proof of (\ref{Ucon})~(i) is similar but easier.\qed

\subsection{Outline} 

In Section~\ref{Wcatsec}, we verify that
$\Wi_{\rm cat}$ is a renormalization class, we prove
Proposition~\ref{identprop}, which connects the renormalization
transformations $F_c$ to the log-Laplace operators $\Ui_\ga$, and we
collect a number of technical properties of the operators $\Ui_\ga$
that will be needed later on. In Section~\ref{convsec} we prove
Theorem~\ref{supercon} about the convergence of the renormalization
branching process to a time-homogeneous limit.

Sections~\ref{S:WFi}--\ref{S:WFiii} are devoted to the super-Wright-Fisher diffusio $\Yi^0$, i.e., the limiting process from Theorem~\ref{supercon}~(b). These sections have been written in such a way that they can be read independently of the rest of this chapter. In fact, we generalize a bit by allowing for an arbitrary positive constant to appear in front of the $u(1-u)$ term in (\ref{cau}). This generatization reveals that the case where this constant is one is in fact a critical case, marking the boundary between two types of long-time behavior. Section~\ref{S:WFi} gives an introduction to the super-Wright-Fisher diffusion, while Sections~\ref{S:WFii}--\ref{S:WFiii} contain proofs. The central tool in these proofs is a weighted superprocess, rather than embedded particle systems which are our main tool for studying the renormalization branching process $\Xc$

In Section~\ref{embsec}, we take up the study of $\Xc$ and its
embedded particle systems. In particular, we prove the statements from
Section~\ref{exex} about extinction versus unbounded growth of
embedded particle systems, with the exception of Lemma~\ref{00lem},
which is proved in Section~\ref{00sec}. In Section~\ref{final},
finally, we combine all results derived by that point to prove our
main theorem.\med

\noi
{\bf\large Acknowledgements} Work sponsored by the DFG. The authors thank Janos Engl\"ander for answering our questions about his work and Jan Seidler for answering questions about the strong Feller property. Achim Klenke, Dmitry Turaev, and Anita Winter are thanked for useful discussions and comments. We than an anonymous referee for comments which lead to an improved exposition. We thank Anton Wakolbinger and Martin M\"ohle for pointing out reference \cite{Ewe04} and the fact that the distribution in (\ref{Gadef}) is a $\bet$-distribution.

\section{The renormalization class $\Wi_{\rm cat}$}\label{Wcatsec}

In this section we prove Theorem~\ref{main}~(a) and
Proposition~\ref{identprop}, as well as
Lemmas~\ref{contlem}--\ref{clusK} from Section~\ref{gensec}, and
Lemma~\ref{L:shape}. The section is organized according to the
techniques used. Section~\ref{gensub} collects some facts that hold
for general renormalization classes on compact sets. In
Section~\ref{coupsec} we use the SDE (\ref{catsde}) to couple
catalytic Wright-Fisher diffusions. In Section~\ref{dualsec} we apply
the moment duality for the Wright-Fisher diffusion to the catalyst and
to the reactant conditioned on the catalyst. In
Section~\ref{concavesec} we prove that monotone concave catalyzing
functions form a preserved class under renormalization.

\subsection{Renormalization classes on compact sets}\label{gensub}

In this section, we prove the lemmas stated in Section~\ref{gensec},
as well as Lemma~\ref{L:shape}. Recall that $D\sub\R^d$ is open,
bounded, and convex, and that $\Wi$ is a prerenormalization class on
$\ov D$, equipped with the topology of uniform convergence.\med

\noi
{\bf Proof of Lemma~\ref{contlem}} To see that
$(x,c,w)\mapsto\nu^{c,w}_x$ is continuous, let $(x_n,c_n,w_n)$ be a
sequence converging in $\ov D\times(0,\infty)\times\Wi$ to a limit
$(x,c,w)$. By the compactness of $\ov D$, the sequence
$(\nu_{x_n}^{c_n,w_n})_{n\geq 0}$ is tight, and each limit point
$\nu^\ast$ satisfies
\be\label{invar}
\li\nu^\ast,A^{c,w}_xf\re=0\qquad(f\in\Ci^{(2)}(D)).
\ee
Therefore, by \cite[Theorem~4.9.17]{EK}, $\nu^\ast$ is an invariant
law for the martingale problem associated with $A^{c,w}_x$. Since we
are assuming uniqueness of the invariant law, $\nu^\ast=\nu^{c,w}_x$
and therefore $\nu^{c_n,w_n}_{x_n}\Rightarrow\nu^{c,w}_x$. The
continuity of $F_cw(x)$ is a simple consequence of the continuity of
$\nu^{c,w}_x$.\qed

\noi
{\bf Proof of Lemma~\ref{schaal}} Formula (\ref{schafo})~(i) follows
from the fact that rescaling the time in solutions $(\y_t)_{t\geq 0}$
to the martingale problem for $A^{c,w}_x$ by a factor $\la$ has no
influence on the invariant law. Formula (\ref{schafo})~(ii) is a
direct consequence of formula (\ref{schafo})~(i).\qed

\noi
{\bf Proof of Lemma~\ref{numom}} This follows by inserting the
functions $f(x)=x_i$ and $f(x)=x_ix_j$ into the equilibrium equation
(\ref{invar}).\qed

\noi
{\bf Proof of Lemma~\ref{effinv}} If $x\in\pa_wD$, then $\y_t:=x$
($t\geq 0$) is a stationary solution to the martingale problem for
$A^{c,w}_x$, and therefore $\nu^{c,w}_x=\de_x$ and $F_cw(x)=w(x)=0$.
On the other hand, if $x\not\in\pa_wD$, then $\y_t:=x$ ($t\geq 0$) is
not a stationary solution to the martingale problem for $A^{c,w}_x$
and therefore $\int_{\ov D}\nu^{c,w}_x(\di y)|y-x|^2>0$. Let
$\tr(w(y)):=\sum_iw_{ii}(y)$ denote the trace of $w(y)$. By
(\ref{moments2})~(ii), $\frac{1}{c}\tr(F_cw)(x)=\frac{1}{c}\int_{\ov
  D}\nu^{c,w}_x(\di y)\tr(w(y))=\int_{\ov D}\nu^{c,w}_x(\di
y)|y-x|^2>0$ and therefore $F_cw(x)\neq 0$.\qed

\noi
{F}rom now on assume that $\Wi$ is a renormalization class. Note that
\be\label{iterK}
K^{w,(n)}=\nu^{c_{n-1},F^{(n-1)}w}\cdots\nu^{c_0,w}\qquad(n\geq 1),
\ee
where we denote the composition of two probability kernels $K,L$ on $\ov D$ by
\be
(K L)_x(\di z):=\int_{\ov D} K_x(\di y)L_y(\di z).
\ee
{\bf Proof of Lemma~\ref{basic}} This is a direct consequence of
Lemmas~\ref{contlem} and \ref{numom}. In particular, the relations
(\ref{moments}) follow by iterating the relations (\ref{moments2}).\qed

\noi
{\bf Proof of Lemma~\ref{clusK}} Recall that $\tr(w(y))$ denotes the
trace of $w(y)$. Formulas (\ref{KF}) and (\ref{moments})~(ii) show that
\be
\int_{\ov D}\! K^{w,(n)}_x(\di y)\,|y-x|^2
=s_n\!\int_{\ov D}\! K^{w,(n)}_x(\di y)\,\tr(w(y)).
\ee
Since $\ov D$ is compact, the left-hand side of this equation is
bounded uniformly in $x\in\ov D$ and $n\geq 1$, and therefore,
since we are assuming $s_n\to\infty$,
\be\label{tranul}
\lim_{n\to\infty}\sup_{x\in D}\int_{\ov D} K^{w,(n)}_x(\di y)\tr(w(y))=0.
\ee
Since $w$ is symmetric and nonnegative definite, $\tr(w(y))$ is
nonnegative, and zero if and only if $y\in\pa_w D$. If $f\in\Ci(\ov
D)$ satisfies $f=0$ on $\pa_w D$, then, for every $\eps>0$, the sets
$C_m:=\{x\in\ov D:|f(x)|\geq\eps+m\,\tr(w(x))\}$ are compact with
$C_m\down\emptyset$ as $m\up\infty$, so there exists an $m$ (depending
on $\eps$) such that $|f|<\eps+m\,\tr(w)$. Therefore,
\be\ba{l}
\dis\limsup_{n\to\infty}\:\sup_{x\in\ov D}
\Big|\int_{\ov D}K^{w,(n)}_x(\di y)f(y)\Big|
\leq\limsup_{n\to\infty}\:\sup_{x\in\ov D}
\int_{\ov D}K^{w,(n)}_x(\di y)|f(y)|\\
\dis\quad\leq\eps+m\limsup_{n\to\infty}\:\sup_{x\in\ov D}
\int_{\ov D}K^{w,(n)}_x(\di y)\tr(w(y))=\eps.
\ec
Since $\eps>0$ is arbitrary, (\ref{concen}) follows.\qed

\noi
{\bf Proof of Lemma~\ref{L:shape}} By (\ref{sFw}), (\ref{sumga}), and
(\ref{gamma}), $w^\ast_{\ga^\ast}=\lim_{n\to\infty}(\ov
F_{\ga^\ast})^nw$ for each $w\in\Wi$. By Lemma~\ref{contlem}~(b), $\ov
F_{\ga^\ast}:\Wi\to\Wi$ is continuous, so $w^\ast_{\ga^\ast}$ is the
unique fixed point of $\ov F_{\ga^\ast}$. This proves part~(a).

Now let $0\neq w\in\Wi$ and assume that $\hat\Wi=\{\la w:\la>0\}$ is a
fixed shape. Then $\hat\Wi\ni s_nF^{(n)}w\asto{n}w^\ast_{\ga^\ast}$
whenever $s_n\to\infty$ and $s_{n+1}/s_n\to1+\ga^\ast$ for some
$0<\ga^\ast<\infty$, which shows that $\hat\Wi=\{\la
w^\ast_{\ga^\ast}:\la>0\}$. Thus, $\Wi$ can contain at most one fixed
shape, and if it does, then the $w^\ast_{\ga^\ast}$ for different
values of $\ga^\ast$ must be constant multiples of each other. This
proves part~(c) and the uniqueness statement in part~(b).

To complete the proof of part~(b), note that if $w^\ast=w^\ast_{\ga^\ast}$ does not depend on $\ga^\ast$, then $w^\ast\in\Wi$ solves (\ref{afp})~(i) for all $0<\ga^\ast<\infty$, hence $F_cw^\ast=(1+\frac{1}{c})^{-1}w^\ast$ for all $c>0$, and therefore, by scaling (Lemma~\ref{schaal}), $F_c(\la w^\ast)=\la F_{c/\la}(w^\ast)=\la(1+\frac{\la}{c})^{-1}w^\ast=(\frac{1}{\la}+\frac{1}{c})^{-1}w^\ast$.\qed

\subsection{Coupling of catalytic Wright-Fisher diffusions}\label{coupsec}

In this section we verify condition (i) of Definition~\ref{defi} for
the class $\Wi_{\rm cat}$, and we prepare for the verification of
conditions (ii)--(iv) in Section~\ref{dualsec}. In fact, we will show
that the larger class $\ov\Wi_{\rm cat}:=\{w^{\al,p}:\al>0,\
p\in\Ci_+[0,1]\}$ is also a renormalization class, and the equivalents
of Theorem~\ref{main}~(a) and Proposition~\ref{identprop} remain true
for this larger class. (We do not know, however, if the convergence
statements in Theorem~\ref{main}~(b) also hold in this larger class;
see the discussion in Section~\ref{probsec}.)

For each $c\geq 0$, $w\in\ov\Wi_{\rm cat}$ and $x\in[0,1]^2$, the
operator $A^{c,w}_x$ is a densely defined linear operator on
$\Ci([0,1]^2)$ that maps the identity function into zero and, as one
easily verifies, satisfies the positive maximum principle. Since
$[0,1]^2$ is compact, the existence of a solution to the martingale
problem for $A^{c,w}_x$, for each $[0,1]^2$-valued initial condition,
now follows from general theory (see \cite{RW87}, Theorem~5.23.5, or
\cite[Theorem~4.5.4 and Remark~4.5.5]{EK}).

We are therefore left with the task of verifying uniqueness of
solutions to the martingale problem for $A^{c,w}_x$. By
\cite[Problem~4.19, Corollary~5.3.4, and Theorem 5.3.6]{EK}, it
suffices to show that solutions to (\ref{catsde}) are pathwise unique.
\bl[Monotone coupling of Wright-Fisher diffusions]\label{leq}
Assume that $0\leq x\leq \ti x\leq 1$, $c\geq 0$ and that
$(P_t)_{t\geq 0}$ is a progressively measurable, nonnegative
process such that $\sup_{t\geq 0,\oo\in\om}P_t(\oo)<\infty$.
Let $\y,\ti\y$ be $[0,1]$-valued solutions to the SDE's
\bc
\di\y_t&=&\dis c\,(x-\y_t)\di t+\sqrt{2P_t\y_t(1-\y_t)}\di B_t,\\
\di\ti\y_t&=&\dis c\,(\ti x-\ti\y_t)\di t+\sqrt{2P_t\ti\y_t(1-\ti\y_t)}\di B_t,
\ec
where in both equations $B$ is the same Brownian motion. If
$\y_0\leq \ti\y_0$ a.s., then
\be
\y_t\leq \ti\y_t\quad\forall t\geq 0\quad\as
\ee
\el
{\bf Proof} This is an easy adaptation of a technique due to
Yamada and Watanabe \cite{YW71}. Since $\int_{0+}\frac{\di x}{x}=\infty$,
it is possible to choose $\rho_n\in\Ci\half$ such that
$\int_0^\infty\rho_n(x)\di x=1$ and
\be
0\leq\rho_n(x)\leq\frac{1}{nx}1_{(0,1]}(x)\qquad\qquad(x\geq 0).
\ee
Define $\phi_n\in\Ci^{(2)}(\R)$ by
\be
\phi_n(x):=\int_0^{x\vee 0}\!\!\di y\int_0^y\!\!\di z\,\rho_n(z).
\ee
One easily verifies that $\phi_n(x)$, $x\phi'_n(x)$, and $x\phi''_n(x)$
are nonnegative and converge, as $n\to\infty$, to $x\vee 0$, $x\vee 0$,
and $0$, respectively. By It\^o's formula:
\be\ba{r@{\,}c@{\,}l@{\qquad}r}
\dis E[\phi_n(\y_t-\ti\y_t)]&=&\dis E[\phi_n(\y_0-\ti\y_0)]&{\rm(i)}\\
&&\dis+c\,(x-\ti x)\int_0^t E[\phi'_n(\y_s-\ti\y_s)]\di s
-c\int_0^tE[(\y_s-\ti\y_s)\phi'_n(\y_s-\ti\y_s)]\di s&{\rm(ii)}\\
&&\dis+\int_0^tE\Big[P_s\Big(\sqrt{\y_s(1-\y_s)}
-\sqrt{\ti\y_s(1-\ti\y_s)}\Big)^2\phi''_n(\y_s-\ti\y_s)\Big]\di s.&{\rm(iii)}
\ec
Here the terms in (ii) are nonpositive, and hence, letting
$n\to\infty$ and using the elementary estimate
\be\label{Hoelder}
|\sqrt{y(1-y)}-\sqrt{\ti y(1-\ti y)}|\leq |y-\ti y|^{\frac{1}{2}}
\qquad(y,\ti y\in[0,1]),
\ee
the properties of $\phi_n$, and the fact that the process $P$ is
uniformly bounded, we find that
\be
E[0\vee(\y_t-\ti\y_t)]\leq E[0\vee(\y_0-\ti\y_0)]=0,
\ee
by our assumption that $\y_0\leq\ti\y_0$. This shows that
$\y_t\leq\ti\y_t$ a.s. for each fixed $t\geq 0$, and by the
continuity of sample paths the statement holds for all
$t\geq 0$ almost surely.\qed

\begin{corollary}[Pathwise uniqueness]\label{uni}
For all $c\geq 0$, $\al>0$, $p\in\Ci_+[0,1]$ and $x\in[0,1]$,
solutions to the SDE (\ref{catsde}) are pathwise unique.
\end{corollary}
{\bf Proof} Let $(\y^1,\y^2)$ and $(\ti\y^1,\ti\y^2)$ be solutions to
(\ref{catsde}) relative to the same pair $(B^1,B^2)$ of Brownian
motions, with $(\y^1_0,\y^2_0)=(\ti\y^1_0,\ti\y^2_0)$. Applying
Lemma~\ref{leq}, with inequality in both directions, we see that
$\y^1=\ti\y^1$ a.s. Applying Lemma~\ref{leq} two more times, this time
using that $\y^1=\ti\y^1$ a.s., we see that also $\y^2=\ti\y^2$
a.s.\qed

\begin{corollary}[Exponential coupling]\label{excoup}
Assume that $x\in[0,1]$, $c\geq 0$, and $\al>0$. Let $\y,\ti\y$
be solutions to the SDE
\be\label{alWF}
\di\y_t=c\,(x-\y_t)\di t+\sqrt{2\al\y_t(1-\y_t)}\di B_t,
\ee
relative to the same Brownian motion $B$. Then
\be\label{exco}
E\big[|\ti\y_t-\y_t|\big]=e^{-ct}E\big[|\ti\y_0-\y_0|\big].
\ee
\end{corollary}
{\bf Proof} If $\y_0=y$ and $\ti\y_0=\ti y$ are deterministic and
$y\leq\ti y$, then by Lemma~\ref{leq} and a simple moment calculation
\be\label{absmom}
E\big[|\ti\y_t-\y_t|\big]=E[\ti\y_t-\y_t]=e^{-ct}|\ti y-y|.
\ee
The same argument applies when $y\geq\ti y$. The general case where
$\y_0$ and $\ti\y_0$ are random follows by conditioning on
$(\y_0,\ti\y_0)$.\qed
\begin{corollary}[Ergodicity]\label{ergo}
The Markov process defined by the SDE (\ref{Yclx}) has a unique
invariant law $\Ga^\ga_x$ and is ergodic, i.e, solutions to
(\ref{Yclx}) started in an arbitrary initial law $\Li(\y_0)$
satisfy $\Li(\y_t)\Asto{t}\Ga^\ga_x$.
\end{corollary}
{\bf Proof} Since our process is a Feller diffusion on a compactum,
the existence of an invariant law follows from a simple time averaging
argument. Now start one solution $\ti\y$ of (\ref{Yclx}) in this
invariant law and let $\y$ be any other solution, relative to the same
Brownian motion. Corollary~\ref{excoup} then gives ergodicity and, in
particular, uniqueness of the invariant law.\qed

\brm{\bf (Density of invariant law)}\label{R:betadis}
It is well-known (see, for example \cite[formula~(5.70)]{Ewe04}) that
$\Ga^\ga_x$ is a $\bet(\al_1,\al_2)$-distribution, where
$\al_1:=x/\ga$ and $\al_2:=(1-x)/\ga$, i.e., $\Ga^\ga_x=\de_x$
$(x\in\{0,1\})$ and
\be\label{Gadef}
\Ga^\ga_x(\di y)=\frac{\Ga(\al_1+\al_2)}{\Ga(\al_1)\Ga(\al_2)}
\,y^{\al_1-1}(1-y)^{\al_2-1}\di y\qquad(x\in(0,1)).
\ee
\erm

\noi
We conclude this section with a lemma that prepares for the verification
of condition (iv) in Definition~\ref{defi} for the class $\Wi_{\rm cat}$.
\bl\label{monot}
{\bf (Monotone coupling of stationary Wright-Fisher diffusions)} Assume
that $c>0$, $\al>0$ and $0\leq x\leq\ti x\leq 1$. Then the pair of equations
\bc\label{kopwf}
\di\y_t&=&\dis c\,(x-\y_t)\di t+\sqrt{2\al\y_t(1-\y_t)}\di B_t,\\
\di\ti\y_t&=&\dis c\,(\ti x-\ti\y_t)\di t+\sqrt{2\al\ti\y_t(1-\ti\y_t)}\di B_t
\ec
has a unique stationary solution $(\y_t,\ti\y_t)_{t\in\R}$. This
stationary solution satisfies
\be\label{leqt}
\y_t\leq\ti\y_t\quad\forall t\in\R\quad\as
\ee
\el
{\bf Proof} Let $(\y_t,\ti\y_t)_{t\geq 0}$ be a solution of
(\ref{kopwf}) and let $(\y'_t,\ti\y'_t)_{t\geq 0}$ be another one,
relative to the same Brownian motion $B$. Then, by Lemma~\ref{excoup},
$E[|\y_t-\y'_t|]\to 0$ and also $E[|\ti\y_t-\ti\y'_t|]\to 0$ as
$t\to\infty$. Hence we may argue as in the proof of
Corollary~\ref{ergo} that (\ref{kopwf}) has a unique invariant law and
is ergodic. Now start a solution of (\ref{kopwf}) in an initial
condition such that $\y_0\leq\ti\y_0$. By ergodicity, the law of this
solution converges as $t\to\infty$ to the invariant law of
(\ref{kopwf}) and using Lemma~\ref{leq} we see that this invariant law
is concentrated on $\{(y,\ti y)\in[0,1]^2:y\leq\ti y\}$. Now consider,
on the whole real time axis, the stationary solution to (\ref{kopwf})
with this invariant law. Applying Lemma~\ref{leq} once more, we see
that (\ref{leqt}) holds.\qed

\subsection{Duality for catalytic Wright-Fisher diffusions}\label{dualsec}

In this section we prove Theorem~\ref{main}~(a) and
Proposition~\ref{identprop}. Moreover, we will show that their
statements remain true if the renormalization class $\Wi_{\rm cat}$ is
replaced by the larger class $\ov\Wi_{\rm cat}:=\{w^{\al,p}:\al>0,\
p\in\Ci_+[0,1]\}$. We begin by recalling the usual moment duality for
Wright-Fisher diffusions.

For $\ga>0$ and $x\in[0,1]$, let $\y$ be a solution to the SDE
\be\label{WFsde}
\di\y(t)=\ffrac{1}{\ga}\,(x-\y(t))\di t+\sqrt{2\y(t)(1-\y(t))}\di B(t),
\ee
i.e., $\y$ is a Wright-Fisher diffusion with a linear drift towards $x$.
It is well-known that $\y$ has a moment dual. To be precise, let
$(\phi,\psi)$ be a Markov process in $\N^2=\{0,1,\ldots\}^2$ that jumps as:
\be\ba{r@{\,}c@{\,}l@{\qquad}l}\label{phipsi}
(\phi_t,\psi_t)&\to&(\phi_t-1,\psi_t)
&\mbox{with rate }\ \phi_t(\phi_t-1)\\
(\phi_t,\psi_t)&\to&(\phi_t-1,\psi_t+1)
&\mbox{with rate }\ \ffrac{1}{\ga}\phi_t.
\ec
Then one has the following {\em duality relation} (see for example
Lemma 2.3 in \cite{Shi80} or Proposition~1.5 in \cite{GKW01})
\be\label{WFdual}
E^y\big[\y_t^nx^m\big]=E^{(n,m)}\big[y^{\phi_t}x^{\psi_t}\big]
\qquad(y\in[0,1],\ (n,m)\in\N^2),
\ee
where $0^0:=1$. The duality in (\ref{WFdual}) has the following
heuristic explanation. Consider a population containing a fixed, large
number of organisms, that come in two genetic types, say I and II.
Each pair of organisms in the population is {\em resampled} with rate
$2$. This means that one organism of the pair (chosen at random) dies,
while the other organism produces one child of its own genetic type.
Moreover, each organism is replaced with rate $\frac{1}{\ga}$ by an
organism chosen from an infinite reservoir where the frequency of type
I has the fixed value $x$. In the limit that the number of organisms
in the population is large, the relative frequency $\y_t$ of type I
organisms follows the SDE (\ref{WFsde}). Now $E[\y_t^n]$ is the
probability that $n$ organisms sampled from the population at time $t$
are all of type I. In order to find this probability, we follow the
ancestors of these organisms back in time. Viewed backwards in time,
these ancestors live for a while in the population, until, with rate
$\frac{1}{\ga}$, they jump to the infinite reservoir. Moreover, due to
resampling, each pair of ancestors coalesces with rate $2$ to one
common ancestor. Denoting the number of ancestors that lived at time
$t-s$ in the population and in the reservoir by $\phi_s$ and $\psi_s$,
respectively, we see that the probability that all ancestors are of
type~I is $E^y[\y_t^n]=E^{(n,0)}[y^{\phi_t}x^{\psi_t}]$. This gives a
heuristic explanation of (\ref{WFdual}).

Since eventually all ancestors of the process $(\phi,\psi)$ end up in
the reservoir, we have $(\phi_t,\psi_t)\to(0,\psi_\infty)$ as
$t\to\infty$ a.s.\ for some $\N$-valued random variable $\psi_\infty$.
Taking the limit $t\to\infty$ in (\ref{WFdual}), we see that the
moments of the invariant law $\Ga^\ga_x$ from Corollary~\ref{ergo} are
given by:
\be\label{Gamom}
\int\Ga^\ga_x(\di y)y^n=E^{(n,0)}[x^{\psi_\infty}]\qquad(n\geq 0).
\ee
It is not hard to obtain an inductive formula for the moments of
$\Ga^\ga_x$, which can then be solved to yield the formula
\be\label{WFmoments}
\int\Ga^\ga_x(\di y)y^n=\prod_{k=0}^{n-1}\frac{x+k\ga}{1+k\ga}\qquad(n\geq 1).
\ee
In particular, it follows that
\be\label{WFfix}
\int\Ga^\ga_x(\di y)y(1-y)=\frac{1}{1+\ga}x(1-x).
\ee
This is the important {\em fixed shape property} of the Wright-Fisher
diffusion (see formula (\ref{fix})).

We now consider catalytic Wright-Fisher diffusions $(\y^1,\y^2)$ as in
(\ref{catsde}) with $p\in\Ci_+[0,1]$ and apply duality to the catalyst
$\y^2$ conditioned on the reactant $\y^1$. Let
$(\y^1_t,\y^2_t)_{t\in\R}$ be a stationary solution to the SDE
(\ref{catsde}) with $c=1/\ga$. Let $(\ti\phi,\ti\psi)$ be a
$\N^2$-valued process, defined on the same probability space as
$(\y^1,\y^2)$, such that conditioned on the past path
$(\y^1_{-t})_{t\leq 0}$, the process $(\ti\phi,\ti\psi)$ is a
(time-inhomogeneous) Markov process that jumps as:
\be\ba{r@{\,}c@{\,}l@{\qquad}l}\label{yphipsi}
(\ti\phi_t,\ti\psi_t)&\to&(\ti\phi_t-1,\ti\psi_t)
&\mbox{with rate }\ p(\y^1_{-t})\ti\phi_t(\ti\phi_t-1),\\
(\ti\phi_t,\ti\psi_t)&\to&(\ti\phi_t-1,\ti\psi_t+1)
&\mbox{with rate }\ \ffrac{1}{\ga}\ti\phi_t.
\ec
Then, in analogy with (\ref{WFdual}),
\be\label{reactdual}
E[(\y^2_0)^nx_2^m|(\y^1_{-t})_{t\leq 0}]
=E^{(n,m)}[(\y^2_{-t})^{\ti\phi_t}x_2^{\ti\psi_t}|(\y^1_{-t})_{t\leq 0}]
\qquad((n,m)\in\N^2,\ t\geq 0).
\ee
We may interpret (\ref{yphipsi}) by saying that pairs of ancestors in
a finite population coalesce with time-dependent rate $2p(\y^1_{-t})$
and ancestors jump to an infinite reservoir with constant rate
$\frac{1}{\ga}$.  Again, eventualy all ancestors end up in the
reservoir, and therefore $(\ti\phi_t,\ti\psi_t)\to(0,\ti\psi_\infty)$
as $t\to\infty$ a.s.\ for some $\N$-valued random variable
$\ti\psi_\infty$. Taking the limit $t\to\infty$ in (\ref{reactdual})
we find that
\be\label{reactdual2}
E[(\y^2_0)^nx_2^m|(\y^1_{-t})_{t\leq 0}]
=E^{(n,m)}[x_2^{\ti\psi_\infty}|(\y^1_{-t})_{t\leq 0}]
\qquad((n,m)\in\N^2,\ t\geq 0).
\ee
\bl{\bf(Uniqueness of invariant law)}\label{unin}
For each $c>0$, $w\in\ov\Wi_{\rm cat}$, and $x\in[0,1]^2$, there exists
a unique invariant law $\nu^{c,w}_x$ for the martingale problem for
$A^{c,w}_x$.
\el
{\bf Proof} Our process being a Feller diffusion on a compactum, the
existence of an invariant law follows from time averaging. We need to
show uniqueness. If $(\y^1,\y^2)=\y^1_t,\y^2_t)_{t\in\R}$ is a
stationary solution, then $\y^1$ is an autonomous process, and
$\Li(\y^1_0)=\Ga^{1/c}_x$, the unique invariant law from
Corollary~\ref{ergo}. Therefore, $\Li((\y^1_t)_{t\in\R})$ is
determined uniquely by the requirement that $(\y^1,\y^2)$ be
stationary. By (\ref{reactdual2}), the conditional distribution of
$\y^2_0$ given $(\y^1_t)_{t\leq 0}$ is determined uniquely, and
therefore the joint distribution of $\y^2_0$ and $(\y^1_t)_{t\leq 0}$
is determined uniquely. In particular,
$\Li(\y^1_0,\y^2_0)=\nu^{c,w}_x$ is determined uniquely.\qed

\brm{\bf (Reversibility)} 
It seems that the invariant law $\nu^{c,w}_x$ from Lemma~\ref{unin} is
reversible. In many cases (densities of) reversible invariant measures
can be obtained in closed form by solving the equations of detailed
balance. This is the case, for example, for the one-dimensional
Wright-Fisher diffusion. We have not attempted this for the catalytic
Wright-Fisher diffusion.
\erm
The next proposition implies Proposition~\ref{identprop} and prepares
for the proof of Theorem~\ref{main}~(a).
\bp{\bf(Extended renormalization class)}\label{exclass}
The set $\ov\Wi_{\rm cat}$ is a renormalization class on $[0,1]^2$, and
\be\label{hutfo}
\hut F_\ga w^{\txt 1,p}=w^{\txt 1,\Ui_\ga p}\qquad(p\in\Ci_+[0,1],\ \ga>0).
\ee
\ep
{\bf Proof} To see that $\ov\Wi_{\rm cat}$ is a renormalization class
we need to check conditions~(i)--(iv) from Definition~\ref{defi}. By
Lemma~\ref{uni}, the martingale problem for $A^{c,w}_x$ is well-posed
for all $c\geq 0$, $w\in\Wi_{\rm cat}$ and $x\in[0,1]^2$. By
Lemma~\ref{unin}, the corresponding Feller process on $[0,1]^2$ has a
unique invariant law $\nu^{c,w}_x$. This shows that conditions~(i) and
(ii) from Definition~\ref{defi} are satisfied. Note that by the
compactness of $[0,1]^2$, any continuous function on $[0,1]^2$ is
bounded, so condition~(iii) is automatically satisfied. Hence $\Wi$ is
a prerenormalization class. As a consequence, for any
$p\in\Ci_+[0,1]$, $\hut F_\ga w^{1,p}$ is well-defined by (\ref{Fc})
and (\ref{hatFc}). We will now first prove (\ref{hutfo}) and then show
that $\ov\Wi_{\rm cat}$ is a renormalization class.

Fix $\ga>0$, $p\in\Ci_+[0,1]$, and $x\in[0,1]^2$. Let
$(\y^1_t,\y^2_t)_{t\in\R}$ be a stationary solution to the SDE
(\ref{catsde}) with $\al=1$ and $c=1/\ga$. Then
\be
\hut F_\ga w^{1,p}_{ij}(x)=(1+\ga)E[w^{1,p}_{ij}(\y^1_0,\y^2_0)]
\qquad(i,j=1,2).
\ee
Since $w^{1,p}_{ij}=0$ if $i\neq j$, it is clear that $\hut F_\ga
w^{1,p}_{ij}(x)=0$ if $i\neq j$. Since $\Li(\y^1_0)=\Ga^\ga_x$ it
follows from (\ref{WFfix}) that $\hut F_\ga
w^{1,p}_{11}(x)=x_1(1-x_1)$. We are left with the task of showing that
\be
\hut F_\ga w^{1,p}_{22}(x)=\Ui_\ga p(x_1)x_2(1-x_2).
\ee
Here, by (\ref{moments2})~(ii),
\bc\label{momtruc}
\dis\hut F_\ga w^{1,p}_{22}(x)&=&\dis(1+\ga)E[p(\y^1_0)\y^2_0(1-\y^2_0)]\\[5pt]
&=&\dis(\ffrac{1}{\ga}+1)E[(\y^2_0-x_2)^2].
\ec
By (\ref{reactdual2}), using the fact that $E[\y^2_0]=x_2$ (which
follows from (\ref{reactdual}) or more elementary from (\ref{moments})~(i)),
we find that
\be\label{coal}
E[(\y^2_0-x_2)^2]=E[(\y^2_0)^2]-(x_2)^2
=E^{(2,0)}[x_2^{\ti\psi_\infty}]-(x_2)^2=P^{(2,0)}[\ti\psi_\infty=1]x_2(1-x_2)
\qquad(t\geq 0).
\ee
Note that $P^{(2,0)}[\ti\psi_\infty=1]$ is the probability that the
two ancestors coalesce before one of them leaves the population. The
probability of {\em noncoalescence} is given by
\be\label{noncoa}
P^{(2,0)}[\ti\psi_\infty=2]
=E\big[\ex{-\int_0^{\frac{1}{2}\tau_\ga}2p(y^1_{-t})\di t}\big],
\ee
where $\tau_\ga$ is an exponentially distributed random variable
with mean $\ga$. Combining this with (\ref{momtruc}) and (\ref{coal})
we find
that
\bc
\dis\hut F_\ga w^{1,p}_{22}(x)&=&\dis(\ffrac{1}{\ga}+1)
E\big[1-\ex{-\int_0^{\tau_\ga}p(y^1_{-t/2})\di t}\big]x_2(1-x_2)\\[5pt]
&=&\dis q_\ga E\big[1-\ex{-\li\Zi^\ga_x,p\re}\big]x_2(1-x_2)\\[5pt]
&=&\dis\Ui_\ga p(x_1)x_2(1-x_2),
\ec
where we have used the definition of $\Ui_\ga$.

We still have to show that $\ov\Wi_{\rm cat}$ satisfies condition~(iv)
from Definition~\ref{defi}. For any $\al>0$ and $p\in\Ci_+[0,1]$, by
scaling (Lemma~\ref{schaal}) and (\ref{hutfo}),
\be\label{wp}
F_cw^{\txt\al,p}=\al F_{\frac{c}{\al}}w^{\txt 1,\frac{p}{\al}}
=\al(1+\frac{\al}{c})^{-1}\ov F_{\frac{c}{\al}}w^{\txt 1,\frac{p}{\al}}
=w^{\txt(\frac{1}{\al}+\frac{1}{c})^{-1},
(\frac{1}{\al}+\frac{1}{c})^{-1}\Ui_{\frac{c}{\al}}(\frac{p}{\al})}.
\ee
By Lemma~\ref{contlem}, this diffusion matrix is continuous, which
implies that $\Ui_{\frac{c}{\al}}(\frac{p}{\al})$ is continuous.\qed

\noi
Our proof of Propostion~\ref{exclass} has a corollary.
\bcor{\bf(Continuity in parameters)}\label{Ugacont}
The map $(x,\ga)\mapsto\Qi_\ga(x,\cdot)$ from $[0,1]\times(0,\infty)$ to
$\Mi_1(\Mi[0,1])$ and the map $(x,\ga,p)\mapsto\Ui_\ga p(x)$ from
$[0,1]\times(0,\infty)\times\Ci_+[0,1]$ to $\R$ are continuous.
\ecor
{\bf Proof} By Lemma~\ref{contlem}, the diffusion matrix in (\ref{wp})
is continuous in $x,\ga$, and $p$, which implies the continuity of
$\Ui_\ga p(x)$. It follows that the map
$(x,\ga)\mapsto\int\Qi_\ga(x,\di\chi)\ex{-\li\chi,f\re}$ is continuous
for all $f\in\Ci_+[0,1]$, so by \cite[Theorem~4.2]{Kal76},
$(x,\ga)\mapsto\Qi_\ga(x,\cdot)$ is continuous.\qed

\noi
{\bf Proof of Theorem~\ref{main}~(a)}\label{aproof} We need to show
that $\Wi_{\rm cat}$ is a renormalization class and that $F_c$ maps
the subclasses $\Wi^{l,r}_{\rm cat}$ into themselves. Since these
classes correspond to the different possible effective boundaries of
diffusion matrices in $\Wi_{\rm cat}$, this latter fact is in fact a
consequence of Lemma~\ref{effinv}. Since in Proposition~\ref{exclass}
it has been shown that $\ov\Wi_{\rm cat}$ is a renormalization class,
we are left with the task to show that $F_c$ maps $\Wi_{\rm cat}$ into
itself. By (\ref{hutfo}) and scaling, it suffices to show that
$\Ui_\ga$ maps $\Hi$ into itself.

Fix $0\leq x\leq\ti x\leq 1$. By Lemma~\ref{monot}, we can couple the
processes $\y^\ga_x$ and $\y^\ga_{\ti x}$ from (\ref{Yclx}) such that
\be\label{leqt2}
\y^\ga_x(t)\leq\y^\ga_{\ti x}(t)\quad\forall t\leq 0\quad\as
\ee
Since the function $z\mapsto 1-e^{-z}$ on $\half$ is Lipschitz continuous
with Lipschitz constant $1$,
\be\ba{l}\label{Lipone}
\dis\big|\Ui_\ga p(\ti x)-\Ui_\ga p(x)\big|\\[5pt]
\dis\quad=\Big|(\ffrac{1}{\ga}+1)
E\big[1-\ex{-\int_0^{\tau_\ga}p(\y^\ga_{\ti x}(-t/2))\di t}\big]
-(\ffrac{1}{\ga}+1)
E\big[1-\ex{-\int_0^{\tau_\ga}p(\y^\ga_x(-t/2))\di t}\big]\Big|\\[5pt]
\dis\quad\leq (\ffrac{1}{\ga}+1)E\Big[\int_0^{\tau_\ga}\big|
p(\y^\ga_{\ti x}(-t/2))-p(\y^\ga_x(-t/2))\big|\di t\Big]\\[5pt]
\dis\quad\leq (\ffrac{1}{\ga}+1)LE\Big[\int_0^{\tau_\ga}\big|
\y^\ga_{\ti x}(-t/2)-\y^\ga_x(-t/2)\big|\di t\Big]\\[5pt]
\dis\quad=(\ffrac{1}{\ga}+1)L\ga(\ti x-x)=L(1+\ga)|\ti x-x|,
\ec
where $L$ is the Lipschitz constant of $p$ and we have used the
same exponentially distributed $\tau_\ga$ for $\y^\ga_x$ and
$\y^\ga_{\ti x}$.\qed

\subsection{Monotone and concave catalyzing functions}\label{concavesec}

In this section we prove that the log-Laplace operators $\Ui_\ga$ from
(\ref{Uiga}) map monotone functions into monotone functions, and
monotone concave functions into monotone concave functions. We do not
know if in general $\Ui_\ga$ maps concave functions into concave
functions.
 \bp{\bf(Preservation of monotonicity and concavity)}\label{moncon}
Let $\ga>0$. Then:\smallskip

\noi
{\bf (a)} If $f\in\Ci_+[0,1]$ is nondecreasing, then $\Ui_\ga f$ is
nondecreasing.\smallskip

\noi
{\bf (b)} If $f\in\Ci_+[0,1]$ is nondecreasing and concave, then
$\Ui_\ga f$ is nondecreasing and concave.
\ep
{\bf Proof} Our proof of Proposition~\ref{moncon} is in part based on
ideas from \cite[Appendix~A]{BCGH97}. The proof is quite long and will
depend on several lemmas. We remark that part~(a) can be proved in a
more elementary way using Lemma~\ref{monot}.

We recall some facts from Hille-Yosida theory. A linear operator $A$
on a Banach space $V$ is closable and its closure $\ov A$ generates a
strongly continuous contraction semigroup $(S_t)_{t\geq 0}$ if and only if
\be\ba{rl}\label{HilYos1}
{\rm (i)}&\Di(A)\mbox{ is dense},\\
{\rm (ii)}&A\mbox{ is dissipative},\\
{\rm (iii)}&\Ri(1-\al A)\mbox{ is dense for some, and hence for all }\al>0.
\ec
\detail{To see that condition (iii) holds for all $\al>0$ if it holds
for some $\al>0$, one needs the elementary fact that $\ov{\Ri(1-\al A)}
=\Ri(1-\al\ov A)$.}
Here, for any linear operator $B$ on $V$, $\Di(B)$ and $\Ri(B)$ denote
the domain and range of $B$, respectively. For each $\al>0$, the operator
$(1-\al\ov A):\Di(\ov A)\to V$ is a bijection and its inverse
$(1-\al\ov A)^{-1}:V\to\Di(\ov A)$ is a bounded linear operator, given by
\be
(1-\al\ov A)^{-1}u=\int_0^\infty\!\! S_tu\;\al^{-1}e^{-t/\al}\di t
\qquad(u\in V,\ \al>0).
\ee
If $E$ is a compact metrizable space and $\Ci(E)$ is the Banach space
of continuous real functions on $E$, equipped with the supremumnorm,
then a linear operator $A$ on $\Ci(E)$ is closable and its closure
$\ov A$ generates a Feller semigroup if and only if (see
\cite[Theorem~4.2.2 and remarks on page~166]{EK})
\be\ba{rl}\label{HilYos2}
{\rm (i)}&1\in\Di(\ov A)\mbox{ and }\ov A1=0,\\
{\rm (ii)}&\Di(A)\mbox{ is dense},\\
{\rm (iii)}&A\mbox{ satisfies the positive maximum principle},\\
{\rm (iv)}&\Ri(1-\al A)\mbox{ is dense for some, and hence for all }\al>0.
\ec
If $\ov A$ generates a Feller semigroup and $g\in\Ci(E)$, then the
operator $\ov A+g$ (with domain $\Di(\ov A+g):=\Di(\ov A)$) generates
a strongly continuous semigroup $(S^g_t)_{t\geq 0}$ on $\Ci(E)$. If
$g\leq 0$ then $(S^g_t)_{t\geq 0}$ is contractive. If $(\xi_t)_{t\geq 0}$
is the Feller process with generator $\ov A$, then one has the
Feynman-Kac representation
\be\label{FeyKac}
S^g_tu(x)=E^x[u(\xi(t))\ex{\int_0^tg(\xi(s))\di s}\big]
\qquad(t\geq 0,\ x\in E,\ g,u\in\Ci(E)).
\ee
Let $\Ci^{(n)}([0,1]^2)$ denote the space of continuous real functions
on $[0,1]^2$ whose partial derivatives up to $n$-th order exist and
are continuous on $[0,1]^2$ (including the boundary), and put
$\Ci^{(\infty)}([0,1]^2):=\bigcap_n\Ci^{(n)}([0,1]^2)$. Define a
linear operator $B$ on $\Ci([0,1]^2)$ with domain
$\Di(B):=\Ci^{(\infty)}([0,1]^2)$ by
\be\label{Bdef}
Bu(x,y):=y(1-y)\diff{y}u(x,y)+\ffrac{1}{\ga}(x-y)\dif{y}u(x,y).
\ee
Below, we will prove:
\bl{\bf(Feller semigroup)}\label{Felgen}
The closure in $\Ci([0,1]^2)$ of the operator $B$ generates a Feller
semigroup on $\Ci([0,1]^2)$.
\el
Write
\be\ba{c@{\,}c@{\,}l}
\Ci_+&:=&\dis\big\{u\in\Ci([0,1]^2):u\geq 0\big\},\\[5pt]
\Ci_{1+}&:=&\dis
\big\{u\in\Ci^{(1)}([0,1]^2):\dif{y}u,\dif{x}u\geq 0\big\},\\[5pt]
\Ci_{2+}&:=&\dis
\big\{u\in\Ci^{(2)}([0,1]^2):\diff{y}u,\difif{x}{y}u,\diff{x}u\geq 0\big\}.
\ec
Let $\ov\Si$ denote the closure of a set $\Si\sub\Ci([0,1]^2)$.
We need the following lemma.
\bl{\bf (Preserved classes)}\label{Lapcon}
Let $g\in\Ci([0,1]^2)$ and let $(S^g_t)_{t\geq 0}$ be the strongly
continuous semigroup with generator $\ov B+g$. Then, for each
$t\geq 0$:\smallskip

\noi
{\bf (a)} If $g\in\ov{\Ci_{1+}}$, then $S^g_t$ maps
$\ov{\Ci_+\cap\Ci_{1+}}$ into itself.\smallskip

\noi
{\bf (b)} If $g\in\ov{\Ci_{1+}\cap\Ci_{2+}}$, then $S^g_t$ maps
$\ov{\Ci_+\cap\Ci_{1+}\cap\Ci_{2+}}$ into itself.
\el
To see why Lemma~\ref{Lapcon} implies Proposition~\ref{moncon}, let
$(\x(t),\y(t))_{t\geq 0}$ denote the Feller process in $[0,1]^2$
generated by $\ov B$. It is easy to see that $\x(t)=\x(0)$ a.s.\ for
all $t\geq 0$. For fixed $\x(0)=x$, the process $(\y(t))_{t\geq 0}$ is
the diffusion given by the SDE (\ref{WFsde}). Therefore, by
Feynman-Kac, for each $g\in\Ci([0,1]^2)$,
\be\label{FK}
E^y\big[\ex{\int_0^tg(x,\y(s))\di s}\big]=S^g_t1(x,y),
\ee
where $1$ denotes the constant function $1\in\Ci([0,1]^2)$. By (\ref{Uiga}),
\be\label{Uiexpr}
\Ui_\ga f(x)=(\ffrac{1}{\ga}+1)\Big(1-\int\Ga^\ga_x(\di y)
E^y\big[\ex{-\int_0^{\tau_\ga}f(\y_x(s))\di s}\big]\Big)
\qquad(f\in\Ci_+[0,1]),
\ee
where $\Ga^\ga_x$ is the invariant law of $(\y(t))_{t\geq 0}$ from
Corollary~\ref{ergo} and $\tau_\ga$ is an exponential time with mean
$\ga$, independent of $(\y(t))_{t\geq 0}$. Setting $g(x,y):=-f(y)$ in
(\ref{FK}), using the ergodicity of $(\y(t))_{t\geq 0}$ (see
Corollary~\ref{ergo}), we find that for each $z\in[0,1]$ and $t\geq
0$,
\bc\label{limla}
\dis\int\Ga^\ga_x(\di y)E^y\big[\ex{-\int_0^tf(\y(s))\di s}\big]
&=&\dis\lim_{r\to\infty}\int P^z[\y(r)\in\di y]\,
E^y\big[\ex{-\int_0^tg(x,\y(s))\di s}\big]\\[7pt]
&=&\dis\lim_{r\to\infty}S^0_rS^g_t1(x,z).
\ec
It follows from Lemma~\ref{Lapcon} that for each fixed $r,t$, and $z$,
the function $x\mapsto S^0_rS^g_t1(x,z)$ is nondecreasing if $f$ is
nonincreasing, and nondecreasing and convex if $f$ is nonincreasing
and concave. Therefore, taking the expectation over the randomness of
$\tau_\ga$, the claims follow from (\ref{Uiexpr}) and
(\ref{limla}).\qed

\noi
We still need to prove Lemmas~\ref{Felgen} and \ref{Lapcon}.\med

\noi
{\bf Proof of Lemma~\ref{Felgen}} It is easy to see that the operator
$B$ from (\ref{Bdef}) is densely defined, satisfies the positive
maximum principle, and maps the constant function $1$ into $0$.
Therefore, by Hille-Yosida (\ref{HilYos2}), we must show that the
range $\Ri(1-\al B)$ is dense in $\Ci([0,1]^2)$ for some, and hence
for all $\al>0$. Let $\Pc_n$ denote the space of polynomials on
$[0,1]^2$ of $n$-th and lower order, i.e., the space of functions
$f:[0,1]^2\to\R$ of the form
\be
f(x,y)=\sum_{k,l\geq 0}a_{kl}\,x^ky^l\quad\mbox{ with $a_{k,l}=0$ for $k+l>n$.}
\ee
Set $\Pc_\infty:=\bigcup_n\Pc_n$. It is easy to see that $B$ maps the
space $\Pc_n$ into itself, for each $n\geq 0$. Since each $\Pc_n$ is
finite-dimensional, a simple argument (see
\cite[Proposition~1.3.5]{EK}) shows that the image of $\Pc_\infty$
under $1-\al B$ is dense in $\Ci([0,1]^2)$ for all but countably many,
and hence for all $\al>0$.\qed

\noi
As a first step towards proving Lemma~\ref{Lapcon}, we prove:
\bl{\bf(Smooth solutions to Laplace equation)}\label{smoothlap}
Let $\al>0$, $g\in\Ci^{(2)}([0,1])$, $g\leq 0$, $v\in\Ci([0,1]^2)$,
and assume that $u\in\Ci^{(\infty)}([0,1]^2)$ solves the Laplace equation
\be\label{Lap}
(1-\al(B+g))u=v.
\ee
{\bf (a)} If $g\in\Ci_{1+}$, then $v\in\Ci_{+}\cap\Ci_{1+}$ implies
$u\in\Ci_{+}\cap\Ci_{1+}$.\smallskip

\noi
{\bf (b)} If $g\in\Ci_{1+}\cap\Ci_{2+}$, then
$v\in\Ci_{+}\cap\Ci_{1+}\cap\Ci_{2+}$ implies
$u\in\Ci_{+}\cap\Ci_{1+}\cap\Ci_{2+}$.
\el
{\bf Proof} Let $u^y:=\dif{y}u$, $u^{xy}:=\difif{x}{y}u$, etc.\ denote
the partial derivatives of $u$ and similarly for $v$ and $g$, whenever
they exist. Set $c:=\frac{1}{\ga}$. Define linear operators $B'$ and
$B''$ on $\Ci([0,1]^2)$ with domains
$\Di(B')=\Di(B''):=\Ci^{(\infty)}([0,1]^2)$ by
\bc
B'&:=&\dis y(1-y)\diff{y}+\big(c(x-y)+2(\ffrac{1}{2}-y)\big)\dif{y},\\[5pt]
B''&:=&\dis y(1-y)\diff{y}+\big(c(x-y)+4(\ffrac{1}{2}-y)\big)\dif{y}.
\ec
Then
\be\ba{r@{\,}c@{\,}l@{\quad}r@{\,}c@{\,}l}\label{Bcommut}
\dif{y}Bu&=&(B'-c)u^y,&\dif{y}B'u&=&(B''-c-2)u^y,\\[5pt]
\dif{x}Bu&=&Bu^x+cu^y,&\dif{x}B'u&=&B'u^x+cu^y.
\ec
Therefore, it is easy to see that
\be\ba{rr@{\,}c@{\,}l}\label{partdi}
{\rm (i)}&(1-\al(B'-c+g))u^y&=&v^y+\al g^yu,\\
{\rm (ii)}&(1-\al(B+g))u^x&=&v^x+\al(cu^y+g^xu),\\
{\rm (iii)}&(1-\al(B''-2c-2+g))u^{yy}&=&v^{yy}+\al(2g^yu^y+g^{yy}u),\\
{\rm (iv)}&(1-\al(B'-c+g))u^{xy}&=&v^{xy}+\al(cu^{yy}+g^yu^x+g^{xy}u+g^xu^y),\\
{\rm (v)}&(1-\al(B+g))u^{xx}&=&v^{xx}+\al(2cu^{xy}+2g^xu^x+g^{xx}u),
\ec
where in (i) and (ii) we assume that $v\in\Ci^{(1)}([0,1]^2)$ and in
(iii)--(v) we assume that $v\in\Ci^{(2)}([0,1]^2)$. By
Lemma~\ref{Felgen}, the closure of the operator $B$ generates a Feller
processes in $[0,1]^2$. Exactly the same proof shows that $B'$ and
$B''$ also generate Feller processes on $[0,1]^2$. Therefore, by
Feynman-Kac, $u$ is nonnegative if $v$ is nonnegative and
$u^y,\ldots,u^{xx}$ are nonnegative if the right-hand sides of the
equations (i)--(v) are well-defined and nonnegative. (Instead of using
Feynman-Kac, this follows more elementarily from the fact that $B,B'$,
and $B''$ satisfy the positive maximum principle.) In particular, if
$g^y,g^x\geq 0$ and $v\in\Ci^{(1)}([0,1]^2)$, $v,v^y,v^x\geq 0$, then
it follows that $u,u^y,u^x\geq 0$. If moreover
$g^{yy},g^{xy},g^{xx}\geq 0$ and $v\in\Ci^{(2)}([0,1]^2)$,
$v^{yy},v^{xy},v^{yy}\geq 0$, then also $u^{yy},u^{xy},u^{yy}\geq
0$.\qed

\noi
In order to prove Lemma~\ref{Lapcon}, based on Lemma~\ref{smoothlap},
we will show that the Laplace equation (\ref{Lap}) has smooth
solutions $u$ for sufficiently many functions $v$. Here `suffiently
many' will mean dense in the topology of uniform convergence of
functions and their derivatives up to second order. To this aim, we
make $\Ci^{(2)}([0,1]^2)$ into a Banach space by equipping it with the
norm
\be\label{2norm}
\|u\|_{(2)}:=\|u\|+\|u^y\|+\|u^x\|+\|u^{yy}\|+2\|u^{xy}\|+\|u^{xx}\|.
\ee
Here, to reduce notation, we denote the supremumnorm by
$\|f\|:=\|f\|_\infty$. Note the factor 2 in the second term
from the right in (\ref{2norm}), which is crucial for the next key lemma.
\bl{\bf(Semigroup on twice diffferentiable functions)}\label{twicegen}
The closure in $\Ci^{(2)}([0,1]^2)$ of the operator $B$ generates a
strongly continuous contraction semigroup on $\Ci^{(2)}([0,1]^2)$.
\el
{\bf Proof} We must check the conditions (i)--(iii) from
(\ref{HilYos1}). It is well-known (see for example
\cite[Proposition~7.1 from the appendix]{EK}) that the space
$\Pc_\infty$ of polynomials is dense in $\Ci^{(2)}([0,1]^2)$.
Therefore $\Di(B)=\Ci^{(\infty)}([0,1]^2)$ is dense, and copying the
proof of Lemma~\ref{Felgen} we see that $\Ri(1-\al B)$ is dense for
all but countably many $\al$. To complete the proof, we must show that
$B$ is dissipative, i.e., that
\be
\|(1-\eps B)u\|_{(2)}\geq\|u\|_{(2)}\qquad
(\eps>0,\ u\in\Ci^{(\infty)}([0,1]^2)).
\ee
Using (\ref{Bcommut}), we calculate
\bc\label{doordif}
\dis\dif{y}(1-\eps B)u&=&\dis(1-\eps(B'-c))u^y,\\[5pt]
\dis\dif{x}(1-\eps B)u&=&\dis(1-\eps B)u^x-\eps cu^y,\\[5pt]
\dis\diff{y}(1-\eps B)u&=&\dis(1-\eps(B''-2c-2))u^{yy},\\[5pt]
\dis\difif{x}{y}(1-\eps B)u&=&\dis(1-\eps(B'-c))u^{xy}-\eps cu^{yy},\\[5pt]
\dis\diff{x}(1-\eps B)u&=&\dis(1-\eps B)u^{xx}-2\eps cu^{xy}.
\ec
Using the disipativity of $B,B'$, and $B''$ with respect to the
supremumnorm (which follows from the positive maximum principle) we
see that $\|(1-\eps(B'-c))u^y\|=(1+\eps c)\|(1-\frac{\eps}{1+\eps c}B)u^y\|
\geq(1+\eps c)\|u^y\|$ etc. We conclude therefore from
(\ref{doordif}) that
\bc
\dis\|(1-\eps B)u\|_{(2)}&\geq &\dis\|(1-\eps B)u\|
+\|(1-\eps(B'-c))u^y\|+\|(1-\eps B)u^x\|-\eps c\|u^y\|\\[5pt]
&&\dis+\|(1-\eps(B''-2c-2))u^{yy}\|+2\|(1-\eps(B'-c))u^{xy}\|
-2\eps c\|u^{yy}\|\\[5pt]
&&\dis+\|(1-\eps B)u^{xx}\|-2\eps c\|u^{xy}\|\\[5pt]
&\geq&\dis\|u\|+(1+\eps c)\|u^y\|+\|u^x\|-\eps c\|u^y\|\\[5pt]
&&\dis+(1+\eps(2c+2))\|u^{yy}\|+2(1+\eps c)\|u^{xy}\|-2\eps c\|u^{yy}\|\\[5pt]
&&\dis+\|u^{xx}\|-2\eps c\|u^{xy}\|\geq\|u\|_{(2)}
\ec
for each $\eps>0$, which shows that $B$ is dissipative with respect to
the norm $\|\cdot\|_{(2)}$.\qed

\noi
{\bf Proof of Lemma~\ref{Lapcon}} Let $g\in\Ci^{(2)}([0,1]^2)$. Then
$u\mapsto gu$ is a bounded operator on both $\Ci([0,1]^2)$ and
$\Ci^{(2)}([0,1]^2)$, so we can choose a $\la>0$ such that
\be
\|gu\|\leq\la\|u\|\quad\mbox{and}\quad\|gu\|_{(2)}\leq\la\|u\|_{(2)}
\ee
for all $u$ in $\Ci([0,1]^2)$ and $\Ci^{(2)}([0,1]^2)$, respectively.
Put $\ti g:=g-\la$. By Lemma~\ref{Felgen}, $\ov B+\ti g$ generates a
strongly continuous contraction semigroup $(S^{\ti g}_t)_{t\geq 0}
=(e^{-\la t}S^g_t)_{t\geq 0}$ on $\Ci([0,1]^2)$. Note that
$\Ri(1-\al(B+\ti g))$ is the space of all $v\in\Ci([0,1]^2)$ for which
the Laplace equation $(1-\al(B+\ti g))u=v$ has a solution
$u\in\Ci^{(\infty)}([0,1]^2)$. Therefore, by Lemma~\ref{smoothlap},
for each $\al>0$:
\be\ba{rl}\label{smoothpres}
{\rm (i)}&\mbox{If $g\in\Ci_{1+}$, then $(1-\al(\ov B+\ti g))^{-1}$
maps $\Ri(1-\al(B+\ti g))\cap\Ci_+\cap\Ci_{1+}$ into
$\Ci_+\cap\Ci_{1+}$.}\\[5pt]
{\rm (ii)}&\mbox{If $g\in\Ci_{1+}\cap\Ci_{2+}$, then
$(1-\al(\ov B+\ti g))^{-1}$ maps
$\Ri(1-\al(B+\ti g))\cap\Ci_+\cap\Ci_{1+}\cap\Ci_{2+}$}\\
&\mbox{into $\Ci_+\cap\Ci_{1+}\cap\Ci_{2+}$.}
\ec
By Lemma~\ref{twicegen}, the restriction of the semigroup
$(S^{\ti g}_t)_{t\geq 0}$ to $\Ci^{(2)}([0,1]^2)$ is strongly continuous and
contractive in the norm $\|\cdot\|_{(2)}$. Therefore, by Hille-Yosida
(\ref{HilYos1}), $\Ri(1-\al(B+\ti g))$ is dense in
$\Ci^{(2)}([0,1]^2)$ for each $\al>0$. It follows that
$\Ri(1-\al(B+\ti g))\cap\Ci_+\cap\Ci_{1+}$ is dense in
$\Ci_+\cap\Ci_{1+}$ and likewise $\Ri(1-\al(B+\ti
g))\cap\Ci_+\cap\Ci_{1+}\cap\Ci_{2+}$ is dense in
$\Ci_+\cap\Ci_{1+}\cap\Ci_{2+}$, both in the norm $\|\cdot\|_{(2)}$.
Note that we need density in the norm $\|\cdot\|_{(2)}$ here: if we
would only know that $\Ri(1-\al(B+\ti g))$ is a dense subset of
$\Ci([0,1]^2)$ in the norm $\|\cdot\|$, then $\Ri(1-\al(B+\ti
g))\cap\Ci_+\cap\Ci_{1+}$ might be empty. By approximation in the norm
$\|\cdot\|_{(2)}$ it follows from (\ref{smoothpres}) that:
\be\ba{rl}\label{C2pres}
{\rm (i)}&\mbox{If $g\in\Ci_{1+}$, then $(1-\al(\ov B+\ti g))^{-1}$
maps $\Ci_+\cap\Ci_{1+}$ into itself.}\\[5pt]
{\rm (ii)}&\mbox{If $g\in\Ci_{1+}\cap\Ci_{2+}$, then
$(1-\al(\ov B+\ti g))^{-1}$ maps $\Ci_+\cap\Ci_{1+}\cap\Ci_{2+}$ into itself.}
\ec
Using also continuity in the norm $\|\cdot\|$ we find that:
\be\ba{rl}\label{ovCpres}
{\rm (i)}&\mbox{If $g\in\Ci_{1+}$, then $(1-\al(\ov B+\ti g))^{-1}$
maps $\ov{\Ci_+\cap\Ci_{1+}}$ into itself.}\\[5pt]
{\rm (ii)}&\mbox{If $g\in\Ci_{1+}\cap\Ci_{2+}$, then
$(1-\al(\ov B+\ti g))^{-1}$ maps $\ov{\Ci_+\cap\Ci_{1+}\cap\Ci_{2+}}$
into itself.}
\ec
For $\eps>0$ let
\be
G_\eps:=\eps^{-1}\big((1-\eps(\ov B+\ti g))^{-1}-1\big)
\ee
be the Yosida approximation to $\ov B+\ti g$. Then
\be
e^{G_\eps t}=e^{-\eps^{-1}t}\sum_{n=0}^\infty\frac{t^n}{n!}
(1-\eps(\ov B+\ti g))^{-n}\qquad(t\geq 0),
\ee
and therefore, by (\ref{ovCpres}), for each $t\geq 0$:
\be\ba{rl}\label{epGep}
{\rm (i)}&\mbox{If $g\in\Ci_{1+}$, then $e^{G_\eps t}$ maps
$\ov{\Ci_+\cap\Ci_{1+}}$ into itself.}\\[5pt]
{\rm (ii)}&\mbox{If $g\in\Ci_{1+}\cap\Ci_{2+}$, then
$e^{G_\eps t}$ maps $\ov{\Ci_+\cap\Ci_{1+}\cap\Ci_{2+}}$ into itself.}
\ec
Finally
\be
e^{-\la t}S^g_tu=S^{\ti g}_tu=\lim_{\eps\to 0}\ex{G_\eps t}u
\qquad(t\geq 0,\ u\in\Ci([0,1]^2)),
\ee
so (\ref{epGep}) implies that for each $t\geq 0$:
\be\ba{rl}
{\rm (i)}&\mbox{If $g\in\Ci_{1+}$, then $S^g_t$ maps
$\ov{\Ci_+\cap\Ci_{1+}}$ into itself.}\\[5pt]
{\rm (ii)}&\mbox{If $g\in\Ci_{1+}\cap\Ci_{2+}$, then
$S^g_t$ maps $\ov{\Ci_+\cap\Ci_{1+}\cap\Ci_{2+}}$ into itself.}
\ec
Using the continuity of $S^g_t$ in $g$ (which follows from
Feynman-Kac (\ref{FeyKac})) we arrive at the statements in
Lemma~\ref{Lapcon}.\qed

\section{Convergence to a time-homogeneous process}\label{convsec}

\subsection{Convergence of certain Markov chains}\label{Markcon}

Section~\ref{convsec} is devoted to the proof of
Theorem~\ref{supercon}. In the present subsection, we start by
formulating a theorem about the convergence of certain Markov chains
to con\-tin\-uous-time processes. In Section~\ref{Brasuper} we
specialize to Poisson-cluster branching processes and superprocesses.
In Section~\ref{superproof}, finally, we carry out the necessary
calculations for the specific processes from Theorem~\ref{supercon}.

Let $E$ be a {\em compact} metrizable space. We equip the space
$\Ci(E)$ of continuous real functions on $E$ with the supremumnorm
$\|\cdot\|_\infty$. By definition, $\Di_E\half$ is the space of cadlag
functions $w:\half\to E$, equipped with the Skorohod topology. Let
$A:\Di(A)\to\Ci(E)$ be an operator defined on a domain
$\Di(A)\sub\Ci(E)$. We say that a process $\y=(\y_t)_{t\geq 0}$ solves
the martingale problem for $A$ if $\y$ has sample paths in
$\Di_E\half$ and for each $f\in\Di(A)$, the process $(M^f_t)_{t\geq
0}$ given by
\be
M^f_t:=f(\y_t)-\int_0^tAf(\y_s)\di s\qquad(t\geq 0)
\ee
is a martingale with respect to the filtration generated by $\y$. We
say that existence (uniqueness) holds for the martingale problem for
$A$ if for each probability measure $\mu$ on $E$ there is at least one
(at most one (in law)) solution $\y$ to the martingale problem for $A$
with initial law $\Li(\y_0)=\mu$. If both existence and uniqueness
hold we say that the martingale problem is well-posed. For each $n\geq
0$, let $X^{(n)}=(X^{(n)}_0,\ldots,X^{(n)}_{m(n)})$ (with $1\leq
m(n)<\infty$) be a (time-inhomogeneous) Markov process in $E$ with
$k$-th step transition probabilities
\be
P_k(x,\di y)=P\big[X^{(n)}_k\in\di y\big|X^{(n)}_{k-1}=x\big]
\qquad(1\leq k\leq m(n)).
\ee
We assume that the $P_k$ are continuous probability kernels on $E$.
Let $(\eps^{(n)}_k)_{1\leq k\leq m(n)}$ be positive constants. Set
\be
A^{(n)}_kf(x):=(\eps^{(n)}_k)^{-1}\Big(\int_EP_k(x,\di y)f(y)-f(x)\Big)
\qquad(1\leq k\leq m(n),\ f\in\Ci(E)).
\ee
Define $t^{(n)}_0:=0$ and
\be\label{tdef}
t^{(n)}_k:=\sum_{l=1}^k\eps^{(n)}_l\qquad(1\leq k\leq m(n)),
\ee
and put
\be\label{kdef}
k^{(n)}(t):=\max\big\{k\;:\;0\leq k\leq m(n),\ t^{(n)}_k\leq t\big\}
\qquad(t\geq 0).
\ee
Define processes $\y^{(n)}=(\y^{(n)}_t)_{t\geq 0}$ with sample paths
in $\Di_E\half$ by
\be
\y^{(n)}_t:=X^{(n)}_{k^{(n)}(t)}\qquad(t\geq 0).
\ee
By definition, a space $\Ai$ of real functions is called an algebra
if $\Ai$ is a linear space and $f,g\in\Ai$ implies $fg\in\Ai$.
\bt{\bf(Convergence of Markov chains)}\label{genth}
Assume that $\Li(X^{(n)}_0)\Rightarrow\mu$ as $n\to\infty$ for some
probability law $\mu$ on $E$. Suppose that there exists at most one
(in law) solution to the martingale problem for $A$ with initial law
$\mu$. Assume that the linear span of $\Di(A)$ contains an algebra
that separates points. Assume that
\be\label{epscon}
{\rm (i)}\ \lim_{n\to\infty}\sum_{k=1}^{m(n)}\eps^{(n)}_k=\infty,
\qquad{\rm (ii)}\ \lim_{n\to\infty}\;\sup_{k:\ t^{(n)}_k\leq T}\eps^{(n)}_k=0,
\ee
and
\be\label{Acon}
\lim_{n\to\infty}\;\sup_{k:\ t^{(n)}_k\leq T}\big\|A^{(n)}_kf-Af\|_\infty=0
\qquad(f\in\Di(A))
\ee
for each $T>0$. Then there exists a unique solution $\y$ to the martingale
problem for $A$ with initial law $\mu$ and moreover
$\Li(\y^{(n)})\Rightarrow\Li(\y)$, where $\Rightarrow$
denotes weak convergence of probability measures on $\Di_E\half$.
\et
{\bf Proof} We apply \cite[Corollary~4.8.15]{EK}. Fix $f\in\Di(A)$.
We start by observing that
\be
f(X^{(n)}_k)-\sum_{i=1}^k\eps^{(n)}_iA^{(n)}_if(X^{(n)}_{i-1})
\qquad(0\leq k\leq m(n))
\ee
is a martingale with respect to the filtration generated by $X^{(n)}$
and therefore,
\be\label{Ymart}
f(\y^{(n)}_t)\;-\!\!\!\!\sum_{i=1}^{k^{(n)}(t)}
\eps^{(n)}_iA^{(n)}_if(\y^{(n)}_{t^{(n)}_{i-1}})\qquad(t\geq 0)
\ee
is a martingale with respect to the filtration generated by $\y^{(n)}$. Put
\be
\lfloor t\rfloor^{(n)}:=t^{(n)}_{k^{(n)}(t)}\qquad(t\geq 0)
\ee
and set
\be
\phi^{(n)}_t:=A^{(n)}_{k^{(n)}(t)+1}f(\y^{(n)}_{\lfloor t\rfloor^{(n)}})
1_{\{t<t^{(n)}_{m(n)}\}}\qquad(t\geq 0)
\ee
and
\be
\xi^{(n)}_t:=f(\y^{(n)}_t)+\int_{\lfloor t\rfloor^{(n)}}^t\phi^{(n)}_s\di s
\qquad(t\geq 0).
\ee
Then we can rewrite the martingale in (\ref{Ymart}) as
\be
\xi^{(n)}_t-\int_0^t\phi^{(n)}_s\di s.
\ee
By \cite[Corollary~4.8.15]{EK} and the compactness of the state space, it
suffices to check the following conditions on $\phi^{(n)}$ and $\xi^{(n)}$:
\be\ba{rl}\label{EKcondi}
{\rm (i)}&\dis\sup_{n\geq N}\;\sup_{t\leq T}E\big[|\xi^{(n)}_t|\big]
<\infty,\\[5pt]
{\rm (ii)}&\dis\sup_{n\geq N}\;\sup_{t\leq T}E\big[|\phi^{(n)}_t|\big]
<\infty,\\[5pt]
{\rm (iii)}&\dis\lim_{n\to\infty}
E\Big[\big(\xi^{(n)}_T-f(\y^{(n)}_T)\big)\prod_{i=1}^rh_i(\y^{(n)}_{s_i})\Big]
=0,\\[5pt]
{\rm (iv)}&\dis\lim_{n\to\infty}E\Big[\big(\phi^{(n)}_T-Af(\y^{(n)}_T)\big)
\prod_{i=1}^rh_i(\y^{(n)}_{s_i})\Big]=0,\\[5pt]
{\rm (v)}&\dis\lim_{n\to\infty}E\Big[\sup_{t\in\Q\cap[0,T]}\big|\xi^{(n)}_t
-f(\y^{(n)}_t)\big|\Big]=0,\\[5pt]
{\rm (vi)}&\dis\sup_{n\geq N}E\big[\|\phi^{(n)}\|_{p,T}\big]<\infty
\qquad\mbox{for some }p\in(1,\infty],
\ec
for some $N\geq 0$ and for each $T>0$, $r\geq 1$,
$0\leq s_1<\cdots<s_r\leq T$, and $h_1,\ldots,h_r\in\Hi\sub\Ci(E)$.
Here $\Hi$ is separating, i.e., $\int h\di\mu=\int h\di\nu$ for all
$h\in\Hi$ implies $\mu=\nu$ whenever $\mu,\nu$ are probability measures
on $E$. In (vi):
\be
\|g\|_{p,T}:=\Big(\int_0^T|g(t)|^p\di t\Big)^{1/p}\qquad(1\leq p<\infty)
\ee
and $\|g\|_{\infty,T}$ denotes the essential supremum of $g$ over $[0,T]$.

The conditions (\ref{EKcondi})~(i)--(vi) are implied by the stronger conditions
\be\ba{rl}\label{strongEK}
{\rm (i)}&\dis\lim_{n\to\infty}\,\sup_{0\leq t\leq T}\big\|\xi^{(n)}_t
-f(\y^{(n)}_t)\big\|_\infty=0,\\[5pt]
{\rm (ii)}&\dis\lim_{n\to\infty}\,\sup_{0\leq t\leq T}\big\|\phi^{(n)}_t
-Af(\y^{(n)}_t)\big\|_\infty=0,
\ec
where we denote the essential supremumnorm of a real-valued random variable
$X$ by $\|X\|_\infty:=\inf\{K\geq 0:|X|\leq K\ \as\}$. Condition
(\ref{strongEK})~(ii) is implied by (\ref{epscon})~(i) and (\ref{Acon}).
To see that also (\ref{strongEK})~(i) holds, set
\be
M_n:=\sup_{0\leq t\leq T}\big\|\phi^{(n)}_t\big\|_\infty,
\ee
and estimate
\be\label{Mnsup}
\sup_{0\leq t\leq T}\big\|\xi^{(n)}_t-f(\y^{(n)}_t)\big\|_\infty
\leq M_n\sup\{\eps^{(n)}_k\,:\,1\leq k\leq m(n),\ t^{(n)}_k\leq T\}.
\ee
Condition (\ref{strongEK})~(ii) implies that $\limsup_nM_n<\infty$
and therefore the right-hand side of (\ref{Mnsup}) tends to zero by
assumption (\ref{epscon})~(ii).\qed

\subsection{Convergence of certain branching processes}\label{Brasuper}

In this section we apply Theorem~\ref{genth} to certain branching
processes and superprocesses.

Throughout this section, $E$ is a compact metrizable space and
$A:\Di(A)\to\Ci(E)$ is a linear operator on $\Ci(E)$ such that the
closure $\ov A$ of $A$ generates a Feller process $\xi=(\xi_t)_{t\geq
  0}$ in $E$ with Feller semigroup $(P_t)_{t\geq 0}$ given by
$P_tf(x):=E^x[f(\xi_t)]$ ($t\geq 0,\ f\in\Ci(E)$).

Let $\al\in\Ci_+(E)$ and $\bet,f\in\Ci(E)$. By definition, a function
$t\mapsto u_t$ from $\half$ into $\Ci(E)$ is a {\em classical} solution
to the semilinear Cauchy problem
\be\label{gencauchy}
\left\{\ba{r@{\,}c@{\,}l}
\dif{t}u_t&=&\ov Au_t+\bet u_t-\al u_t^2\qquad(t\geq 0),\\[5pt]
u_0&=&f
\ea\right.
\ee
if $t\mapsto u_t$ is continuously differentiable (in $\Ci(E)$),
$u_t\in\Di(\ov A)$ for all $t\geq 0$, and (\ref{gencauchy}) holds. We
say that $u$ is a {\em mild} solution to (\ref{gencauchy}) if
$t\mapsto u_t$ is continuous and
\be\label{mild}
u_t=P_tf+\int_0^tP_{t-s}(\bet u_s-\al u_s^2)\di s\qquad(t\geq 0).
\ee
\bl{\bf(Mild and classical solutions)}\label{L:semlin}
Equation (\ref{gencauchy}) has a unique $\Ci_+(E)$-valued mild
solution $u$ for each $f\in\Ci_+(E)$, and $f>0$ implies that $u_t>0$
for all $t\geq 0$. If moreover $f\in\Di(\ov A)$ then $u$ is a
classical solution. For each $t\geq 0$, $u_t$ depends continuously on
$f\in\Ci_+(E)$.
\el
{\bf Proof} It follows from \cite[Theorems~6.1.2, 6.1.4, and
6.1.5]{Paz83} that for each $f\in\Ci(E)$, (\ref{gencauchy}) has a
unique solution $(u_t)_{0\leq t<T}$ up to an explosion time $T$, and
that this is a classical solution if $f\in\Di(\ov A)$. Moreover, $u_t$
depends continuously on $f$. Using comparison arguments based on the
fact that $\ov A$ satisfies the positive maximum principle (which
follows from Hille-Yosida (\ref{HilYos2})) one easily proves the other
statements; compare \cite[Lemmas~23 and 24]{FStrim}.\qed

\noi
We denote the (mild or classical) solution of (\ref{gencauchy}) by
$\Ui_tf:=u_t$; then $\Ui_t:\Ci_+(E)\to\Ci_+(E)$ are continuous
operators and $\Ui=(\Ui_t)_{t\geq 0}$ is a (nonlinear) semigroup on
$\Ci_+(E)$.

Since $E$ is compact, the spaces $\{\mu\in\Mi(E):\mu(E)\leq M\}$ are
compact for each $M\geq 0$. In particular, $\Mi(E)$ is locally
compact. We denote its one-point compactification by
$\Mi(E)_\infty=\Mi(E)\cup\{\infty\}$. We define functions
$F_f\in\Ci(\Mi(E)_{\infty})$ by $F_f(\infty):=0$ and
\be\label{Ffdef}
F_f(\mu):=\ex{-\li\mu,f\re}\qquad(f\in\Ci_+(E),\ f>0,\ \mu\in\Mi(E)).
\ee
We introduce an operator $\Gi$ with domain
\be
\Di(\Gi):=\{F_f:f\in\Di(A),\ f>0\},
\ee
given by $\Gi F_f(\infty):=0$ and
\be\label{Gidef}
\Gi F_f(\mu):=-\li\mu,Af+\bet f-\al f^2\re\,\ex{-\li\mu,f\re}
\qquad(\mu\in\Mi(E)).
\ee
Note that $\Gi F_f\in\Ci(\Mi(E)_{\infty})$ for all $F_f\in\Di(\Gi)$.
\bp{\bf($(\ov A,\al,\bet)$-superprocesses)}\label{superdef}
The martingale problem for the operator $\Gi$ is well-posed. The
solutions to this martingale problem define a Feller process
$\Yi=(\Yi_t)_{t\geq 0}$ in $\Mi(E)_\infty$ with continuous sample
paths, called the $(\ov A,\al,\bet)$-superprocess. If $\Yi_0=\infty$
then $\Yi_t=\infty$ for all $t\geq 0$. If $\Yi_0=\mu\in\Mi(E)$ then
\be\label{logLa}
E^\mu\big[\ex{-\li\Yi_t,f\re}\big]=\ex{-\li\mu,\Ui_tf\re}\qquad(f\in\Ci_+(E)).
\ee
\ep
{\bf Proof} Results of this type are well-known, see for example
\cite[Theorem~9.4.3]{EK}, \cite{Fit88}, and
\cite[Th\'eor\`eme~7]{ER91}. Since, however, it is not completely
straightforward to derive the proposition above from these references,
we give a concise autonomous proof of most of our statements. Only for
the continuity of sample paths we refer the reader to
\cite[Corollary~(4.7)]{Fit88} or \cite[Corollaire~9]{ER91}.

We are going to extend $\Gi$ to an operator $\hat\Gi$ that is linear
and satisfies the conditions of the Hille-Yosida Theorem
(\ref{HilYos2}). For any $\ga\in\Ci_+(E)$ and $\mu\in\Mi(E)$, let
${\rm Clust}_\ga(\mu)$ denote a random measure such that on
$\{\ga=0\}$, ${\rm Clust}_\ga(\mu)$ is equal to $\mu$, and on
$\{\ga>0\}$, ${\rm Clust}_\ga(\mu)$ is a Poisson cluster measure with
intensity $\frac{1}{\ga}\mu$ and cluster mechanism
$\Qi(x,\cdot)=\Li(\tau_{\ga(x)}\de_x)$, where $\tau_{\ga(x)}$ is
exponentially distributed with mean $\ga(x)$. It is not hard to see
that
\be
E\big[\ex{-\li{\rm Clust}_\ga(\mu),f\re}\big]=\ex{-\li\mu,\Vi_\ga f\re}
\qquad(f\in\Ci(E),\ f>0),
\ee
where $\Vi_\ga f(x):=(\frac{1}{f(x)}+\ga(x))^{-1}$.
\detail{Indeed, by (\ref{randclust}),
\[\ba{l}
\dis\Vi_\ga f(x)=\ga(x)^{-1}E\big[1-\ex{-\li\tau_{\ga(x)}\de_x,f\re}\big]
=\ga(x)^{-1}\Big(1-\int_0^\infty e^{-f(x)s}\ffrac{1}{\ga(x)}
e^{s/\ga(x)}\di s\Big)\\[5pt]
\quad\dis=\ga(x)^{-1}\Big(1-\ffrac{1}{\ga(x)}\int_0^\infty 
e^{-(f(x)+1/\ga(x))s}\di s\Big)=\ga(x)^{-1}
\Big(1-\frac{\frac{1}{\ga(x)}}{f(x)+\frac{1}{\ga(x)}}\Big)\\[5pt]
\quad\dis=\frac{f(x)\ga(x)^{-1}}{f(x)+\ga(x)^{-1}}.
\ea\]
}
Note that since $\Vi_\ga 1$ is bounded, the previously mentioned
Poisson cluster measure mentioned above is well-defined. By definition,
we put ${\rm Clust}_\ga(\infty):=\infty$.

Define a linear operator $\Gi_\al$ on $\Ci(\Mi(E))_\infty)$ by
\be\label{G1}
\Gi_\al F(\mu):=\lim_{\eps\to 0}\eps^{-1}
\big(E[F({\rm Clust}_{\eps\al}(\mu))]-F(\mu)\big)
\ee
with as domain $\Di(\Gi_\al)$ the space of all $F\in\Ci(\Mi(E)_\infty)$
for which the limit exists. Define a linear operator $\Gi_\bet$ by
\be\label{G2}
\Gi_\bet F(\mu):=\lim_{\eps\to 0}\eps^{-1}\big(F((1+\eps\bet)\mu)-F(\mu)\big)
\ee
with domain $\Di(\Gi_\bet):=\Ci(\Mi(E))_\infty)$. Define
$P^\ast_t:\Mi(E)_\infty\to\Mi(E)_\infty$ by
$\li P^\ast_t\mu,f\re:=\li\mu,P_tf\re$ $(t\geq 0,\ f\in\Ci(E),\ \mu\in\Mi(E))$
and $P^\ast_t\infty:=\infty$ ($t\geq 0$). Finally, let $\Gi_{\ov A}$
be the linear operator on $\Ci(\Mi(E))_\infty)$ defined by
\be\label{G3}
\Gi_{\ov A}F(\mu):=\lim_{\eps\to 0}\eps^{-1}\big(F(P^\ast_\eps\mu)-F(\mu)\big),
\ee
with as domain $\Di(\Gi_{\ov A})$ the space of all $F$ for which the
limit exists. Define an operator $\hat\Gi$ by
\be
\hat\Gi:=\Gi_\al+\Gi_\bet+\Gi_{\ov A},
\ee
with domain $\Di(\hat\Gi):=\Di(\Gi_\al)\cap\Di(\Gi_{\ov A})$.
If $f\in\Di(\ov A)$, $f>0$, and $F_f$ is as in (\ref{Ffdef}),
then it is not hard to see that $\hat\Gi F_f(\infty)=0$ and
\be
\hat\Gi F_f(\mu):=-\li\mu,\ov Af+\bet f-\al f^2\re\,\ex{-\li\mu,f\re}
\qquad(\mu\in\Mi(E)).
\ee
\detail{$\dif{t}\Vi_{t\al}f(x)=\dif{t}(\frac{1}{f(x)}+t\al(x))^{-1}
=-\al(x)(\frac{1}{f(x)}+t\al(x))^{-2}=-\al(x)\Vi_{t\al}f(x)$.}
In particular, $\hat\Gi$ extends the operator $\Gi$ from
(\ref{Gidef}). Since $\Di(\ov A)$ is dense in $\Ci(E)$, it is easy to
see that $\{F_f:f\in\Di(\ov A),\ f>0\}$ is dense in
$\Ci(\Mi(E)_\infty)$. Hence $\Di(\hat\Gi)$ is dense. Using
(\ref{G1})--(\ref{G3}) it is not hard to show that $\hat\Gi$ satisfies
the positive maximum principle. Moreover, by Lemma~\ref{L:semlin}, for
$f\in\Di(\ov A)$ with $f>0$, the function $t\mapsto F_{\Ui_tf}$ from
$\half$ into $\Ci(\Mi(E)_\infty)$ is continuously differentiable,
satisfies $F_{\Ui_tf}\in\Di(\hat\Gi)$ for all $t\geq 0$, and
\be\label{FC}
\dif{t}F_{\Ui_tf}=\hat\Gi F_{\Ui_tf}\qquad(t\geq 0).
\ee
{F}rom this it is not hard to see that $\hat\Gi$ also satisfies
condition~(\ref{HilYos2})~(ii), so the closure of $\hat\Gi$ generates
a Feller semigroup $(S_t)_{t\geq 0}$ on $\Ci(\Mi(E)_\infty)$. It is
easy to see that $S_tF_f=F_{\Ui_tf}$ $(t\geq 0)$. By
\cite[Theorem~4.2.7]{EK}, this semigroup corresponds to a Feller
process $\Yi$ with cadlag sample paths in $\Mi(E)_\infty$. This means
that $E^\mu[F_f(\Yi_t)]=F_{\Ui_tf}(\mu)$ for all $f\in\Di(\ov A)$ with
$f>0$. If $\mu=\infty$ this shows that $\Yi_t=\infty$ for all $t\geq
0$. If $\mu\in\Mi(E)$ we obtain (\ref{logLa}) for $f\in\Di(\ov A)$,
$f>0$; the general case follows by approximation.\qed

\noi
Now let $(q_\eps)_{\eps>0}$ be continuous weight functions and let
$(\Qi_\eps)_{\eps>0}$ be continuous cluster mechanisms on $E$. Assume that
\be
Z_\eps(x):=\int\Qi_\eps(x,\di\chi)\li\chi,1\re<\infty\qquad(x\in E)
\ee
and define probability kernels $K_\eps$ on $E$ by
\be
\int K_\eps(x,\di y)f(y):=\frac{1}{Z_\eps(x)}
\int\Qi_\eps(x,\di\chi)\li\chi,f\re\qquad(f\in B(E)).
\ee
For each $n\geq 0$, let $(\eps^{(n)}_k)_{1\leq k\leq m(n)}$ (with
$1\leq m(n)<\infty$) be positive constants. Let
$\Xc^{(n)}=(\Xc^{(n)}_0,\ldots,\Xc^{(n)}_{m(n)})$ be a Poisson-cluster
branching process with weight functions
$q_{\eps^{(n)}_1},\ldots,q_{\eps^{(n)}_{m(n)}}$ and cluster mechanisms
$\Qi_{\eps^{(n)}_1},\ldots,\Qi_{\eps^{(n)}_{m(n)}}$. Define
$t^{(n)}_k$ and $k^{(n)}(t)$ as in (\ref{tdef})--(\ref{kdef}). Define
processes $\Yi^{(n)}$ by
\be
\Yi^{(n)}_t:=\Xc^{(n)}_{k^{(n)}(t)}\qquad(t\geq 0).
\ee
\bt{\hspace{1pt}\bf(Convergence of Poisson-cluster branching processes)}\label{bratosup}
\hspace{3pt}Assume that $\Li(\Xc^{(n)}_0)\Rightarrow\rho$ as $n\to\infty$ for some
probability law $\rho$ on $\Mi(E)$. Suppose that the constants
$\eps^{(n)}_k$ fulfill (\ref{epscon}). Assume that
\be\ba{rr@{\,}c@{\,}l}\label{Qlim}
{\rm (i)}&\dis q_\eps(x)\int\Qi_\eps(x,\di\chi)\li\chi,1\re&
=&\dis 1+\eps\bet(x)+o(\eps),\\[5pt]
{\rm (ii)}&\dis q_\eps(x)\int\Qi_\eps(x,\di\chi)\li\chi,1\re^2&
=&\dis\eps\,2\al(x)+o(\eps),\\[5pt]
{\rm (iii)}&\dis q_\eps(x)\int
\Qi_\eps(x,\di\chi)\li\chi,1\re^21_{\{\li\chi,1\re>\de\}}&=&o(\eps)
\ec
for each $\de>0$, and
\be\label{Klim}
\int K_\eps(x,\di y)f(y)=f(x)+\eps Af(x)+o(\eps)
\ee
for each $f\in\Di(A)$, uniformly in $x$ as $\eps\to 0$.
Then $\Li(\Yi^{(n)})\Rightarrow\Li(\Yi)$, where $\Yi$ is the
$(\ov A,\al,\bet)$-superprocess with initial law $\rho$.
\et
Here $\Rightarrow$ denotes weak convergence of probability
measures on $\Di_{\Mi(E)}\half$.\med

\noi
{\bf Proof} We apply Theorem~\ref{genth} to the operator $\Gi$, where
we use the fact that if we view $\Mi_1(\Di_{\Mi(E)}\half)$ as a
subspace of $\Mi_1(\Di_{\Mi(E)_\infty}\half)$ (note the
compactification), equipped with the topology of weak convergence,
then the induced topology on $\Mi_1(\Di_{\Mi(E)}\half)$ is again the
topology of weak convergence.

By Proposition~\ref{superdef}, solutions to the martingale problem for
$\Gi$ are unique. Since $F_fF_g=F_{f+g}$ and $\Di(A)$ is a linear
space, the linear span of the domain of $\Gi$ is an algebra. Using the
fact that $\Di(A)$ is dense in $\Ci(E)$ we see that this algebra
separates points. Therefore, we are left with the task to check
(\ref{Acon}).

Define $\Ui_\eps:\Ci_+(E)\to\Ci_+(E)$ by
\be
\Ui_\eps f(x):=q_\eps(x)\int\Qi_\eps(x,\di\chi)
\big(1-\ex{-\li\chi,f\re}\big)\qquad(x\in E,\ f\in\Ci_+[0,1],\ f>0,\ \eps>0),
\ee
and define transition probabilities $P_\eps(\mu,\di\nu)$ on
$\Mi(E)_\infty$ by $P_\eps(\infty,\,\cdot\,):=\de_\infty$ and
\be\label{pepstr}
\int P_\eps(\mu,\di\nu)\ex{-\li\nu,f\re}=\ex{-\li\mu,\Ui_\eps f\re}.
\ee
We will show that 
\be\label{Uicon}
\lim_{\eps\to 0}\big\|\eps^{-1}(\Ui_\eps f-f)
-(Af+\bet f-\al f^2)\big\|_\infty=0\qquad(f\in\Di(A),\ f>0).
\ee
Together with (\ref{pepstr}) this implies that
\be\label{Gicon}
\int P_\eps(\mu,\di\nu)F_f(\nu)=F_f(\mu)+\eps\Gi F_f(\mu)+o(\eps)
\qquad(f\in\Di(A),\ f>0),
\ee
uniformly in $\mu\in\Mi(E)_\infty$ as $\eps\to 0$. Therefore, the
result follows from Theorem~\ref{genth}.

It remains to prove (\ref{Uicon}). Set $g(z):=1-z+\frac{1}{2}z^2-e^{-z}$
$(z\geq 0)$ and write
\be\label{123}
\Ui_\eps f(x)=q_\eps(x)\int\Qi_\eps(x,\di\chi)\big(\li\chi,f\re
-\ffrac{1}{2}\li\chi,f\re^2+g(\li\chi,f\re)\big).
\ee
Since
\be
g(z)=\int_0^z\di y\int_0^y\di x\int_0^x\di t\,e^{-t}\qquad(z\geq 0),
\ee
it is easy to see that $g$ is nondecreasing on $\half$ and
(since $0\leq e^{-t}\leq 1$ and $\int_0^x\di t\,e^{-t}\leq 1$)
\be
0\leq g(z)\leq\ffrac{1}{2}z^2\wedge\ffrac{1}{6}z^3\qquad(z\geq 0).
\ee
Using these facts and (\ref{Qlim})~(ii) and (iii), we find that
\be\ba{l}
\dis q_\eps(x)\int\Qi_\eps(x,\di\chi)g(\li\chi,f\re)\\[5pt]
\dis\quad\leq\|f\|_\infty q_\eps(x)\Big\{\int\Qi_\eps(x,\di\chi)
g(\li\chi,1\re)1_{\{\li\chi,1\re\leq\de\}}+\int\Qi_\eps(x,\di\chi)
g(\li\chi,1\re)1_{\{\li\chi,1\re>\de\}}\Big\}\\[5pt]
\dis\qquad\leq\|f\|_\infty q_\eps(x)\Big\{\ffrac{1}{6}\de
\int\Qi_\eps(x,\di\chi)\li\chi,1\re^21_{\{\li\chi,1\re\leq\de\}}
+\ffrac{1}{2}\int\Qi_\eps(x,\di\chi)
\li\chi,1\re^21_{\{\li\chi,1\re>\de\}}\Big\}\\[10pt]
\dis\qquad=\ffrac{1}{6}\de\|f\|_\infty\big(\eps\,2\al(x)+o(\eps)\big)+o(\eps).
\ec
Since this holds for any $\de>0$, we conclude that
\be\label{restest}
q_\eps(x)\int\Qi_\eps(x,\di\chi)g(\li\chi,f\re)=o(\eps)
\ee
uniformly in $x$ as $\eps\to 0$. By (\ref{Qlim})~(i) and (\ref{Klim}),
\be\ba{l}\label{1est}
\dis q_\eps(x)\int\Qi_\eps(x,\di\chi)\li\chi,f\re=\Big(q_\eps(x)
\int\Qi_\eps(x,\di\chi)\li\chi,1\re\Big)
\Big(\int K_\eps(x,\di y)f(y)\Big)\\[5pt]
\dis\qquad=\big(1+\eps\bet(x)+o(\eps)\big)\big(f(x)
+\eps Af(x)+o(\eps)\big)\\[5pt]
\dis\qquad=f(x)+\eps\bet(x)f(x)+\eps Af(x)+o(\eps).
\ec
Finally, write
\be\ba{l}\label{pre2est}
\dis q_\eps(x)\int\Qi_\eps(x,\di\chi)\li\chi,f\re^2\\[5pt]
\dis\qquad=q_\eps(x)\int\Qi_\eps(x,\di\chi)\big(\li\chi,f(x)\re^2
+2\li\chi,f(x)\re\li\chi,f-f(x)\re+\li\chi,f-f(x)\re^2\big).
\ec
Then, by (\ref{Qlim})~(ii),
\be\label{fkwa}
q_\eps(x)\int\Qi_\eps(x,\di\chi)\li\chi,f(x)\re^2=f(x)^2
\big(\eps\,2\al(x)+o(\eps)\big).
\ee
We will prove that
\be\label{fmin1}
q_\eps(x)\int\Qi_\eps(x,\di\chi)\li\chi,f-f(x)\re^2=o(\eps).
\ee
Then, by H\"older's inequality, (\ref{Qlim})~(ii), and (\ref{fmin1}),
\be\ba{l}\label{f1f}
\dis\big|q_\eps(x)\int\Qi_\eps(x,\di\chi)\li\chi,f
-f(x)\re\li\chi,f(x)\re\big|\\[5pt]
\dis\qquad\leq\Big(q_\eps(x)
\int\Qi_\eps(x,\di\chi)\li\chi,f-f(x)\re^2\Big)^{1/2}
\Big(q_\eps(x)\int\Qi_\eps(x,\di\chi)\li\chi,f(x)\re^2\Big)^{1/2}\\[5pt]
\dis\qquad\leq\big(o(\eps)(2\al(x)\eps+o(\eps))\big)^{1/2}=o(\eps).
\ec
Inserting (\ref{fkwa}), (\ref{fmin1}) and (\ref{f1f}) into (\ref{pre2est})
we find that
\be\label{2est}
q_\eps(x)\int\Qi_\eps(x,\di\chi)\li\chi,f\re^2=\eps\,2\al(x)f(x)^2+o(\eps).
\ee
Inserting (\ref{restest}), (\ref{1est}) and (\ref{2est}) into (\ref{123}),
we arrive at (\ref{Uicon}). We still need to prove (\ref{fmin1}). To this
aim, we estimate, using (\ref{1est}),
\be\ba{l}
\dis q_\eps(x)\int\Qi_\eps(x,\di\chi)\li\chi,f-f(x)\re^2
1_{\{\li\chi,1\re\leq\de\}}\\[5pt]
\dis\qquad\leq\de\|f-f(x)\|_\infty q_\eps(x)
\int\Qi_\eps(x,\di\chi)\li\chi,f-f(x)\re\\[5pt]
\dis\qquad=\de\|f-f(x)\|_\infty\big(\eps Af(x)+o(\eps)\big)
\ec
and, using (\ref{Qlim})~(iii),
\be\ba{l}
\dis q_\eps(x)\int\Qi_\eps(x,\di\chi)\li\chi,f-f(x)\re^2
1_{\{\li\chi,1\re>\de\}}\\[5pt]
\dis\qquad\leq\|f-f(x)\|_\infty q_\eps(x)
\int\Qi_\eps(x,\di\chi)\li\chi,1\re^21_{\{\li\chi,1\re>\de\}}=o(\eps).
\ec
It follows that
\be
q_\eps(x)\int\Qi_\eps(x,\di\chi)\li\chi,f-f(x)\re^2\leq
\de\eps\|f-f(x)\|_\infty Af(x)+o(\eps)
\ee
for any $\de>0$. This implies (\ref{fmin1}) and completes
the proof of (\ref{Uicon}).\qed

\subsection{Application to the renormalization branching
process}\label{superproof}

{\bf Proof of Theorem~\ref{supercon}~(a)} For any
$f_0,\ldots,f_k\in\Ci_+[0,1]$ one has
\be\ba{l}\label{exexpl}
E\big[\ex{-\li\Xc_{-n},f_0\re}\cdots\ex{-\li\Xc_{-n+k},f_k\re}\big]\\[5pt]
\dis\quad=E\big[\ex{-\li\Xc_{-n},f_0\re}\cdots
\ex{-\li\Xc_{-n+k-1},f_{k-1}+\Ui_{\ga_{n-k}}f_k\re}\big]\\[5pt]
\dis\quad=\cdots=\dis\quad E\big[\ex{-\li\Xc_{-n},g_k\re}\big],
\ec
where we define inductively
\be
g_0:=f_k\quad\mbox{and}\quad g_{m+1}:=f_{k-m-1}+\Ui_{\ga_{n-k+m}}g_m.
\ee
By the compactness of $[0,1]$ and Corollary~\ref{Ugacont}, the map
$(\ga,f)\mapsto\Ui_\ga f$ from
$(0,\infty)\times\Ci_+[0,1]$ to $\Ci_+[0,1]$ (equipped with the
supremumnorm) is continuous. Using this fact and
(\ref{exexpl}) we find that
\be
E\big[\ex{-\li\Xc_{-n},f_0\re}\cdots\ex{-\li\Xc_{-n+k},f_k\re}\big]\asto{n}
E\big[\ex{-\li\Yi^{\ga^\ast}_{-n},f_0\re}\cdots
\ex{-\li\Yi^{\ga^\ast}_{-n+k},f_k\re}\big].
\ee
Since $f_1,\ldots,f_k$ are arbitrary, (\ref{disco}) follows.\qed

\noi
{\bf Proof of Theorem~\ref{supercon}~(b)} We apply
Theorem~\ref{bratosup} to the weight functions $q_\ga$ and cluster
mechanisms $\Qi_\ga$ from (\ref{Qqdef}) and to $A_{\rm WF}
=x(1-x)\diff{x}$ with domain $\Di(A_{\rm WF})=\Ci^{(2)}[0,1]$,
and $\al=\bet=1$. It is well-known that $\ov A_{\rm WF}$ generates a
Feller semigroup \cite[Theorem~8.2.8]{EK}. We observe that
\be
\int\Qi_\ga(x,\di\chi)\li\chi,f\re
=E\big[2\int_0^{\tau_\ga}\!f(\y^\ga_x(-t))\big]
=2E[\tau_\ga]E\big[f(\y^\ga_x(0))\big]
=\ga\int\Ga^\ga_x(\di y)f(y),
\ee
where $\Ga^\ga_x$ is the equilibrium law of the process $\y^\ga_x$
from Corollary~\ref{ergo}. It follows from (\ref{WFmoments}) that
\be\ba{rr@{\,}c@{\,}l}
{\rm (i)}&\dis\int\Ga^\ga_x(\di y)(y-x)&=&0,\\[5pt]
{\rm (ii)}&\dis\int\Ga^\ga_x(\di y)(y-x)^2&=&
\dis\frac{\ga x(1-x)}{1+\ga},\\[5pt]
{\rm (iii)}&\dis\int\Ga^\ga_x(\di y)(y-x)^4&=&\dis O(\ga^2),
\ec
uniformly in $x$ as $\ga\to 0$.
\detail{Indeed, writing $(y-x)^4=y^4-4xy^3+6x^2y^2-4x^3y+x^4$, we find that
\[\ba{r@{\,}c@{\,}l}
\dis\int\Ga^\ga_x(\di y)(y-x)^4
&=&\dis\big\{x(x+\ga)(x+2\ga)(x+3\ga)-4x^2(x+\ga)(x+2\ga)(1+3\ga)\\
&&\dis\phantom{\big\{}+6x^3(x+\ga)(1+2\ga)(1+3\ga)
-3x^4(1+\ga)(1+2\ga)(1+3\ga)\big\}\\
&&\dis\phantom{\big\{}\cdot\big((1+\ga)(1+2\ga)(1+3\ga)\big)^{-1}.
\ea\]
Collecting the terms with different powers of $\ga$, we find $0\cdot\ga^0$ and
\[\ba{l}
\dis\ga\big\{6x^3-4(3x^3+3x^4)+6(x^3+5x^4)-3\cdot 6x^4\big\}\\
\dis\quad=\ga\big\{(6-12+6)x^3+(-12+30-24+6)x^4\big\}=0\cdot\ga^1.
\ea\]}
Therefore, for any $\de>0$,
\be\ba{rr@{\,}c@{\,}l}
{\rm (i)}&\dis\int\Ga^\ga_x(\di y)(y-x)&=&0,\\[8pt]
{\rm (ii)}&\dis\int\Ga^\ga_x(\di y)(y-x)^2&=&\dis\ga x(1-x)+o(\ga),\\[8pt]
{\rm (iii)}&\dis\int\Ga^\ga_x(\di y)1_{\{|y-x|>\de\}}&=&\dis o(\ga),
\ec
uniformly in $x$ as $\ga\to 0$. Consequently, a Taylor expansion of
$f$ around $x$ yields
\be
\int\Ga^\ga_x(\di y)f(x)=f(x)+\ga\ffrac{1}{2}x(1-x)\diff{x}f(x)
+o(\ga)\qquad(f\in\Ci^{(2)}[0,1]),
\ee
uniformly in $x$ as $\ga\to 0$. (For details, in particular the
uniformity in $x$, see for example \cite[Proposition~B.1.1]{Swa99}.)
This shows that condition (\ref{Klim}) is satisfied. Moreover,
\be\ba{l}
\dis\int\Qi_\ga(x,\di\chi)\li\chi,1\re=E[2\tau_\ga]=\ga,\\[5pt]
\dis\int\Qi_\ga(x,\di\chi)\li\chi,1\re^2=E[(2\tau_\ga)^2]
=\int_0^\infty z^2\ffrac{1}{\ga}e^{-z/\ga}\di z=2\ga^2,\\[5pt]
\dis\int\Qi_\ga(x,\di\chi)\li\chi,1\re^3=E[(2\tau_\ga)^3]
=\int_0^\infty z^3\ffrac{1}{\ga}e^{-z/\ga}\di z=6\ga^3,
\ec
which, using the fact that $q_\ga=(\frac{1}{\ga}+1)$, gives
\be\ba{l}
\dis q_\ga\int\Qi_\ga(x,\di\chi)\li\chi,1\re=1+\ga,\\[10pt]
\dis q_\ga\int\Qi_\ga(x,\di\chi)\li\chi,1\re^2=2\ga+o(\ga),\\[10pt]
\dis q_\ga\int\Qi_\ga(x,\di\chi)\li\chi,1\re^3=o(\ga).
\ec
This shows that (\ref{Qlim}) is fulfilled. In particular,
\be
q_\ga\int\Qi_\ga(x,\di\chi)\li\chi,1\re^21_{\{\li\chi,1\re>\de\}}
\leq\de^{-1}q_\ga\int\Qi_\ga(x,\di\chi)\li\chi,1\re^3=o(\ga)
\ee
for all $\de>0$.\qed

\section{The super-Wright-Fisher diffusion: introduction}\label{S:WFi}

\newcommand{\fixi}{0}
\newcommand{\fixii}{p_{0,0}}
\newcommand{\fixiii}{p_{1,0}}
\newcommand{\fixiv}{p_{0,1}}
\newcommand{\fixv}{p_{1,1}}

\subsection{Superprocesses and binary splitting particle systems}\label{supsplit}

Let $E$ be a compact metrizable space, $G$ the generator of a Feller process $\xi=(\xi_t)_{t\geq 0}$ in $E$, and $\al\in\Ci_+(E)$, $\bet\in\Ci(E)$. Then, for each $f\in B_+(E)$, the semilinear Cauchy problem in $B_+(E)$
\be\left\{\:\ba{r@{\,}c@{\,}l}\label{cauchy}
\dif{t}u_t&=&Gu_t+\bet u_t-\al u_t^2\qquad(t\geq 0),\\
u_0&=&f,
\ea\right.\ee
has a unique mild solution $u_t=:\Ui_tf$. Moreover, there exists a unique (in law) Markov process $\Ys$ with continuous sample paths in the space $\Mi(E)$ of finite measures on $E$, defined by its Laplace functionals
\be\label{laplace}
E^{\txt\,\mu}[\ex{-\li\Ys_t,f\re}]=\ex{-\li\mu,\Ui_tf\re}\qquad(t\geq 0,\ \mu\in\Mi(E),\ f\in B_+(E)).
\ee
The process $\Ys$ is called the {\em superprocess} in $E$ with {\em underlying motion generator} $G$, {\em activity} $\al$ and {\em growth parameter} $\bet$ (the last two terms are our terminology), or in short the $(G,\al,\bet)$-{\em superprocess}. The operators $(\Ui_t)_{t\geq 0}=\Ui=\Ui(G,\al,\bet)$ form a semigroup, called the {\em log-Laplace semigroup} of $\Ys$. 

The process $\Ys$ can be constructed in several ways and is nowadays standard. We outlined one such construction in Section~\ref{Brasuper}; see also, e.g., \cite{Fit88,Fit91,Fit92}. We can think of $\Ys$ as describing a population where mass flows with generator $G$, and during a time interval $\di t$ a bit of mass $\di m$ at position $x$ produces offspring with mean $(1+\bet(x)\di t)\di m$ and finite variance $2\al(x)\di t\,\di m$. For basic facts on superprocesses we refer to \cite{Daw93,Eth00,Dyn02}.

Similarly, when $G$ is (again) the generator of a Feller process on a compact metrizable space $E$ and $\gal\in\Ci_+(E)$, then, for any $f\in B_{[0,1]}(E)$, the semilinear Cauchy problem
\be\left\{\:\ba{r@{\,}c@{\,}l}\label{cauchyga}
\dif{t}u_t&=&Gu_t+\gal u_t(1-u_t)\qquad(t\geq 0),\\
u_0&=&f,
\ea\right.\ee
has a unique mild solution $u_t=:U_tf$ in $B_{[0,1]}(E)$. Moreover, there exists a unique Markov process $\Yp$ with cadlag sample paths in the space $\Ni(E)$ of finite counting measures on $E$, defined by its generating functionals
\be\label{momgen}
E^{\txt\nu}\big[(1-f)^{\txt \Yp_t}\big]=(1-U_tf)^{\txt\nu}\qquad(t\geq 0,\ \nu\in\Ni(E),\ f\in B_{[0,1]}(E)).
\ee
Here if $\nu=\sum_{i=1}^n\de_{x_i}$ is a finite counting measure and $g\in B_{[0,1]}(E)$, then $g^\nu:=\prod_{i=1}^ng(x_i)$. We call $\Yp$ the {\em binary splitting particle system} in $E$ with {\em underlying motion generator} $G$ and {\em splitting rate} $\gal$, or in short the {\em $(G,\gal)$-bin-split-process}. The semigroup $(U_t)_{t\geq 0}=U=U(G,\gal)$ is called the {\em generating semigroup} of $\Yp$. The process $\Yp$ consists of particles that independently move according to the generator $G$, and additionally split with local rate $\gal$ into two new particles, created at the position of the old one.

\subsection{Statement of the problem and motivation}

Let $\ov A$ be the closure in $\Ci[0,1]$ (equipped with the supremum norm) of the operator
\be\label{AdefSWF}
A=\ffrac{1}{2}x(1-x)\diff{x}.
\ee
It is well-known that $\ov A$ is the generator of a Feller process $\xi$ on $[0,1]$, called the (standard) {\em Wright-Fisher diffusion}, see \cite[Theorem~8.2.8]{EK}. We are interested in mild solutions to the Cauchy equation
\be\left\{\:\ba{r@{\,}c@{\,}l}\label{cauchyons}
\dif{t}u_t&=&\ov Au_t+\gal u_t(1-u_t)\qquad(t\geq 0),\\
u_0&=&f,
\ea\right.\ee
where $\gal>0$ is a constant. We wish to find all fixed points of (\ref{cauchyons}) and determine their domains of attraction.

For $f\in B_+[0,1]$, the mild solution of (\ref{cauchyons}) is given by $u_t=\Ui_tf$, where $\Ui=\Ui(\ov A,\gal,\gal)$ is the log-Laplace semigroup of a superprocess $\Ys$ in $[0,1]$ with underlying motion generator $G=\ov A$, and activity and growth parameter both equal to $\gal$. We call $\Ys$ the {\em super-Wright-Fisher diffusion} (with activity and growth parameter $\gal>0$).\footnote{More generally, if $\Zi$ is the $(\ov A,\al',\gal)$-superprocess, with $\al',\gal>0$ constants, then $\frac{\gal}{\al'}\Zi=\Ys$ in law, and therefore this more general case can be reduced to the case $\al'=\gal$.}

Our main interest is in the case $\gal=1$. In this case, we have proved in Theorem~\ref{supercon}~(b) above that a suitably rescaled version of the renormalization branching process converges to $\Ys$. In particular, we will need Proposition~\ref{Ucor} below for $\gal=1$ in our proof of Lemmas~\ref{00lem} and \ref{01lem} (see Propositions~\ref{00hom}~(b) and \ref{01hom}~(b) below). We will generalize a bit and treat general $\gal>0$. This will not be much more work and will give a more complete picture. In particular, we will see that the case $\gal=1$ is a critical case, since $\Ys$ dies out on the interior if and only if $\gal\leq 1$, and the weighted process $\Ys^v$ from (\ref{Xv}) is critical for $\gal=1$.

If $f\in B_{[0,1]}[0,1]$, then the solution of (\ref{cauchyons}) is also given by $u_t=U_tf$, where $U=U(\ov A,\gal)$ is the generating semigroup of a system $\Yp$ of {\em binary splitting Wright-Fisher diffusions}, with splitting rate $\gal$. The process $\Yp$ can be obtained from $\Ys$ by Poissonization with the constant function $1$ (compare Proposition~\ref{Poisprop}). In fact, $\Yp$ is the {\em trimmed tree} of $\Ys$, i.e., the particles in $\Yp$ correspond to those infinitesimal bits of mass in $\Ys$, that have offspring at all later times. For a precise statement of this fact we refer the reader to \cite{FStrim}.

\begin{figure}
\begin{center}
\includegraphics[width=13.5cm,height=4.5cm]{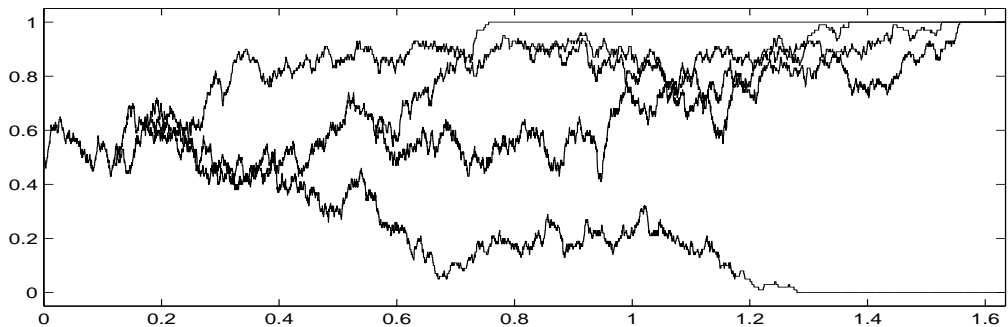}
\end{center}
\caption{A system of binary splitting Wright-Fisher diffusions with splitting rate $\gal=1$.}\label{brapla}
\end{figure}

See Figure~\ref{brapla} for a simulation of $\Yp$ for $\gal=1$. The points $0,1$ are accessible traps for the Wright-Fisher diffusion, and therefore a natural question is whether eventually all particles of $\Yp$ end up in $0$ or $1$. This question will be answered for all $\gal>0$ in Proposition~\ref{fix2} below.

Binary splitting Wright-Fisher diffusions have been studied before in \cite{GKW01}. In particular, the authors of that paper investigated the function $p$, which is defined in terms of the system $\Yp$ of binary splitting Wright-Fisher diffusions with splitting rate $\gal=1$, as
\be\label{interpret}
p(x):=\lim_{t\to\infty}P^{\de_x}[\Yp_t(\{1\})>0]=\lim_{t\to\infty}P^{\de_x}[\Yp_t((0,1])>0]\qquad(x\in[0,1]).
\ee
In order to show that the two expressions for $p$ in (\ref{interpret}) are identical, in \cite{GKW01} the authors note that both expressions correspond to a fixed point $p$ of the generating semigroup $U(\ov A,1)$ with boundary conditions $p(0)=0$ and $p(1)=1$. Assuming that $p$ is sufficiently smooth, the fixed point property means that $p$ solves the equation
\be\label{pdif}
\ffrac{1}{2}x(1-x)\diff{x}p(x)+\gal p(x)(1-p(x))=0\qquad(x\in[0,1]).
\ee
Though stated only for the case $\gal=1$, the proof of Lemma~1.13 in \cite{GKW01} shows that equation (\ref{pdif}) has at most one solution with boundary conditions $p(0)=0$ and $p(1)=1$ when $\gal<z_0^2/8\cong1.836$, where $z_0$ is the smallest non-trivial zero of the Bessel function of the first kind with parameter 1. The authors do not answer the question whether solutions to (\ref{pdif}) with these boundary condions are unique for $\gal\geq z_0^2/8$, or what solutions may exist for other boundary conditions. Proposition~\ref{Ucor} below settles these questions. We show moreover that all fixed points of $U(\ov A,\gal)$ are smooth, a fact tacitly assumed in \cite{GKW01}.

\subsection{Results}

\begin{figure}
\begin{center}
\includegraphics[width=5cm,height=5cm]{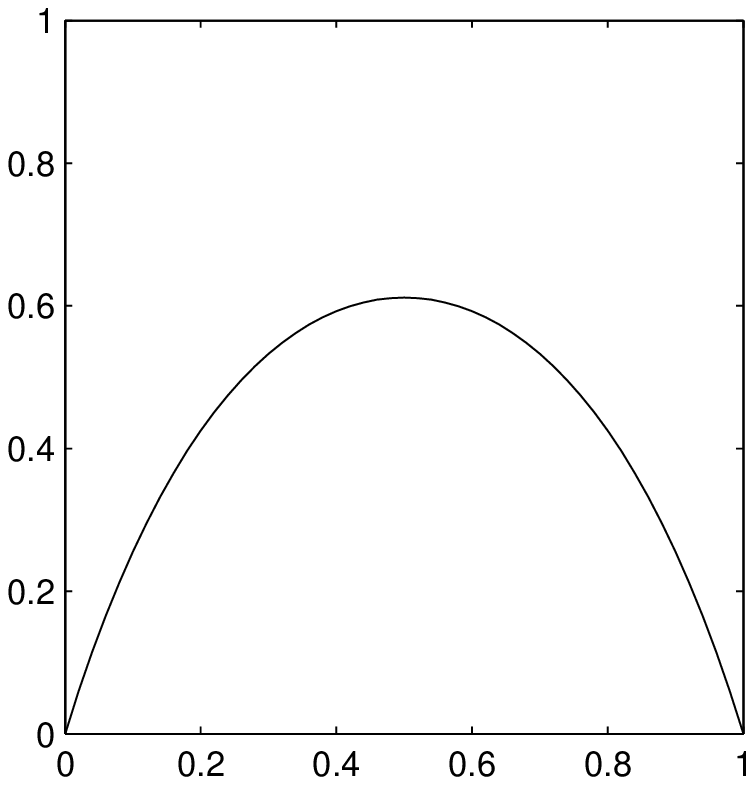}
\hspace{1.2cm}  
\includegraphics[width=5cm,height=5cm]{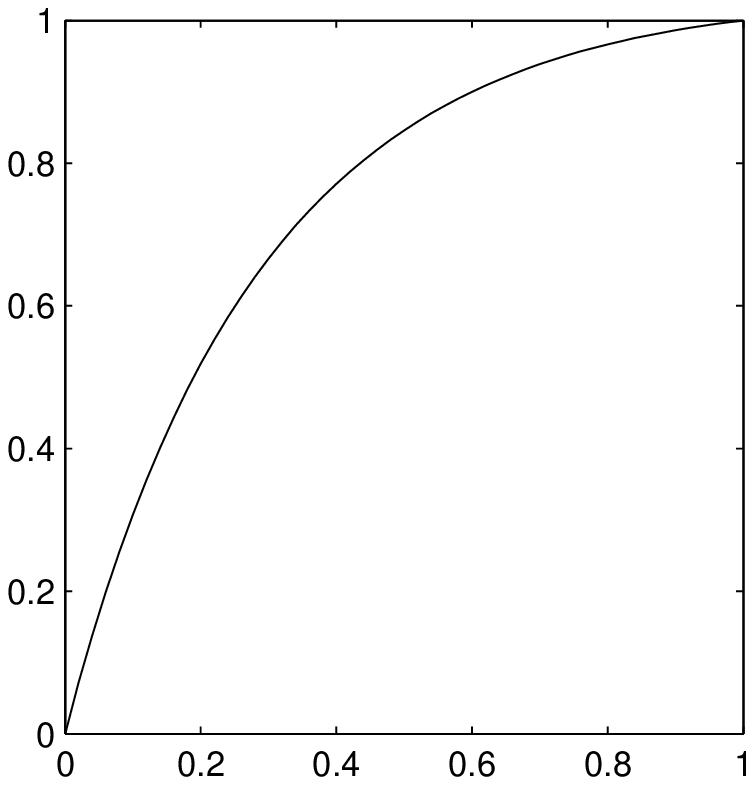}
\end{center}
\caption{Two solutions to the differential equation $\ffrac{1}{2}x(1-x)\diff{x}p(x)+2p(x)(1-p(x))=0$.}\label{ppla}
\setlength{\unitlength}{1cm}
\begin{picture}(0,0)(0,0)
\put(3.3,5.7){$\fixii$}
\put(9.7,5.7){$\fixiv$}
\end{picture}
\end{figure}

The following theorem is our main result. We write `eventually' behind an event, depending on $t$, to denote the existence of a (random) time $\tau<\infty$ such that the event holds for all $t\geq\tau$.
\bt{\bf (Long-time behavior of the super-Wright-Fisher diffusion)}\label{mainSWF}
Let $\Ys$ be the super-Wright-Fisher diffusion with activity and growth parameter equal to the same constant $\gal>0$, started in $\mu\in\Mi[0,1]$. Set
\be\label{vdef}
v(x):=6x(1-x)\qquad(x\in[0,1]).
\ee
Then there exist nonnegative random variables $W_0$, $W_1$, $W_{(0,1)}$ (depending on $\mu$) such that
\be\ba{rlll}\label{WWdef}
{\rm (i)}&\dis\lim_{t\to\infty}e^{-\gal t}\li\Ys_t,1_{\{r\}}\re=W_r\quad&{\rm a.s.}\qquad&(r=0,1),\\[5pt]
{\rm (ii)}&\dis\lim_{t\to\infty}e^{-(\gal-1)t}\li\Ys_t,v\re=W_{(0,1)}\quad&{\rm a.s.}
\ec
and
\be\ba{rlll}\label{Wext}
{\rm (i)}&\dis\{W_r=0\}=\{\Ys_t(\{r\})=0\mbox{ eventually}\}\quad&{\rm a.s.}\qquad&(r=0,1),\\[5pt]
{\rm (ii)}&\dis\{W_{(0,1)}=0\}=\{\Ys_t((0,1))=0\mbox{ eventually}\}\quad&{\rm a.s.}
\ec
Moreover,
\be\label{WWW}
\{W_{(0,1)}>0\}\sub\{W_0>0\}\cap\{W_1>0\}\quad{\rm a.s.}
\ee
If $\gal\leq 1$, then
\be\label{sterf}
W_{(0,1)}=0\qquad{\rm a.s.}
\ee
If $\gal>1$, then $W_{(0,1)}$ satisfies
\be\label{Wmom}
E^\mu(W_{(0,1)})=\li\mu,v\re\quad\mbox{and}\quad\var^\mu(W_{(0,1)})\leq3\ffrac{\gal}{\gal-1}\li\mu,v\re
\ee
as well as
\be\label{explSWF}
\lim_{t\to\infty}E^\mu\Big[\big|e^{-(\gal-1)t}\li\Ys_t,vf\re-W_{(0,1)}\li\ell,vf\re\big|^2\Big]=0\quad\forall f\in B[0,1],
\ee
where $\ell$ denotes the Lebesgue measure on $(0,1)$.
\et
Except for the statement about smoothness (of the functions $\fixii,\ldots,\fixv$ below) and the uniformity of the limit in (\ref{Up5}), the following result about the log-Laplace semigroup $\Ui(\ov A,\gal,\gal)$ is an immediate consequence of Theorem~\ref{mainSWF}.
\bp{\boldmath\bf(Long-time behavior of $\Ui(\ov A,\gal,\gal)$)}\label{Ucor}
Let $\Ys$, $W_0,W_1,W_{(0,1)}$ be as in Theorem~\ref{mainSWF} and let $\Ui=\Ui(\ov A,\gal,\gal)$. Then, for all $f\in B_+[0,1]$, uniformly on $[0,1]$,
\be\label{Up5}
\lim_{t\to\infty}\Ui_t f=\left\{\ba{ll}
\fixi&\mbox{if}\quad f(0)=f(1)=\li\ell,f\re=0,\\
\fixii&\mbox{if}\quad f(0)=f(1)=0,\ \li\ell,f\re>0,\\
\fixiii&\mbox{if}\quad f(0)>0,\ f(1)=0,\\
\fixiv&\mbox{if}\quad f(0)=0,\ f(1)>0,\\
\fixv&\mbox{if}\quad f(0)>0,\ f(1)>0,\ea\right.
\ee
where the constant function $0$ and
\be\left.\ba{r@{\,}c@{\,}l}\label{p1}
\fixii(x)&:=&-\log P^{\de_x}[W_{(0,1)}=0],\\
\fixiii(x)&:=&-\log P^{\de_x}[W_0=0]=P^{\de_x}[W_0=W_{(0,1)}=0],\\
\fixiv(x)&:=&-\log P^{\de_x}[W_1=0]=P^{\de_x}[W_1=W_{(0,1)}=0],\\
\fixv(x)&:=&-\log P^{\de_x}[W_0=W_1=0]=P^{\de_x}[W_0=W_1=W_{(0,1)}=0]
\ea\quad\right\}\quad(x\in[0,1]).
\ee
are all fixed points of the log-Laplace semigroup $\Ui(\ov A,\gal,\gal)$. Here $\fixii=0$ if $\gal\leq 1$, and $\fixii>0$ on $(0,1)$ if $\gal>1$. The functions $p_{l,r}$ $(l,r\in\{0,1\}$ satisfy $p_{l,r}(0)=l$ and $p_{l,r}(1)=r$, are twice continuously differentiable on $[0,1]$, and solve (\ref{pdif}).
\ep
Since conversely, every nonnegative twice continuously differentiable solution to (\ref{pdif}) is a fixed point of $\Ui(\ov A,\gal,\gal)$, we see that (\ref{pdif}) has precisely four solutions when $\gal\leq 1$ and precisely five solutions when $\gal>1$. The functions $\fixii,\ldots,\fixv$ are $[0,1]$-valued and therefore fixed points of the generating semigroup $U(\ov A,\gal)$ as well. Our final result describes $\fixii,\ldots,\fixv$ in terms of the system $\Yp$ of binary splitting Wright-Fisher diffusions with splitting rate $\gal$.
\bp{\bf (Fixed points of $U(\ov A,\gal)$)}\label{fix2}
The functions $\fixii,\ldots,\fixv$ in (\ref{p1}) satisfy
\be\left.\ba{r@{\,}c@{\,}l}\label{p2}
\fixii(x)&=&P^{\de_x}[\Yp_t((0,1))>0\mbox{ eventually}],\\
\fixiii(x)&=&P^{\de_x}[\Yp_t(\{0\})>0\mbox{ eventually}]=P^{\de_x}[\Yp_t([0,1))>0\mbox{ eventually}],\\
\fixiv(x)&=&P^{\de_x}[\Yp_t(\{1\})>0\mbox{ eventually}]=P^{\de_x}[\Yp_t((0,1])>0\mbox{ eventually}],\\
\fixv(x)&=&1
\ea\!\!\right\}\,(x\in[0,1]).
\ee
\ep
See Figure~\ref{ppla} for a plot of the functions $\fixii$ and $\fixiv$ (for $\gal=2$).

\subsection{Methods and related work}

An essential tool in the proof of Theorem~\ref{mainSWF} is the {\em weighted super-Wright-Fisher diffusion} $\Ys^v$, defined as
\be\label{Xv}
\Ys^v_t(\di x):=v(x)\Ys_t(\di x)\qquad(t\geq 0),
\ee
where $v$ is defined in (\ref{vdef}). Note that $v$ is an eigenfunction of the operator $\ov A$, with eigenvalue $-1$. For convenience, we have normalized $v$ such that $\li\ell,v\re=1$.

When a superprocess is weighted with a sufficiently smooth density, the result is a new superprocess, with a new activity and growth parameter and a new underlying motion, which is a compensated {\it h}-transform of the old one. For the case that the underlying motion is a locally uniformly elliptic diffusion on a open domain $D\sub\R^d$, weighted superprocesses were developed by \cite{EP99}. In our case, where uniform ellepticity does not hold, the following can be proved without too much effort.
\bl{\bf\hspace{-1.2pt}(Weighted super-Wright-Fisher diffusion)}\label{weilem}
Let $\Ys$ be the super-Wright-Fisher diffusion with $\gal>0$ and let $\Ys^v$ be defined as in (\ref{Xv}). Then $\Ys^v$ is the $(\ov{A^v},\gal v,\gal-1)$-superprocess in $[0,1]$, where $\ov{A^v}$ is the closure of the operator
\be\label{Avdef}
A^v:=\ffrac{1}{2}x(1-x)\diff{x}+2(\ffrac{1}{2}-x)\dif{x}.
\ee
\el
Indeed, $\ov{A^v}$ generates a Feller process $\xi^v$ in $[0,1]$, see \cite[Theorem~8.2.1]{EK}. The diffusion $\xi^v$ is a compensated {\it h}-transform (with $h=v$) of the Wright-Fisher diffusion $\xi$. This compensated {\it v}-transformed Wright-Fisher diffusion $\xi^v$ is ergodic with invariant law $v\ell$ (Lemma~\ref{vergo} below). For $\gal>1$, the $(\ov{A^v},\gal v,\gal-1)$-superprocess is supercritical, and in this case one expects $e^{-(\gal-1)t}\Ys^v_t$ to converge, in some way, to a random multiple of $v\ell$. This is the idea behind formula (\ref{explSWF}).

Recently, \cite{ET02}, have shown for a certain class of superdiffusions $\Ys$ in $\R^d$ with underlying motion generator $G$, growth parameter $\bet$ and activity $\al$, the convergence in law
\be\label{ETres}
e^{-\la_ct}\li\Ys,g\re\Rightarrow W\li\rho,g\re\quad\mbox{as }t\to\infty,
\ee
where $W$ is a nonnegative random variable, $\la_c$ is the generalized principal eigenvalue of $G+\bet$ (which is assumed to be positive), $\rho$ is a measure on $\R^d$, defined in terms of $G+\bet$, and $g$ is any compactly supported continuous function on $\R^d$. In their work, the weighted superprocess $\Ys^\phi_t(\di x):=\phi(x)\Ys_t(\di x)$ plays a central role, where $\phi$ is the principal eigenfunction of the operator $G+\bet$. Their dynamical system methods are based on a result on the existence of an invariant curve of the log-Laplace semigroup of their superprocess. Using this invariant curve, they give an expression for the Laplace-transform of the law of the random variable $W$ in (\ref{ETres}). Their results are in line with our results for the super-Wright-Fisher diffusion restricted to $(0,1)$, where in our case $\la_c=\gal-1$ and $\phi=v$. However, their methods use in an essential way the fact that their underlying space is $\R^d$ (and not an open subset of $\R^d$, like $(0,1)$), and therefore their results are not applicable to our situation. It is stated as an open problem by \cite{ET02} whether the random variable $W$ in (\ref{ETres}) in general satisfies $P[W=0]=P[\Ys_t=0\mbox{ eventually}]$. For a recent result on {\em local} extinction versus {\em local} exponential growth of superdiffusions on open domains $D\sub\R^d$, we refer to \cite{EK02}.

In our set-up, we can prove that $\{W_{(0,1)}=0\}=\{\Ys_t((0,1))=0\mbox{ eventually}\}$ because of the following property of the weighted super-Wright-Fisher diffusion $\Ys^v$.
\bl{\bf (Finite ancestry)}\label{Wfa}
For all $\gal>0$, the weighted super-Wright-Fisher diffusion $\Ys^v$ satisfies
\be\label{Wfin}
\inf_{x\in[0,1]}P^{\de_x}[\Ys^v_t=0]>0\qquad\forall t>0.
\ee
\el
Formula (\ref{Wfin}) has been called the {\em finite ancestry property} (of $\Ys^v$); for a justification of this terminology we refer the reader to \cite{FStrim}. A sufficient condition for a superprocess to enjoy the finite ancestry property is that the activity be bounded away from zero (see Lemma~\ref{exti} below). This condition is not necessary. In fact, the activity of $\Ys^v$ is $\gal v$, which is zero on $\{0,1\}$. Our proof of Lemma~\ref{Wfa} is quite long. It is not clear whether the weighted superprocesses $\Ys^\phi$ occurring in \cite{ET02} will in general satisfy a formula of the form (\ref{Wfin}). Therefore, we mention as an open problem:
\begin{quote}
How to check, in a practical way, whether a given superprocess has the finite ancestry property (\ref{Wfin})?
\end{quote}
Another problem that is left open in here, is whether the $L_2$-convergence in (\ref{explSWF}) can be replaced by almost sure convergence. In fact, we suspect that (\ref{explSWF}) can be strengthened to
\be
\lim_{t\to\infty}e^{-(\gal-1)t}\li\Ys_t,1_{(0,1)}f\re=W_{(0,1)}\li\ell,f\re\quad\forall f\in B[0,1]\quad{\rm a.s.},
\ee
but we do not have a proof.\vc

\noi
The following sections are organized as follows. Sections~\ref{prepar} and \ref{fasec} contain some general facts about $(G,\al,\bet)$-superprocesses and on $(G,\al,\bet)$-superprocesses enjoying the finite ancestry property, respectively. After some preparatory work in Sections~\ref{smosec} and \ref{absec}, we prove Lemmas~\ref{weilem} and \ref{Wfa} in Section~\ref{trafosec}. In Sections~\ref{ergsec} and \ref{Longwei} we derive some properties of the weighted super-Wright-Fisher diffusion $\Ys^v$, culminating in the proof of Theorem~\ref{mainSWF} in Section~\ref{Longsup}. Finally, Sections~\ref{LongU}--\ref{Smosec} contain the proofs of Propositions~\ref{Ucor} and \ref{fix2}.

\section{The super-Wright-Fisher diffusion: preparatory results}\label{S:WFii}

\subsection{Some general facts about log-Laplace semigroups}\label{prepar}

Let $E$ be a compact metrizable space and let $\Ci(E)$ be the space of continuous real functions on $E$, equipped with the supremum norm $\|\cdot\|_\infty$. Let $\xi=(\xi_t)_{t\geq 0}$ be a Feller process in $E$ with semigroup $S_tf(x):=E^x[f(\xi_t)]$ ($t\geq 0,\ x\in E, f\in B(E)$). By definition, the (full) generator $G$ of $\xi$ is the linear operator on $\Ci(E)$ given by $Gf:=\lim_{t\to 0}t^{-1}(S_tf-f)$ where the domain $\Di(G)$ of $G$ is the space of all functions $f\in\Ci(E)$ for which the limit exists in $\Ci(E)$.

Let $\al\in\Ci_+(E)$, $\bet\in\Ci(E)$, and $f\in\Ci_+(E)$. By definition, we call $u$ a {\em classical solution} of the Cauchy problem (\ref{cauchy}) if $u:\half\to\Ci_+(E)\cap\Di(G)$ is continuously differentiable in $\Ci(E)$ (i.e., the derivative $\dif{t}u_t:=\lim_{s\to t}s^{-1}(u_{t+s}-u_t)$ exists in $\Ci(E)$ for all $t\geq 0$ and the map $\dif{t}u:[0,\infty)\to\Ci(E)$ is continuous) and (\ref{cauchy}) holds. A measurable function $u:\half\times E\to\half$ is called a {\em mild solution} of (\ref{cauchy}) if $u$ is bounded on finite time intervals and solves (pointwise)
\be\label{int}
u_t=S_tf+\int_0^t\!S_{t-s}\big(\bet u_s-\al u_s^2\big)\di s\qquad(t\geq 0).
\ee
Equation (\ref{cauchy}) has a unique mild solution for all $f\in B_+(E)$, see \cite{Fit88} and this solution is a classical solution if $f\in\Ci_+(E)\cap\Di(G)$. (See \cite{Paz83}, Theorems~6.1.4 and 6.1.5. The fact that $f$ is nonnegative and $\al\geq 0$ implies that solutions cannot explode. Our definition of a classical solution is slightly stronger than the one used in \cite{Paz83}, since we require $u$ to be continuously differentiable on $\half$ instead of $(0,\infty)$. However, the proof of Theorem~6.1.5 in \cite{Paz83} shows that $u$ is continuously differentiable on $\half$ if $f\in\Ci_+(E)\cap\Di(G)$.)

The $(G,\al,\bet)$-superprocess $\Ys$ is defined as the unique strong Markov process with continuous sample paths in $\Mi(E)$, equipped with the topology of weak convergence, such that (\ref{laplace}) holds for all $f\in B_+(E)$; see \cite{Fit88,Fit91,Fit92}.

Note the following elementary properties of the log-Laplace semigroup $\Ui(G,\al,\bet)$. Here, we write $\bplim_{n\to\infty}f_n=f$ if $f$ is the bounded pointwise limit of the sequence $(f_n)_{n\geq 0}$.
\bl{\bf (Continuity and monotonicity of log-Laplace semigroups)}\label{Uico}
For each $t\geq 0$, $\Ui_t:\Ci_+(E)\to\Ci_+(E)$ is continuous. Moreover, if $\bplim_{n\to\infty}f_n=f$ for some sequence $f_n\in B_+(E)$, then $\bplim_{n\to\infty}\Ui_t f_n=\Ui_tf$. Finally, $f\leq g$ implies $\Ui_tf\leq\Ui_tg$ $(f,g\in B_+(E))$.
\el
{\bf Proof} The continuity of $\Ui_t:\Ci_+(E)\to\Ci_+(E)$ follows from \cite[Theorem~6.1.2]{Paz83} and the fact that solutions do not explode. Continuity of $\Ui_t$ with respect to bounded pointwise limits is obvious from (\ref{laplace}), and the same formula also makes clear that $\Ui_t:B_+(E)\to B_+(E)$ is monotone.\qed

\noi
Recall that (\ref{cauchy}) has a classical solution for $f\in\Ci_+(E)\cap\Di(G)$. Because of the following, for many purposes it suffices to work with classical solutions.
\bl{\bf (Closure and bp-closure)}\label{bpclos}
For $t\geq 0$ fixed, $\{(f,\Ui_tf):f\in\Ci_+(E)\}$ is the closure in $\Ci(E)$ of $\{(f,\Ui_tf):f\in\Ci_+(E)\cap\Di(G)\}$, and $\{(f,\Ui_tf):f\in B_+(E)\}$ is the bp-closure of $\{(f,\Ui_tf):f\in\Ci_+(E)\}$.
\el
Here, the bp-closure of a set $B$ is the smallest set $\ov B$ such that $B\sub\ov B$ and $f\in\ov B$ whenever $\bplim_{n\to\infty}f_n=f$ for some sequence $f_n\in\ov B$.\vc

\noi
{\bf Proof of Lemma~\ref{bpclos}} It follows from the Hille-Yosida Theorem, see \cite[Theorem~1.2.6]{EK} that $\Di(G)$ is dense in $\Ci(E)$. Since $\Di(G)$ is a linear space and $1\in\Di(G)$, it is not hard to see that $\Ci_+(E)\cap\Di(G)$ is dense in $\Ci_+(E)$. The fact that $\{(f,\Ui_tf):f\in\Ci_+(E)\}$ is the closure in $\Ci(E)$ of $\{(f,\Ui_tf):f\in\Ci_+(E)\cap\Di(G)\}$ now follows from the continuity of $\Ui_t:\Ci_+(E)\to\Ci_+(E)$.

In \cite[Proposition 3.4.2]{EK}, it is proved that $\Ci(E)$ is bp-dense in $B(E)$; the argument can easily be adapted to show that $\Ci_+(E)$ is bp-dense in $B_+(E)$. Therefore Lemma~\ref{bpclos} follows from the continuity of $\Ui_t$ with respect to bounded pointwise limits.\qed

\noi
$\Ui_tf$ may be defined unambiguously such that (\ref{laplace}) holds also for functions $f$ that are not bounded, or even infinite.
\bl{\bf (Extension of $\Ui$ to unbounded functions)}\label{unbound}
For each measurable $f:E\to[0,\infty]$ and $t\geq 0$ there exists a unique measurable $\Ui_tf:E\to[0,\infty]$ such that (\ref{laplace}) holds for all $\mu\in\Mi(E)$, where we put $e^{-\infty}:=0$.
\el
{\bf Proof} Define $\Ui_tf$ by $\Ui_tf(x):=-\log E^{\de_x}[e^{-\li\Ys_t,f\re}]$ where $\log 0:=-\infty$. To see that (\ref{laplace}) holds again for all $\mu\in\Mi(E)$, choose $B_+(E)\ni f_n\up f$, note that $\Ui_t f_n\up\Ui_t f$, and take the limit in (\ref{laplace}).\qed

\noi
We will often need the following comparison result, compare \cite[Theorem~10.1]{Smo83}.
\bl{\bf (Sub- and supersolutions)}\label{subsol}
Assume that $T>0$ and that $\ti u:[0,T]\to\Ci_+(E)\cap\Di(G)$ is continuously differentiable in $\Ci(E)$ and solves
\be
\dif{t}\ti u_t\leq G\ti u_t+\bet \ti u_t-\al\ti u_t^2\qquad(t\in[0,T]).
\ee
Then $\ti u_T\leq\Ui_T\ti u_0$. The same holds with both inequality signs reversed. 
\el
{\bf Proof} Let $g:[0,T]\to\Ci_+(E)$ be defined by the formula
\be
\dif{t}\ti u_t=G\ti u_t+\bet \ti u_t-\al\ti u_t^2-g_t\qquad(t\in[0,T]).
\ee
Set $u_t:=\Ui_t\ti u_0$. Then $u:[0,T]\to\Ci_+(E)$ is the classical solution of
\be\left\{\:\ba{r@{\,}c@{\,}l}\label{cauchyeq}
\dif{t}u_t&=&Gu_t+\bet u_t-\al u_t^2\qquad(t\in[0,T]),\\
u_0&=&\ti u_0.
\ea\right.\ee
Put $\De_t:=u_t-\ti u_t$ ($t\in[0,T]$). Then $\De$ solves
\be\left\{\:\ba{r@{\,}c@{\,}l}\label{difde}
\dif{t}\De_t&=&G\De_t+\bet\De_t-\al\,(u_t+\ti u_t)\De_t+g_t\qquad(t\in[0,T]),\\
\De_0&=&0.
\ea\right.\ee
The generator $G$ satisfies the positive maximum principle, see \cite[Theorem~4.2.2]{EK} and therefore (\ref{difde}) implies that $\De\geq 0$. For imagine that $\De_t(x)<0$ somewhere on $[0,T]\times E$. Let $R$ be a constant such that $\bet-\al\,(u_t+\ti u_t)+R<0$. Then $\ti\De_t:=e^{Rt}\De_t$ solves
\be\left\{\:\ba{r@{\,}c@{\,}l}\label{diftide}
\dif{t}\ti\De_t&=&G\ti\De_t+\{\bet-\al\,(u_t+\ti u_t)+R\}\ti\De_t+g_te^{Rt}\qquad(t\in[0,T]),\\[5pt]
\ti\De_0&=&0.
\ea\right.\ee
If $\ti\De_t(x)<0$ for some $(t,x)\in[0,T]\times E$, then $\ti\De$ must assume a negative minimum over $[0,T]\times E$ in some point $(s,y)$, with $s>0$ since $\ti\De_0=0$. But in such a point one would have $\dif{s}\ti\De_s(y)\leq 0$ while $G\ti\De_s(y)+\{\bet(y)-\al(y)\,(u_s(y)+\ti u_s(y))+R\}\ti\De_s(y)+g_s(y)e^{Rs}>0$, in contradiction with (\ref{diftide}).

The same argument applies when both inequality signs are reversed.\qed

\noi
Lemma~\ref{subsol} has the following application.
\bl{\bf (Bounds on log-Laplace semigroups)}\label{exti}
Let $\Ui=\Ui(G,\al,\bet)$, $\ov\Ui=\Ui(G,\un\al,\ov\bet)$, where $\al,\un\al\in\Ci_+(E)$ and $\bet,\ov\bet\in\Ci(E)$ satisfy
\be
\al\geq\un\al\quad\mbox{and}\quad\bet\leq\ov\bet.
\ee
Then
\be\label{UU}
\Ui_tf\leq\ov\Ui_tf\mbox{ for all measurable }f:E\to[0,\infty]\quad(t\geq 0).
\ee
In particular, if $\un\al,\ov\bet$ are constants and $\un\al>0$, then, for $t>0$,
\be\label{ovU}
\ov\Ui_t\infty=\frac{\ov\bet}{\un\al\,(1-e^{-\ov\bet t})}\quad(\ov\bet\neq 0)\quad\mbox{ and }\quad\ov\Ui_t\infty=\frac{1}{\un\al\, t}\quad(\ov\bet=0),
\ee
and (\ref{UU}) with $f=\infty$ gives
\be\label{sterfkans}
P^\mu[\Ys_t=0]\geq\ex{-\li\mu,\ov\Ui_t\infty\re}\qquad(t>0).
\ee
\el
{\bf Proof} For each $f\in\Ci_+(E)\cap\Di(G)$, the function $\ti u_t:=\Ui_tf$ solves
\be
\dif{t}\ti u_t= G\ti u_t+\bet \ti u_t-\al\ti u_t^2\leq G\ti u_t+\ov\bet \ti u_t-\un\al\ti u_t^2\qquad(t\geq 0),
\ee
and therefore $\Ui_tf=\ti u_t\leq\ov\Ui_tf$ by Lemma~\ref{subsol}. Using Lemmas~\ref{bpclos} and \ref{unbound} this is easily extended to measurable $f:E\to[0,\infty]$, giving (\ref{UU}). Define $\ov u$ by the right-hand side of the equations in (\ref{ovU}). Then it is easy to check that $\ov u$ solves $\dif{t}\ov u_t=\ov\bet\ov u_t-\un\al\ov u_t^2$ ($t>0$) with $\lim_{t\to 0}\ov u_t=\infty$, and therefore (\ref{sterfkans}) follows from the fact that
\be\label{Uinf}
P^\mu[\Ys_t=0]=E^\mu[e^{-\li\Ys_t,\infty\re}]=\ex{-\li\mu,\Ui_t\infty\re}\qquad(t\geq 0,\ \mu\in\Mi(E)),
\ee
and a little approximation argument.\qed

\subsection{Some consequences of the finite ancestry property}\label{fasec}

Let $\Ys$ be a $(G,\al,\bet)$-superprocess as in the last section. In line with Lemma~\ref{Wfa}, we say that $\Ys$ has the {\em finite ancestry property} if
\be\label{fa}
\inf_{x\in E}P^{\de_x}[\Ys_t=0]>0\qquad(t>0).
\ee
Note that by (\ref{Uinf}), property (\ref{fa}) is equivalent to $\|\Ui_t\infty\|_\infty<\infty$ ($t>0$). In this section we prove three simple consequences of the finite ancestry property.
\bl{\bf (Extinction versus unbounded growth)}\label{exfa}
Assume that the $(G,\al,\bet)$-super\-pro\-cess $\Ys$ has the finite ancestry property. Then, for any $\mu\in\Mi(E)$,
\be\label{exexSWF}
P^\mu\big[\Ys_t=0\mbox{ eventually}\mbox{\ \ or\ \ }\lim_{t\to\infty}\li\Ys_t,1\re=\infty\big]=1.
\ee
\el
{\bf Proof} We use a general fact about tail events of strong Markov processes, the statement and proof of which can be found in Section~\ref{S:zeone}. Consider the tail event $A:=\{\Ys_t=0\mbox{ eventually}\}$. By Lemma~\ref{tail} below,
\be\label{survWF}
\lim_{t\to\infty}P^{\Ys_t}(A)=1_A\qquad{\rm a.s.}
\ee
For any fixed $T>0$, by (\ref{Uinf}),
\be
P^\mu(A)\geq P^\mu[\Ys_T=0]=\ex{-\li\mu,\Ui_T\infty\re}\geq\ex{-\li\mu,1\re\|\Ui_T\infty\|_\infty}\qquad(\mu\in\Mi(E)).
\ee
Hence (\ref{survWF}) implies that
\be
\liminf_{t\to\infty}\ex{-\li\Ys_t,1\re\|\Ui_T\infty\|_\infty}\leq 1_A\qquad{\rm a.s.}
\ee
By the finite ancestry property, $\|\Ui_T\infty\|_\infty<\infty$ and therefore $\lim_{t\to\infty}\li\Ys_t,1\re=\infty$ a.s.\ on $A^\co$.\qed

\noi
The following is a simple consequence of Lemma~\ref{exfa}.
\bl{\bf (Extinction of (sub-) critical processes)}\label{subcri}
Assume that the $(G,\al,\bet)$-super\-process $\Ys$ has the finite ancestry property and that $\bet\leq 0$. Then, for any $\mu\in\Mi(E)$,
\be
P^\mu\big[\Ys_t=0\mbox{ eventually}\big]=1.
\ee
\el
{\bf Proof} Since $E^\mu[\li\Ys_t,1\re]\leq\li\mu,1\re$, $P^\mu[\lim_{t\to\infty}\li\Ys_t,1\re=\infty]=0$. Now the claim follows from Lemma~\ref{exfa}.\qed

\noi
Our final result of this section is the following.
\bl{\bf (Extinction versus exponential growth)}\label{versusSWF}
Assume that the $(G,\al,\bet)$-super\-process $\Ys$ has the finite ancestry property and that $\bet>0$ is a constant. Then, for any $\mu\in\Mi(E)$, there exists a nonnegative random variable $W$, depending on $\mu$, such that
\be\ba{rl}\label{vier}
{\rm (i)}&\dis\lim_{t\to\infty}e^{-\bet t}\li\Ys_t,1\re=W\qquad P^\mu{\rm -a.s.},\\[5pt]
{\rm (ii)}&\dis\lim_{t\to\infty}E^\mu\big[|e^{-\bet t}\li\Ys_t,1\re-W|^2\big]=0,\\[5pt]
{\rm (iii)}&\dis E^\mu(W)=\li\mu,1\re,\\[5pt]
{\rm (iv)}&\dis\var^\mu(W)\leq 2\bet^{-1}\|\al\|_\infty\li\mu,1\re,\\[5pt]
{\rm (v)}&\dis\{W=0\}=\{\Ys_t=0\mbox{ eventually}\}\qquad P^\mu{\rm -a.s.}
\ec
\el
{\bf Proof} Put $\Vi_tf:=e^{\bet t}S_t$. The mean and covariance of $\Ys$ are given by the following formulas, see, for example, \cite{Fit88}:
\be\left.\ba{rr@{\,}c@{\,}l}\label{momSWF}
{\rm (i)}&\dis E^\mu[\li\Ys_t,f\re]&=&\dis\li\mu,\Vi_tf\re\\
{\rm (ii)}&\dis\cov^\mu(\li\Ys_t,f\re,\li\Ys_t,g\re)&=&\dis 2\!\int_0^t\!\di s\,\li\mu,\Vi_s(\al\,(\Vi_{t-s}f)(\Vi_{t-s}g))\re\\
\ea\ \right\}\quad(t\geq 0,\ f,g\in B(E)).
\ee
Therefore,
\be\label{vmean}
E^\mu[\li\Ys_t,f\re]=e^{\bet t}\li\mu,S_tf\re\qquad(t\geq 0,\ f\in B(E)),
\ee
and
\bc\label{timean}
\dis\var^\mu(\li\Ys_t,f\re)
&=&\dis 2\int_0^t\!\!\di s\, e^{\bet s}e^{2\bet(t-s)}\li\mu,S_s(\al(S_{t-s}f)^2)\re\\[5pt]
&\leq&\dis 2\|\al\|_\infty\|f\|_\infty^2\li\mu,1\re e^{\bet t}\int_0^t\!\!\di s\, e^{\bet(t-s)}\\[10pt]
&\leq &\dis2\bet^{-1}\|\al\|_\infty\|f\|_\infty^2\li\mu,1\re e^{2\bet t}\qquad(t\geq 0,\ f\in B(E)).
\ec
Let $(\Fi_t)_{t\geq 0}$ be the filtration generated by $\Ys$ and put
\be
\ti\Ys_t:=e^{-\bet t}\Ys_t\qquad(t\geq 0).
\ee
Then (\ref{vmean}) and (\ref{timean}) show that for any $0\leq s\leq t$ and $f\in B(E)$,
\be\ba{rr@{\,}c@{\,}l}\label{timi}
{\rm (i)}&\dis E^\mu\big[\li\ti\Ys_t,f\re\big|\Fi_{s}\big]&=&\dis\li\ti\Ys_{s},S_{t-s}f\re\quad{\rm a.s.},\\[5pt]
{\rm (ii)}&\dis\var^\mu\big[\li\ti\Ys_t,f\re\big|\Fi_{s}\big]&\leq&\dis2\bet^{-1}\|\al\|_\infty\|f\|_\infty^2\li\ti\Ys_s,1\re e^{-\bet s}\quad{\rm a.s.}
\ec
Since $S_{t-s}1=1$, formula (\ref{timi})~(i) shows that $(\li\ti\Ys_t,1)\re_{t\geq 0}$ is a nonnegative martingale, and hence there exists a nonnegative random variable $W$ such that (\ref{vier})~(i) holds. Setting $s=0$ in (\ref{timi})~(ii), we see that
\be\label{varbo}
\var^\mu\big[\li\ti\Ys_t,1\re\big]\leq\dis 2\bet^{-1}\|\al\|_\infty\li\mu,1\re\qquad(t\geq 0).
\ee
This implies (\ref{vier})~(ii), and, using Fatou, (\ref{vier})~(iv). Moreover, by (\ref{varbo}) the random variables $\li\Ys_t,1\re_{t\geq 0}$ are uniformly integrable, and therefore (\ref{vier})~(iii) holds.

We are left with the task to prove (\ref{vier})~(v). The inclusion $\supset$ is trivial. Formulas (\ref{vier})~(iii) and (\ref{vier})~(iv) imply that
\be\label{truc}
\li\mu,1\re^2P^\mu[W=0]\leq\var^\mu(W)\leq 2\bet^{-1}\|\al\|_\infty\li\mu,1\re,
\ee
and therefore
\be\label{grkans}
P^\mu[W>0]\geq 1-2\bet^{-1}\|\al\|_\infty\li\mu,1\re^{-1}\qquad(\mu\neq 0).
\ee
Note that $\{W>0\}$ is a tail event. Thus, by Lemma~\ref{tail},
\be\label{limW}
\lim_{t\to\infty}P^{\Ys_t}[W>0]=1_{\{W>0\}}\quad{\rm a.s.}
\ee
Formula (\ref{grkans}) shows that
\be\label{liminfW}
\liminf_{t\to\infty}P^{\Ys_t}[W>0]\geq 1_{\{\lim_{t\to\infty}\li\Ys_t,1\re=\infty\}}.
\ee
Combining Lemma~\ref{exfa} with formulas (\ref{limW}) and (\ref{liminfW}) we see that $\{\Ys_t=0\mbox{ eventually}\}^\co\sub\{\lim_{t\to\infty}\li\Ys_t,1\re=\infty\}\sub\{W>0\}$ a.s.\qed

\subsection{Smoothness of two log-Laplace semigroups}\label{smosec}

We return to the special situation $E=[0,1]$ and $G=\ov A$ or $G=\ov{A^v}$, where $\ov A$ and $\ov{A^v}$ are the closures in $\Ci(E)$ of the operators $A$ in (\ref{AdefSWF}) and $A^v$ in (\ref{Avdef}), respectively, with domains $\Di(A)=\Di(A^v):=\Ci^{(2)}[0,1]$, the space of real functions on $[0,1]$ that are twice continuously differentiable. Let $\Ui=\Ui(\ov A,\gal,\gal)$ and $\Ui^v=\Ui(\ov{A^v},\gal v,\gal-1)$ denote the log-Laplace semigroups of the super-Wright-Fisher diffusion $\Ys$ and the weighted super-Wright-Fisher diffusion $\Ys^v$, respectively, where $\gal>0$ is constant. In this section we prove:
\bl{\bf (Smoothing property of $\Ui$ and $\Ui^v$)}\label{UconSWF}
One has $\Ui_t(B_+[0,1])\sub\Ci_+[0,1]$ and $\Ui^v_t(B_+[0,1])\sub\Ci_+[0,1]$ for all $t>0$. Moreover, if $\bplim_{n\to\infty}f_n=f$ for some $f_n,f\in B_+[0,1]$, then $\lim_{n\to\infty}\|\Ui_tf_n-\Ui_tf\|_\infty=0$ and $\lim_{n\to\infty}\|\Ui^v_tf_n-\Ui^v_tf\|_\infty=0$ for all $t>0$.
\el
To prepare for the proof, we start with the following elementary property of the semigroups $S$ and $S^v$ generated by $\ov A$ and $\ov{A^v}$, respectively (recall (\ref{AdefSWF}) and (\ref{Avdef})).
\bl{\bf (Strong Feller property)}\label{Smo}
The semigroups $S$ and $S^v$ have the strong Feller property, i.e., $S_t(B[0,1])\sub\Ci[0,1]$ and $S^v_t(B[0,1])\sub\Ci[0,1]$ for all $t>0$.
\el
{\bf Proof} Couple two realizations $\xi^x,\xi^y$ of the process with generator $\ov A$, started in $x,y\in[0,1]$, in such a way that $\xi^x$ and $\xi^y$ move independently up to the random time $\tau:=\inf\{t\geq 0:\xi^x_t=\xi^y_t\}$, and such that $\xi^x_t=\xi^y_t$ for all $t\geq\tau$. (Here the superscript in $\xi^x$ refers to the initial condition, and not, like elsewhere, to a compensated {\it h}-transform.) Then it is not hard to see that
\be\label{indcoup}
P[\xi^y_t=\xi^x_t]\to 1\quad\mbox{as}\quad y\to x\quad\forall t>0.
\ee
In particular, (\ref{indcoup}) holds also for $x\in\{0,1\}$ since the boundary is attainable. Since $|S_tf(x)-S_tf(y)|\leq 2\|f\|_\infty P[\xi^x_t\neq\xi^y_t]$, formula (\ref{indcoup}) shows that $S_tf\in\Ci[0,1]$ for all $f\in B[0,1]$ and $t>0$. For the process with generator $\ov{A^v}$ the argument is similar but easier, since in this case $\{0,1\}$ is an entrance boundary.\qed

\noi
{\bf Proof of Lemma~\ref{UconSWF}} For each $f\in B[0,1]$, the function $u_t:=\Ui_tf$ is a mild solution of (\ref{cauchyons}), i.e., (see (\ref{int}))
\be\label{int3}
\Ui_tf=S_tf+\int_0^t\!S_{t-s}\big(\gal\Ui_sf(1-\Ui_sf)\big)\,\di s\qquad(t\geq 0).
\ee
By the strong Feller property of $(S_t)_{t\geq 0}$ (Lemma~\ref{Smo}), the functions $S_tf$ and $S_{t-s}(\gal\Ui_sf(1-\Ui_sf))$ are continuous for each $0\leq s<t$, and therefore $\Ui_tf$ is continuous.

Now let $f_n\to f$ in a bounded pointwise way for some $f_n,f\in B_+[0,1]$, and let $t>0$. By Lemma~\ref{Uico}, $\Ui_tf_n\to\Ui_t f$ in a bounded pointwise way. By the strong Feller property of $(S_t)_{t\geq 0}$ and \cite[Prop.~1.5.8 and Thm.~1.5.9]{Rev84}, $S_tf_n$ converges uniformly to $S_tf$ and the function $(x,s)\mapsto S_{t-s}\big(\gal\Ui_sf_n(1-\Ui_sf_n)\big)(x)$ converges uniformly on $[0,1]\times[0,t-\eps]$ to $S_{t-s}\big(\gal\Ui_sf(1-\Ui_sf)\big)(x)$, for all $\eps>0$. By (\ref{int3}), it follows that $\Ui_tf_n\to\Ui_t f$ uniformly on $[0,1]$.

The same arguments apply to $\Ui^v_tf$.\qed

\subsection{Bounds on the absorption probability}\label{absec}

Let $\Ui=\Ui(\ov A,\gal,\gal)$. Since the points $0,1$ are traps for the Wright-Fisher diffusion, $f(r)=0$ implies $\Ui_tf(r)=0$ ($r=0,1$). We have already seen (Lemma~\ref{UconSWF}) that $\Ui_tf$ is continuous for each $t>0$. The following lemma shows that if $f(r)=0$, then $\Ui_tf$ has a finite slope at $r=0,1$, for all $t>0$. By symmetry, it suffices to consider the case $r=0$.
\bl{\bf(Absorption of the super-Wright-Fisher diffusion)}\label{absorb}
Let $\Ui=\Ui(\ov A,\gal,\gal)$, with $\gal>0$. Then
\be\label{lui}
\Ui_t(\infty1_{(0,1]})(x)\leq K_t\,x\qquad(t>0,\ x\in[0,1]),
\ee
with
\be
K_t:=\frac{e^{\gal t/2}}{1-e^{-\gal t/2}}\Big(\frac{8}{t}+2\Big)\qquad(t>0).
\ee
\el
Note that (\ref{lui}) implies that
\be\label{supabs}
P^{\de_x}[\Ys_t((0,1])>0]\leq 1-\ex{-K_t\,x}\leq K_t\,x\qquad(t>0,\ x\in[0,1]).
\ee
We begin with a preparatory lemma.
\bl{\bf (Absorption of the Wright-Fisher diffusion)}\label{exbo}
For the Wright-Fisher diffusion $\xi$,
\be\label{noabs}
P^x[\xi_t>0]\leq\Big(\frac{4}{t}+2\Big)x\qquad(t>0,\ x\in[0,1]).
\ee
\el
{\bf Proof} For $x\geq 0$ put 
\be\label{fdef}
f_0(x):=1_{\{0\}}(x)\quad\mbox{and}\quad f_t(x):=(1-2x)\ex{-\frac{4x}{t}}1_{[0,\frac{1}{2}]}(x)\qquad(t>0).
\ee
A little calculation shows that for $t>0$ and $x\geq 0$,
\bc
\dis\dif{t}f_t(x)&=&\dis 4x(1-2x)t^{-2}e^{-\frac{4x}{t}}1_{[0,\frac{1}{2}]}(x)\\[5pt]
\dis\ffrac{1}{2}x(1-x)D^2_xf_t(x)&=&\dis\big(8x(1-x)(1-2x)t^{-2}e^{-\frac{4x}{t}}+8x(1-x)t^{-1}e^{-\frac{4x}{t}}\big)1_{[0,\frac{1}{2}]}(x)\\[5pt]
&&\dis+2e^{-\frac{2}{t}}\de_{\frac{1}{2}}(x),
\ec
where $D^2_x$ denotes the generalized second derivative with respect to $x$ and $\de_{\frac{1}{2}}$ is the delta-function at $\frac{1}{2}$. Since $4x\leq 8x(1-x)$ for all $x\in[0,\frac{1}{2}]$, it follows that
\be\label{subf}
\dif{t}f_t(x)\leq\ffrac{1}{2}x(1-x)D^2_xf_t(x)\qquad(t>0,\ x\geq 0).
\ee
If $f_t$ were contained in $\Di(\ov A)$, then (\ref{subf}) would mean that $\dif{t}f_t\leq\ov Af_t$ for $t>0$, and a standard argument (compare Lemma~\ref{subsol}) would tell us that $f_t\leq S_tf_0$, where $S$ is the semigroup of $\xi$. In the present case, we need a little approximation argument.

Let $\phi_n\geq 0$ ($n\geq 0$) denote $\Ci^{(\infty)}$-functions defined on $\half$ with support contained in $[0,\ffrac{1}{3}]$, say, such that $\phi_n(x)\di x$ are probability measures converging weakly to the $\de$-measure $\de_0$ as $n\to\infty$. Put
\be
f^n_t(x):=\int_0^\infty\!\di y\,\phi_n(y)f_t(x+y)=:\phi_n\ast f_t(x)\qquad(t>0,\ x\geq 0).
\ee
Then
\bc
\dis\dif{t}f^n_t(x)&=&\dis\phi_n\ast\dif{t}f_t(x)\\[5pt]
\dis\diff{x}f^n_t(x)&=&\dis\phi_n\ast D^2_xf_t(x),
\ec
and therefore (\ref{subf}) shows that
\be\label{subfn}
\dif{t}f^n_t(x)\leq\ffrac{1}{2}x(1-x)\diff{x}f^n_t(x)\qquad(t>0,\ x\geq 0,\ n\geq 0).
\ee
Since $f^n_t\in\Di(\ov A)$ for all $t>0$, the argument mentioned above gives
\be
f^n_{t+\eps}\leq S_tf^n_\eps\qquad(t\geq 0,\ \eps>0).
\ee
Letting $n\to\infty$ and afterwards $\eps\to 0$ we find that
\be\label{fbound}
f_t(x)\leq S_tf_0(x)=P^x[\xi_t=0]\qquad(t\geq 0,\ x\in[0,1]).
\ee
Note that $\dif{x}(1-f_t(x))=(1-2x)4t^{-1}e^{-\frac{4x}{t}}+2e^{-\frac{4x}{t}}\leq(\frac{4}{t}+2)$ for $x\in[0,\ffrac{1}{2}]$. Therefore (\ref{fbound}) implies (\ref{noabs}). (Note that (\ref{noabs}) is trivial for $x\in[\ffrac{1}{2},1]$.)\qed

\noi
{\bf Proof of Lemma~\ref{absorb}} Fix $f\in B_+[0,1]$ satisfying $f(0)=0$ and write $\Ui_tf=\Ui_{t/2}\Ui_{t/2}f$. By (\ref{sterfkans}) from Lemma~\ref{exti}, $\Ui_{t/2}f\leq (1-e^{-\gal t/2})^{-1}$. Since moreover $\Ui_{t/2}f(0)=0$ because of absorption at zero, we have
\be
\Ui_tf\leq\Ui_{t/2}((1-e^{-\gal t/2})^{-1}1_{(0,1]})\qquad(t>0).
\ee
Using (\ref{UU}) from Lemma~\ref{exti}, we may estimate $\Ui(\ov A,\gal,\gal)$ in terms of $\Ui(\ov A,0,\gal)$, which is just the linear semigroup $(e^{\gal t}S_t)_{t\geq 0}$. Thus, by Lemma~\ref{exbo},
\bc
\dis\Ui_tf(x)&\leq&\dis e^{\gal t/2}S_{t/2}((1-e^{-\gal t/2})^{-1}1_{(0,1]})(x)\\[5pt]
&\leq&\dis e^{\gal t/2}(1-e^{-\gal t/2})^{-1}(\ffrac{8}{t}+2)x\qquad(t>0,\ x\in[0,1]).
\ec
Letting $f\up\infty$, by monotonicity we arrive at (\ref{lui}).\qed

\subsection{The weighted super-Wright-Fisher diffusion}\label{trafosec}

In this section we prove Lemmas~\ref{weilem} and \ref{Wfa}. Recall that $\xi,\xi^v$ are the diffusions in $[0,1]$ with generators $\ov A,\ov{A^v}$ defined in (\ref{AdefSWF}) and (\ref{Avdef}), and associated semigroups $S,S^v$, respectively, and that $\Ui=\Ui(\ov A,\gal,\gal)$ and $\Ui^v=\Ui(\ov{A^v},\gal v,\gal-1)$.
\bl{\bf({\it v}-transformed log-Laplace semigroup)}
If $f\in\Di(\ov{A^v})$, then $vf\in\Di(\ov A)$ and
\be\label{aatra}
\ov A(vf)=v\,(\ov{A^v}-1)f.
\ee
Moreover,
\be\label{utra}
\Ui_t(vf)=v\,\Ui^v_tf\qquad(t\geq 0,\ f\in B_+[0,1]).
\ee
\el
{\bf Proof} For any $f\in\Ci^{(2)}[0,1]$, it is easy to check that
\be\label{atra}
A(vf)=v\,(A^v-1)f.
\ee
Fix $f\in\Di(\ov{A^v})$ and choose $f_n\in\Ci^{(2)}[0,1]$ such that $f_n\to f$ in $\Ci[0,1]$. Then (\ref{atra}) shows that $A(vf_n)\to v\,(\ov{A^v}-1)f$, which implies that $vf\in\Di(\ov A)$ and that (\ref{aatra}) holds.

Now fix $f\in\Ci_+[0,1]\cap\Di(\ov{A^v})$ and put $u^v_t:=\Ui^v_tf$ ($t\geq 0$). Then $u^v$ is the classical solution of the Cauchy equation
\be\left\{\:\ba{r@{\,}c@{\,}l}\label{cauchy4}
\dif{t}u^v_t&=&\ov{A^v}u^v_t+(\gal-1)u^v_t-\gal v\,(u^v_t)^2\qquad(t\geq 0),\\
u^v_0&=&f.
\ea\right.\ee
It follows from (\ref{aatra}) that
\be\ba{r@{\,}c@{\,}l}
\dif{t}vu^v_t&=&v\dif{t}u^v_t=v\ov{A^v}u^v_t+(\gal-1)vu^v_t-\gal\,(vu^v_t)^2\\[5pt]
&=&\ov A(vu^v_t)+\gal vu^v_t-\gal\,(vu^v_t)^2\qquad(t\geq 0),
\ec
i.e., $u_t:=vu^v_t$ is the classical solution to the Cauchy equation
\be\left\{\:\ba{r@{\,}c@{\,}l}\label{cauchy5}
\dif{t}u_t&=&\ov Au_t+\gal u_t-\gal u_t^2\qquad(t\geq 0),\\
u_0&=&vf.
\ea\right.\ee
This proves that $\Ui_t(vf)=u_t=vu^v_t=v\Ui^v_tf$ for all $f\in\Ci_+[0,1]\cap\Di(\ov{A^v})$. The general case follows from Lemma~\ref{bpclos} and the fact that the class of $f\in B_+[0,1]$ for which (\ref{utra}) holds is closed under bounded pointwise limits.\qed

\noi
{\bf Proof of Lemma~\ref{weilem}} Set $\Fi_t:=\sig(\Ys_s:0\leq s\leq t)$. Then by (\ref{utra}), for all $0\leq s\leq t$ and $f\in B_+[0,1]$,
\be\ba{l}
E\big[\ex{-\li v\Ys_t,f\re}\big|\Fi_s\big]=E\big[\ex{-\li\Ys_t,vf\re}\big|\Fi_s\big]=\ex{-\li\Ys_s,\Ui_{t-s}(vf)\re}\\[5pt]
=\ex{-\li\Ys_s,v\Ui^v_{t-s}f\re}=\ex{-\li v\Ys_s,\Ui^v_{t-s}f\re}.
\ec
It follows that $(v\Ys_t)_{t\geq 0}$ is a Markov process and that its transition probabilities coincide with those of the $(\ov{A^v},\gal v,\gal-1)$-superprocess. Since $\Ys$ has continuous sample paths, so has $v\Ys$.\qed

\noi
{\bf Proof of Lemma~\ref{Wfa}} We need to prove (\ref{Wfin}), which by (\ref{Uinf}) is equivalent to the statement that $\|\Ui^v_t\infty\|_\infty<\infty$ for all $t>0$. Assume that $f\in B_+[0,1]$ satisfies $f(0)=f(1)=0$. By Lemma~\ref{absorb}, $\Ui_tf(x)\leq K_t\,x$ for the constant $K_t$ mentioned there. By symmetry, one also has $\Ui_tf(x)\leq K_t\,(1-x)$ and, since $x\wedge(1-x)\leq\frac{1}{3}v(x)$, $\Ui_tf(x)\leq\frac{1}{3}K_t\,v(x)$. Let $g\in B_+[0,1]$. By formula (\ref{utra}) and the fact that $(vg)(0)=(vg)(1)=0$, we see that $\Ui^v_tg(x)=\frac{1}{v(x)}\Ui_t(vg)(x)\leq\frac{1}{3}K_t$ for all $x\in(0,1)$. By Lemma~\ref{Smo}, $\Ui^v_tg$ is continuous on $[0,1]$ and therefore $\Ui^v_tg(x)\leq\frac{1}{3}K_t$ holds also for $x=0,1$. Taking the limit $g\up\infty$ we see that $\|\Ui^v_t\infty\|_\infty\leq\frac{1}{3}K_t<\infty$ for all $t>0$.\qed

\subsection{A zero-one law for Markov processes}\label{S:zeone}

Let $E$ be a Polish space and let $(P^x)^{x\in E}$ be a family of probability measures on $\Di_E\half$ (the space of cadlag functions $w:\half\to E$) such that under $(P^x)^{x\in E}$, the coordinate projections $\{w\mapsto w_t=:\xi_t(w):t\geq 0\}$ form a Borel right process in the sense of \cite{Sha88}. This is true, for example, if $(P^x)^{x\in E}$ are the laws of a Feller process on a locally compact Polish space, or a $(G,\al,\bet)$-superprocess as introduced in Section~\ref{prepar}, see \cite{Fit88}. Let $\Ti:=\bigcap_{t\geq 0}\sig(\xi_s:s\geq t)$ denote the tail-\si-field of $\xi$. Let $(\tet_tw)_s:=w_{t+s}$ ($t,s\geq 0$) be the time-shift on $\Di_E\half$. Then the following holds.
\bl{\bf (Zero-one law for Markov processes)}\label{tail}
Assume that $A\in\Ti$. Then for each $x\in E$,
\be\label{asli}
\lim_{t\to\infty}P^{\xi_t}(\tet_t^{-1}(A))=1_A\quad P^x{\rm -a.s.}
\ee
\el
{\bf Proof} Let $\Fi_t:=\sig(\xi_s:0\leq s\leq t)$ ($t\geq 0$) be the filtration generated by $\xi$ and set $\Fi_\infty:=\sig(\xi_s:s\geq 0)$. Since $\xi$ is a Markov process, $P^{\xi_t}(\tet_t^{-1}(A))=P[A|\Fi_t]$ a.s. For any sequence of times $t_n\up\infty$ one has $\Fi_{t_n}\up\Fi_\infty$ and therefore $P[A|\Fi_{t_n}]\to P[A|\Fi_\infty]=1_A$ a.s., see \cite[\S~29, Complement~10~(b)]{Loe63}. Since $\xi$ is a right process, the function $t\mapsto P^{\xi_t}(\tet_t^{-1}(A))$ is a.s.\ right-continuous, see \cite[Theorem~(7.4.viii)]{Sha88}, and we conclude that (\ref{asli}) holds.\qed

\section{The super-Wright-Fisher diffusion: long-time behavior}\label{S:WFiii}

\subsection{Ergodicity of the compensated {\it v}-transformed Wright-Fisher diffusion}\label{ergsec}

Recall that $\xi^v$ is the diffusion on $[0,1]$ with generator $\ov{A^v}$ defined in (\ref{Avdef}) and associated semigroup $S^v$. As in Theorem~\ref{mainSWF}, $\ell$ denotes the Lebesgue measure on $(0,1)$ and $v$ is defined by (\ref{vdef}). In this section we prove:
\bl{\bf (Ergodicity of the compensated {\it v}-transformed Wright-Fisher diffusion)}\label{vergo}
The Markov process $\xi^v$ has the unique invariant law $v\ell$ and is ergodic:
\be\label{fB}
\lim_{t\to\infty}\|S^v_tf-\li v\ell,f\re\|_\infty=0\qquad\forall f\in B[0,1].
\ee
\el
{\bf Proof} Since
\be
\dif{x}\big[\ffrac{1}{2}x(1-x)v(x)\big]=2(\ffrac{1}{2}-x)v(x)\qquad(x\in[0,1]),
\ee
$v\ell$ is a (reversible) invariant law for the process with generator $\ov{A^v}$, see \cite[Proposition~4.9.2]{EK}. Fix $x\in[0,1]$. Let $\xi^v$ be the process started in $x$ and let $\ti\xi^v$ be the process started in the invariant law $v\ell$. Then $\xi^v,\ti\xi^v$ may represented as solutions to the SDE
\be
\di\xi^v_t=2(\ffrac{1}{2}-\xi^v_t)\di t+\sqrt{\xi^v_t(1-\xi^v_t)}\di B_t,
\ee
relative to the same Brownian motion $B$. Using the technique of Yamada \& Watanabe (see \cite{YW71} or, for example, \cite[Theorem~5.3.8]{EK}), it is easy to prove that
\be
E[|\xi^v_t-\ti\xi^v_t|]=e^{-2t}E[|\xi^v_0-\ti\xi^v_0|]\leq e^{-2t}\qquad(t\geq 0).
\ee
It follows that for any function $f$ satisfying $|f(y)-f(z)|\leq |y-z|\quad (y,z\in[0,1]$),
\be
\big|E[f(\xi^v_t)]-\li v\ell,f\re\big|\leq E[|f(\xi^v_t)-f(\ti\xi^v_t)|]\leq e^{-2t}.
\ee
This implies that the function $x\mapsto \Li^x(\xi^v_t)$ from $[0,1]$ into the space $\Mi_1[0,1]$ of probability measures on $[0,1]$, converges as $t\to\infty$ uniformly to the constant function $v\ell$. This shows that (\ref{fB}) holds for all $f\in\Ci[0,1]$. Since $\xi^v$ has the strong Feller property (Lemma~\ref{Smo}), (\ref{fB}) holds for all $f\in B[0,1]$.\qed

\subsection{Long-time behavior of the weighted super-Wright-Fisher diffusion}\label{Longwei}

\newcommand{\q}{(\gal-1)}

The following lemma prepares for the proof of formula~(\ref{explSWF}) in Theorem~\ref{mainSWF}.
\bl{\bf (Mean square convergence)}\label{mesq}
Assume that $\gal>1$. Let $\Ys^v$ be the $(\ov{A^v},\gal v,\gal-1)$-superprocess started in $\Ys^v_0=\mu\in\Mi[0,1]$. Then there exists a nonnegative random variable $W$, depending on $\mu$, such that
\be\ba{rl}\label{toW}
{\rm (i)}&\dis\lim_{t\to\infty}e^{-(\gal-1)t}\li\Ys^v_t,1\re=W\quad{\rm a.s.}\\
{\rm (ii)}&\dis\lim_{t\to\infty}E^\mu\Big[\big|e^{-(\gal-1)t}\li\Ys^v_t,f\re-W\li v\ell,f\re\big|^2\Big]=0\quad\forall f\in B[0,1].
\ec
Moreover,
\be\label{WW}
E^\mu(W)=\li\mu,1\re\quad\mbox{and}\quad\var^\mu(W)\leq 3\ffrac{\gal}{\gal-1}\li\mu,1\re,
\ee
and
\be\label{wext}
\{W=0\}=\{\Ys^v_t=0\mbox{ eventually}\}\qquad{\rm a.s.}
\ee
\el
{\bf Proof} Except for formula (\ref{toW})~(ii), all statements are direct consequences of the fact that $\Ys^v$ has the finite ancestry property (Lemma~\ref{Wfa}) and of Lemma~\ref{versusSWF} (note that $\|\gal v\|_\infty=\frac{3}{2}\gal$).

Fix $f\in B[0,1]$. Let $(\Fi_t)_{t\geq 0}$ be the filtration generated by $\Ys^v$ and put $\ti\Ys^v_t:=e^{-(\gal-1)t}\Ys^v_t$ ($t\geq 0$). Pick $1\leq s_n\leq t_n$ such that $s_n\to\infty$ and $t_n-s_n\to\infty$. Then, by (\ref{timi}),
\be\label{coqua}
E^\mu\Big[\big|\li\ti\Ys^v_{t_n},f\re-\li\ti\Ys^v_{s_n},S^v_{t_n-s_n}f\re\big|^2\Big|\Fi_{s_n}\Big]\leq3\ffrac{\gal}{\gal-1}\|f\|_\infty^2\li\ti\Ys^v_{s_n},1\re e^{-(\gal-1)s_n}\quad{\rm a.s.}
\ee
Taking expectations on both sides in (\ref{coqua}), one finds that
\be\label{emu}
E^\mu\Big[\big|\li\ti\Ys^v_{t_n},f\re-\li\ti\Ys^v_{s_n},S^v_{t_n-s_n}f\re\big|^2\Big]\leq3\ffrac{\gal}{\gal-1}\|f\|_\infty^2\li\mu,1\re e^{-(\gal-1)s_n}.
\ee
By (\ref{vier})~(ii),
\be\label{masqua}
\lim_{t\to\infty}E^\mu\big[|\li\ti\Ys^v_t,1\re-W|^2\big]=0.
\ee
Using Lemma~\ref{vergo} (about the ergodicity of $\xi^v$) and (\ref{masqua}), it is easy to show that
\be
\lim_{n\to\infty}E^\mu\Big[\big|\li\ti\Ys^v_{s_n},S^v_{t_n-s_n}f\re-W\li v\ell,f\re\big|^2\Big]=0.
\ee
Combining this with (\ref{emu}), we see that
\be
\lim_{n\to\infty}E^\mu\Big[\big|\li\ti\Ys^v_{t_n},f\re-W\li v\ell,f\re\big|^2\Big]=0.
\ee
Since this is true for any $t_n\to\infty$, (\ref{toW})~(ii) follows.\qed

\subsection{Long-time behavior of the super-Wright-Fisher diffusion}\label{Longsup}

{\bf Proof of Theorem~\ref{mainSWF}} Using Lemma~\ref{weilem}, we can translate our results on the weighted super-Wright-Fisher diffusion $\Ys^v$ to the super-Wright-Fisher diffusion $\Ys$. Thus, Lemma~\ref{mesq} proves formulas (\ref{WWdef})~(ii), (\ref{Wext})~(ii), and (\ref{Wmom})--(\ref{explSWF}), where $W_{(0,1)}$ is the random variable $W$ from Lemma~\ref{mesq}. Formula (\ref{sterf}) follows from Lemma~\ref{subcri}. To finish the proof of Theorem~\ref{mainSWF}, it suffices to prove (\ref{WWdef})~(i), (\ref{Wext})~(i) and (\ref{WWW}).\vc

\noi
{\bf\boldmath $1^\circ$.\ Proof of formula (\ref{WWdef})~(i)}
One has $E^\mu[\li\Ys_t,f\re]=e^{\gal t}\li\mu,S_tf\re$ for all $t\geq 0$, $f\in B[0,1]$ by (\ref{vmean}). Since the points $r=0,1$ are traps for the Wright-Fisher diffusion, $E^\mu[\li\Ys_t,1_{\{r\}}\re]=e^{\gal t}\li\mu,S_t1_{\{r\}}\re\geq e^{\gal t}\li\mu,1_{\{r\}}\re$ for all $t\geq 0$, $r=0,1$. Thus, the processes $(e^{-\gal t}\li\Ys_t,1_{\{r\}}\re)_{t\geq 0}$ ($r=0,1$) are nonnegative submartingales, and hence there exist random variables $W_r$ ($r=0,1$) such that (\ref{WWdef})~(i) holds.\vc

\noi
{\bf\boldmath $2^\circ$.\ Proof of formula (\ref{WWW})} For $\gal\leq 1$ the statement is trivial by (\ref{sterf}), so assume $\gal>1$. By symmetry it suffices to consider the case $r=0$. From the $L_2$-convergence formula (\ref{explSWF}) we have, for any $K>0$,
\be\label{Kgeq}
\{W_{(0,1)}>0\}\sub\big\{\forall T<\infty\ \exists t\geq T\mbox{ such that }\Ys_t([\ffrac{1}{4},\ffrac{1}{3}])\geq K\big\}\quad{\rm a.s.}
\ee
Assume for the moment that for some $t>0$ and (sufficiently large) $K$,
\be\label{intrand}
\inf_{\mu:\,\mu([\frac{1}{4},\frac{1}{3}])\geq K}P^\mu[W_0>0]>0.
\ee
Then we see from (\ref{Kgeq}) and (\ref{intrand}) that
\be
\{W_{(0,1)}>0\}\sub\big\{\lim_{t\to\infty}P^{\Ys_t}[W_0>0]=0\big\}^\co\sub\{W_0>0\}\quad{\rm a.s.},
\ee
where the second inclusion follows from the fact that, by Lemma~\ref{tail},
\be\label{zeroone}
\lim_{t\to\infty}P^{\Ys_t}[W_0>0]=1_{\{W_0>0\}}\quad{\rm a.s.}
\ee
Thus, we are done if we can prove (\ref{intrand}). By the branching property, it suffices to prove (\ref{intrand}) for measures $\mu$ that are concentrated on $[\ffrac{1}{4},\ffrac{1}{3}]$. Fix any $t>0$. Formulas~(\ref{vmean}) and (\ref{timean}) give
\be\ba{rr@{\,}c@{\,}l}\label{timir}
{\rm (i)}&\dis E^\mu\big[\li\Ys_t,1_{\{0\}}\re\big]&=&\dis\li\mu,S_t1_{\{0\}}\re e^{\gal t},\\[5pt]
{\rm (ii)}&\dis\var^\mu\big[\li\Ys_t,1_{\{0\}}\re\big]&\leq&\dis2\li\mu,1\re e^{2\gal t}.
\ec
It follows from formula (\ref{fbound}) (recall (\ref{fdef})) that
\be
\inf_{x\in[\frac{1}{4},\frac{1}{3}]}S_t1_{\{0\}}(x)>0.
\ee
Denoting the infimum by $\eps$, we get the bounds
\be\ba{rr@{\,}c@{\,}l}\label{timir2}
{\rm (i)}&\dis E^\mu\big[\li\Ys_t,1_{\{0\}}\re\big]&\geq&\dis\eps\li\mu,1\re e^{\gal t},\\[5pt]
{\rm (ii)}&\dis\var^\mu\big[\li\Ys_t,1_{\{0\}}\re\big]&\leq&\dis2\li\mu,1\re e^{2\gal t}.
\ec
These formulas show that for large $\li\mu,1\re$, the standard deviation of $\li\Ys_t,1_{\{0\}}\re$ is small compared to its mean. Therefore, using Chebyshev's inequality, it is easy to show that for every $M>0$ there exists a $K>0$ such that
\be\label{MK}
\inf_{\mu\in\Mi[\frac{1}{4},\frac{1}{3}]:\,\li\mu,1\re\geq K}P^\mu[\li\Ys_t,1_{\{0\}}\re\geq M]>0.
\ee
Hence, by the Markov property, in order to prove (\ref{intrand}) it suffices to show that for $M$ sufficiently large,
\be\label{muM}
\inf_{\mu:\,\mu(\{0\})\geq M}P^\mu[W_0>0]>0.
\ee
By the branching property, it suffices to prove (\ref{muM}) for measures $\mu$ that are concentrated on $\{0\}$. In that case, $\Ys_t(\{0\})_{t\geq 0}$ is an autonomous supercritical Feller's branching diffusion (a superprocess in a single-point space is just a Feller's branching diffusion). Applying Lemma~\ref{versusSWF} to this Feller's branching diffusion, again using Chebyshev, it is not hard to prove (\ref{muM}). Since the arguments are very similar to those we have already seen, we skip the details.\vc

\noi
{\bf\boldmath $3^\circ$.\ Proof of formula (\ref{Wext})~(i)} The inclusion $\{W_r=0\}\supset\{\Ys_t(\{r\})=0\mbox{ eventually}\}$ a.s.\ is trivial. By (\ref{WWW}) and (\ref{Wext})~(ii), $\{W_r=0\}\sub\{W_{(0,1)}=0\}\sub\{\Ys_t((0,1))=0\mbox{ eventually}\}$ a.s. Therefore, by the strong Markov property, it suffices to prove $\{W_r=0\}\sub\{\Ys_t(\{r\})=0\mbox{ eventually}\}$ a.s.\ for the process started in $\mu$ with $\mu((0,1))=0$. In this case, $(\Ys_t(\{r\}))_{t\geq 0}$ is an autonomous supercritical Feller's branching diffusion, and the statement is easy (see the previous parapraph).\qed

\subsection{Long-time behavior of the log-Laplace semigroup}\label{LongU}

{\bf Proof of Proposition~\ref{Ucor}} We start by proving that for all $\mu\in\Mi[0,1]$ and $f\in B_+[0,1]$,
\be\ba{l}\label{limU}
\dis\lim_{t\to\infty}\ex{-\li\mu,\Ui_tf\re}\\[5pt]
\dis=P^{\,\txt\mu}\Big[\big\{f(0)=0\mbox{ or }W_0=0\big\}\cap\big\{f(1)=0\mbox{ or }W_1=0\big\}\cap\big\{\li\ell,f\re=0\mbox{ or }W_{(0,1)}=0\big\}\Big]\\[5pt]
=\left\{\!\ba{l@{\ }l}
\dis 1&\mbox{if}\quad f(0)=f(1)=\li\ell,f\re=0,\\
\dis P^\mu\big[W_{(0,1)}=0\big]&\mbox{if}\quad f(0)=f(1)=0,\ \li\ell,f\re>0,\\
\dis P^\mu\big[W_0=0\big]=P^\mu\big[W_0=W_{(0,1)}=0\big]&\mbox{if}\quad f(0)>0,\ f(1)=0,\\
\dis P^\mu\big[W_1=0\big]=P^\mu\big[W_1=W_{(0,1)}=0\big]&\mbox{if}\quad f(0)=0,\ f(1)>0,\\
\dis P^\mu\big[W_0=W_1=0\big]=P^\mu\big[W_0=W_1=W_{(0,1)}=0\big]&\mbox{if}\quad f(0)>0,\ f(1)>0,\\
\ea\right.
\ec
where $P^\mu[W_{(0,1)}=0]<1$ if and only if $\gal>1$ and $\li\mu,v\re>0$.

Indeed, by formula (\ref{laplace}),
\be\label{prod}
\ex{-\li\mu,\Ui_tf\re}=E^\mu\big[\ex{-f(0)\Ys_t(\{0\})}\ex{-f(1)\Ys_t(\{1\})}\ex{-\li\Ys_t,1_{(0,1)}f\re}\big].
\ee
By (\ref{WWdef})~(i) and (\ref{Wext})~(i) in Theorem~\ref{mainSWF},
\be\label{raco}
\lim_{t\to\infty}\ex{-f(r)\Ys_t(\{r\})}=1_{\{f(r)=0\mbox{\ \scriptsize or }W_r=0\,\}}\quad{\rm a.s.}\quad(r=0,1).
\ee
Now, if $\li\ell,f\re=0$ for some $f\in B_+[0,1]$, then $e^{-\li\Ys_t,1_{(0,1)}f\re}=1$ a.s.\ for each $t>0$. To see this, note that by (\ref{momSWF}), $E^{\de_x}[\li\Ys_t,1_{(0,1)}f\re]=e^{\gal t}\li\de_x,S_t1_{(0,1)}f\re=e^{\gal t}E^x[1_{(0,1)}(\xi_t)f(\xi_t)]$ where $\xi$ is the Wright-Fisher diffusion. Since the law of the Wright-Fisher diffusion at any time $t>0$ (started in an arbitrary initial condition) on $(0,1)$ is absolutely continuous with respect to Lebesgue measure, we see that $E^{\de_x}[\li\Ys_t,1_{(0,1)}f\re]=0$ and hence $\li\Ys_t,1_{(0,1)}f\re=0$ $P^{\de_x}$-a.s. (Actually, since $\Ys$ is a one-dimensional superprocess, one can prove that $\Ys_t$, restricted to $(0,1)$, for $t>0$ is almost surely absolutely continuous with respect to Lebesgue measure.)

On the other hand, if $\li\ell,f\re>0$, then by formulas (\ref{WWdef})~(ii), (\ref{Wext})~(ii), (\ref{sterf}), and (\ref{explSWF}) in Theorem~\ref{mainSWF},
\be
\ex{-\li\Ys_t,1_{(0,1)}f\re}\stackrel{{\rm P}}{\longrightarrow}1_{\{W_{(0,1)}=0\}}.
\ee
Hence, for general $f\in B_+[0,1]$,
\be\label{mico}
\ex{-\li\Ys_t,1_{(0,1)}f\re}\stackrel{{\rm P}}{\longrightarrow}1_{\{\li\ell,f\re=0\mbox{\ \scriptsize or }W_{(0,1)}=0\}},
\ee
where $\stackrel{{\rm P}}{\longrightarrow}$ denotes convergence in probability. Inserting (\ref{raco}) and (\ref{mico}) into (\ref{prod}) we arrive at the first equality in (\ref{limU}). Using formula (\ref{WWW}) and checking the eight possibilities for $f(0),f(1),\li\ell,f\re$ to be zero or positive, we find the second equality in (\ref{limU}).

In particular, setting $\mu=\de_x$ in (\ref{limU}) we see that $\Ui_t f$ converges in a bounded pointwise way to $\fixi$ or to one of the functions $\fixii,\ldots,\fixv$ from (\ref{p1}), where $\fixii=0$ if $\gal\leq 1$ and $\fixii>0$ on $(0,1)$ otherwise. It follows from Lemma~\ref{UconSWF} that the convergence in (\ref{Up5}) is in fact uniform.

The fact that $p_{l,r}(0)=l$ and $p_{l,r}(1)=r$ will follow from Proposition~\ref{fix2}. The statements about smoothness of fixed points will be proved in Section~\ref{Smosec} below.\qed

\noi
{\bf Proof of Proposition~\ref{fix2}} By Proposition~\ref{Ucor}, for the functions $\fixii,\ldots,\fixv$ from (\ref{p1}),
\be\left.\ba{r@{\,}c@{\,}l}\label{p3}
\fixii(x)&=&\lim_{t\to\infty}\Ui_t1_{(0,1)}(x),\\
\fixiii(x)&=&\lim_{t\to\infty}\Ui_t1_{\{0\}}(x)=\lim_{t\to\infty}\Ui_t1_{[0,1)}(x),\\
\fixiv(x)&=&\lim_{t\to\infty}\Ui_t1_{\{1\}}(x)=\lim_{t\to\infty}\Ui_t1_{(0,1]}(x),\\
\fixv(x)&=&\lim_{t\to\infty}\Ui_t1
\ea\quad\right\}\quad(x\in[0,1]).
\ee
Since by formula (\ref{momgen}), for each Borel measurable $B\sub[0,1]$, $P^{\de_x}[\Yp_t(B)>0]=U_t1_B=\Ui_t1_B(x)$ ($t\geq 0,\ x\in[0,1]$), we can rewrite the expressions in the right-hand side of (\ref{p3}) as in (\ref{p2}).\qed

\subsection{Smoothness of fixed points}\label{Smosec}

In order to finish the proof of Proposition~\ref{Ucor} we need to show that the functions $\fixii,\ldots,\fixv$ occurring there are twice continuously differentiable on $[0,1]$. We begin with the following.
\bl{\bf\hspace{-2.3pt}(Smoothness of fixed points)\hspace{-2.2pt}}\label{smooth}
If $p\in B_+[0,1]$ is a fixed point under $\Ui(\ov A,\gal,\gal)$, then $p\in\Di(\ov A)$ and $\ov A p+\gal\, p(1-p)=0$.
\el
{\bf Proof} For any $t\geq0$, Lemma~\ref{UconSWF} implies that $p=\Ui_tp\in\Ci_+[0,1]$. Moreover, since $u_t:=p$ ($t\geq 0$) is a mild solution of (\ref{cauchyons}) (recall (\ref{int3})),
\be
p=S_tp+\int_0^t\!S_s\big(\gal p(1-p)\big)\di s\qquad(t\geq 0).
\ee
Hence
\be
\ov Ap:=\dis\lim_{t\to 0}\,t^{-1}(S_tp-p)=-\lim_{t\to 0}\,t^{-1}\!\!\int_0^t\!S_s\big(\gal\, p(1-p)\big)\di s=-\gal\, p(1-p),
\ee
where the limit exists in $\Ci[0,1]$.\qed

\noi
In this one-dimensional situation, the domain of $\ov A$ is known explicitly. One has, see \cite[Theorem~8.1.1]{EK}
\be\label{Dov}
\Di(\ov A)=\Big\{f\in\Ci[0,1]\cap\Ci^{(2)}(0,1):\lim_{x\to r}\ffrac{1}{2}x(1-x)\diff{x}f(x)=0\ \ (r=0,1)\Big\}.
\ee
Here $\Ci[0,1]\cap\Ci^{(2)}(0,1)$ denotes the class of continuous real functions on $[0,1]$ that are twice continuously differentiable on $(0,1)$.\vc

\noi
{\bf Proof of the smoothness of fixed points} It suffices to show that $\fixii$ and $\fixiv$ are twice continuously differentiable on $[0,1]$ and solve (\ref{pdif}). The statement for $\fixiii$ then follows by symmetry, while for the constant functions $0$ and $\fixv=1$ (see Proposition~\ref{fix2}), the claim is obvious. Since $\fixii,\fixiv$ are fixed points under $\Ui(\ov A,\gal,\gal)$, it follows from Lemma~\ref{smooth} and formula (\ref{Dov}) that $\fixii,\fixiv$ are continuous on $[0,1]$, twice continuously differentiable on $(0,1)$, and solve equation (\ref{pdif}) on $(0,1)$. We are done if we can show that their first and second derivatives can be extended to continuous functions on $[0,1]$. (If $f$ is twice continuously differentiable on $(0,1)$ and the limits $\lim_{x\to r}\dif{x}f(x)$ and $\lim_{x\to r}\diff{x}f(x)$ exists ($r=0,1$), then these limits coincide with the one-sided derivatives on the boundary. This follows, for example, from Corollary 6.3 in the appendix of \cite{EK}.)

Proposition~\ref{fix2} shows that $\fixii,\fixiv\leq 1$ and therefore, since they solve (\ref{pdif}) on $(0,1)$, $\fixii$ and $\fixiv$ are concave. Proposition~\ref{fix2} also shows that $\fixii(0)=\fixii(1)=0$ and $\fixiv(0)=0$, $\fixiv(1)=1$. (See Figure~\ref{ppla} as an illustration.) Since $\fixii$ is concave, $\dif{x}\fixii(x)$ increases to a limit in $(-\infty,\infty]$ as $x\down 0$. Lemma~\ref{absorb} implies that this limit is finite, and therefore $\dif{x}\fixii(x)$ is continuous at $x=0$. Since $\fixii$ solves (\ref{pdif}) on $(0,1)$,
\be\label{ddp2}
\lim_{x\to 0}\diff{x}\fixii(x)=-\lim_{x\to 0}\frac{2\gal \fixii(x)(1-\fixii(x))}{x(1-x)}=-2\gal\dif{x}\fixii(x)\big|_{x=0},
\ee
which proves that $\diff{x}\fixii(x)$ is continuous at $x=0$. The same argument proves that $\dif{x}\fixii(x)$ and $\diff{x}\fixii(x)$ are continuous at $x=1$, and that $\dif{x}\fixiv(x)$ and $\diff{x}\fixiv(x)$ are continuous at $x=0$. Since $\fixiv$ is concave, $\dif{x}\fixiv(x)$ decreases to a limit in $[-\infty,\infty)$ as $x\up 1$. Since $\fixiv(1)=1$ and $\fixiv\leq 1$, $\dif{x}\fixiv(x)\big|_{x=1}\geq 0$. Since $\fixiv$ solves (\ref{pdif}) on $(0,1)$ and $\dif{x}[\fixiv(x)(1-\fixiv(x))]\big|_{x=1}=-\dif{x}\fixiv(x)\big|_{x=1}$,
\be\label{ddp4}
\lim_{x\up 1}\diff{x}\fixiv(x)=-\lim_{x\up 1}\frac{2\gal \fixiv(x)(1-\fixiv(x))}{x(1-x)}=-2\gal\dif{x}\fixiv(x)\big|_{x=1},
\ee
which proves that $\dif{x}\fixiv(x)$ and $\diff{x}\fixiv(x)$ are continuous at $x=1$.\qed

\section[Renormalization branching process: embedded particles]{The renormalization branching process: embedded particle systems}\label{embsec}

In this section we use embedded particle systems to prove
Proposition~\ref{Uit}. An essential ingredient in the proofs
is Proposition~\ref{00hom}~(a), which will be proved in the Section~\ref{00sec}.

\subsection{Weighting and Poissonization}\label{Poissec}

{\bf Proof of Proposition~\ref{weightprop}} Obviously
$q^h_k\in\Ci_+(E^h)$ for each $k=1,\ldots,n$. Since $h\in\Ci_+(E)$ and
$h$ is bounded, it is easy to see that the map $\mu\mapsto h\mu$ from
$\Mi(E)$ into $\Mi(E^h)$ is continuous, and therefore the cluster
mechanisms defined in (\ref{Qih}) are continuous. Since
\be
\Ui^h_kf(x)=\frac{q_k(x)}{h(x)}E\big[1-\ex{-\li h\Zi_x,f\re}\big]
=\frac{\Ui_k(hf)(x)}{h(x)}\qquad(x\in E^h,\ f\in B_+(E^h)),
\ee
formula (\ref{htrafo1}) holds on $E^h$. To see that (\ref{htrafo1})
holds on $E\beh E^h$, note that by assumption $\Ui_kh\leq Kh$ for some
$K<\infty$, so if $x\in E\beh E^h$, then $\Ui_kh(x)=0$. By
monotonicity also $\Ui_k(hf)(x)=0$, while $h\Ui^h_kf(x)=0$ by
definition. Since $\sup_{x\in E^h}\Ui^h_k1(x)=\sup_{x\in E^h}
\frac{\Ui_kh(x)}{h(x)}\leq K<\infty$, the log-Laplace operators
$\Ui^h_k$ satisfy (\ref{fincon}). If $\Xc$ is started in an initial
state $\Xc_0$, then the Poisson-cluster branching process $\Xc^h$ with
log-Laplace operators $\Ui^h_1,\ldots,\Ui^h_n$ started in
$\Xc^h_0=h\Xc_0$ satisfies
\bc
\dis E\big[\ex{-\li h\Xc_k,f\re}\big]&
=&\dis E\big[\ex{-\li\Xc_0,\Ui_1\circ\cdots\circ\Ui_k(hf)\re}\big]\\[5pt]
&=&\dis E\big[\ex{-\li\Xc_0,h\Ui^h_1\circ\cdots\circ\Ui^h_k(f)\re}\big]
=E\big[\ex{-\li\Xc^h_k,f\re}\big]\qquad(f\in B_+(E^h)),
\ec
which proves (\ref{weighting}).\qed

\noi
{\bf Proof of Proposition~\ref{Poisprop}} We start by noting that by
(\ref{Vk}),
\be\label{Uiintu}
\Ui_kf(x)=q(x)E\big[1-\ex{-\li\Zi^k_x,f\re}\big]
=q_k(x)P[\Pois(f\Zi^k_x)\neq 0]\qquad(x\in E,\ f\in B_+(E)).
\ee
Into (\ref{Zxdef}), we insert
\be\ba{l}\label{condpois}
\dis P\big[\Pois(h\Zi^k_x)\in\cdot\,\big]\\[5pt]
\dis\quad=P\big[\Pois(h\Zi^k_x)\in\cdot\,\big|\,\Pois(h\Zi^k_x)\neq 0\big]
P[\Pois(h\Zi^k_x)\neq 0]+\de_0P[\Pois(h\Zi^k_x)=0].
\ec
Here and in similar formulas below, if in a conditional probability
the symbol $\Pois(\,\cdot\,)$ occurs twice with the same argument,
then it always refers to the same random variable (and not to
independent Poisson point measures with the same intensity, for
example). Using moreover (\ref{Uiintu}) we can rewrite (\ref{Zxdef})
as
\be\label{Zxfo}
Q^h_k(x,\,\cdot\,)=\frac{\Ui_k h(x)}{h(x)}
P\big[\Pois(h\Zi^k_x)\in\cdot\,\big|\,\Pois(h\Zi^k_x)\neq 0\big]
+\frac{h(x)-\Ui_k h(x)}{h(x)}\de_0(\,\cdot\,).
\ee
In particular, since we are assuming that $h$ is $\Ui_k$-subharmonic,
this shows that $Q^h_k(x,\,\cdot\,)$ is a probability measure. Let
$X^h$ be the branching particle system with offspring mechanisms
$Q^h_1,\ldots,Q^h_k$. Let $Z^{h,k}_x$ be random variables such that
$\Li(Z^{h,k}_x)=Q^h_k(x,\,\cdot\,)$. Then, by (\ref{VV}),
(\ref{Zxdef}), (\ref{elrel}), and (\ref{Uiintu}),
\be\ba{l}
\dis U^h_kf(x)=P[\Thin_f(Z^{h,k}_x)\neq 0]
=\frac{q_k(x)}{h(x)}P[\Thin_f(\Pois(h\Zi^k_x))\neq 0]\\[5pt]
\dis\qquad=\frac{q_k(x)}{h(x)}P[\Pois(hf\Zi^k_x)\neq 0]
=\frac{1}{h(x)}\Ui_k(hf)(x)\qquad(x\in E^h).
\ec
If $x\in E\beh E^h$, then $\Ui_k(hf)(x)\leq\Ui_k(h)(x)\leq
h(x)=0=:h\Ui^h(f)(x)$. This proves (\ref{htrafo}). To see that $Q^h_k$
is a {\em continuous} offspring mechanism, by
\cite[Theorem~4.2]{Kal76} it suffices to show that $x\mapsto\int
Q^h_k(x,\di\nu)\ex{-\li\nu,g\re}$ is continuous for all bounded
$g\in\Ci_+(E^h)$. Indeed, setting $f:=1-e^{-g}$, one has $\int
Q^h_k(x,\di\nu)\ex{-\li\nu,g\re}=\int
Q^h_k(x,\di\nu)(1-f)^\nu=1-\Ui^h_kf(x)=1-\Ui_k(hf)(x)/h(x)$ which is
continuous on $E^h$ by the continuity of $q_k$ and $\Qi_k$.

To see that also (\ref{Poissonization}) holds, just note that by
(\ref{V2}), (\ref{htrafo}), and (\ref{Vi2}),
\be\ba{l}
\dis P^{\Li(\Pois(h\mu))}[\Thin_f(X^h_n)=0]
=P[\Thin_{U^h_1\circ\cdots\circ U^h_nf}(\Pois(h\mu))=0]\\[5pt]
\dis\quad=P[\Pois((hU^h_1\circ\cdots\circ U^h_nf)\mu)=0]
=P[\Pois((\Ui_1\circ\cdots\circ\Ui_n(hf))\mu)=0]\\[5pt]
\dis\quad=P^\mu[\Pois(hf\Xc_n)=0]=P^\mu[\Thin_f(\Pois(h\Xc_n))=0].
\ec
Here $P^{\Li(\Pois(h\mu))}$ denotes the law of the process started
with initiallaw $\Li(\Pois(h\mu))$. Since this formula holds for all
$f\in B_{[0,1]}(E^h)$, formula (\ref{Poissonization}) follows.\qed

\brm{\bf(Boundedness of $h$)}
Propositions~\ref{weightprop} and \ref{Poisprop} generalize to the
case that $h$ is unbounded, except that in this case the cluster
mechanism in (\ref{Qih}) and the offspring mechanism in (\ref{Zxdef})
need in general not be continuous. Here, in order for (\ref{htrafo1})
and (\ref{htrafo}) to be well-defined, one needs to extend the
definition of $\Ui_kf$ to unbounded functions $f$, which can always
be done unambiguously (see Lemma~\ref{unbound}).
\erm

\subsection{Sub- and superharmonic functions}

This section contains a number of pivotal calculations involving the
log-Laplace operators $\Ui_\ga$ from (\ref{Uiga}). In particular, we
will prove that the functions $h_{1,1}$, $h_{0,0}$, and $h_{0,1}$ from
Lemmas~\ref{11lem}, \ref{00lem}, and \ref{01lem}, respectively, are
$\Ui_\ga$-superharmonic.

We start with an observation that holds for general log-Laplace operators.
\bl{\bf(Constant multiples)}\label{rhlem}
Let $\Ui$ be a log-Laplace operator of the form (\ref{Vk}) satisfying
(\ref{fincon}) and let $f\in B_+(E)$. Then $\Ui(rf)\leq r\Ui f$ for all
$r\geq 1$, and $\Ui(rf)\geq r\Ui f$ for all $0\leq r\leq 1$. In particular,
if $f$ is $\Ui$-superharmonic then $rf$ is $\Ui$-superharmonic for each
$r\geq 1$, and if $f$ is $\Ui$-subharmonic then $rf$ is $\Ui$-superharmonic
for each $0\leq r\leq 1$.
\el
{\bf Proof} If $\Xc$ is a branching process and $\Ui$ is the log-Laplace
operator of the transition law from $\Xc_0$ to $\Xc_1$ then, using Jensen's
inequality, for all $r\geq 1$,
\be\label{rhsup}
\ex{-\li\mu,\Ui(rf)\re}=E^\mu\big[\ex{-\li\Xc_1,rf\re}\big]
=E^\mu\big[\big(\ex{-\li\Xc_1,f\re}\big)^r\big]
\geq\big(E^\mu\big[\ex{-\li\Xc_1,f\re}\big]\big)^r=\ex{-\li\mu,r\Ui f\re}.
\ee
Since this holds for all $\mu\in\Mi(E)$, it follows that
$\Ui(rf)\leq r\Ui f$. The proof of the statements for $0\leq r\leq 1$
is the same but with the inequality signs reversed.\qed

\noi
We next turn our attention to the functions $h_{1,1}$ and $h_{0,0}$.
\bl{\bf(The catalyzing function $h_{1,1}$)}\label{sup11}
One has
\be\label{cofo}
\Ui_\ga(rh_{1,1})(x)=\frac{1+\ga}{\frac{1}{r}+\ga}\qquad(\ga,r>0,\ x\in[0,1]).
\ee
In particular, $h_{1,1}$ is $\Ui_\ga$-harmonic for each $\ga>0$.
\el
{\bf Proof} Recall (\ref{Zidef})--(\ref{Uiga}). Let $\sig_{1/r}$ be an
exponentially distributed random variable with mean $1/r$, independent
of $\tau_\ga$. Then
\be
\Ui_\ga(rh_{1,1})(x)=(\ffrac{1}{\ga}+1)
E\big[1-\ex{-\int_0^{\tau_\ga}r\di t}\big]
=(\ffrac{1}{\ga}+1)P[\sig_{1/r}<\tau_\ga]
=(\ffrac{1}{\ga}+1)\frac{\ga}{\frac{1}{r}+\ga},
\ee
which yields (\ref{cofo}).\qed

\bl{\bf(The catalyzing function $h_{0,0}$)}\label{sup00}
One has $\Ui_\ga(rh_{0,0})\leq rh_{0,0}$ for each $\ga,r>0$.
\el
{\bf Proof} Let $\Ga^\ga_x$ be the invariant law from Corollary~\ref{ergo}.
Then, for any $\ga>0$ and $f\in B_+[0,1]$,
\bc\label{linest}
\Ui_\ga f(x)&=&\dis(\ffrac{1}{\ga}+1)E\big[1-\ex{-\li\Zi^\ga_x,f\re}\big]
\leq(\ffrac{1}{\ga}+1)E[\li\Zi^\ga_x,f\re]\\[5pt]
&=&\dis(\ffrac{1}{\ga}+1)E\big[\int_0^{\tau_\ga}\!\!\!
f(\y^\ga_x(-t/2))\,\di t\big]=(1+\ga)\li\Ga^\ga_x,f\re\qquad(x\in[0,1]),
\ec
where we have used that $\tau_\ga$ is independent of $\y^\ga_x$ and has
mean $\ga$. In particular, setting $f=rh_{0,0}$ and using (\ref{WFfix})
we find that $\Ui_\ga(rh_{0,0})\leq rh_{0,0}$.\qed

\noi
The aim of the remainder of this section is to derive various bounds on
$\Ui_\ga f$ for $f\in\Hi_{0,1}$. We start with a formula for $\Ui_\ga f$
that holds for general $[0,1]$-valued functions $f$.
\bl{\bf(Action of $\Ui_\ga$ on $[0,1]$-valued functions)}\label{Ui01}
 Let $\y^\ga_x$ be the stationary solution to (\ref{Yclx}) and let
$\tau_{\ga/2}$ be an independent exponentially distributed random
variable with mean $\ga/2$. Let $(\bet_i)_{i\geq 1}$ be independent
exponentially distributed random variables with mean $\frac{1}{2}$,
independent of $\y^\ga_x$ and $\tau_{\ga/2}$, and let
$\sig_k:=\sum_{i=1}^k\bet_i$ $(k\geq 0)$. Then
\be\label{U1}
1-\Ui_\ga f(x)=E\Big[\prod_{k\geq 0:\;\sig_k<\tau_\ga}
\big(1-f(\y^\ga_x(-\sig_k))\big)\Big]
\qquad(\ga>0,\ f\in B_{[0,1]}[0,1],\ x\in[0,1]).
\ee 
\el
{\bf Proof} By Lemma~\ref{sup11}, the constant function
$h_{1,1}(x):=1$ satisfies $\Ui_\ga h_{1,1}=h_{1,1}$ for all $\ga>0$.
Therefore, by Proposition~\ref{Poisprop}, Poissonizing the
Poisson-cluster branching process $\Xc$ with the density $h_{1,1}$
yields a branching particle system
$X^{h_{1,1}}=(X^{h_{1,1}}_{-n},\ldots,X^{h_{1,1}}_0)$ with generating
operators $U^{h_{1,1}}_{\ga_{n-1}},\ldots,U^{h_{1,1}}_{\ga_0}$, where
\be
U^{h_{1,1}}_\ga f=\Ui_\ga f\qquad(f\in B_{[0,1]}[0,1],\ \ga>0).
\ee
By (\ref{VV}) and (\ref{Zxfo}),
\be
U^{h_{1,1}}_\ga f(x)=1-E\big[(1-f)^{\txt\Pois(\Zi^\ga_x)}\,\big|
\,\Pois(\Zi^\ga_x)\neq 0\big]\quad(f\in B_{[0,1]}[0,1],\ x\in[0,1],\ \ga>0).
\ee
Therefore, (\ref{U1}) will follow provided that
\be\label{conpos}
P\big[\Pois(\Zi^\ga_x)\in\cdot\,\big|\,\Pois(\Zi^\ga_x)\neq 0\big]
=\Li\Big(\sum_{k\geq 0:\;\sig_k<\tau_{\ga/2}}\de_{\y^\ga_x(-\sig_k)}\Big).
\ee
Indeed, it is not hard to see that
\be\label{unconpos}
\Pois(\Zi^\ga_x)\isd\sum_{k>0:\;\sig_k<\tau_{\ga/2}}\de_{\y^\ga_x(-\sig_k)}.
\ee
This follows from the facts that $\Zi^\ga_x=2\int_0^{\tau_{\ga/2}}
\de_{\y^\ga_x(-s)}\di s$ and
\be
\sum_{k>0:\;\sig_k<\tau_{\ga/2}}\de_{-\sig_k}
\isd\Pois(2\,1_{(-\tau_{\ga/2},0]}).
\ee
Conditioning $\Pois(2\,1_{(-\tau_{\ga/2},0]})$ on being nonzero means
conditioning on $\tau_{\ga/2}>\sig_1$. Since $\tau_{\ga/2}-\sig_1$,
conditioned on being nonnegative, is exponentially distributed with
mean $\ga/2$, using the stationarity of $\y^\ga_x$, we arrive at
(\ref{conpos}).\qed

\noi
The next lemma generalizes the duality (\ref{WFdual}) to mixed moments
of the Wright-Fisher diffusion $\y$ at multiple times. We can
interpret the left-hand side of (\ref{mixdua}) as the probability that
$m_1,\ldots,m_n$ organisms sampled from the population at times
$t_1,\ldots,t_n$ are all of the genetic type~I.
\bl{\bf(Sampling at multiple times)}\label{multidual}
Fix $0\leq t_1<\cdots<t_n=t$ and nonnegative integers $m_1,\ldots,m_n$.
Let $\y$ be the diffusion in (\ref{WFsde}). Then
\be\label{mixdua}
E^y\Big[\prod_{k=1}^n\y_{t_k}^{m_k}\Big]=E\big[y^{\phi_t}x^{\psi_t}\big],
\ee
where $(\phi_s,\psi_s)_{s\in[0,t]}$ is a Markov process in $\N^2$ started
in $(\phi_0,\psi_0)=(m_n,0)$, that jumps deterministically as
\be
(\phi_s,\psi_s)\to(\phi_s+m_k,\psi_s)\quad\mbox{at time}\quad t-t_k\quad(k<n),
\ee
and between these deterministic times jumps with rates as in (\ref{phipsi}).
\el
{\bf Proof} Induction, with repeated application of (\ref{WFdual}).\qed

\noi
For any $m\geq 1$, we put
\be
h_m(x):=1-(1-x)^m\qquad(x\in[0,1]).
\ee
The next lemma shows that we have particular good control on the action
of $\Ui_\ga$ on the functions $h_m$.
\bl{\bf(Action of $\Ui_\ga$ on the functions $h_m$)}\label{L:Uihm}
Let $m\geq 1$ and let $\tau_\ga$ be an exponentially distributed random
variable with mean $\ga$. Conditional on $\tau_\ga$, let
$(\phi'_t,\psi'_t)_{t\geq 0}$ be a Markov process in $\N^2$, started in
$(\phi'_0,\psi'_0)=(m,0)$ that jumps at time $t$ as:
\be\ba{r@{\,}c@{\,}l@{\qquad}l}\label{1phipsi}
(\phi'_t,\psi'_t)&\to&(\phi'_t-1,\psi'_t)
&\mbox{with rate}\ \phi'_t(\phi'_t-1),\\
(\phi'_t,\psi'_t)&\to&(\phi'_t-1,\psi'_t+1)
&\mbox{with rate}\ \ffrac{1}{\ga}\phi'_t,\\
(\phi'_t,\psi'_t)&\to&(\phi'_t+m,\psi'_t)
&\mbox{with rate}\ 1_{\{\tau_{\ga/2}<t\}}.
\ec
Then the limit $\lim_{t\to\infty}\psi'_t=:\psi'_\infty$ exists a.s., and
\be\label{suphest}
\Ui_\ga h_m(x)=E^{(m,0)}\big[1-(1-x)^{\psi'_\infty}\big]
\qquad(m\geq 1,\ x\in[0,1]).
\ee
\el
{\bf Proof} Let $\y^\ga_x$, $\tau_{\ga/2}$, and $(\sig_k)_{k\geq 0}$
be as in Lemma~\ref{Ui01}. Then, by (\ref{U1}),
\be\label{preUh01}
\Ui_\ga h_m(x)=1-E\Big[\prod_{k\geq 0:\;\sig_k<\tau_{\ga/2}}
\big(1-\y^\ga_x(-\sig_k)\big)^m\Big].
\ee
Let $(\phi',\psi')=(\phi'_t,\psi'_t)_{t\geq 0}$ be a $\N^2$-valued
process started in $(\phi'_0,\psi'_0)=(m,0)$ such that conditioned on
$\tau_\ga$ and $(\sig_k)_{k\geq 0}$, $(\phi',\psi')$ is a Markov
process that jumps deterministically as
\be
(\phi'_t,\psi'_t)\to(\phi'_t+m,\psi'_s)\quad\mbox{at time}
\quad\sig_k\quad(k\geq 1:\;\sig_k<\tau_{\ga/2})
\ee
and between these times jumps with rates as in (\ref{phipsi}). Then
$(\phi'_t,\psi'_t)\to(0,\psi'_\infty)$ as $t\to\infty$ a.s.\ for some
$\N$-valued random variable $\psi'_\infty$, and (\ref{suphest})
follows from Lemma~\ref{multidual}, using the symmetry
$y\leftrightarrow 1-y$. Since $\sig_{k+1}-\sig_k$ are independent
exponentially distributed random variables with mean one,
$(\phi',\psi')$ is the Markov process with jump rates as in
(\ref{1phipsi}).\qed

\noi
The next result is a simple application of Lemma~\ref{L:Uihm}.
\bl{\bf(The catalyzing function $h_1$)}\label{subx}
The function $h_1(x):=x$ $(x\in[0,1])$ is $\Ui_\ga$-subharmonic
for each $\ga>0$.
\el
{\bf Proof} Since $\psi'_\infty\geq 1$ a.s., one has
$1-(1-x)^{\psi'_\infty}\geq x$ a.s.\ $(x\in[0,1])$ in
(\ref{suphest}). In particular, setting $m=1$ yields
$\Ui_\ga h_1\geq h_1$.\qed

\noi
We now set out to prove that $h_7$, which is the function $h_{0,1}$
from Lemma~\ref{01lem}, is $\Ui_\ga$-super\-harmonic. In order to do
so, we will derive upper bounds on the expectation of $\psi'_\infty$.
We derive two estimates: one that is good for small $\ga$ and one that
is good for large $\ga$.

In order to avoid tedious formal arguments, it will be convenient to
recall the interpretation of the process $(\phi',\psi')$ and
Lemma~\ref{multidual}. Recall from the discussion following
(\ref{WFdual}) that $(\y^\ga_x(t))_{t\in\R}$ describes the equilibrium
frequency of genetic type $I$ as a function of time in a population
that is in genetic exchange with an infinite reservoir. From this
population we sample at times $-\sig_k$ ($k\geq 0$,
$\sig_k<\tau_{\ga/2}$) each time $m$ individuals, and ask for the
probability that they are not all of the genetic type~II. In order to
find this probability, we follow the ancestors of the sampled
individuals back in time. Then $\phi'_t$ and $\psi'_t$ are the number
of ancestors that lived at time $-t$ in the population and the
reservoir, respectively, and $E[1-(1-x)^{\psi'_\infty}]$ is the
probability that at least one ancestor is of type~I.
\bl{\bf(Bound for small $\ga$)}\label{smaga}
For each $\ga\in(0,\infty)$ and $m\geq 1$,
\be\label{smag}
\frac{1}{m}E^{(m,0)}[\psi'_\infty]\leq\frac{1}{m}
\sum_{i=0}^{m-1}\frac{1+\ga}{1+i\ga}=:\chi_m(\ga).
\ee
The function $\chi_m$ is concave and satisfies $\chi_m(0)=1$
for each $m\geq 1$.
\el
{\bf Proof} Note that
\be
E\big[\big|\{k\geq 0:\;\sig_k<\tau_{\ga/2}\}\big|\big]=1+\ga.
\ee
We can estimate $(\phi',\psi')$ from above by a process where
ancestors from individuals sampled at different times cannot
coalesce. Therefore,
\be\label{aces}
E^{(m,0)}[\psi'_\infty]\leq(1+\ga)E^{(m,0)}[\psi_\infty],
\ee
where $(\phi,\psi)$ is the Markov process in (\ref{phipsi}).
Note that if $(\phi,\psi)$ is in the state $(m+1,0)$, then the
next jump is to $(m,1)$ with probability
\be
\frac{\frac{1}{\ga}(m+1)}{\frac{1}{\ga}(m+1)+m(m+1)}=\frac{1}{1+m\ga}
\ee
and to $(m,0)$ with one minus this probability. Therefore,
\bc
\dis E^{(m+1,0)}[\psi_\infty]&=&\dis\frac{1}{1+m\ga}
E^{(m,1)}[\psi_\infty]+\Big(1-\frac{1}{1+m\ga}\Big)
E^{(m,0)}[\psi_\infty]\\[5pt]
&=&\dis\frac{1}{1+m\ga}\Big(E^{(m,0)}[\psi_\infty]+1\Big)
+\Big(1-\frac{1}{1+m\ga}\Big)E^{(m,0)}[\psi_\infty]\\[5pt]
&=&\dis E^{(m,0)}[\psi_\infty]+\frac{1}{1+m\ga}.
\ec
By induction, it follows that
\be\label{psim}
E^{(m,0)}[\psi_\infty]=\sum_{i=0}^{m-1}\frac{1}{1+i\ga}.
\ee
Inserting this into (\ref{aces}) we arrive at (\ref{smag}). Finally, since
\be
\Diff{\ga}\frac{1+\ga}{1+i\ga}=\frac{2i(i-1)}{(1+i\ga)^3}\geq 0
\qquad(i\geq 0,\ \ga\geq 0),
\ee
the function $\chi_m$ is convex.\qed

\bl{\bf(Bound for large $\ga$)}\label{laga}
For each $\ga\in(0,\infty)$ and $m\geq 1$,
\be\label{lag}
E^{(m,0)}[\psi'_\infty]\leq(\ffrac{1}{\ga}+1)\sum_{k=1}^m\frac{1}{k}
+\frac{3}{2}.
\ee
\el
{\bf Proof} We start by observing that $\dif{t}E[\psi'_t]
=\frac{1}{\ga}E[\phi'_t]$, and therefore
\be\label{Epsi}
E[\psi'_\infty]=\ffrac{1}{\ga}\int_0^\infty E[\phi'_t]\di t.
\ee
Unlike in the proof of the last lemma, this time we cannot fully
ignore the coalescence of ancestors sampled at different times. In
order to deal with this we use a trick: at time zero we introduce an
extra ancestor that can only jump to the reservoir when
$t\geq\tau_\ga$ and there are no other ancestors left in the
population. We further assume that all other ancestors do not jump to
the reservoir on their own. Let $\xi_t$ be one as long as this extra
ancestor is in the population and zero otherwise, and let $\phi''_t$
be the number of other ancestors in the population according to these
new rules. Then we have at a Markov process $(\xi,\phi'')$ started in
$(\xi_0,\phi''_0)=(1,m)$ that jumps as:
\be\ba{r@{\,}c@{\,}l@{\qquad}l}\label{2phipsi}
(\xi_t,\phi''_t)&\to&(\xi_t,\phi''_t-1)
&\mbox{with rate}\ (\phi''_t+1)\phi''_t,\\
(\xi_t,\phi''_t)&\to&(\xi_t,\phi''_t+m)
&\mbox{with rate}\ 1_{\{\tau_{\ga/2}<t\}},\\
(\xi_t,\phi''_t)&\to&(\xi_t-1,\phi''_t)
&\mbox{with rate}\ \frac{1}{\ga}1_{\{\tau_{\ga/2}\geq t\}}1_{\{\phi''_t=0\}}.
\ec
It is not hard to show that $(\xi,\phi'')$ and $\phi'$ can be coupled
such that $\xi_t+\phi''_t\geq\phi'_t$ for all $t\geq 0$. We now
simplify even further and ignore all coalescence between ancestors
belonging to the process $\phi''$ that are introduced at different
times. Let $\phi^{(k)}_t$ be the number of ancestors in the population
that were introduced at the time $\sig_k$ $(k\geq 0)$. Thus, for
$t<\sig_k$ one has $\phi^{(k)}_t=0$, for $t=\sig_k$ one has
$\phi^{(k)}_t=m$, while for $t>\sig_k$, the process $\phi^{(k)}_t$
jumps from $n$ to $n-1$ with rate $(n+1)n$. Then it is not hard to see
that, for an appropriate coupling, $\phi''_t\leq\sum_{k\geq 0
:\sig_k<\tau_{\ga/2}}\phi^{(k)}_t$ for all $t\geq 0$. We let $\xi'$
be a process such that $\xi'_0=1$ and $\xi'_t$ jumps to zero with rate
\be
\frac{1}{\ga}1_{\{\tau_{\ga/2}\geq t\}}
\prod_{k\geq 0:\sig_k<\tau_{\ga/2}}1_{\{\phi^{(k)}_t=0\}}.
\ee
Then for an appropriate coupling $\xi'_t\geq\xi_t$ $(t\geq 0)$.
Thus, we can estimate
\be\label{phint}
\int_0^\infty E[\phi'_t]\di t\leq\int_0^\infty E[\xi'_t]\di t
+\int_0^\infty E\Big[\sum_{k\geq 0:\sig_k<\tau_{\ga/2}}\phi^{(k)}_t\Big]\di t.
\ee
Set $\rho:=\inf\{t\geq\tau_{\ga/2}:\phi^{(k)}_t=0\ \forall k\geq 0
\mbox{ with }\sig_k<\tau_{\ga/2}\}$ and $\pi:=\inf\{t\geq 0:\xi'_t=0\}$.
Then
\be\label{xint}
\int_0^\infty E[\xi'_t]\di t=E[\tau_{\ga/2}]+E[\rho-\tau_{\ga/2}]
+E[\pi-\rho]=\frac{3}{2}\ga+E[\rho-\tau_{\ga/2}].
\ee
Since
\bc
\dis E[\rho-\tau_{\ga/2}]&\leq&\dis\int_0^\infty
E\Big[1_{\{\sum_{k\geq 0:\sig_k<\tau_{\ga/2}}
\phi^{(k)}_t\neq 0\}}\Big]\di t\\[10pt]
&\leq&\dis\int_0^\infty E\Big[\sum_{k\geq 0:\sig_k<\tau_{\ga/2}}
1_{\{\phi^{(k)}_t\neq 0\}}\Big]\di t,
\ec
using moreover (\ref{phint}) and (\ref{xint}), we can estimate
\be
\int_0^\infty E[\phi'_t]\di t\leq\frac{3}{2}\ga
+\int_0^\infty E\Big[\sum_{k\geq 0:\sig_k<\tau_{\ga/2}}
(\phi^{(k)}_t+1_{\{\phi^{(k)}_t\neq 0\}})\Big]\di t.
\ee
Since $E\big[\big|\{k\geq 0:\;\sig_k<\tau_{\ga/2}\}\big|\big]=1+\ga$, we obtain
\be\label{obta}
\int_0^\infty E[\phi'_t]\di t\leq\frac{3}{2}\ga
+(1+\ga)\int_0^\infty E[\phi^{(0)}_t+1_{\{\phi^{(0)}_t\neq 0\}}]\di t.
\ee
Since $\phi^{(0)}_t$ jumps from $n$ to $n-1$ with rate $(n+1)n$, the
expected total time
that $\phi^{(0)}_t=n$ equals $1/((n+1)n)$, and therefore
\be
\int_0^\infty E[\phi^{(0)}_t+1_{\{\phi^{(0)}_t\neq 0\}}]\di t
=\sum_{n=1}^m\frac{1}{(n+1)n}(n+1_{\{n\neq 0\}})
=\sum_{n=1}^m\frac{1}{n}.
\ee
Inserting this into (\ref{obta}), using (\ref{Epsi}), we arrive
at (\ref{lag}).\qed

\bl{\bf(The catalyzing function $h_{0,1}$)}\label{sup01}
One has $\Ui_\ga(h_{0,1})\leq h_{0,1}$ for each $\ga>0$. Moreover,
for each $r>1$ and $\ga>0$,
\be\label{strongsup}
\sup_{x\in(0,1]}\frac{\Ui_\ga(rh_{0,1})(x)}{rh_{0,1}(x)}<1.
\ee
\el
{\bf Proof} Recall that $h_{0,1}(x)=h_7(x)=1-(1-x)^7$ $(x\in[0,1])$.
We will show that
\be
E^{(7,0)}[\psi'_\infty]<7
\ee
for each $\ga\in(0,\infty)$. The function $\chi_m(\ga)$ from
Lemma~\ref{smaga} satisfies
\be
\chi_m(1)=\frac{1}{m}\sum_{n=1}^m\frac{2}{n}<1\qquad(m\geq 5).
\ee
Since $\chi_m(\ga)$ is concave in $\ga$ and satisfies $\chi_m(0)=1$,
it follows that $\chi_m(\ga)<1$ for all $0<\ga\leq 1$ and $m\geq 5$.
By Lemma~\ref{laga}, for all $\ga\geq 1$,
\be
E^{(m,0)}[\psi'_\infty]\leq 2\sum_{k=1}^m\frac{1}{k}+\frac{3}{2}<m
\qquad(m\geq 7).
\ee
Therefore, if $m\geq 7$, then $m':=E^{(m,0)}[\psi'_\infty]<m$. It
follows by (\ref{suphest}) and Jensen's inequality applied to the
concave function $z\mapsto 1-(1-x)^z$ that
\be\label{hmjens}
\Ui_\ga h_m(x)\leq 1-(1-x)^{E^{(m,0)}[\psi'_\infty]}=1-(1-x)^{m'}
\leq h_m(x)\qquad(x\in[0,1],\ \ga>0).
\ee
This shows that $h_m$ is $\Ui_\ga$-superharmonic for each $\ga>0$.
By Lemma~\ref{rhlem}, for each $r>1$,
\be\label{leftb}
\frac{\Ui_\ga(rh_m)(x)}{rh_m(x)}\leq\frac{r\Ui_\ga(h_m)(x)}{rh_m(x)}
\leq\frac{1-(1-x)^{m'}}{1-(1-x)^m}\qquad(x\in(0,1]).
\ee
By Lemma~\ref{sup11} and the monotonicity of $\Ui_\ga$,
\be\label{rightb}
\frac{\Ui_\ga(rh_m)(x)}{rh_m(x)}\leq\frac{\Ui_\ga(r)(x)}{rh_m(x)}
\leq\frac{1+\ga}{1+r\ga}\frac{1}{(1-(1-x)^m)}\qquad(x\in(0,1]).
\ee
Since the right-hand side of (\ref{leftb}) is smaller than $1$ for
$x\in(0,1)$ and tends to $m'/m<1$ as $x\to 0$, since the right-hand
side of (\ref{rightb}) is smaller than $1$ for $x$ in an open
neighborhood of $1$, and since both bounds are continuous,
(\ref{strongsup}) follows.\qed

\subsection{Extinction versus unbounded growth}\label{versus}

In this section we show that Lemmas~\ref{11lem}--\ref{01lem} are
equivalent to Proposition~\ref{P:exgr}. (This follows from the
equivalence of conditions (i) and (ii) in Lemma~\ref{exgro} below.) We
moreover prove Lemmas~\ref{11lem} and \ref{01lem} and prepare for the
proof of Lemma~\ref{00lem}. We start with some general facts about
log-Laplace operators and branching processes.

For the next lemma, let $E$ be a separable, locally compact,
metrizable space. For $n\geq 0$, let $q_n\in\Ci_+(E)$ be continuous
weight functions, let $\Qi_n$ be continuous cluster mechanisms on $E$,
and assume that the associated log-Laplace operators $\Ui_n$ defined
in (\ref{Vk}) satisfy (\ref{fincon}). Assume that $0\neq h\in\Ci_+(E)$
is bounded and $\Ui_n$-superharmonic for all $n$, let $E^h:=\{x\in
E:h(x)>0\}$, and define generating operators $U^h_n:B_{[0,1]}(E^h)\to
B_{[0,1]}(E)$ as in (\ref{htrafo}). For each $n\geq 0$, let
$(\Xc^{(n)}_0,\Xc^{(n)}_1)$ be a one-step Poisson cluster branching
process with log-Laplace operator $\Ui_n$, and let
$(X^{(n),h}_0,X^{(n),h}_1)$ be the one-step branching particle system
with generating operator $U^h_n$. (In a typical application of this
lemma, the operators $\Ui_n$ will be iterates of other log-Laplace
operators, and $\Xc^{(n)}_0,\Xc^{(n)}_1$ will be the initial and final
state, respectively, of a Poisson cluster branching process with many
time steps.)
\bl{\bf(Extinction versus unbounded growth)}\label{exgro}
Assume that $\rho\in\Ci_{[0,1]}(E^h)$ and put
\be
p(x):=\left\{\ba{ll}h(x)\rho(x)\quad&\mbox{if}\ x\in E^h,\\
0\quad&\mbox{if}\ x\in E\beh E^h.\ea\right.
\ee
Then the following statements are equivalent:
\[\ba{rl}
{\rm (i)}&\dis P^{\de_x}\big[|X^{(n),h}_1|\in\cdot\,\big]\Asto{n}
\rho(x)\de_\infty+(1-\rho(x))\de_0\\
&\hspace{1.5cm}\mbox{locally uniformly for }x\in E^h,\\[5pt]
{\rm (ii)}&\dis P^{\de_x}\big[\li\Xc^{(n)}_1,h\re\in\cdot\,\big]\Asto{n}
\ex{-p(x)}\de_0+\big(1-\ex{-p(x)}\big)\de_\infty\\
&\hspace{1.5cm}\mbox{locally uniformly for }x\in E,\\[5pt]
{\rm (iii)}&\Ui_n(\la h)(x)\asto{n}p(x)\\
&\hspace{1.5cm}\mbox{locally uniformly for }x\in E\quad\forall\la>0,\\[5pt]
{\rm (iv)}&\exists0<\la_1<\la_2<\infty:\quad\Ui_n(\la_i h)(x)\asto{n}p(x)\\
&\hspace{1.5cm}\mbox{locally uniformly for }x\in E\qquad(i=1,2).
\ea\]
\el
{\bf Proof of Lemma~\ref{exgro}} It is not hard to see that (i) is
equivalent to
\be
P^{\de_x}[\Thin_\la(X^{(n),h}_1)\neq 0]\asto{n}\rho(x)
\ee
locally uniformly for $x\in E^h$, for all $0<\la\leq 1$. It follows from
(\ref{V2}) and (\ref{htrafo}) that $h(x)P^{\de_x}[\Thin_\la(X^{(n),h}_1)\neq 0]
=hU^h(\la)(x)=\Ui(\la h)(x)$ $(x\in E)$, so (i) is equivalent to
\[\ba{l}
{\rm (i)'}\quad\Ui_n(\la h)(x)\asto{n}p(x)\\
\hspace{1.5cm}\mbox{locally uniformly for }x\in E\quad\forall0<\la\leq 1.
\ea\]
By (\ref{Vi}), condition~(ii) implies that
\be
\ex{-\Ui_n(\la h)(x)}=E^{\de_x}\big[\ex{-\la\li\Xc_1,h\re}\big]\asto{n}
\ex{-p(x)}
\ee
locally uniformly for $x\in E$ for all $\la>0$, and therefore (ii)
implies (iii). Obviously (iii)$\volgt{\rm (i)'}\volgt$(iv) so we are
done if we show that (iv)$\volgt$(ii). Indeed, (iv) implies that
\be
E^{\de_x}\big[\ex{-\la_1\li\Xc^{(n)}_1,h\re}
-\ex{-\la_2\li\Xc^{(n)}_1,h\re}\big]\asto{n}0
\ee
locally uniformly for $x\in E$, which shows that
\be
P^{\de_x}\big[c<\li\Xc^{(n)}_1,h\re<C\big]\asto{n}0
\ee
for all $0<c<C<\infty$. Using (iv) once more we arive at (ii).\qed

\noi
Our next lemma gives sufficient conditions for the $n$-th iterates of
a single log-Laplace operator $\Ui$ to satisfy the equivalent
conditions of Lemma~\ref{exgro}. Let $E$ (again) be separable, locally
compact, and metrizable. Let $q\in\Ci_+(E)$ be a weight function,
$\Qi$ a continuous cluster mechanism on $E$, and assume that the
associated log-Laplace operator $\Ui$ defined in (\ref{Vk}) satisfies
(\ref{fincon}). Let $\Xc=(\Xc_0,\Xc_1,\ldots)$ be the Poisson-cluster
branching process with log-Laplace operator $\Ui$ in each step, let
$0\neq h\in\Ci_+(E)$ be bounded and $\Ui$-superharmonic, and let
$X^h=(X^h_0,X^h_1,\ldots)$ denote the branching particle system on
$E^h$ obtained from $\Xc$ by Poissonization with a $\Ui$-superharmonic
function $h$, in the sense of Proposition~\ref{Poisprop}.
\bl{\bf(Sufficient condition for extinction versus unbounded
growth)}\label{L:sufexgro}
Assume that
\be\label{strsup}
\sup_{x\in E^h}\frac{\Ui h(x)}{h(x)}<1.
\ee
Then the process $X^h$ started in any initial law $\Li(X^h_0)\in\Mi_1(E^h)$
satisfies
\be\label{extgrow}
\lim_{k\to\infty}|X^h_k|=\infty\quad\mbox{or}\quad\exists k\geq 0
\mbox{ s.t.\ }X^h_k=0\qquad\as
\ee
Moreover, if the function $\rho:E^h\to[0,1]$ defined by
\be
\rho(x):=P^{\de_x}[X^h_n\neq 0\quad\forall n\geq 0]\qquad(x\in E^h)
\ee
satisfies $\inf_{x\in E^h}\rho(x)>0$, then $\rho$ is continuous.
\el
{\bf Proof of Lemma~\ref{L:sufexgro}} Let $\Ai$ denote the tail event
$\Ai=\{X^h_n\neq 0\ \forall n\geq 0\}$ and let $(\Fi_k)_{k\geq 0}$ be
the filtration generated by $X^h$. Then, by the Markov property and
continuity of the conditional expectation with respect to increasing
limits of \si-fields (see Complement 10(b) from
\cite[Section~29]{Loe63} or \cite[Section~32]{Loe78})
\be
P[X^h_n\neq 0\ \forall n\geq 0|X_k]=P(\Ai|\Fi_k)\asto{k}1_\Ai\qquad\as
\ee
In particular, this implies that a.s.\ on the event $\Ai$ one must
have $P[X^h_{k+1}=0|X^h_k]\to 0$ a.s. By (\ref{V2}) and
(\ref{htrafo}), $P^{\de_x}[X^h_1\neq 0]=U^h 1(x)=(\Ui h(x))/h(x)$,
which is uniformly bounded away from one by (\ref{strsup}). Therefore,
$P[X^h_{k+1}=0|X^h_k]\to 0$ a.s.\ on $\Ai$ is only possible if the
number of particles tends to infinity.

The continuity of $\rho$ can be proved by a straightforward adaptation
of the proof of \cite[Proposition~5 (d)]{FStrim} to the present
setting with discrete time and noncompact space $E$. An essential
ingredient in the proof, apart from (\ref{strsup}), is the fact that
the map $\nu\mapsto P^\nu[X^h_n\in\cdot\,]$ from $\Ni(E)$ to
$\Mi_1(\Ni(E))$ is continuous, which follows from the continuity of
$Q^h$.\qed

\noi
We now turn our attention more specifically to the renormalization
branching process $\Xc$. In the remainder of this section,
$(\ga_k)_{k\geq 0}$ is a sequence of positive constants such that
$\sum_n\ga_n=\infty$ and $\ga_n\to\ga^\ast$ for some
$\ga^\ast\in\half$, and $\Xc=(\Xc_{-n},\ldots,\Xc_0)$ is the Poisson
cluster branching process on $[0,1]$ defined in Section~\ref{RBP}. We
put $\Ui^{(n)}:=\Ui_{\ga_{n-1}}\circ\cdots\circ\Ui_{\ga_0}$. If $0\neq
h\in\Ci[0,1]$ is $\Ui_{\ga_k}$-superharmonic for all $k\geq 0$, then
$\Xc^h$ and $X^h$ denote the branching process and the branching
particle system on $\{x\in[0,1]:h(x)>0\}$ obtained from $\Xc$ by
weighting and Poissonizing with $h$ in the sense of
Propositions~\ref{weightprop} and \ref{Poisprop}, respectively.\med

\noi
{\bf Proof of Lemma~\ref{11lem}} By induction, it follows from
Lemma~\ref{sup11} that
\be
\Ui^{(n)}(\la h_{1,1})=\frac{\prod_{k=0}^{n-1}(1+\ga_k)}
{\prod_{k=0}^{n-1}(1+\ga_k)-1+\frac{1}{\la}}\qquad(\la>0).
\ee
It is not hard to see (compare the footnote at (\ref{sumga})) that
\be
\prod_{k=0}^\infty(1+\ga_k)=\infty\quad\mbox{if and only if}
\quad\sum_{k=0}^\infty\ga_k=\infty.
\ee
Therefore, since we are assuming that $\sum_n\ga_n=\infty$,
\be
\Ui^{(n)}(\la h_{1,1})\asto{n}h_{1,1},
\ee
uniformly on $[0,1]$ for all $\la>0$. The result now follows from
Lemma~\ref{exgro} (with $h=h_{1,1}$ and $\rho(x)=1$ $(x\in[0,1])$).\qed

\brm{\bf(Conditions on $(\ga_n)_{n\geq 0}$)}
Our proof of Lemma~\ref{11lem} does not use that $\ga_n\to\ga^\ast$
for some $\ga^\ast\in\half$. On the other hand, the proof shows that
$\sum_n\ga_n=\infty$ is a necessary condition for (\ref{expl}).
\erm
We do not know if the assumption that $\ga_n\to\ga^\ast$ for some
$\ga^\ast\in\half$ is needed in Lemma~\ref{00lem}. We guess that
it can be dropped, but it will greatly simplify proofs to have it around.

We will show that in order to prove Lemmas~\ref{00lem} and
\ref{01lem}, it suffices to prove their analogues for embedded
particle systems in the time-homogeneous processes $\Yi^{\ga^\ast}$
($\ga^\ast\in\half$). More precisely, we will derive
Lemmas~\ref{00lem} and \ref{01lem} from the following two results.
Below, $(\Ui^0_t)_{t\geq 0}$ is the log-Laplace semigroup of the
super-Wright-Fisher diffusion $\Yi^0$, defined in (\ref{cau}). The
functions $p^\ast_{0,1,\ga^\ast}$ ($\ga^\ast\in\half$) are defined in
(\ref{pdef}).
\bp{\bf(Time-homogeneous embedded particle system with
$h_{0,0}$)}\label{00hom}\smallskip

\noi
{\bf (a)} For any $\ga^\ast>0$, one has $\dis(\Ui_{\ga^\ast})^n h_{0,0}
\asto{n}0$ uniformly on $[0,1]$.\smallskip

\noi
{\bf (b)} One has $\dis\Ui^0_th_{0,0}\asto{t}0$ uniformly on $[0,1]$.
\ep

\bp{\bf(Time-homogeneous embedded particle system with $h_{0,1}$)}\label{01hom}

\noi
{\bf (a)} For any $\ga^\ast>0$, one has
$(\Ui_{\ga^\ast})^n(\la h_{0,1})\asto{n}p^\ast_{0,1,\ga^\ast}$
uniformly on $[0,1]$, for all $\la>0$.\med

\noi
{\bf (b)} One has
$\Ui^0_t(\la h_{0,1})\asto{t}p^\ast_{0,1,0}$ uniformly on $[0,1]$,
for all $\la>0$.
\ep
Propositions~\ref{00hom}~(b) and \ref{01hom}~(b) follow from Proposition~\ref{Ucor}. Proposition~\ref{00hom}~(a) will be proved in Section~\ref{002sec}.\med

\noi
{\bf Proof of Proposition~\ref{01hom}~(a)} By formula (\ref{strongsup})
from Lemma~\ref{sup01}, for each $r>1$ the function $rh_{0,1}$ satisfies
condition (\ref{strsup}) from Lemma~\ref{L:sufexgro}. Set $\rho(x)
:=P^{\de_x}[Y^{\ga^\ast,rh_{0,1}}_n\neq 0\ \forall n]$. Then, by
(\ref{V2}) and (\ref{htrafo}),
\bc\label{rholow}
\dis\rho(x)&=&\dis\lim_{n\to\infty}P^{\de_x}[Y^{\ga^\ast,rh_{0,1}}_n\neq 0]
=\lim_{n\to\infty}(U^{rh_{0,1}}_{\ga^\ast})^n1(x)\\[5pt]
&=&\dis\lim_{n\to\infty}\frac{(\Ui_{\ga^\ast})^n(rh_{0,1})(x)}{rh_{0,1}(x)}
\geq\frac{h_1(x)}{rh_{0,1}(x)}\qquad(x\in(0,1]),
\ec
where $h_1(x)=x$ $(x\in[0,1])$ is the $\Ui_{\ga^\ast}$-subharmonic function
from Lemma~\ref{subx}. It follows that $\inf_{x\in(0,1]}\rho(x)>0$ and
therefore, by Lemma~\ref{L:sufexgro}, $\rho$ is continuous in $x$.

By Lemma~\ref{L:sufexgro}, we see that the Poissonized particle system
$X^{rh_{0,1}}$ exhibits extinction versus unbounded growth in the sense
of Lemma~\ref{exgro}, which implies the statement in
Proposition~\ref{01hom}~(a).\qed

\noi
We now show that Propositions~\ref{00hom} and \ref{01hom} imply
Lemmas~\ref{00lem} and \ref{01lem}, respectively.\med

\noi
{\bf Proof of Lemma~\ref{00lem}} We start with the proof that the
embedded particle system $X^{h_{0,0}}$ is critical. For any
$f\in B_+[0,1]$ and $k\geq 1$, we have, by Poissonization
(Proposition~\ref{Poisprop}) and the definition of $\Xc$,
\be\ba{l}\label{critcalc}
\dis h_{0,0}(x)E^{-k,\de_x}[\li X^{h_{0,0}}_{-k+1},f\re]
=E^{-k,\Li(\Pois(h_{0,0}\de_x))}[\li X^{h_{0,0}}_{-k+1},f\re]
=E^{-k,\de_x}[\li\Pois(h_{0,0}\Xc_{-k+1}),f\re]\\[5pt]
\dis\qquad=E^{-k,\de_x}[\li\Xc_{-k+1},h_{0,0}f\re]
=(\ffrac{1}{\ga}+1)E[\li\Zi^\ga_x,h_{0,0}f\re]
=(\ffrac{1}{\ga}+1)\li\Ga^{\ga_{k-1}}_x,h_{0,0}f\re,
\ec
where $\Ga^\ga_x$ is the invariant law of $\y^\ga_x$ from
Corollary~\ref{ergo}. In particular, setting $f=1$ gives
$h_{0,0}(x)E^{-k,\de_x}[|X^{h_{0,0}}_{-k+1}|]=h_{0,0}(x)$ by (\ref{WFfix}).

To prove (\ref{ext}), by Lemma~\ref{exgro} it suffices to show that
\be\label{TB}
\Ui^{(n)}(\la h_{0,0})\asto{n}0
\ee
uniformly on $[0,1]$ for all $0<\la\leq 1$. We first treat the case
$\ga^\ast>0$. Then, by Theorem~\ref{supercon}~(a), for each fixed
$l\geq 1$ and $f\in\Ci_+[0,1]$,
\be\label{Uitohom}
\Ui_{\ga_{n-1}}\circ\cdots\circ\Ui_{\ga_{n-l}}f\asto{n}(\Ui_{\ga^\ast})^lf
\ee
uniformly on $[0,1]$. Therefore, by a diagonal argument, we can find
$l(n)\to\infty$ such that
\be
\big\|(\Ui_{\ga^\ast})^{l(n)}h_{0,0}-\Ui_{\ga_{n-1}}\circ\cdots\circ
\Ui_{\ga_{n-l(n)}}h_{0,0}\big\|_\infty\asto{n}0.
\ee
Using the fact that the function $h_{0,0}$ is $\Ui_\ga$-superharmonic
for each $\ga>0$ and the monotonicity of the operators $\Ui_\ga$, we
derive from Proposition~\ref{00hom}~(a) that
\be
\Ui^{(n)}(\la h_{0,0})\leq\Ui_{\ga_{n-1}}\circ\cdots\circ
\Ui_{\ga_{n-l(n)}}h_{0,0}\asto{n}0
\ee
uniformly on $[0,1]$ for all $0<\la\leq 1$. This proves (\ref{TB})
in the case $\ga^\ast>0$.

The proof in the case $\ga^\ast=0$ is similar. In this case, by
Theorem~\ref{supercon}~(b), for each fixed $t>0$ and $f\in\Ci_+[0,1]$,
\be\label{Uitohom2}
\Ui_{\ga_{n-1}}\circ\cdots\circ\Ui_{\ga_{k_n(t)}}f(x_n)\asto{n}
\Ui^0_tf(x)\qquad\forall x_n\to x\in[0,1],
\ee
which shows that $\Ui_{\ga_{n-1}}\circ\cdots\circ\Ui_{\ga_{k_n(t)}}f$
converges to $\Ui^0_tf$ uniformly on $[0,1]$. By a diagonal argument,
we can find $t(n)\to\infty$ such that
\be
\big\|\Ui^0_t(h_{0,0})-\Ui_{\ga_{n-1}}\circ\cdots\circ
\Ui_{\ga_{k_n(t(n))}}(h_{0,0})\big\|_\infty\asto{n}0,
\ee
and the proof proceeds in the same way as before.\qed

\noi
{\bf Proof of Lemma~\ref{01lem}} By Lemma~\ref{exgro} and the monotonicity
of the operators $\Ui_\ga$ it suffices to show that
\be\ba{rl}\label{TB2}
{\rm (i)}&\dis\limsup_{n\to\infty}\Ui^{(n)}(h_{0,1})
\leq p^\ast_{0,1,\ga^\ast},\\[5pt]
{\rm (ii)}&\dis\liminf_{n\to\infty}\Ui^{(n)}(\ffrac{1}{2}h_{0,1})
\geq p^\ast_{0,1,\ga^\ast},
\ec
uniformly on $[0,1]$. We first consider the case $\ga^\ast>0$.
By (\ref{Uitohom}) and a diagonal argument, we can find $l(n)\to\infty$
such that
\be
\big\|(\Ui_{\ga^\ast})^{l(n)}h_{0,1}-\Ui_{\ga_{n-1}}\circ\cdots\circ
\Ui_{\ga_{n-l(n)}}h_{0,1}\big\|_\infty\asto{n}0.
\ee
Therefore, by Proposition~\ref{01hom}~(a), the fact that $h_{0,1}$
is $\Ui_{\ga_k}$-superharmonic for each $k\geq 0$, and the monotonicity
of the operators $\Ui_\ga$, we find that
\be
\Ui^{(n)}h_{0,1}\leq\Ui_{\ga_{n-1}}\circ\cdots\circ
\Ui_{\ga_{n-l(n)}}h_{0,1}\asto{n}p^\ast_{0,1,\ga^\ast},
\ee
uniformly on $[0,1]$. This proves (\ref{TB2})~(i). To prove also
(\ref{TB2})~(ii) we use the $\Ui_\ga$-{\em sub}harmonic (for each
$\ga>0$) function $h_1$ from Lemma~\ref{subx}. By Lemma~\ref{rhlem}
also $\ffrac{1}{2}h_1$ is $\Ui_\ga$-subharmonic. By bounding
$\ffrac{1}{2}h_1$ from above and below with multiples of $h_{0,1}$ it
is easy to derive from Proposition~\ref{01hom}~(a) that
\be\label{Ughl}
(\Ui_{\ga^\ast})^n(\ffrac{1}{2}h_1)\asto{n}p^\ast_{0,1,\ga^\ast}
\ee
uniformly on $[0,1]$. Arguing as before, we can find $l(n)\to\infty$ such that
\be
\big\|(\Ui_{\ga^\ast})^{l(n)}(\ffrac{1}{2}h_1)-\Ui_{\ga_{n-1}}
\circ\cdots\circ\Ui_{\ga_{n-l(n)}}(\ffrac{1}{2}h_1)\big\|_\infty\asto{n}0.
\ee
Therefore, by (\ref{Ughl}) and the facts that $\ffrac{1}{2}h_1$ is
$\Ui_{\ga_k}$-subharmonic for each $k\geq 0$ and
$\ffrac{1}{2}h_1\leq\ffrac{1}{2}h_{0,1}$,
\be
\Ui^{(n)}(\ffrac{1}{2}h_{0,1})\geq\Ui_{\ga_{n-1}}\circ\cdots\circ
\Ui_{\ga_{n-l(n)}}(\ffrac{1}{2}h_1)\asto{n}p^\ast_{0,1,\ga^\ast},
\ee
uniformly on $[0,1]$, which proves (\ref{TB2})~(ii). The proof
of (\ref{TB2}) in case $\ga^\ast=0$ is completely analogous.\qed

\section[Renormalization branching process: extinction on interior]{The renormalization branching process: extinction on the interior}\label{00sec}

\subsection{Basic facts}

In this section we prove Proposition~\ref{00hom}~(a). To simplify
notation, throughout this section $h$ denotes the function $h_{0,0}$.
We fix $0<\ga^\ast<\infty$, we let $Y^h:=Y^{\ga^\ast,h}$ denote the
branching particle system on $(0,1)$ obtained from
$\Yi^{\ga^\ast}=(\Yi^{\ga^\ast}_0,\Yi^{\ga^\ast}_1,\ldots)$ by
Poissonization with $h$ in the sense of Proposition~\ref{Poisprop},
and we denote its log-Laplace operator by $U^h_{\ga^\ast}$. We will
prove that
\be\label{extin}
\rho(x):=P^{\de_x}\big[Y^h_n\neq 0\ \forall n\geq 0\big]=0\qquad(x\in(0,1)).
\ee
Since for each $n$ fixed, $x\mapsto\rho_n(x):=P^{\de_x}[Y^h_n\neq 0]$
is a continuous function that decreases to $\rho(x)$, (\ref{extin})
implies that $\rho_n(x)\to 0$ locally uniformly on $(0,1)$, which, by
an obvious analogon of Lemma~\ref{exgro}, yields
Proposition~\ref{00hom}~(a).

As a first step, we prove:
\bl{\bf (Continuous survival probability)}\label{rholem}
One has either $\rho(x)=0$ for all $x\in(0,1)$ or there exists a
continuous function $\ti\rho:(0,1)\to[0,1]$ such that
$\rho(x)\geq\ti\rho(x)>0$ for all $x\in(0,1)$.
\el
{\bf Proof} Put $p(x):=h(x)\rho(x)$. We will show that either
$p=0$ on $(0,1)$ or there exists a continuous function
$\ti p:(0,1)\to(0,1]$ such that $p\geq\ti p$ on $(0,1)$. Indeed,
\be\ba{l}
\dis p(x)=h(x)P^{\de_x}\big[Y^h_n\neq 0\ \forall n\geq 0\big]
=\lim_{n\to\infty}h(x)P^{\de_x}\big[Y^h_n\neq 0\big]\\[5pt]
\dis\quad=h(x)\lim_{n\to\infty}(U^h_{\ga^\ast})^n1(x)
=\lim_{n\to\infty}(\Ui_{\ga^\ast})^nh(x)\qquad(x\in(0,1)),
\ec
where we have used (\ref{V2}) and (\ref{htrafo}) in the last two steps.
Using the continuity of $\Ui_{\ga^\ast}$ with respect to decreasing
sequences, it follows that
\be\label{fixp}
\Ui_{\ga^\ast}p=p.
\ee

We claim that for any $f\in B_{[0,1]}[0,1]$, one has the bounds
\be\label{Gast}
\li\Ga^\ga_x,f\re\leq\Ui_\ga f(x)\leq(1+\ga)\li\Ga^\ga_x,f\re
\qquad(\ga>0,\ x\in[0,1]).
\ee
Indeed, by Lemma~\ref{Ui01}, $\Ui_\ga f(x)\geq 1-E[(1-f(\y^\ga_x(0)))]
=\li\Ga^\ga_x,f\re$, while the upper bound in (\ref{Gast}) follows
from (\ref{linest}).

By Remark~\ref{R:betadis}, $(0,1)\ni x\mapsto\li\Ga^\ga_x,f\re$ is
continuous for all $f\in B_{[0,1]}[0,1]$. Moreover,
$\li\Ga^\ga_x,f\re=0$ for some $x\in(0,1)$ if and only if $f=0$ almost
everywhere with respect to Lebesgue measure.

Applying these facts to $f=p$ and $\ga=\ga^\ast$, using (\ref{fixp}),
we see that there are two possibilities. Either $p=0$ a.s.\ with
respect to Lebesgue measure, and in this case $p=0$ by the upper bound
in (\ref{Gast}), or $p$ is not almost everywhere zero with respect to
Lebesgue measure, and in this case the function $x\mapsto\ti
p(x):=\li\Ga^\ga_x,f\re$ is continuous, positive on $(0,1)$, and
estimates $p$ from below by the lower bound in (\ref{Gast}).\qed

\subsection{A representation for the Campbell law}\label{002sec}

(Local) extinction properties of critical branching processes are
usually studied using Palm laws. Our proof of formula (\ref{extin}) is
no exception, except that we will use the closely related Campbell
laws. Loosely speaking, Palm laws describe a population that is
size-biased at a given position, plus `typical' particle sampled from
that position, while Campbell laws describe a population that is
size-biased as a whole, plus a `typical' particle sampled from a
random position.

Let $\Pc$ be a probability law on $\Ni(0,1)$ with
$\int_{\Ni(0,1)}\Pc(\di\nu)|\nu|=1$. Then the {\em size-biased law}
$\Pc_{\rm size}$ associated with $\Pc$ is the probability law on
$\Ni(0,1)$ defined by
\be
\Pc_{\rm size}(\,\cdot\,)
:=\int_{\Ni(0,1)}\Pc(\di\nu)\,|\nu|1_{\txt\{\nu\in\cdot\,\}}.
\ee
The {\em Campbell law} associated with $\Pc$ is the probability
law on $(0,1)\times\Ni(0,1)$ defined by
\be\label{Campdef}
\Pc_{\rm Camp}(A\times B):=\int_{\Ni(0,1)}\Pc(\di\nu)\,\nu(A)
1_{\txt\{\nu\in B\}}
\ee
for all Borel-measurable $A\sub(0,1)$ and $B\sub\Ni(0,1)$. If $(v,V)$
is a $(0,1)\times\Ni(0,1)$-valued random variable with law
$\Pc_{\rm Camp}$, then $\Li(V)=\Pc_{\rm size}$, and $v$ is the position
of a `typical' particle chosen from $V$.

Let
\be
\Pc^{x,n}(\,\cdot\,):=P^{\de_x}\big[Y^h_n\in\cdot\,]
\ee
denote the law of $Y^h$ at time $n$, started at time $0$ with one
particle at position $x\in(0,1)$. Note that by criticality,
$\int_{\Ni(0,1)}\Pc^{x,n}(\di\nu)|\nu|=1$. Using again criticality, it
is easy to see that in order to prove the extinction formula
(\ref{extin}), it suffices to show that
\be\label{Sexpl}
\lim_{n\to\infty}\Pc^{x,n}_{\rm size}\big(\{1,\ldots,N\}\big)=0
\qquad(x\in(0,1),\ N\geq 1).
\ee
In order to prove (\ref{Sexpl}), we will write down an expression for
$\Pc^{x,n}_{\rm Camp}$. Let $Q^h$ denote the offspring mechanism of
$Y^h$, and, for fixed $x\in(0,1)$, let $Q^h_{\rm Camp}(x,\,\cdot\,)$
denote the Campbell law associated with $Q^h(x,\,\cdot\,)$. The next
proposition is a time-inhomogeneous version of Kallenberg's famous
backward tree technique; see \cite[Satz~8.2]{Lie81}.
\bp{\bf(Representation of Campbell law)}\label{P:Camprep}
Let $(\vb_k,V_k)_{k\geq 0}$ be the Markov process in $(0,1)\times\Ni(0,1)$
with transition laws
\be
P\big[(\vb_{k+1},V_{k+1})\in\cdot\,\big|\,(\vb_k,V_k)=(x,\nu)\big]
=Q^h_{\rm Camp}(x,\,\cdot\,)\qquad((x,\nu)\in(0,1)\times\Ni(0,1)),
\ee
started in $(\vb_0,V_0)=(\de_x,0)$. Let $(Y^{h,(k)})^{k\geq 1}$ be
branching particle systems with offspring mechanism $Q^h$, conditionally
independent given $(\vb_k,V_k)_{k\geq 0}$, started in
$Y^{h,(k)}_0=V_k-\de_{\vb_k}$. Then
\be\label{Camprep}
\Pc^{x,n}_{\rm Camp}=\Li\Big(\vb_n,\de_{\vb_n}
+\sum_{k=1}^n Y^{h,(k)}_{n-k}\Big).
\ee
\ep
Formula (\ref{Camprep}) says that the Campbell law at time $n$ arises
in such a way, that an `immortal' particle at positions
$\vb_0,\ldots,\vb_n$ sheds off offspring
$V_1-\de_{\vb_1},\ldots,V_n-\de_{\vb_n}$, distributed according to the
size-biased law with one `typical' particle taken out, and this
offspring then evolve under the usual forward dynamics till time $n$.
Note that the position of the immortal particle $(\vb_k)_{k\geq 0}$ is
an autonomous Markov chain.

We need a bit of explicit control on $Q^h_{\rm Camp}$.
\bl{\bf(Campbell law)}\label{L:Campexpl}
One has
\be\label{Campexpl}
Q^h_{\rm Camp}(x,A\times B)=\frac{\frac{1}{\ga^\ast}+1}{h(x)}
\int P[\Pois(h\Zi^{\ga^\ast}_x)\in\di\chi]\chi(A)1_{\{\chi\in A\}},
\ee
where the random measures $\Zi^{\ga^\ast}_x$ are defined in (\ref{Zidef}).
\el
{\bf Proof} By the definition of the Campbell law (\ref{Campdef}),
and (\ref{Zxdef}),
\bc
\dis Q^h_{\rm Camp}(x,A\times B)
&=&\dis\int Q^h(x,\di\chi)\chi(A)1_{\{\chi\in B\}}\\[5pt]
&=&\dis\frac{\frac{1}{\ga^\ast}+1}{h(x)}
\int P[\Pois(h\Zi^{\ga^\ast}_x)\in\di\chi]\chi(A)1_{\{\chi\in B\}}
+\Big(1-\frac{\frac{1}{\ga^\ast}+1}{h(x)}\Big)\cdot 0.
\ec
\qed

\noi
Recall that by (\ref{Zidef}),
\be
\Zi^{\ga^\ast}_x:=\int_0^{\tau_{\ga_\ast}}\de_{\y^{\ga^\ast}_x(-t/2)}\di t,
\ee
where $(\y^{\ga^\ast}_x(t))_{t\in\R}$ is a stationary solution to the
SDE (\ref{Yclx}) with $\ga=\ga^\ast$. By Lemma~\ref{L:Campexpl}, the
transition law of the Markov chain $(\vb_k)_{k\geq 0}$ from
Proposition~\ref{P:Camprep} is given by
\be\label{vbtr}
P[\vb_{k+1}\in\di y|\vb_k=x]=\frac{\frac{1}{\ga^\ast}+1}{h(x)}
E[\Pois(h\Zi^{\ga^\ast}_x)(\di y)]
=\frac{1+\ga^\ast}{h(x)}h(y)\Ga^{\ga^\ast}_x(\di y),
\ee
where $\Ga^{\ga^\ast}_x$ is the invariant law of $\y^{\ga^\ast}_x$
from Corollary~\ref{ergo}. In the next section we will prove the
following lemma.
\bl{\bf(Immortal particle stays in interior)}\label{stayin}
The Markov chain $(\vb_k)_{k\geq 0}$ started in any $\vb_0=x\in(0,1)$ satisfies
\be
(\vb_k)_{k\geq 0}\mbox{ has a cluster point in }(0,1)\quad\as
\ee
\el
We now show that Lemma~\ref{stayin}, together with our previous results,
implies Proposition \ref{00hom}~(a).\med

\noi
{\bf Proof of Proposition~\ref{00hom}~(a)} We need to prove (\ref{extin}).
By our previous analysis, it suffices to prove (\ref{Sexpl}) under the
assumption that $\rho\neq 0$. By Proposition~\ref{P:Camprep},
\be\label{sizerep}
\Pc^{x,n}_{\rm size}=\Li\Big(\de_{\vb_n}+\sum_{k=1}^n Y^{h,(k)}_{n-k}\Big).
\ee
Conditioned on $(\vb_k,V_k)_{k\geq 0}$, the $(Y^{h,(k)}_{n-k})_{k=1,\ldots,n}$
are independent random variables with
\be
P\big[Y^{h,(k)}_{n-k}\neq 0\big]
\geq P\big[Y^{h,(k)}_m\neq 0\ \forall m\geq 0\big]
=P[\Thin_\rho(V_k-\de_{\vb_k})\neq 0].
\ee
Therefore, (\ref{Sexpl}) will follow by Borel-Cantelli provided
that we can show that
\be\label{BC}
\sum_{k=1}^\infty P[\Thin_\rho(V_k-\de_{\vb_k})\neq 0|\vb_{k-1}]=\infty\quad\as
\ee
Define $f(x):=P[\Thin_\rho(V_k-\de_{\vb_k})\neq 0|\vb_{k-1}=x]$
$(x\in(0,1))$. We need to show that $\sum_{k=1}^\infty f(x)=\infty$
a.s. Using Lemma~\ref{rholem} and Lemma~\ref{L:Campexpl} we can estimate
\be\label{notnil}
f(x)\geq P[\Thin_{\ti\rho}(V_k-\de_{\vb_k})\neq 0|\vb_{k-1}=x]
=\int_{\Ni(0,1)}Q^h_{\rm Camp}(x,\di y,\di\nu)
\{1-(1-\ti\rho)^{\nu-\de_y}\big\}>0
\ee
for all $x\in(0,1)$. Since $\Qi_{\ga^\ast}$, defined in (\ref{Qqdef}),
is a continuous cluster mechanism, also $Q^h_{\rm Camp}(x,\cdot)$ is
continuous as a function of $x$, hence the bound in (\ref{notnil}) is
locally uniform on $(0,1)$, hence Lemma~\ref{stayin} implies that
there is an $\eps>0$ such that
\be
P[\Thin_\rho(V_k-\de_{\vb_k})\neq 0|\vb_{k-1}]\geq\eps
\ee
at infinitely many times $k-1$, which in turn implies (\ref{BC}).\qed

\subsection{The immortal particle}

{\bf Proof of Lemma~\ref{stayin}} Let $K(x,\di y)$ denote the transition
kernel (on $(0,1)$) of the Markov chain $(\vb_k)_{k\geq 0}$, i.e., by
(\ref{vbtr}),
\be
K(x,\di y)=(1+\ga^\ast)\frac{y(1-y)}{x(1-x)}\Ga^{\ga^\ast}_x(\di y).
\ee
It follows from (\ref{WFmoments}) that
\be\label{Kyy}
\int K(x,\di y)y(1-y)
=\frac{x(1-x)+\ga^\ast(1+\ga^\ast)}{(1+2\ga^\ast)(1+3\ga^\ast)}.
\ee
\detail{Indeed,
\[\ba{l}
\dis\int K(x,\di y)y(1-y)=\frac{1+\ga}{x(1-x)}
\int\Ga^\ga_x(\di y)\{y^2-2y^3+y^4\}\\[5pt]
\dis\quad=\frac{1+\ga}{x(1-x)}\Big\{\left(\frac{x+0}{1+0}\right)
\left(\frac{x+\ga}{1+\ga}\right)-2\left(\frac{x+0}{1+0}\right)
\left(\frac{x+\ga}{1+\ga}\right)\left(\frac{x+2\ga}{1+2\ga}\right)\\
\dis\qquad\phantom{=}+\left(\frac{x+0}{1+0}\right)
\left(\frac{x+\ga}{1+\ga}\right)\left(\frac{x+2\ga}{1+2\ga}\right)
\left(\frac{x+3\ga}{1+3\ga}\right)\Big\}\\[5pt]
\dis\quad=\frac{x(x+\ga)}{x(1-x)(1+2\ga)(1+3\ga)}
\Big\{(1+2\ga)(1+3\ga)-2(x+2\ga)(1+3\ga)+(x+2\ga)(x+3\ga)\Big\}\\[5pt]
\dis\quad=\frac{x+\ga}{(1-x)(1+2\ga)(1+3\ga)}
\Big\{(1+2\ga)(1+3\ga)-2(((x-1)+1+2\ga)(1+3\ga)\\
\dis\quad\phantom{=\frac{x+\ga}{(1-x)(1+2\ga)(1+3\ga)}
\Big\{}+((x-1)+1+2\ga)((x-1)+1+3\ga)\Big\}\\[5pt]
\dis\quad=\frac{x+\ga}{(1-x)(1+2\ga)(1+3\ga)}
\Big\{(x-1)\big\{-2(1+3\ga)+(1+2\ga)+(1+3\ga)\big\}+(x-1)^2\Big\}\\[5pt]
\dis\quad=\frac{(x+\ga)\{\ga+(1-x)\}}{(1+2\ga)(1+3\ga)}
=\frac{x(1-x)+\ga(1+\ga)}{(1+2\ga)(1+3\ga)}.
\ea\]}
Set
\be
g(x):=\int K(x,\di y)y(1-y)-x(1-x)\qquad(x\in(0,1)).
\ee
Then
\be
M_n:=\vb_n(1-\vb_n)-\sum_{k=0}^{n-1}g(\vb_k)\qquad(n\geq 0)
\ee
defines a martingale $(M_n)_{n\geq 0}$. Since $g>0$ in an open
neighborhood of $\{0,1\}$,
\be
P[(\vb_k)_{k\geq 0}\mbox{ has no cluster point in }(0,1)]\leq
P[\lim_{n\to\infty}M_n=-\infty]=0,
\ee
where in the last equality we have used that $(M_n)_{n\geq 0}$
is a martingale.\qed

\section{Proof of the main result}\label{final}

{\bf Proof of Theorem~\ref{main}} Part~(a) has been proved in
Section~\ref{dualsec}. It follows from (\ref{sumga}), (\ref{gamma}),
(\ref{hutform}), and (\ref{renbra}) that part~(b) is equivalent to
the following statement. Assuming that
\be\label{varsga}
{\rm (i)}\quad\sum_{n=1}^\infty\ga_n=\infty\qquad\mbox{and}
\qquad{\rm (ii)}\quad\ga_n\asto{n}\ga^\ast
\ee
for some $\ga^\ast\in\half$, one has, uniformly on $[0,1]$,
\be\label{varwast}
\Ui_{\ga_{n-1}}\circ\cdots\circ\Ui_{\ga_0}(p)\asto{n}p^\ast_{l,r,\ga^\ast},
\ee
where $p^\ast_{l,r,\ga^\ast}$ is the unique solution in $\Hi_{l,r}$ of
\be\ba{rr@{\,}c@{\,}ll}\label{varwiga}
{\rm (i)}&\Ui_{\ga^\ast}p^\ast&=&p^\ast\quad
&\mbox{if }0<\ga^\ast<\infty,\\[5pt]
{\rm (ii)}&\ffrac{1}{2}x(1-x)\diff{x}p^\ast(x)-p^\ast(x)(1-p^\ast(x))
&=&0\quad(x\in[0,1])\quad&\mbox{if }\ga^\ast=0.
\ec
It follows from Proposition~\ref{Uit} that the left-hand side of
(\ref{varwast}) converges uniformly to a limit $p^\ast_{l,r,\ga^\ast}$
which is given by (\ref{pdef}). We must show $1^\circ$ that
$p^\ast_{l,r,\ga^\ast}\in\Hi_{l,r}$ and $2^\circ$ that
$p^\ast_{l,r,\ga^\ast}$ is the unique solution in this class to
(\ref{varwiga}). We first treat the case $\ga^\ast>0$.

$1^\circ$ Since $p^\ast_{0,0,\ga^\ast}\equiv 0$ and
$p^\ast_{1,1,\ga^\ast}\equiv 1$, it is obvious that
$p^\ast_{0,0,\ga^\ast}\in\Hi_{0,0}$ and
$p^\ast_{1,1,\ga^\ast}\in\Hi_{1,1}$. Therefore, by symmetry, it
suffices to show that $p^\ast_{0,1,\ga^\ast}\in\Hi_{0,1}$. By
Lemmas~\ref{subx} and \ref{sup01}, $x\leq p\leq 1-(1-x)^7$ implies
$x\leq\Ui_{\ga_k}p\leq 1-(1-x)^7$ for each $k$. Iterating this
relation, using (\ref{varwast}), we find that
\be\label{lowup}
x\leq p^\ast_{0,1,\ga^\ast}(x)\leq 1-(1-x)^7.
\ee
By Proposition~\ref{moncon}, the left-hand side of (\ref{varwast}) is
nondecreasing and concave in $x$ if $p$ is, so taking the limit we
find that $p^\ast_{0,1,\ga^\ast}$ is nondecreasing and concave.
Combining this with (\ref{lowup}) we conclude that
$p^\ast_{0,1,\ga^\ast}$ is Lipschitz continuous. Moreover
$p^\ast_{0,1,\ga^\ast}(0)=0$ and $p^\ast_{0,1,\ga^\ast}(1)=1$ so
$p^\ast_{0,1,\ga^\ast}\in\Hi_{0,1}$.

$2^\circ$ Taking the limit $n\to\infty$ in
$(\Ui_{\ga^\ast})^{n}p=\Ui_{\ga^\ast}(\Ui_{\ga^\ast})^{n-1}p$, using
the continuity of $\Ui_{\ga^\ast}$ (Corollary~\ref{Ugacont}) and
(\ref{varwast}), we find that
$\Ui_{\ga^\ast}p^\ast_{l,r,\ga^\ast}=p^\ast_{l,r,\ga^\ast}$. It
follows from (\ref{varwast}) that $p^\ast_{l,r,\ga^\ast}$ is the only
solution in $\Hi_{l,r}$ to this equation.

For $\ga^\ast=0$, it has been shown in \cite[Proposition~3]{FSsup}
that $p^\ast_{l,r,0}$ is the unique solution in $\Hi_{l,r}$ to
(\ref{varwiga})~(ii). In particular, it has been shown there that
$p^\ast_{0,1,0}$ is twice continuously differentiable on $[0,1]$
(including the boundary). This proves parts~(b) and (c) of the
theorem.\qed

\chapter{Branching-coalescing particle systems.}\label{C:braco}

\section{Introduction and main results}

\subsection{Introduction}

In this chapter we study systems of particles subject to a stochastic
dynamics with the following description. $1^\circ$ Each particle
moves independently of the others according to a continuous time
Markov process on a lattice $\La$, which jumps from site $i$ to site
$j$ with rate $a(i,j)$. $2^\circ$ Each particle splits with rate
$b\geq 0$ into two new particles, created on the position of the old
one. $3^\circ$ Each pair of particles, present on the same site,
coalesces with rate $2c$ (with $c\geq 0$) to one particle. $4^\circ$
Each particle dies with rate $d\geq 0$. Throughout this chapter, we make
the following assumptions.
\begin{enumerate}
\item $\La$ is a finite or countably infinite set.
\item The transition rates $a(i,j)$ are irreducible, i.e., if $\De\sub\La$
is neither $\La$ nor $\emptyset$,\\ then there exist $i\in\De$ and
$j\in\La\beh\De$ such that $a(i,j)>0$ or $a(j,i)>0$.\label{irred}
\item $\sup_i\sum_ja(i,j)<\infty$.
\item\label{doub} $\sum_ja^\dgg(i,j)=\sum_ja(i,j)$,
where $a^\dgg(i,j):=a(j,i)$.
\item $b,c$, and $d$ are nonnegative constants.
\end{enumerate}
Here and elsewhere sums and suprema over $i,j$ always run over $\La$,
unless stated otherwise. Assumption~(iv) says that the counting
measure is an invariant \si-finite measure for the Markov process with
jump rates $a$. With respect to this invariant measure, the
time-reversed process jumps from $i$ to $j$ with rate $a^\dgg(i,j)$.

Let $X_t(i)$ denote the number of particles present at site $i\in\La$
and time $t\geq 0$. Then $X=(X_t)_{t\geq 0}$, with
$X_t=(X_t(i))_{i\in\La}$, is a Markov process with formal generator
\bc\label{GdefC3}
Gf(x)&:=&\dis\sum_{ij}a(i,j)x(i)\{f(x+\de_j-\de_i)-f(x)\}
+b\sum_ix(i)\{f(x+\de_i)-f(x)\}\\
&&\dis+c\sum_ix(i)(x(i)-1)\{f(x-\de_i)-f(x)\}
+d\sum_ix(i)\{f(x-\de_i)-f(x)\},
\ec
where $\de_i(j):=1$ if $i=j$ and $\de_i(j):=0$ otherwise. The process
$X$ can be defined for finite initial states and also for some infinite
initial states in an appropriate Liggett-Spitzer space (see
Section~\ref{prelim}). We call $(X_t)_{t\geq 0}$ a branching
coalescing particle system with underlying motion $(\La,a)$,
branching rate $b$, coalescence rate $c$ and death rate $d$, or
shortly the $(a,b,c,d)$-braco-process.

Some typical examples of underlying motions we have in mind are
nearest neighbour random walk on $\La=\Z^d$ and on $\La=\T^d$, the
homogeneous tree of degree $d+1$. We will not restrict ourselves to
symmetric underlying motions (i.e., $a=a^\dgg$) but also allow
$a(i,j)=1_{\{j=i+1\}}$ on $\Z$, for example. The reason why we do not
restrict ourselves to graphs, is that we also want to include the case
$\La=\om_d$, the hierarchical group with freedom $d$, i.e.,
\be
\om_d:=\{i=(i_0,i_1,\ldots):i_\al\in\{0,\ldots,d-1\}\ \forall
 \al\geq 0,\ i_\al\neq 0\mbox{ finitely often}\,\},
\ee
equipped with componentwise addition modulo $n$. On $\om_d$, one
typically chooses transition rates $a(i,j)$ that depend only on the
hierarchical distance $|i-j|:=\min\{\al\geq 0:i_\bet=j_\bet\
\forall\bet\geq\al\}$. The hierarchical group has found widespread
applications in population biology and is therefore a natural choice
for the underlying space.

\subsection{Motivation}

Our motivation for studying branching-coalescing particle systems
comes from three directions.

{\em Reaction diffusion models, Schl\"ogl's first model.}
Branching-coalescing particle systems are known in the physics
literature as a reaction diffusion models. More precisely, our model
is a special case of Schl\"ogl's first model \cite{Sch72}, where in
the latter there is an additional rate with which particles are
spontaneously created. For $d=0$, our model is known as the
autocatalytic reaction. Reaction diffusion models have been studied
intensively by physicists and more recently also by probabilists
\cite{DDL90,Mou92,Neu90}. All work that we are aware of is restricted
to the case $\La=\Z^d$.

{\em Population dynamics, the contact process.} Branching-coalescing
particle systems may be thought of as a more or less realistic model
for the spread and growth of a population of organisms. Here, the
underlying motion models the migration of organisms, births and deaths
have their obvious interpretations, while coalescence of particles
should be thought of as additional deaths, caused by local
overpopulation. In this respect, our model is similar to the contact
process. The latter is often referred to as a model for the spread of
an infection, but in fact it is a reasonable model for the population
dynamics of many organisms, from trees in a forest to killer
bees. There are two striking differences between the contact process
and branching-coalescing particle systems.  First, whereas the total
population at one site is subject to a rigid bound in the contact
process (namely one), it may reach arbitrarily high values in a
branching-coalescing system.  However, when the local population is
high, the coalescence (which grows quadratically in the number of
organisms) dominates the branching (which grows linearly), and in this
way the population is reduced. A second difference is that in the
contact process, if one site infects its neighbor, the original site
is still infected. As opposed to this, even when the death rate is
zero, it is possible that a branching coalescing particle system goes
to local extinction due to migration only. Thus, we can say that the
gain from infection is guaranteed in the contact process, whereas the
reward for migration is uncertain in a branching-coalescing particle
system.

{\em Resampling with selection and negative mutations.} Our third
motivation also comes from population dynamics, but from a different
perspective. Assume that at each site $i\in\La$ there lives a large,
fixed number of organisms, and that each of these organisms carries a
gene that comes in two types: a healthy and a defective one. Let us
model the evolution of the population as follows. $1^\circ$ with rate
$a(i,j)$, we let an organism at site $i$ migrate to site $j$.
$2^\circ$ to model the effect of natural selection, we let each
organism with rate $b$ choose another organism, living on the same
site. If the first organism carries a healthy gene and the second
organism a defective gene, then the latter is replaced by an organism
with a healthy gene. $3^\circ$ to model the effect of random mating,
we resample each pair of organisms living at the same site with rate
$2c$, i.e., we choose one of the two at random and replace it by an
organism with the type of the other one. $4^\circ$ with rate $d$, we
let a healthy gene mutate into a defective gene. In the limit that the
number of organisms at each site is large, the frequencies $\Xc_t(i)$
of healthy organisms at site $i$ and time $t$ are described by the
unique pathwise solution to the infinite dimensional stochastic
differential equation (SDE) (see \cite{SU86}):
\bc\label{sde}
\di \Xc_t(i)&=&\dis\sum_ja(j,i)(\Xc_t(j)-\Xc_t(i))\,\di t
+b\Xc_t(i)(1-\Xc_t(i))\,\di t-d\Xc_t(i)\,\di t\\
&&\dis+\sqrt{2c\Xc_t(i)(1-\Xc_t(i))}\,\di B_t(i)
\qquad\qquad(t\geq 0,\ i\in\La).
\ec
We call the $[0,1]^\La$-valued process $\Xc=(\Xc_t)_{t\geq 0}$ the
resampling-selection process with underlying motion $(\La,a)$,
selection rate $b$, resampling rate $c$ and mutation rate $d$, or
shortly the $(a,b,c,d)$-resem-process (the letters in `resem' standing
for \underline{re}sampling, \underline{se}lection and
\underline{m}utation).

It is known that branching-coalescing particle systems are dual to
resampling-selection processes. To be precise, for any
$\phi\in[0,1]^\La$ and $x\in\N^\La$, write
\be\label{mux}
\phi^x:=\prod_i\phi(i)^{x(i)},
\ee
where $0^0:=1$. Let $\Xc$ be the $(a,b,c,d)$-resem-process and let $X^\dgg$
be the $(a^\dgg,b,c,d)$-braco-process. Then (see Theorem~\ref{duth}~(a)
below)
\be\label{sdedual}
E^\phi[(1-\Xc_t)^x]=E^x[(1-\phi)^{X^\dgg_t}].
\ee
Formula (\ref{sdedual}) has the following interpretation:
$E^\phi[(1-\Xc_t)^x]$ is the probability that $x$ organisms, sampled
from the population at time $t$, all have defective genes. If we want
to calculate this probability, we must follow back in time those
organisms that could possibly be healthy ancestors of these $x$
organisms. In this way we end up with a system of branching coalescing
$a^\dgg$-random walks, which die when a mutation occurs, coalesce when
two potential ancestors descend from the same ancestor, and branch
when a selection event takes place. If we end up with at least one
healthy potential ancestor at time zero, then we know that not all the
$x$ particles have defective genes.

Resampling-selection processes of the form (\ref{sde}) are also known
as {\em stepping stone models} (with selection and one type of
mutation). These were studied by Shiga and Uchiyama in \cite{SU86}, a
paper similar in spirit to ours. The duality (\ref{sdedual}) is a
special case of Lemma~2.1 \cite{SU86}. Moment duals for genetic
diffusions in a more general but non-spatial context go back to
\cite{Shi81}. The idea of incorporating selection in resampling models
by introducing branching into the usual coalescent dual seems to have
been independently reinvented in \cite{KN97}. They were probably the
first to interpret the duality (\ref{sdedual}) in terms of potential
ancestors. For some recent versions of this duality, see also
\cite{DK99,DG99,BES02}. A SDE that is dual to branching-{\em
annihilating} random walks occurs in \cite[Lemma~2.1]{BEM03}. A SPDE
version of (\ref{sde}) (with $d=0$) has been derived as the rescaled
limit of long-range biased voter models in \cite[Theorem~2]{MT95}.

Note that for $c=0$, the process $\Xc$ is deterministic. In this case,
the semigroup $(U_t)_{t\geq 0}$ defined by $U_t\phi:=\Xc_t$ ($t\geq
0$), where $\Xc$ is the deterministic solution of (\ref{sde}) with
initial state $\Xc_0=\phi\in[0,1]^\La$, is called the generating
semigroup of the branching particle system $X^\dgg$. (For this
terminology, see for example \cite{FStrim}.) Thus, the duality relation
(\ref{sdedual}) says that, loosely speaking, branching-coalescing
particle systems have a random generating semigroup. The SDE
(\ref{sde}) will be our main tool for studying branching-coalescing
particle systems.

\subsection{Preliminaries}\label{prelim}

In this section we introduce the notation and definitions that we will
use throughout the chapter.\medskip

\noi
{\bf(Inner product and norm notation)} For
$\phi,\psi\in[-\infty,\infty]^\La$, we write
\be
\li\phi,\psi\re:=\sum_i\phi(i)\psi(i)\qquad\mbox{and}
\qquad|\phi|:=\sum_i|\phi(i)|,
\ee
whenever the infinite sums are defined.\medskip

\noi
{\bf (Poisson measures)} If $\phi$ is a $\half^\La$-valued random
variable, then by definition a Poisson measure with random intensity
$\phi$ is an $\N^\La$-valued random variable $\Pois(\phi)$ whose law
is uniquely determined by
\be\label{Poisdef}
E[(1-\psi)^{\Pois(\phi)}]=E[\ex{-\li\phi,\psi\re}]\qquad(\psi\in[0,1]^\La).
\ee
In particular, when $\phi$ is nonrandom, then the components
$(\Pois(\phi)(i))_{i\in\La}$ are independent Poisson distributed
random variables with intensity $\phi(i)$.\medskip

\noi
{\bf (Thinned point measures)} If $x$ and $\phi$ are
random variables taking values in $\N^\La$ and $[0,1]^\La$,
respectively, then by definition a $\phi$-thinning of $x$ is an
$\N^\La$-valued random variable $\Thin_\phi(x)$ whose law is uniquely
determined by
\be\label{Thindef}
E[(1-\psi)^{\Thin_\phi(x)}]=E[(1-\phi\psi)^x]\qquad(\psi\in[0,1]^\La).
\ee
In particular, when $x$ and $\phi$ are nonrandom, and
$x=\sum_{n=1}^m\de_{i_n}$, then a $\phi$-thinning of $x$ can be
constructed as $\Thin_\phi(x):=\sum_{n=1}^m\chi_n\de_{i_n}$ where the
$\chi_n$ are independent $\{0,1\}$-valued random variables with
$P[\chi_n=1]=\phi(i_n)$.

If $\phi$ and $x$ are both random, then it will always be understood
that they are independent. Thus, $\Li(\Thin_\phi(x))$ depends on the
laws $\Li(\phi)$ and $\Li(x)$ alone, and it is only the map
$(\Li(\phi),\Li(x))\mapsto\Li(\Thin_\phi(x))$ that is of interest to
us. We have chosen the present notation in terms of random variables
instead of their laws to keep things simple if $\phi$ and $x$ are
nonrandom.

We leave it to the reader to check the elementary relations
\be\label{nothard}
\Thin_\psi(\Thin_\phi(x))\isd\Thin_{\psi\phi}(x)\quad\mbox{and}
\quad\Thin_\psi(\Pois(\phi))\isd\Pois(\psi\phi),
\ee
where $\isd$ denote equality in distribution.\medskip

\noi
{\bf (Weak convergence)} We let $\ov\N=\N\cup\{\infty\}$ denote the
one-point compactification of $\N$, and equip $\ov\N^\La$ with the
product topology.  We say that probability measures $\nu_n$ on
$\ov\N^\La$ converge weakly to a limit $\nu$, denoted as
$\nu_n\Rightarrow\nu$, when $\int\nu_n(\di x)f(x)\to\int\nu(\di
x)f(x)$ for every $f\in\Ci(\ov\N^\La)$, the space of continuous real
functions on $\ov\N^\La$.  One has $\nu_n\Rightarrow\nu$ if and only
if $\nu_n(\{x:x(i)=y(i) \ \forall i\in\De\})\to\nu(\{x:x(i)=y(i)\
\forall i\in\De\})$ for all finite $\De\sub\La$ and $y\in\N^\De$.

\detail{Let $\Fi$ be the space of functions $x\mapsto 1_{\{x(i)=y(i)\
\forall i\in\De\}}$ with $\De\sub\La$ finite and $y\in\N^\La$ (the
fact that we don't write $\ov\N^\La$ is not a typoe!). Then
$\Fi\sub\Ci(\ov\N^\La)$, $\Fi$ is closed under multiplication,
contains the identity, and separates points. (Note that the functions
of the form $x\mapsto 1_{\{x(i)=y(i)\}}$ with $i\in\La$ and $y\in\N$
already separate points.) By the Stone-Weierstra\ss\ theorem, $\Fi$ is
dense in $\Ci(\ov\N^\La)$. Since $\ov\N^\La$ is compact, every
sequence $\nu_n$ of probability measures on $\ov\N^\La$ is tight. If
$\int\nu_n(\di x)f(x)\to\int\nu(\di x)f(x)$ for every $f\in\Fi$ then
each weak cluster point $\nu^\ast$ of the $\nu_n$ satisfies
$\int\nu^\ast(\di x)f(x)=\int\nu(\di x)f(x)$ for all $f\in\Fi$ and
therefore $\nu^\ast=\nu$. This shows that $\nu_n\Rightarrow\nu$.}

We equip the space $[0,1]^\La$ with the product topology,
and we say that probability measures $\mu_n$ on $[0,1]^\La$ converge
weakly to a limit $\mu$, denoted as $\mu_n\Rightarrow\mu$, when
$\int\mu_n(\di\phi)f(\phi)\to\int\mu(\di\phi)f(\phi)$ for every
$f\in\Ci([0,1]^\La)$.\medskip

\noi
{\bf (Monotone convergence)} If $\nu_1,\nu_2$ are probability measures
on $\ov\N^\La$, then we say that $\nu_1$ and $\nu_2$ are
stochastically ordered, denoted as $\nu_1\leq\nu_2$, if
$\ov\N^\La$-valued random variables $Y_1,Y_2$ with laws
$\Li(Y_i)=\nu_i$ ($i=1,2$) can be coupled such that $Y_1\leq Y_2$. We
say that a sequence of probability measures $\nu_n$ on $\N^\La$
decreases (increases) stochastically to a limit $\nu$, denoted as
$\nu_n\down\nu$ ($\nu_n\up\nu$), if random variables $Y_n,Y$ with laws
$\Li(Y_n)=\nu_n$ and $\Li(Y)=\nu$ can be coupled such that $Y_n\down
Y$ ($Y_n\up Y$). It is not hard to see that $\nu_n\down\nu$
($\nu_n\up\nu$) implies $\nu_n\Rightarrow\nu$.
Stochastic ordering and monotone convergence of
probability measures on $[0,1]^{\La}$ are defined in the same
way.\medskip

\noi
{\bf (Finite systems)} We denote the set of finite particle
configurations by $\Ni(\La):=\{x\in\N^\La:|x|<\infty\}$ and let
\be\label{Sidef}
\Si(\Ni(\La))
:=\{f:\Ni(\La)\to\R:|f(x)|\leq K|x|^k+M\mbox{ for some }K,M,k\geq 0\}
\ee
denote the space of real functions on $\Ni(\La)$ satisfying a
polynomial growth condition. For finite initial conditions, the
$(a,b,c,d)$-braco-process $X$ is well-defined as a Markov process in
$\Ni(\La)$ (in particular, $X$ does not explode), $f(X_t)$ is
absolutely integrable for each $f\in\Si(\Ni(\La))$ and $t\geq 0$, and
the semigroup
\be\label{Stdef}
S_tf(x):=E^x[f(X_t)]\qquad(t\geq 0,\ x\in\Ni(\La),\ f\in\Si(\Ni(\La)))
\ee
maps $\Si(\Ni(\La))$ into itself (see Proposition~\ref{finmart}
below).\medskip

\noi
{\bf (Liggett-Spitzer space)} Set $a_{\rm s}(i,j):=a(i,j)+a^\dgg(i,j)$.
It follows from our assumptions on $a$ that there exist (strictly)
positive constants $(\ga_i)_{i\in\La}$ such that
\detail{To see this, put $|a_{\rm s}|:=\sup_i\sum_ja_{\rm s}(i,j)$
which is finite by our assumptions (iii) and (iv). Define
$P(i,j):=\frac{1}{|a|}a_{\rm s}(i,j)$ and for any bounded function
$\phi:\La\to\R$ put $P\phi(i):=\sum_jP(i,j)\phi(j)$. Fix $\phi\geq 0$,
$\phi\neq 0$. It is clear that $\|P\phi\|_\infty\leq\|\phi\|_\infty$
and therefore, for any $\eps>0$, the infinite sum
\be
\ga:=\sum_{k=0}^\infty e^{-\eps k}P^k\phi
\ee
converges. One has
\be
P\ga=e^\eps\sum_{k=0}^\infty e^{-\eps(k+1)}P^{k+1}\phi=e^\eps(\ga-\phi),
\ee
and therefore $\sum_ja_{\rm s}(i,j)\ga_j\leq|a_{\rm
s}|e^\eps\ga_i$. The fact that the $\ga_i$ are all positive follows
from the irreducibility of $a$ (assumption (ii)).}
\be\label{gacon}
\sum_i\ga_i<\infty\quad\mbox{and}\quad\sum_ja_{\rm s}(i,j)\ga_j\leq K\ga_i
\quad(i\in\La)
\ee
for some $K<\infty$.
We fix such $(\ga_i)_{i\in\La}$ throughout the chapter
and define the Liggett-Spitzer space (after \cite{LS81})
\be
\Ei_\ga(\La):=\{x\in\N^\La:\|x\|_\ga<\infty\},
\ee
where for $x\in\Z^\La$ we put
\be\label{ganorm}
\|x\|_\ga:=\sum_i\ga_i|x(i)|.
\ee
We let $\Ci_{\rm Lip}(\Ei_\ga(\La))$ denote the class of Lipschitz
functions on $\Ei_\ga(\La)$, i.e., $f:\Ei_\ga(\La)\to\R$ such that
$|f(x)-f(y)|\leq L\|x-y\|_\ga$ for some $L<\infty$.\medskip

\noi
{\bf (Infinite systems)} It is known (\cite{Che87}, see also
Proposition~\ref{Xcon} below) that for each $f\in\Ci_{\rm
Lip}(\Ei_\ga(\La))$ and $t\geq 0$, the function $S_tf$ defined in
(\ref{Stdef}) can be extended to a unique Lipschitz function on
$\Ei_\ga(\La)$, also denoted by $S_tf$. Moreover, there exists a
time-homogeneous Markov process $X$ in $\Ei_\ga(\La)$ (also called
$(a,b,c,d)$-braco-process) with transition laws given by
\be
E^x[f(X_t)]=S_tf(x)\qquad(f\in\Ci_{\rm Lip}(\Ei_\ga(\La)),
\ x\in\Ei_\ga(\La),\ t\geq 0).
\ee
\detail{To see that the restriction of a function $f\in\Ci_{\rm
Lip}(\Ei_\ga(\La))$ to $\Ni(\La)$ yields a function in
$\Si(\Ni(\La))$, so that $S_tf$ is well-defined, note that for such a
function $|f(x)|\leq|f(0)|+L\|x\|_\ga\leq|f(0)|+L\|\ga\|_\infty|x|$.}
We will show (in Proposition~\ref{Xcon} below) that $X$ has a
modification with cadlag sample paths, a fact that may seem
obvious but to our knowledge has not been proved before.\medskip

\noi
{\bf (Survival and extinction)} We say that the
$(a,b,c,d)$-braco-process survives if
\be
P^x[X_t\neq 0\ \forall t\geq 0]>0\quad\mbox{for some}\quad x\in\Ni(\La).
\ee
If $X$ does not survive we say that $X$ dies out.  Note that the
process with death rate $d=0$ survives, since the number of particles
can no longer decrease once only one particle is left. If $\La$ is
finite then the $(a,b,c,d)$-braco-process survives if and only if
$d=0$, but for infinite $\La$ survival often holds also for some
$d>0$. For $\La=\Z^d$ and $b$ sufficiently large survival has been
proved in \cite[Theorem~3.1]{SU86}. We plan to study sufficient
conditions for survival in more detail in a forthcoming paper.\medskip

\noi
{\bf (Nontrivial measures)} We say that a probability measure $\nu$ on
$\ov\N^\La$ is nontrivial if $\nu(\{0\})=0$, where $0\in\ov\N^\La$
denotes the zero configuration. Likewise, we say that a probability
measure $\mu$ on $[0,1]^\La$ is nontrivial if $\mu(\{0\})=0$.\medskip

\noi
{\bf (Homogeneous lattices)} By definition, an {\em automorphism} of
$(\La,a)$ is a bijection $g:\La\to\La$ such that $a(gi,gj)=a(i,j)$ for
all $i,j\in\La$. We denote the group of all automorphisms of $(\La,a)$
by ${\rm Aut}(\La,a)$. We say that a subgroup $G\sub{\rm Aut}(\La,a)$
is {\em transitive} if for each $i,j\in\La$ there exists a $g\in G$
such that $gi=j$. We say that $(\La,a)$ is {\em homogeneous} if ${\rm
Aut}(\La,a)$ is transitive. We define shift operators
$T_g:\N^\La\to\N^\La$ by
\be\label{shift}
T_gx(j):=x(g^{-1}j)\qquad(i\in\La,\ x\in\N^\La,\ g\in{\rm Aut}(\La,a)).
\ee
If $G$ is a subgroup of ${\rm Aut}(\La,a)$, then we say that a
probability measure $\nu$ on $\N^\La$ is {\em $G$-homogeneous} if
$\nu\circ T_g^{-1}=\nu$ for all $g\in G$. For example, if $\La=\Z^d$
and $a(i,j)=1_{\{|i-j|=1\}}$ (nearest-neighbor random walk), then the
group $G$ of translations $i\mapsto i+j$ ($j\in\La$) form a transitive
subgroup of ${\rm Aut}(\La,a)$ and the $G$-homogeneous probability
measures are the translation invariant probability measures. Shift
operators and $G$-homogeneous measures on $[0,1]^\La$ are defined
analogously.

\subsection{Main results}

Our first result is a tool that we exploit substantially towards the
main result. Part~(a) is known \cite[Lemma~2.1]{SU86}, but we
are not aware of parts~(b) and (c) occuring anywhere in the literature.
\bt{\bf (Dualities and Poissonization)}\label{duth}
Let $X$ and $\Xc$ be the $(a,b,c,d)$-braco-process and the
$(a,b,c,d)$-resem-process, respectively, and let $\Xc^\dgg$ denote the
$(a^\dgg,b,c,d)$-resem-process. Then the following holds:\newline
{\bf(a) (Duality)} 
\be\label{du}
P^x[\Thin_\phi(X_t)=0]=P^\phi[\Thin_{\Xc^\dgg_t}(x)=0]
\qquad(t\geq 0,\ \phi\in[0,1]^\La,\ x\in\Ei_\ga(\La)).
\ee
{\bf(b) (Self-duality)} Assume $c>0$, then
\be\label{sdu}
P^\phi[\Pois(\ffrac{b}{c}\Xc_t\psi)=0]
=P^\psi[\Pois(\ffrac{b}{c}\phi\Xc^\dgg_t)=0]
\qquad(t\geq 0,\ \phi,\psi\in[0,1]^\La).
\ee
{\bf(c) (Poissonization)} Assume $c>0$, then
\be\label{Poisfo}
P^{\Li(\Pois(\frac{b}{c}\phi))}[X_t\in\cdot\,]
=P^\phi[\Pois(\ffrac{b}{c}\Xc_t)\in\cdot\,]
\qquad(t\geq 0,\ \phi\in[0,1]^\La),
\ee
i.e., if $X$ is started in the initial law
$\Li(\Pois(\frac{b}{c}\phi))$ and $\Xc$ is started in $\phi$, then
$X_t$ and $\Pois(\frac{b}{c}\Xc_t)$ are equal in law.
\et
Note that $P[\Thin_\phi(x)=0]=(1-\phi)^x$. Therefore,
Theorem~\ref{duth}~(a) is just a reformulation of the duality relation
(\ref{sdedual}).  Theorem~\ref{duth}~(b) says that
resampling-selection processes are in addition dual with respect to
each other. In particular, if the underlying motion is symmetric,
i.e., $a=a^\dgg$, then this is a self-duality. Since
$P[\Pois(\phi)=0]=e^{-|\phi|}$, formula (\ref{sdu}) can be rewritten
as
\be\label{selfo}
E^\phi\big[\ex{-\frac{b}{c}\li\Xc_t,\psi\re}\big]
=E^\psi\big[\ex{-\frac{b}{c}\li\phi,\Xc^\dgg_t\re}\big]
\qquad(t\geq 0,\ \phi,\psi\in[0,1]^\La).
\ee
We note that by \cite[Lemma~15.5.1]{Kal83}, for $b>0$, the
distribution of $\Xc_t$ is determined uniquely by all
$E[e^{-\frac{b}{c}\li\Xc_t,\psi\re}]$ with $\psi\in[0,1]^\La$. To
convince the reader that the notation in (\ref{du}) and (\ref{sdu}),
which may feel a little uneasy in the beginning, is convenient, we
give here the proof of the Poissonization formula
(\ref{Poisfo}).\medskip

\noi
{\bf Proof of Theorem~\ref{duth}~(c)} By (\ref{nothard}) and the
duality relations (\ref{du}) and (\ref{sdu}),
\be\ba{l}
P^{\Li(\Pois(\frac{b}{c}\phi))}[\Thin_\psi(X_t)=0]
=P^\psi[\Thin_{\Xc^\dgg_t}(\Pois(\ffrac{b}{c}\phi))=0]\\[5pt]
\dis\qquad=P^\psi[\Pois(\ffrac{b}{c}\Xc^\dgg_t\phi)=0]
=P^\phi[\Pois(\ffrac{b}{c}\psi\Xc_t)=0]
=P^\phi[\Thin_\psi(\Pois(\ffrac{b}{c}\Xc_t))=0].
\ec
Since this is true for all $\psi\in[0,1]^\La$, the random variables
$X_t$ and $\Pois(\ffrac{b}{c}\Xc_t)$ are equal in distribution.\qed

\noi
Our next result shows that it is possible to start the
$(a,b,c,d)$-braco-process with infinitely many particles at each
site. This result (except for parts~(b) and (f))
has been proved for branching-coalescing particle systems with more
general branching and coalescing mechanisms on $\Z^d$ in
\cite{DDL90}. Their methods are not restricted to the case $\La=\Z^d$,
but we give an independent proof using duality, which has the
additional appeal of yielding the explicit bound in part~(b).

\bt{\bf(The maximal branching-coalescing process)}\label{maxbraco}
Assume that $c>0$. Then there exists an $\Ei_\ga(\La)$-valued process
$X^{(\infty)}=(X^{(\infty)}_t)_{t>0}$ with the following properties:
\medskip

\noi
{\bf (a)} For each $\eps>0$, $(X^{(\infty)}_t)_{t\geq\eps}$ is the
$(a,b,c,d)$-braco-process starting in $X^{(\infty)}_\eps$.\medskip

\noi
{\bf (b)} Set $r:=b-d+c$. Then
\be\label{explicit}
E[X^{(\infty)}_t(i)]\leq\left\{\ba{cl}\frac{r}{c(1-e^{-rt})}\quad&
\mbox{if }r\neq 0,\\[5pt] \frac{1}{ct}\quad&\mbox{if }r=0\ea\right.
\qquad(i\in\La,\ t>0).
\ee

\noi
{\bf (c)} If $X^{(n)}$ are $(a,b,c,d)$-braco-processes starting in
initial states $x^{(n)}\in\Ei_\ga(\La)$ such that
\be
x^{(n)}(i)\up\infty\quad\mbox{as }n\up\infty\qquad(i\in\La),
\ee
then
\be
\Li(X^{(n)}_t)\up\Li(X^{(\infty)}_t)\quad\mbox{as }n\up\infty\qquad(t>0).
\ee
{\bf (d)} There exists an invariant measure $\ov\nu$ of the
$(a,b,c,d)$-braco-process such that
\be
\Li(X^{(\infty)}_t)\down\ov\nu\quad\mbox{as }t\up\infty.
\ee
{\bf (e)} If $\nu$ is another invariant measure for the
$(a,b,c,d)$-braco-process, then $\nu\leq\ov\nu$.\vspace{5pt}

\noi
{\bf (f)} The measure $\ov\nu$ is uniquely characterised by
\be\label{extinct}
\int\ov\nu(\di x)(1-\phi)^x=P^\phi[\exists t\geq 0\mbox{ such that }
\Xc^\dgg_t=0]\qquad(\phi\in[0,1]^\La),
\ee
where $\Xc^\dgg$ denotes the $(a^\dgg,b,c,d)$-resem-process.
\et
We call $X^{(\infty)}$ the maximal $(a,b,c,d)$-braco process and we
call $\ov\nu$ the upper invariant measure. To see why
Theorem~\ref{maxbraco}~(f) holds, note that by Theorem~\ref{duth}~(a)
and Theorem~\ref{maxbraco}~(c),
\be\label{fikst}
P[\Thin_\phi(X^{(\infty)}_t)=0]
=\lim_{n\up\infty}P^\phi[\Thin_{\Xc^\dgg}(x^{(n)})=0]
=P^\phi[\Xc^\dgg_t=0]\qquad(\phi\in[0,1]^\La,\ t>0).
\ee
Now $0$ is an absorbing state for the $(a,b,c,d)$-resem-process, and
therefore $P^\phi[\Xc^\dgg_t=0] =P^\phi[\exists s\leq t\mbox{ such
that }\Xc^\dgg_s=0]$.  Therefore, taking the limit $t\up\infty$ in
(\ref{fikst}) we arrive at (\ref{extinct}).

The $(a,b,c,d)$-resem process has an upper invariant measure too. Of
our next theorem, parts~(a)--(c) are simple, but part~(d) lies
somewhat deeper.
\bt{\bf\hspace{-2.65pt}(The maximal resampling-selection process)\hspace{-1.5pt}}\label{maxresem}
Let $\Xc^1$ denote the $(a,b,c,d)$-resem-process started in
$\Xc^1_0(i)=1$ $(i\in\La)$. Then the following holds.\newline
{\bf (a)} There exists an invariant measure $\ov\mu$ of the
$(a,b,c,d)$-resem process such that
\be
\Li(\Xc^1_t)\down\ov\mu\quad\mbox{as }t\up\infty.
\ee
{\bf (b)} If $\mu$ is another invariant measure, then $\mu\leq\ov\mu$.\medskip

\noi
{\bf (c)} Let $X^\dgg$ denote the $(a^\dgg,b,c,d)$-braco-process. Then
\be\label{sdemom}
\int\ov\mu(\di\phi)(1-\phi)^x=P^x[\exists t\geq 0\mbox{ such that }X^\dgg_t=0]
\qquad(x\in\Ni(\La)),
\ee
and the measure $\ov\mu$ is nontrivial if and only if the
$(a^\dgg,b,c,d)$-braco-process survives.\medskip

\noi
{\bf (d)} Assume that $c>0$ and that $\La$ is infinite. If $\Yi$ is a
random variable such that $\ov\mu=\Li(\Yi)$, then the upper
invariant measure of the $(a,b,c,d)$-braco-process is given by
$\ov\nu=\Li(\Pois(\frac{b}{c}\Yi))$. If $\ov\mu$ is nontrivial then
so is $\ov\nu$.
\et
Note that $\int\ov\mu(\di\phi)(1-\phi)^x$ is the probability that $x$
individuals, sampled from a population with resampling and selection
in the equilibrium measure $\ov\mu$, all have defective genes.

The following is our main result.
\bt{\bf (Convergence to the upper invariant measure)}\label{homconv}
Assume that $(\La,a)$ is infinite and homogeneous, $G$ is a transitive
subgroup of ${\rm Aut}(\La,a)$, and $c>0$.\medskip

\noi
{\bf (a)} Let $X$ be the $(a,b,c,d)$-braco process started in a
$G$-homogeneous nontrivial initial law $\Li(X_0)$. Then
$\Li(X_t)\Rightarrow\ov\nu$ as $t\to\infty$, where $\ov\nu$ is the
upper invariant measure.\medskip

\noi
{\bf (b)} Let $\Xc$ be the $(a,b,c,d)$-resem process started in a
$G$-homogeneous nontrivial initial law $\Li(\Xc_0)$. Then
$\Li(\Xc_t)\Rightarrow\ov\mu$ as $t\to\infty$, where $\ov\mu$ is the
upper invariant measure.
\et
Shiga and Uchiyama \cite[Theorems~1.3 and 1.4]{SU86} proved
Theorem~\ref{homconv}~(b) under the additional assumptions that
$\La=\Z^d$ and that $a$ satisfies a first moment condition in case the
death rate $d$ is zero. As we will show below
Theorem~\ref{homconv}~(b) can be derived from
Theorem~\ref{homconv}~(a) by Poissonization, but not vice versa.

\subsection{Methods}\label{methods}

A key ingredient in the proofs of Theorem~\ref{maxresem}~(d) and
Theorem~\ref{homconv} is the following property of
resampling-selection processes, which is of some interest on its own.
\bl{\bf (Extinction versus unbounded growth)}\label{exgroC3}
Assume that $c>0$. Let $\Xc$ be the $(a,b,c,d)$-resem-process
starting in an initial state $\phi\in[0,1]^\La$ with
$|\phi|<\infty$. Then $e^{-\frac{b}{c}|\Xc_t|}$ is a submartingale,
and a martingale if $d=0$. If moreover $\La$ is infinite, then
\be
\Xc_t=0\quad\mbox{for some }t\geq 0\quad\mbox{or}
\quad\lim_{t\to\infty}|\Xc_t|=\infty\quad{\rm a.s.}
\ee 
\el
Note that by Theorem~\ref{duth}~(b),
\be
E^\phi\big[\ex{-\frac{b}{c}\li\Xc_t,1\re}\big]
=E^1\big[\ex{-\frac{b}{c}\li\phi,\Xc^\dgg_t\re}\big]
\geq \ex{-\frac{b}{c}\li\phi,1\re}\qquad(\phi\in[0,1]^\La),
\ee
with equality if $d=0$, since $1$ is a stationary state for the
$(a^\dgg,b,c,0)$-resem-process. This shows that $e^{-\frac{b}{c}|\Xc_t|}$
is a submartingale, and a martingale if $d=0$. By submartingale
convergence, $|\Xc_t|$ converges a.s.\ to a limit in $[0,\infty]$. All
the hard work of Lemma~\ref{exgroC3} consists of proving that this limit
is a.s.\ either $0$ or $\infty$, and that $\Xc$ gets extinct in finite
time if the limit is zero.

Once Lemma~\ref{exgroC3} is established the proof of
Theorem~\ref{maxresem}~(d) is simple.\medskip

\noi
{\bf Proof of Theorem~\ref{maxresem}~(d)} Let $\Yi$ be a random
variable such that $\ov\mu=\Li(\Yi)$ and let $Y$ be a random variable
such that $\ov\nu=\Li(Y)$. By (\ref{nothard}), Theorem~\ref{duth}~(b),
and Theorem~\ref{maxbraco}~(f)
\be\ba{l}\label{uit}
\dis P[\Thin_\phi(\Pois(\ffrac{b}{c}\Yi))=0]
=\lim_{t\to\infty}P^1[\Pois(\ffrac{b}{c}\phi\Xc_t)=0]
=\lim_{t\to\infty}P^\phi[\Pois(\ffrac{b}{c}\Xc^\dgg_t)=0]\\[5pt]
\dis\qquad\stackrel{!}{=}P^\phi[\exists t\geq 0\mbox{ such that }\Xc^\dgg_t=0]
=P[\Thin_\phi(Y)=0],
\ec
where we have used Lemma~\ref{exgroC3} in the equality marked with `!'.
Since (\ref{uit}) holds for all $\phi\in[0,1]^\La$, the random
variables $\Pois(\frac{b}{c}\Yi)$ and $Y$ are equal in
distribution. By Lemma~\ref{exgroC3}, $|\Yi|\in\{0,\infty\}$ a.s. and
therefore if $\ov\mu$ is nontrivial then $\Li(\Pois(\frac{b}{c}\Yi))$
is nontrivial.\qed

\noi
In view of Theorem~\ref{maxresem}~(d), it is natural to ask if for
infinite lattices, every invariant law of the
$(a,b,c,d)$-braco-process is the Poissonization of an invariant law of
the $(a,b,c,d)$-resem-process. We do not know the answer to this
question.

In order to give a very short proof of Theorem~\ref{homconv}, we need
one more lemma.
\bl{\bf(Systems with particles everywhere)}\label{omnipres}
Assume that $(\La,a)$ is infinite and homogeneous and that $G$ is a
transitive subgroup of ${\rm Aut}(\La,a)$. Let $X$ be the
$(a,b,c,d)$-braco process started in a $G$-homogeneous nontrivial
initial law $\Li(X_0)$. Then, for any $t>0$
\be\label{posthin}
\lim_{n\to\infty}P[\Thin_{\phi_n}(X_t)=0]=0,
\ee
for all $\phi_n\in[0,1]^\La$ satisfying $|\phi_n|\to\infty$.
\el
{\bf Proof of Theorem~\ref{homconv} (a)} Let $\Xc^\dgg$ denote the
$(a^\dgg,b,c,d)$-resem-process started in $\phi$. By Theorem~\ref{duth}~(a),
Lemmas~\ref{exgroC3} and \ref{omnipres}, and Theorem~\ref{maxbraco}~(f),
\be\ba{l}
\dis\lim_{t\to\infty}P[\Thin_\phi(X_t)=0]
=\lim_{t\to\infty}P[\Thin_{\Xc^\dgg_{t-1}}(X_1)=0]\\[5pt]
\dis\qquad=P[\exists t\geq 0\mbox{ such that }\Xc^\dgg_t=0]
=\int\!\!\ov\nu(\di x)\,(1-\phi)^x.
\ec
Since this holds for all $\phi\in[0,1]^\La$, it follows that
$\Li(X_t)\Rightarrow\ov\nu$.\qed

\noi
{\bf\boldmath Proof of Theorem~\ref{homconv}~(b)} Let $X_\infty$ and
$\Xc_\infty$ be random variables with laws $\ov\nu$ and $\ov\mu$,
respectively. Let $\Xc$ be the $(a,b,c,d)$-resem-process started in a
$G$-homogeneous nontrivial initial law $\Li(\Xc_0)$. Let $X$ be the
$(a,b,c,d)$-braco-process started in
$\Li(X_0):=\Li(\Pois(\frac{b}{c}\Xc_0))$. Then by
Theorem~\ref{homconv}~(a), $\Li(X_t)\Rightarrow\Li(X_\infty)$ as
$t\to\infty$. Therefore, by Poissonization (Theorem~\ref{duth} (c))
and by Theorem~\ref{maxresem} (d),
$\Li(\Pois(\frac{b}{c}\Xc_t))\Rightarrow\Li(X_\infty)
=\Li(\Pois(\frac{b}{c}\Xc_\infty))$. It
follows that
\be\ba{l}
\dis P\big[\ex{-\frac{b}{c}\li\Xc_t,\phi\re}\big]
=P[\Thin_\phi(\Pois(\ffrac{b}{c}\Xc_t))=0]\\
\dis\hspace{1.7cm}\Longrightarrow P[\Thin_\phi(\Pois(\ffrac{b}{c}\Xc_\infty))
=0]=P\big[\ex{-\frac{b}{c}\li\Xc_\infty,\phi\re}\big]
\qquad\mbox{as }t\to\infty.
\ec
Since this holds for all $\phi\in[0,1]^\La$, we conclude that
$\Li(\Xc_t)\Rightarrow\Li(\Xc_\infty)$.\qed

\noi
Note that there is no easy way to convert the last argument: if
$\Li(X_0)$ is homogeneous and nontrivial then we cannot in general
find a random variable $\Xc_0$ such that
$\Li(X_0)=\Li(\Pois(\frac{b}{c}\Xc_0))$. For example, this is the case
if $X_0(i)\leq 1$ for each $i\in\La$ a.s. Therefore,
Theorem~\ref{homconv}~(a) is stronger than
Theorem~\ref{homconv}~(b).

Summarizing, all the hard work for getting Theorem~\ref{homconv} is in
proving Lemmas~\ref{exgroC3} and \ref{omnipres}, as well as the more
basic Theorems~\ref{duth} and \ref{maxbraco}. The heart of the proof
of Theorem~\ref{maxbraco} is the bound in part~(b). We derive this
bound using a `duality' relation with a nonnegative error term,
between the $(a,b,c,d)$-braco-process and a super random walk
(Proposition~\ref{subdup}). We call this relation a
subduality. Theorem~\ref{maxbraco}~(b) yields a lower bound on the
finite time extinction probabilities of the $(a,b,c,d)$-resem-process
started with small initial mass (Lemma~\ref{finex}, in particular
formula (\ref{linbo})), which plays a key role in the proof of
Lemma~\ref{exgroC3}.

Our methods are similar to those of Shiga and Uchiyama
\cite{SU86}. Since they prove a version of our
Theorem~\ref{homconv}~(b), while our main focus is on proving the
stronger Theorem~\ref{homconv}~(a), the roles of $X$ and $\Xc$ are
interchanged in their work. Their Lemma~3.2 and Theorem~4.2 are
analogues for the $(a,b,c,d)$-braco-process $X$ of our
Lemma~\ref{exgroC3}. The proof of the latter is considerably more
involved, however. This is because of the fact that we do not want to
use spatial homogeneity and we have to prove that $|\Xc_t|\to 0$
implies $\Xc_t=0$ for some $t\geq 0$, which is obvious for the
$(a,b,c,d)$-braco-process $X$. On the other hand, we can use the
submartingale property of $e^{-\frac{b}{c}|\Xc_t|}$, a very useful
fact that has no analogue for the particle system. Lemma~2.5 in
\cite{SU86} is the analogue for the $(a,b,c,d)$-resem-process $\Xc$ of
our Lemma~\ref{omnipres}. By adapting elements of their proof to our
situation, we were able to simplify and considerably shorten our
original proof of Lemma~\ref{omnipres}.

Our original proof of Lemma~\ref{omnipres} assumed that $\La$ has a
group structure, and used an $L^2$ spatial ergodic theorem for general
countable groups that need not be amenable.

\subsection{Discussion}

Generalizing our model, let $X$ be a process in a Liggett-Spitzer
subspace of $\N^\La$, with local jump rates
\be\ba{ll}\label{RD}
x\mapsto x+\de_j-\de_i\quad&\mbox{with rate }a(i,j)\\[5pt]
x\mapsto x+\de_i\quad&\mbox{with rate }\sum_{n=0}^kb_nx^{(n)},\\[5pt]
x\mapsto x-\de_i\quad&\mbox{with rate }\sum_{n=1}^{k+1}c_nx^{(n)},
\ec
where $x^{(0)}:=1$ and $x^{(n)}:=x(x-1)\cdots(x-n+1)$ ($n\geq 1$). In
particular, the $(a,b,c,d)$-braco-process corresponds to the case
$k=1$, $b_0=0$, $b_1=b$, $c_1=d$, and $c_2=c$. Processes with jump
rates as in (\ref{RD}) are known as reaction-diffusion systems. It has
been known for a long time that if the coefficients satisfy
\be\label{revcond}
a=a^\dgg\quad\mbox{and}\quad b_n=\la c_n\quad\mbox{for some}\quad\la\geq 0,
\ee
then $\Li(\Pois(\la))$ is a reversible equilibrium for the
corresponding reaction-diffusion system. Note that the
$(a,b,c,d)$-braco-process satisfies (\ref{revcond}) if and only if
$a=a^\dgg$ and $d=0$.

The ergodic behavior of reaction-diffusion systems on $\La=\Z^d$
satisfying the reversibility condition (\ref{revcond}) was studied by
Ding, Durrett and Liggett in \cite{DDL90}. For our model with
$a=a^\dgg$ and $d=0$ on $\Z^d$, they show that all homogeneous
invariant measures are convex combinations of $\de_0$ and
$\Li(\Pois(\frac{b}{c}))$. Their proof uses the fact that for a large
block in $\Z^d$, surface terms are small compared to volume terms,
i.e., $\Z^d$ is amenable. Such arguments typically fail on nonamenable
lattices such as trees, and therefore it is not immediately obvious if
their methods can be generalized to such lattices. Our
Theorem~\ref{homconv}~(a) shows that all homogeneous invariant measures of
the $(a,b,c,d)$-braco-process are convex combinations of $\de_0$ and
$\ov\nu$, also in the non-reversible case $d>0$ and for nonamenable
lattices. Thus, neither reversibility nor amenability are essential
here.

On the other hand, we believe that amenability is essential for more
subtle ergodic properties of reaction-diffusion processes. In analogy
with the contact process, let us say that a reaction-diffusion process
with $b_0=0$ exhibits complete convergence, if
\be
P^x[X_t\in\cdot\,]\Rightarrow\rho(x)\ov\nu+(1-\rho(x))\de_0
\quad\mbox{as}\quad t\to\infty\qquad(x\in\Ni(\La)),
\ee
where $\rho(x):=P^x[X_t\neq 0\ \forall t\geq 0]$ denotes the survival
probability. It has been shown by Mountford \cite{Mou92} that complete
convergence holds for reaction-diffusion systems on $\La=\Z^d$
satisfying the reversibility condition (\ref{revcond}), $b_0=0$, and a
first moment condition on $a$. We conjecture that complete convergence
holds more generally if $a=a^\dgg$ and $\La$ is amenable, but not in
general on nonamenable lattices. As a motivation for this conjecture,
we note that complete convergence holds for the contact process on
$\Z^d$ but not in general on $\T^d$; see Liggett \cite{Lig99}.

The self-duality of resampling-selection processes
(Theorem~\ref{duth}~(b)) is reminiscent of the self-duality of the
contact process. It is an interesting question whether our methods can
be adapted to the contact process, to show that the upper invariant
measure of the contact process on a countable group is the limit
started from any homogeneous nontrivial initial law.

Other interesting processes that some of our techniques might be applied to
are multitype bran\-ching-coalescing particle systems. For example, it
seems natural to color the particles in a branching-coalescing
particle system in two (or more) colors, with the rule that in
coalescence of differently colored particles, the newly created
particle chooses the color of one of its parents with equal
probabilities (neutral selection) or with a prejudice towards one
color (positive selection). More difficult questions refer to what
happens when the two colors have different parameters $b,c,d$ or even
different underlying motions $a$.

One also wonders whether the techniques in this chapter can be
generalized to reaction-diffusion processes with higher-order
branching and coalescence as in (\ref{RD}). It seems that at least
some of these systems have some sort of a resampling-selection dual too,
now with `resampling' and `selection' events involving three and more
particles.

We conclude with an intriguing question. Does survival of the
$(a,b,c,d)$-braco-process $X$ imply survival of the
$(a^\dgg,b,c,d)$-braco-process $X^\dgg$? If $X$ survives, then
Theorem~\ref{maxresem}~(c) and (d) and Theorem~\ref{homconv}~(a) show
that the upper invariant measure of $X^\dgg$ is nontrivial, which
suggests that $X^\dgg$ should survive. Survival of $X^\dgg$ is obvious
if $(\La,a)$ and $(\La,a^\dgg)$ are isomorfic, as is the case if
$a=a^\dgg$, or if $\La$ is an Abelian group, with group action denoted
by $+$, and $a(i,j)$ depends only on $j-i$. However, even when
$(\La,a)$ is homogeneous, $(\La,a)$ and $(\La,a^\dgg)$ need in general
not be isomorphic, and in this case we don't know the answer to our
question.

\subsection{Outline}

We start in Section~\ref{mapo} with a few generalities about
martingale problems that will be needed in our proofs. In
Section~\ref{cosec} we construct $(a,b,c,d)$-braco-processes and
$(a,b,c,d)$-resem-processes and prove some of their elementary
properties, such as comparison, approximation with finite systems,
moment estimates and martingale problems. Section~\ref{dusec} contains
the proof of Theorem~\ref{duth} and of the subduality between
branching-coalescing particle systems and super random walks. In
Section~\ref{maxsec} we prove Theorems~\ref{maxbraco} and
\ref{maxresem}. In Section~\ref{consec}, finally, we prove
Lemma~\ref{exgroC3} and Lemma~\ref{omnipres}, thereby completing the
proof of Theorem~\ref{homconv}.\vc

\noi 
{\bf\large Acknowledgements} We thank Klaus Fleischmann who played a
stimulating role during the early stages of this project and answered
a question about Laplace functionals, Claudia Neuhauser for answering
questions about branching-coalescing processes, Olle H\"aggstr\"om for
answering questions on nonamenable groups, and Tokuzo Shiga for
answering our questions about his work. We thank the referee for
drawing our attention to the reference \cite{SU86}. Part of this work
was carried out during the visits of Siva Athreya to the Weierstrass
Institute for Applied Analysis and Stochastics, Berlin and to the
Friedrich-Alexander University Erlangen-Nuremberg, and of Jan Swart to
the Indian Statistical Institute, Delhi. We thank all these places for
their kind hospitality.

\section{Martingale problems}\label{mapo}

\subsection{Definitions}

If $E$ be a metrizable space, we denote by $M(E),B(E)$ the spaces of
real Borel measurable and bounded real Borel measurable functions on
$E$, respectively. If $A$ is a linear operator from a domain
$\Di(A)\sub M(E)$ into $M(E)$ and $X$ is an $E$-valued process, then
we say that $X$ solves the martingale problem for $A$ if $X$ has
cadlag sample paths and for each $f\in\Di(A)$,
\be
E\big[|f(X_t)|\big]<\infty\quad\mbox{and}\quad
\int_0^tE\big[|Af(X_s)|\big]\di s<\infty\qquad(t\geq 0),
\ee
and the process $(M_t)_{t\geq 0}$ defined by
\be\label{duM}
M_t:=f(X_t)-\int_0^t\!Af(X_s)\di s\qquad(t\geq 0)
\ee
is a martingale with respect to the filtration generated by $X$.

\subsection{Duality with error term}

For later use in Section~\ref{dusec}, we formulate a theorem giving
sufficient conditions for two martingale problems to be dual to each
other up to a possible error term. Although the techniques for proving
Theorem~\ref{errorth} below are well-known (see, for example,
\cite[Section~4.4]{EK}), we don't know a good reference for the
theorem as is formulated here.
\bt{\bf(Duality with error term)}\label{errorth}
Assume that $E_1,E_2$ are metrizable spaces and that for $i=1,2$,
$A_i$ is a linear operator from a domain $\Di(A_i)\sub B(E_i)$ into
$M(E_i)$. Assume that $\Psi\in B(E_1\times E_2)$ satisfies
$\Psi(\cdot,x_2)\in\Di(A_1)$ and $\Psi(x_1,\cdot)\in\Di(A_2)$ for each
$x_1\in E_1$ and $x_2\in E_2$, and that
\be
\Phi_1(x_1,x_2):=A_1\Psi(\cdot,x_2)(x_1)\quad\mbox{and}\quad\Phi_2(x_1,x_2)
:=A_2\Psi(x_1,\cdot)(x_2)\qquad(x_1\in E_1,\ x_2\in E_2)
\ee
are jointly measurable in $x_1$ and $x_2$. Assume that $X^1$ and $X^2$
are independent solutions to the martingale problems for $A_1$ and
$A_2$, respectively, and that
\be\label{Phint2}
\int_0^T\!\!\!\di s\int_0^T\!\!\!\di t\;E\big[|\Phi_i(X^1_s,X^2_t)|\big]<\infty
\qquad(T\geq 0,\ i=1,2).
\ee
Then
\be\label{apdual}
E[\Psi(X^1_T,X^2_0)]-E[\Psi(X^1_0,X^2_T)]
=\int_0^T\!\!\!\di t\;E[R(X^1_t,X^2_{T-t})]\qquad(T\geq 0),
\ee
where $R(x_1,x_2):=\Phi_1(x_1,x_2)-\Phi_2(x_1,x_2)
\quad(x_1\in E_1,\ x_2\in E_2)$.
\et
{\bf Proof} Put
\be
F(s,t):=E[\Psi(X^1_s,X^2_t)]\qquad(s,t\geq 0).
\ee
Then, for each $T>0$,
\bc\label{splitC3}
\dis\int_0^T\!\!\di t\,\big\{F(t,0)-F(0,t)\big\}
&=&\dis\int_0^T\!\!\di t\,\big\{F(T-t,t)-F(0,t)-F(T-t,t)+F(t,0)\big\}\\[5pt]
&=&\dis\int_0^T\!\!\di t\,\big\{F(T-t,t)-F(0,t)\big\}
-\int_0^T\!\!\di t\,\big\{F(t,T-t)-F(t,0)\big\},
\ec
where we have subsituted $t\mapsto T-t$ in the term $-F(T-t,t)$.
Since $X^1$ solves the martingale problem for $A_1$,
\be
E\big[\Psi(X^1_{T-t},x_2)\big]-E\big[\Psi(X^1_0,x_2)\big]
=\int_0^{T-t}\!\!\!\!\di s\;E\big[\Phi_1(X^1_s,x_2)\big]\qquad(x_2\in E_2),
\ee
and therefore, integrating the $x_2$-variable with respect to the law of
$X^2_t$, using the independence of $X^1$ and $X^2$ and (\ref{Phint2}), we
find that
\be\ba{l}
\dis\int_0^T\!\!\!\di t\,\big\{F(T-t,t)-F(0,t)\big\}
=\int_0^T\!\!\!\di t\,\big\{E\big[\Psi(X^1_{T-t},X^2_t)\big]
-E\big[\Psi(X^1_0,X^2_t)\big]\big\}\\[5pt]
\dis\qquad=\int_0^T\!\!\!\di t\,\int_0^{T-t}\!\!\!\!\di s
\;E\big[\Phi_1(X^1_s,X^2_t)\big]
=\int_0^T\!\!\!\di t\int_0^t\!\!\!\di s\;E\big[\Phi_1(X^1_{t-s},X^2_s)\big].
\ec
Treating the second term in the right-hand side of (\ref{splitC3}) in the
same way, we find that
\be
\int_0^T\!\!\!\di t\,\big\{F(t,0)-F(0,t)\big\}
=\int_0^T\!\!\!\di t\int_0^t\!\!\!\di s\;E\big[\Phi_1(X^1_{t-s},X^2_s)\big]
-\int_0^T\!\!\!\di t\int_0^t\!\!\!\di s\;E\big[\Phi_2(X^1_{t-s},X^2_s)\big].
\ee
Differentiating with respect to $T$ we arrive at (\ref{apdual}).\qed

\section{Construction and comparison}\label{cosec}

\subsection{Finite branching-coalescing particle systems}

For finite initial conditions, the $(a,b,c,d)$-braco-process $X$ can
be constructed explicitly using exponentially distributed random
variables. The only thing one needs to check is that $X$ does not
explode. This is part of the next proposition. Recall the definitions
of $\Ni(\La)$ and $\Si(\Ni(\La))$ from (\ref{Sidef}) and of $G$ from
(\ref{GdefC3}).
\bp{\bf(Finite braco-processes)}\label{finmart}
Let $X$ be the $(a,b,c,d)$-braco-process started in a finite state
$x$. Then $X$ does not explode. Moreover, with $z^{\li k\re}
:=z(z+1)\cdots(z+k-1)$, one has
\be\label{kmom}
E^x\big[|X|^{\li k\re}_t\big]\leq |x|^{\li k\re}e^{kbt}
\qquad(k=1,2,\ldots,\ t\geq 0).
\ee
For each $f\in\Si(\Ni(\La))$, one has $Gf\in\Si(\Ni(\La))$ and $X$ solves the martingale problem for the operator $G$ with domain $\Si(\Ni(\La))$.
\ep
{\bf Proof} Introduce stopping times $\tau_N:=\inf\{t\geq 0:|X_t|\geq N\}$.
Put $f^k_t(x):=|x|^{\li k\re}e^{-kbt}$. It is easy to see that
\be\label{subha}
\{G+\dif{t}\}f^k_t(x)\leq kb|x|^{\li k\re}e^{-kbt}-kb|x|^{\li k\re}e^{-kbt}=0.
\ee
\detail{Write $g^k(x):=|x|^{\li k\re}$. Then
\be
Gg^k(x)\leq kb|x|^{\li k\re}.
\ee
To see this, note that
\bc
\dis Gg^k(x)&=&\dis 0+b\sum_ix(i)\{(|x|+1)^{\li k\re}-|x|^{\li k\re}\}\\
&&\dis+c\sum_i\ffrac{1}{2}x(i)(x(i)-1)\{(|x|-1)^{\li k\re}-|x|^{\li k\re}\}
+d\sum_ix(i)\{(|x|-1)^{\li k\re}-|x|^{\li k\re}\},
\ec
where
\be\ba{l}
\dis(z+1)^{\li k\re}-z^{\li k\re}=(z+1)(z+2)\cdots(z+k)-z(z+1)\cdots(z+k-1)\\
\dis\quad=\{(z+k)-z\}(z+1)(z+2)\cdots(z+k-1),
\ec
and therefore
\be
z\{(z+1)^{\li k\re}-z^{\li k\re}\}=kz(z+1)(z+2)\cdots(z+k-1)=kz^{\li k\re}.
\ee}
The stopped process $(X_{t\wedge\tau_N})_{t\geq 0}$ is a jump process
in $\{x\in\N^\La:|x|\leq N\}$ with bounded jump rates, and therefore
standard theory tells us that the process $(M_t)_{t\geq 0}$ given by
\be
M_t:=f^k_{t\wedge\tau_N}(X_{t\wedge\tau_N})
-\int_0^{t\wedge\tau_N}\!\!\big(\{G+\dif{s}\}f^k_s\big)(X_s)\,\di s
\qquad(t\geq 0)
\ee
is a martingale. By (\ref{subha}), it follows that
$E^x\big[|X_{t\wedge\tau_N}|^{\li k\re}e^{-kb(t\wedge\tau_N)}\big]
\leq|x|^{\li k\re}$ and therefore
\be\label{expk}
E^x\big[|X_{t\wedge\tau_N}|^{\li k\re}\big]\leq |x|^{\li k\re}e^{kbt}
\qquad(k=1,2,\ldots,\ t\geq 0).
\ee
In particular, setting $k=1$, we see that
\be
NP^x[\tau_N\leq t]\leq E^x\big[|X_{t\wedge\tau_N}|\big]
\leq|x|e^{bt}\qquad(t\geq 0),
\ee
which shows that $\lim_{N\to\infty}P^x[\tau_N\leq t]=0$ for all $t\geq 0$,
i.e., the process does not explode. Taking the limit $N\up\infty$ in
(\ref{expk}), using Fatou, we arrive at (\ref{kmom}).

If $f\in\Si(\Ni(\La))$ then $f$ is bounded on sets of the form
$\{x\in\N^\La:|x|\leq N\}$, and therefore $Gf$ is well-defined.
By standard theory, the processes
$(M^N_t)_{t\geq 0}$ given by
\be\label{MNt}
M^N_t:=f(X_{t\wedge\tau_N})-\int_0^{t\wedge\tau_N}\!\! Gf(X_s)\di s
\qquad(t\geq 0)
\ee
are martingales. It is easy to see that $f\in\Si(\Ni(\La))$ implies
$Gf\in\Si(\Ni(\La))$, and therefore $\int_0^tE[|Gf(X_s)|\di s<\infty$
for all $t\geq 0$ by (\ref{kmom}). Using (\ref{expk}), one can now
check that for fixed $t\geq 0$, the random variables $\{M^N_t\}_{N\geq
1}$ are uniformly integrable. Taking the pointwise limit in
(\ref{MNt}), one can now check that $X$ solves the martingale problem
for $G$ with domain $\Si(\Ni(\La))$.\qed

\detail{To see that the $\{M^N_t\}_{N\geq 1}$ are uniformly
integrable, choose $K,k$ such that $|f(x)|,|Gf(x)|\leq
K(|x|+1)^k$. Then
\be
E^x[|f(X_{t\wedge\tau_N})|^{\frac{k+1}{k}}]
\leq KE^x[(|X_{t\wedge\tau_N}|+1)^{k+1}],
\ee
where, by (\ref{expk}), the right-hand side can be bounded uniformly
in $N$. This shows that the random variables
$\{f(X_{t\wedge\tau_N})\}_{N\geq 1}$ are uniformly
integrable. Likewise,
\be\ba{l}
\dis E\Big[\Big|\int_0^{t\wedge\tau_N}
\!\! Gf(X_s)\di s\Big|^{\frac{k+1}{k}}\Big]
=E\Big[|t\wedge\tau_N|^{\frac{k+1}{k}}\Big||t\wedge\tau_N|^{-1}
\int_0^{t\wedge\tau_N}\!\! Gf(X_s)\di s\Big|^{\frac{k+1}{k}}\Big]\\[5pt]
\dis\leq E\Big[|t\wedge\tau_N|^{\frac{k+1}{k}}|t\wedge\tau_N|^{-1}
\int_0^{t\wedge\tau_N}\!\!\big| Gf(X_s)\big|^{\frac{k+1}{k}}\di s\Big]\\[5pt]
\dis\leq KE\Big[|t\wedge\tau_N|^{\frac{1}{k}}
\int_0^{t\wedge\tau_N}\!\!(|X_s|+1)^{k+1}\di s\Big]\\[5pt]
\dis\leq Kt^{\frac{1}{k}}
\int_0^t\!\!E\big[(|X_{s\wedge\tau_N}|+1)^{k+1}\big]\di s,
\ec
which can again be bounded  uniformly in $N$ using (\ref{expk}).}

\detail{To see that for $f\in\Si(\Ni(\La))$ the process
$(M_t)_{t\geq 0}$ given by
\be\label{Mt}
M_t:=f(X_t)-\int_0^tGf(X_s)\di s\qquad(t\geq 0)
\ee
is a martingale with respect to the filtration generated by $X$,
recall that uniform integrability and almost sure convergence
imply convergence in $L^1$-norm. Now we are using the following fact:
If $M^{(n)}$ are martingales and $M^{(n)}_t\to M_t$ a.s.\ and in $L^1$
norm for each fixed $t\geq 0$, then $M$ is a martingale. To see that
this is true, note that the fact that $M^{(n)}$ is a martingale is
equivalent to the statement that $E[|M^{(n)}_t|]<\infty$ for all
$t\geq 0$ and
\be\label{martco}
E[M^{(n)}_tF(M^{(n)}_{s_1},\ldots,M^{(n)}_{s_n})]
=E[M^{(n)}_sF(M^{(n)}_{s_1},\ldots,M^{(n)}_{s_n})]
\ee
for each $0\leq s_1\leq\cdots\leq s_n\leq s\leq t$ and measurable
$F:\R^n\to[0,1]$. Since $M^{(n)}_t\to M_t$ a.s.\ and in $L^1$ norm, it
follows that $E[|M_t|]<\infty$ for all $t\geq 0$. One has
\be\ba{l}\label{MF}
\dis\big|E[M^{(n)}_tF(M^{(n)}_{s_1},\ldots,M^{(n)}_{s_n})]
-E[M_tF(M_{s_1},\ldots,M_{s_n})]\big|\\[5pt]
\dis\leq E[|M^{(n)}_t-M_t|F(M^{(n)}_{s_1},\ldots,M^{(n)}_{s_n})]
+E\big[|M_t|\,\big|F(M^{(n)}_{s_1},\ldots,M^{(n)}_{s_n})
-F(M_{s_1},\ldots,M_{s_n})\big|\big]\\[5pt]
\dis\leq E[|M^{(n)}_t-M_t|]
+E\big[|M_t|\,\big|F(M^{(n)}_{s_1},\ldots,M^{(n)}_{s_n})
-F(M_{s_1},\ldots,M_{s_n})\big|\big].
\ec
Taking the limit $n\to\infty$ in (\ref{martco}), using (\ref{MF}) and
dominated convergence, we find that
\be
E[M_tF(M_{s_1},\ldots,M_{s_n})]=E[M_sF(M_{s_1},\ldots,M_{s_n})].
\ee
Thus, $M$ is a martingale. It is easy to see that the process in
(\ref{Mt}) is adapted to the filtration generated by $X$ and therefore
it is also a martingale with respect to this filtration.}

\subsection{Monotonicity and subadditivity}

In this section we present two simple comparison results for finite
branching-coalescing particle systems.
\bl{\bf(Comparison of branching-coalescing particle systems)}\label{complem}
Let $X$ and $\ti X$ be the $(a,b,c,d)$-braco-process and the
$(a,\ti b,\ti c,\ti d)$-braco-process started in finite initial states
$x$ and $\ti x$, respectively. Assume that
\be
x\leq\ti x,\quad b\leq\ti b,\quad c\geq\ti c,\quad d\geq\ti d.
\ee
Then $X$ and $\ti X$ can be coupled in such a way that
\be\label{comp}
X_t\leq\ti X_t\qquad(t\geq 0).
\ee
\el
{\bf Proof} We will construct a bivariate process $(B,W)$, say of
black and white particles, such that $X=B$ are the black particles and
$\ti X=B+W$ are the black and white particles together. To this aim,
we let the particles evolve in such a way that black and white
particles branch with rates $b$ and $\ti b$, respectively, and
additionally black particles give birth to white particles with rate
$\ti b-b$. Moreover, all pairs of particles coalesce with rate $2\ti
c$, where the new particle is black if at least one of its parents is
black, and additionally each pair of black particles is with rate
$2c-2\ti c$ replaced by a pair consisting of one black and one white
particle. Finally, all particles die with rate $\ti d$, and
additionally, black particles change into white particles with rate
$d-\ti d$. It is easy to see that with these rules, $X$ and $\ti X$
are the $(a,b,c,d)$-braco-process and the $(a,\ti b,\ti c,\ti
d)$-braco-process, respectively.\qed

\noi
The next lemma has been proved for $\La=\Z^d$ in
\cite[Lemma~2.2]{SU86}. It can be proved (with particles in three
colors) in a similar way as the previous lemma.
\bl{\bf (Subadditivity)}\label{subad}
Let $X,Y,Z$ be $(a,b,c,d)$-braco-processes started in finite initial states
$x,y$, and $x+y$, respectively. Then $X,Y,Z$ may be coupled in such a
way that $X$ and $Y$ are independent and
\be
Z_t\leq X_t+Y_t\qquad(t\geq 0).
\ee
\el

\detail{{\bf Proof} Consider a trivariate process $(B,W,R)$, say of black,
white and red particles, with initial condition $(x,y,0)$, such that
each color evolves as an autonomous $(a,b,c,d)$-braco-processes, and
additionally, pairs of black and white particles are with rate $2c$
replaced by a black and a red particle, and pairs of white and red
particles are with rate $2c$ replaced by one white particle. It is not
hard to see that $X:=B$ (black), $Y:=W+R$ (white + red) and $Z:=B+W$
(black + white) are $(a,b,c,d)$-braco-systems, and that $X$ and $Y$
are independent.\qed}

\subsection{Infinite branching-coalescing particle systems}
\label{appsec}

In this section we carry out the construction of branching-coalescing
particle systems for infinite initial conditions. We will also derive
two results on the approximation of infinite systems with finite
systems, that are needed later on. Except for the statement about
sample paths, the next proposition has been proved in \cite{Che87},
but we give a proof here for the sake of completeness.
\bp{\bf (Construction of branching-coalescing particle systems)}\label{Xcon}
For each $f\in\Ci_{\rm Lip}(\Ei_\ga(\La))$ and $t\geq 0$, the function
$S_tf$ defined in (\ref{Stdef}) can be extended to a unique Lipschitz
function on $\Ei_\ga(\La)$, also denoted by $S_tf$. There exists a
unique (in distribution) time-homogeneous Markov process with cadlag
sample paths in the space $\Ei_\ga(\La)$ equipped with the norm
$\|\cdot\|_\ga$, such that
\be
E^x[f(X_t)]=S_tf(x)\qquad(f\in\Ci_{\rm Lip}(\Ei_\ga(\La)),
\ x\in\Ei_\ga(\La),\ t\geq 0).
\ee
\ep
We start with the following lemma.
\bl{\bf(Action of the semigroup on Lipschitz functions)}\label{Liplem}
If $f:\Ni(\La)\to\R$ is Lipschitz continuous in the norm
$\|\cdot\|_\ga$ from (\ref{ganorm}), with Lipschitz constant $L$, and
$K$ is the constant from (\ref{gacon}), then
\be\label{Lipest}
|S_tf(x)-S_tf(y)|\leq Le^{(K+b-d)t}\|x-y\|_\ga
\qquad(x,y\in\Ni(\La),\ t\geq 0).
\ee
\el
{\bf Proof} It follows from Propostion~\ref{finmart} that
$\dif{t}E[f(X_t)]=E[Gf(X_t)]$ for all $f\in\Si(\Ni(\La))$, $t\geq
0$. Applying this to the function $f(x):=\|x\|_\ga$ we see that
\bc
\dis\dif{t}E^x[\|X_t\|_\ga]&=&\dis\sum_{ij}a(i,j)(\ga_j-\ga_i)E[X_t(i)]
+(b-d)E^x[\|X_t\|_\ga]\\[5pt]
&&\dis-c\sum_i\ga_iE[X_t(i)(X_t(i)-1)]\;\leq\;
\dis(K+b-d)E[\|X\|_\ga],
\ec
and therefore
\be\label{Kbd}
E^x[\|X_t\|_\ga]\leq e^{(K+b-d)t}\|x\|_\ga\qquad(x\in\Ni(\La)).
\ee
\detail{To see this, set $f(x):=\|x\|_\ga$ and let $G$ be the
generator of the $(a,b,c,d)$-braco-process. Then
\bc
Gf(x)&=&\dis\sum_{ij}a(i,j)x(i)(\ga_j-\ga_i)
+(b-d)\sum_ix(i)\ga_i-c\sum_ix(i)(x(i)-1)\ga_i\\[5pt]
&\leq &\dis\sum_ix(i)K\ga_i+b\sum_ix(i)\ga_i=(K+b)\|x\|_\ga,
\ec
which implies that $\dif{t}E^x[\|\ti X_t\|_\ga]
\leq(K+b)E^x[\|\ti X_t\|_\ga]$.}
Let $X^x$ denote the $(a,b,c,d)$-braco-process started in
$x$. By Lemma~\ref{complem}, we can couple $X^x$, $X^y$, $X^{x\wedge y}$,
and $X^{x\vee y}$ such that $X^{x\wedge y}_t\leq X^x_t,X^y_t\leq
X^{x\vee y}_t$ for all $t\geq 0$. It follows that
\be\label{Kbe}
E[\|X^x_t-X^y_t\|_\ga]\leq E[\|X^{x\vee y}_t-X^{x\wedge y}_t\|_\ga].
\ee
By Lemma~\ref{subad}, we can couple $X^{x\wedge y}$ and $X^{x\vee y}$
to the process $X^{|x-y|}$ such that $X^{x\vee y}_t\leq X^{x\wedge
y}_t+X^{|x-y|}_t$ for all $t\geq 0$. Therefore, by (\ref{Kbe}) and
(\ref{Kbd}),
\be\label{Ega}
E[\|X^x_t-X^y_t\|_\ga]\leq E[\|X^{|x-y|}_t\|_\ga]\leq\|x-y\|_\ga e^{(K+b-d)t},
\ee
which implies that
\be
|S_tf(x)-S_tf(y)|\leq E[|f(X^x_t)-f(X^y_t)|]\leq LE[\|X^x_t-X^y_t\|_\ga]
\leq L\|x-y\|_\ga e^{(K+b-d)t},
\ee
as required.\qed

\noi
Since Lipschitz functions on $\Ni(\La)$ have a unique Lipschitz
extension to $\Ei_\ga(\La)$, Lemma~\ref{Liplem} implies that $S_tf$
can be uniquely extended to a function in $\Ci_{\rm
Lip}(\Ei_\ga(\La))$ for each $f\in\Ci_{\rm Lip}(\Ei_\ga(\La))$.
\bl{\bf(Construction of the process for fixed times)}\label{fixt}
Let $X^{(n)}$ be $(a,b,c,d)$-braco-processes started in initial states
$x^{(n)}\in\Ni(\La)$ such that $x^{(n)}\up x$ for some
$x\in\Ei_\ga(\La)$. Then the $X^{(n)}$ may be coupled such that
$X^{(n)}_t\up X_t$ $(t\geq 0)$ for some $\ov\N^\La$-valued process
$X=(X_t)_{t\geq 0}$. The process $X$ satisfies $X_t\in\Ei_\ga(\La)$
a.s.\ $\forall t\geq 0$ and $X$ is a Markov process with semigroup
$(S_t)_{t\geq 0}$.
\el
{\bf Proof} It follows from Lemma~\ref{complem} that the $X^{(n)}$ can
be coupled such that $X^{(n)}_t\leq X^{(n+1)}_t$ ($t\geq 0$), and
therefore $X^{(n)}_t\up X_t$ ($t\geq 0$) for some $\ov\N^\La$-valued
random variables $X_t$. By (\ref{Ega}),
\be\label{hiero}
E\big[\|X_t-X^{(n)}_t\|_\ga\big]
=\lim_{m\up\infty}E\big[\|X^{(m)}_t-X^{(n)}_t\|_\ga\big]
\leq\|x-x^{(n)}\|_\ga e^{(K+b-d)t}.
\ee
This shows in particular that $E[\|X_t\|_\ga]<\infty$ and therefore
$X_t\in\Ei_\ga(\La)$ a.s.\ $\forall t\geq 0$. If $f\in\Ci_{\rm
Lip}(\Ei_\ga(\La))$ has Lipschitz constant $L$, then by (\ref{hiero}),
\be\ba{l}
\dis|E[f(X_t)]-E[f(X^{(n)}_t)]|\leq E[|f(X_t)-f(X^{(n)}_t)|]\\[5pt]
\dis\qquad\leq LE[\|X_t-X^{(n)}_t\|_\ga]\leq L\|x-x^{(n)}\|_\ga e^{(K+b-d)t},
\ec
and therefore
\be
E[f(X_t)]=\lim_{n\up\infty}E[f(X^{(n)}_t)]=\lim_{n\up\infty}S_tf(x^{(n)})
=S_tf(x).
\ee
This proves that for each $x\in\Ei_\ga(\La)$ and $t\geq 0$ there exists
a probability measure $P_t(x,\cdot)$ on $\Ei_\ga(\La)$ such that $\int
P_t(x,\di y)f(y)=S_tf(x)$ for all $f\in\Ci_{\rm
Lip}(\Ei_\ga(\La))$. We need to show that $X$ is the Markov process
with transition probabilities $P_t(x,\di y)$. Let $\Ci_{\rm
Lip,b}(\Ei_\ga(\La))$ denote the class of bounded Lipschitz functions
on $\Ei_\ga(\La)$. Then $\Ci_{\rm Lip,b}(\Ei_\ga(\La))$ is closed
under multiplication and $S_t$ maps $\Ci_{\rm Lip,b}(\Ei_\ga(\La))$
into itself. Therefore, for all $0\leq t_0<\cdots<t_k$ and
$f_1,\ldots,f_k\in\Ci_{\rm Lip,b}(\Ei_\ga(\La))$, one has
\detail{Indeed, if $f,g$ are bounded and Lipschitz with Lipschitz
constants $L_f,L_g$, repectively, then $fg$ is bounded and Lipschitz
with Lipschitz constant $\|g\|_\infty L_f+\|f\|_\infty L_g$.}
\be\label{Lico}
E\big[f_1(X^{(n)}_{t_1})\cdots f_k(X^{(n)}_{t_k})\big]
=S_{t_1}f_1S_{t_2-t_1}f_2\cdots S_{t_k-t_{k-1}}f_k(x^{(n)}).
\ee
It follows from (\ref{hiero}) that
\be\label{prodli}
\big|E\big[f_1(X_{t_1})\cdots f_k(X_{t_k})\big]
-E\big[f_1(X^{(n)}_{t_1})\cdots f_k(X^{(n)}_{t_k})\big]\big|
\leq\|x-x^{(n)}\|_\ga\sum_{i=1}^kL_ie^{(K+b-d)t_k}\prod_{j\neq i}\|f_j\|_\infty,
\ee
where $L_i$ is the Lipschitz constant of $f_i$. Taking the limit
$n\up\infty$ in (\ref{Lico}), using (\ref{prodli}), we see that
\be\label{SfS}
E\big[f_1(X_{t_1})\cdots f_k(X_{t_k})\big]
=S_{t_1}f_1S_{t_2-t_1}f_2\cdots S_{t_k-t_{k-1}}f_k(x),
\ee
i.e., $X$ is the Markov process with semigroup $(S_t)_{t\geq 0}$.\qed

\noi
{\bf Proof of Proposition~\ref{Xcon}} We need to show that the process
$X$ from Lemma~\ref{fixt} satisfies $X_t\in\Ei_\ga(\La)$ $\forall
t\geq 0$ a.s.\ (and not just for fixed times) and that
$(X_t)_{t\geq 0}$ has cadlag sample paths with respect to the norm
$\|\cdot\|_\ga$. It suffices to prove these facts on the time interval
$[0,1]$. We will do this by constructing an $\Ei_\ga(\La)$-valued
process $Z$ such that $Z$ makes only upward jumps, and the number of
upward jumps of $Z$ dominates the number of upward jumps of $X$.

Couple the process $X^{(n)}$ from Lemma~\ref{fixt} to a process
$Y^{(n)}$ such that the joint process $(X^{(n)},Y^{(n)})$ is the
Markov process in $\Ni(\La)\times\Ni(\La)$ with generator
\be\ba{l}\label{GXYdef}
G_{X,Y}f(x,y):=\\[5pt]
\dis\sum_{ij}a(i,j)x(i)\{f(x+\de_j-\de_i,y+\de_i)-f(x,y)\}
+\sum_{ij}a(i,j)y(i)\{f(x,y+\de_j)-f(x,y)\}\\
\dis+b\sum_ix(i)\{f(x+\de_i,y)-f(x,y)\}+b\sum_iy(i)\{f(x,y+\de_i)-f(x,y)\}\\
\dis+c\sum_ix(i)(x(i)-1)\{f(x-\de_i,y+\de_i)-f(x,y)\}
+d\sum_ix(i)\{f(x-\de_i,y+\de_i)-f(x,y)\}.
\ec
and initial state $(X^{(n)}_0,Y^{(n)}_0)=(x^{(n)},0)$. Indeed, it is not
hard to see that the first component of the process with generator
$G_{X,Y}$ is the $(a,b,c,d)$-braco-process, and that
$Z^{(n)}:=X^{(n)}+Y^{(n)}$ is the Markov process in $\Ni(\La)$ with
generator
\be\label{GZdef}
G_Zf(z):=\sum_{ij}a(i,j)z(i)\{f(z+\de_j)-f(z)\}
+b\sum_iz(i)\{f(z+\de_i)-f(z)\}
\ee
and initial state $Z^{(n)}_0=x^{(n)}$.  In analogy with (\ref{Kbd}) it is
easy to check that
\detail{Set $f(z):=\|z\|_\ga$. Then
\be
G_Zf(z)=\sum_{ij}a(i,j)\ga_j z(i)+b\sum_i\ga_i z(i)\leq(K+b)\sum+i\ga_iz(i).
\ee}
\be\label{Zexp}
E^z[\|Z^{(n)}_t\|_\ga]\leq\|x^{(n)}\|_\ga e^{(K+b)t}\qquad(z\in\Ni(\La),\ t\geq 0).
\ee
$Z^{(n)}$ makes only upward jumps and $Z^{(n)}(i)$ makes at least as
many upward jumps as $X^{(n)}(i)$. Since $X^{(n)}(i)$ cannot become
negative, it follows that
\be
|\{t\in[0,1]:X^{(n)}_{t-}(i)\neq X^{(n)}_t(i)\}|\leq x^{(n)}(i)+2Z^{(n)}_1(i).
\ee
Summing with respect to the $\ga_i$, taking expectations, using
(\ref{Zexp}), we see that
\be\label{finjump}
\sum_i\ga_i\,E\big[|\{t\in[0,1]:X^{(n)}_{t-}(i)\neq X^{(n)}_t(i)\}|\big]
\leq\|x^{(n)}\|_\ga(1+2e^{K+b}).
\ee
Let $Z$ be the increasing limit of the processes $Z^{(n)}$. It follows from
(\ref{Zexp}) that $Z_1\in\Ei_\ga(\La)$ a.s. Now
\be
X_t,X_{t-}\leq Z_t\leq Z_1\qquad\forall t\in[0,1]\quad {\rm a.s.},
\ee
and therefore $X_t,X_{t-}\in\Ei_\ga(\La)$ $\forall t\in[0,1]$ a.s.
Since a.s.\ all jumps occur at different times,
\be
|\{t\in[0,1]:X^{(n)}_{t-}(i)\neq X^{(n)}_t(i)\}|
\up|\{t\in[0,1]:X_{t-}(i)\neq X_t(i)\}|\quad\mbox{as }n\up\infty.
\ee
Thus, taking the limit $n\up\infty$ in (\ref{finjump}) we see that
\be\label{prejumpbound}
\sum_i\ga_i\,E\big[|\{t\in[0,1]:X_{t-}(i)\neq X_t(i)\}|\big]
\leq\|x\|_\ga(1+2e^{K+b}).
\ee
This proves that $X$ has a.s.\ componentwise cadlag sample paths. If
$1\geq t_n\down t$, then $X_{t_n}\to X_t$ pointwise and
$|X_{t_n}-X_t|\leq 2Z_1$, and therefore, by dominated convergence,
\be
\|X_{t_n}-X_t\|_\ga=\sum_i\ga_i|X_{t_n}(i)-X_t(i)|\to 0.
\ee
The same argument shows that $X_{t_n}\to X_{t-}$ for $t_n\up t\leq 1$,
i.e., $X$ has cadlag sample paths with respect to the norm
$\|\cdot\|_\ga$.\qed

\noi
The proof of Proposition~\ref{Xcon} yields a useful corollary.
\bcor{\bf(Locally finite number of jumps)}\label{Finjump}
The $(a,b,c,d)$-braco-process $X$ satisfies
\be\label{jumpbound}
\sum_i\ga_i\,E^x\big[|\{t\in[0,1]:X_{t-}(i)\neq X_t(i)\}|\big]
\leq\|x\|_\ga(1+2e^{K+b}).
\ee
\ecor
We can now prove two approximation lemmas.
\bl{\bf(Convergence of finite dimensional distributions)}\label{fddcon}
Let $X^{x_n},X^x$ be the $(a,b,c,d)$-braco-process started in initial
states $x_n,x\in\Ei_\ga(\La)$, respectively, such that
\be
\lim_{n\to\infty}\|x_n-x\|_\ga=0.
\ee
Then, for all $0\leq t_1<\cdots<t_k$, one has
\be\label{weto}
(X^{(n)}_{t_1},\ldots,X^{(n)}_{t_k})
\Rightarrow(X_{t_1},\ldots,X_{t_k})\qquad\mbox{as }n\to\infty.
\ee
\el
{\bf Proof} Use (\ref{SfS}) for $x_n$ and then let $n\to\infty$.\qed

\detail{In particular, we may take for $f_1,\ldots,f_k$ functions
of the form
\be
x\mapsto 1_{\{x(i)=y(i)\ \forall i\in\De\}},
\ee
where $\De\sub\La$ is finite and $y\in\N^\De$, and in this way we see that
\be
(X^{(n)}_{t_1},\ldots,X^{(n)}_{t_k})\Rightarrow(X_{t_1},\ldots,X_{t_k}).
\ee
(To see that the function $x\mapsto 1_{\{x(i)=y(i)\ \forall
i\in\De\}}$ is Lipschitz continuous, choose $\eps>0$ such that
$\|x-y\|_\ga\leq\eps$ implies $x(i)=y(i)$ for all $i\in\De$. Then
$x\mapsto 1_{\{x(i)=y(i)\forall i\in\De\}}$ is Lipschitz continuous
with Lipschitz constant $\eps^{-1}$.)}

\bl{\bf(Monotonicities for infinite systems)}\label{mocolem}
Lemmas~\ref{complem} and \ref{subad} also hold for infinite initial
states. If $X^x,X^{x_n}$ are $(a,b,c,d)$-braco-process started in
initial states $x,x_n\in\Ei_\ga(\La)$, such that $x_n\up x$, then
$X^x,X^{x_n}$ may be coupled such that
\be\label{moco}
X^{x_n}_t(i)\up X^x_t(i)\quad\mbox{as }n\up\infty
\quad\forall i\in\La,\ t\geq 0\quad{\rm a.s.}
\ee
\el
{\bf Proof} The proof of Proposition~\ref{Xcon} shows that
(\ref{moco}) holds if the $x_n$ are finite. To generalize
Lemma~\ref{complem} to infinite initial states $x,\ti x$, it therefore
suffices to note that if $x\leq\ti x$, then there exist finite
$x_n\leq\ti x_n$ such that $x_n\up x$ and $\ti x_n\up\ti x$, and then
take the limit $n\up\infty$ in (\ref{comp}) using
(\ref{moco}). Lemma~\ref{subad} can be generalized to infinite $x,y$
by approximation with finite $x_n,y_n$ in the same way. Finally, to
see that (\ref{moco}) remains valid if the $x_n$ are infinite, note
that by Lemma~\ref{complem} (which has now been proved in the infinite
case), the processes $X^{x_n}$ can be coupled such that
$X^{x_n}_t(i)\leq X^{x_{n+1}}_t(i)$ for all $i\in\La$ and $t\geq
0$. Denote the increasing limit of the $X^{x_n}$ by
$X^x$. Lemma~\ref{fddcon} shows that $X^x$ has the same finite
dimensional distributions as the $(a,b,c,d)$-braco-process started in
$x$ and it follows from Corollary~\ref{Finjump} that $X^x$ has
componentwise cadlag sample paths, so $X^x$ is a version of
the $(a,b,c,d)$-braco-process started in $x$.\qed

\subsection{Construction and comparison of resampling-selection processes}
\label{resemcon}

We equip the space $[0,1]^\La$ with the product topology and let
$\Ci([0,1]^\La)$ denote the space of continuous real functions on
$[0,1]^\La$, equipped with the supremum norm. By $\Ci^2_{\rm
fin}([0,1]^\La)$ we denote the space of $\Ci^2$ functions on $[0,1]^\La$
depending on finitely many coordinates. By definition, $\Ci^2_{\rm
sum}([0,1]^\La)$ is the space of continuous functions $f$ on
$[0,1]^\La$ such that the partial derivatives $\dif{\phi(i)}f(\phi)$ and
$\difif{\phi(i)}{\phi(j)}f(\phi)$ exist for each $x\in(0,1)^\La$ and such that
the functions
\be\ba{l}
\phi\mapsto\big(\dif{\phi(i)}f(\phi)\big)_{i\in\La}
\quad\mbox{and}\quad
\phi\mapsto\big(\difif{\phi(i)}{\phi(j)}f(\phi)\big)_{i,j\in\La}
\ec
can be extended to continuous functions from $[0,1]^\La$ into the
spaces $\ell^1(\La)$ and $\ell^1(\La^2)$ of absolutely summable
sequences on $\La$ and $\La^2$, respectively, equipped with the
$\ell^1$-norm. Define an operator $\Gi:\Ci^2_{\rm
sum}([0,1]^\La)\to\Ci([0,1]^\La)$ by
\bc\label{GidefC3}
\Gi f(\phi)&:=&\dis\sum_{ij}a(j,i)(\phi(j)-\phi(i))\dif{\phi(i)}f(\phi)
+b\sum_i\phi(i)(1-\phi(i))\dif{\phi(i)}f(\phi)\\[2pt]
&&\dis+c\sum_i\phi(i)(1-\phi(i))\diff{\phi(i)}f(\phi)
-d\sum_i\phi(i)\dif{\phi(i)}f(\phi)\qquad(\phi\in[0,1]^\La).
\ec
One can check that for $f\in\Ci^2_{\rm sum}([0,1]^\La)$, the infinite
sums converge in the supremumnorm and the result does not depend on
the summation order \cite[Lemma~3.4.4]{Swa99}. If a $[0,1]^\La$-valued
process $\Xc$ solves the martingale problem for $\Gi$ with domain
$\Ci_{\rm fin}([0,1]^\La)$, then also for the larger domain $\Ci_{\rm
sum}([0,1]^\La)$ (see \cite[Lemma~3.4.5]{Swa99}).

Let $\Ci_{[0,1]^\La}\half$ denote the space of continuous functions
from $\half$ into $[0,1]^\La$, equipped with the topology of uniform
convergence on compacta.  If $\Xc^{(n)},\Xc$ are
$\Ci_{[0,1]^\La}\half$-valued random variables, then we say that
$\Xc^{(n)}$ converges in distribution to $\Xc$, denoted as
$\Xc^{(n)}\Rightarrow\Xc$, when $\Li(\Xc^{(n)})$ converges weakly to
$\Li(\Xc)$. Convergence in distribution implies convergence of the
finite-dimensional distributions (see \cite[Theorem~3.7.8]{EK}). The
fact that a $\Ci_{[0,1]^\La}\half$-valued random variable $\Xc$ solves
the martingale problem for $\Gi$ is a property of the law of $\Xc$
only. Standard results from \cite{EK} yield the following (for the
details, see for example Lemma~4.1 in \cite{Swa00}):
\bl{\bf(Existence and compactness of solutions to the martingale problem)}
\label{exislem}
For each $\phi\in[0,1]^\La$, there exists a solution $\Xc$ to the
martingale problem for $\Gi$ with initial state $\Xc_0=\phi$, and each
solution to the martingale problem for $\Gi$ has continuous sample
paths. Moreover, the space $\{\Li(\Xc):\Xc\mbox{ solves the martingale
problem for }\Gi\}$ is compact in the topology of weak convergence.
\el
If $\Xc$ solves the SDE (\ref{sde}), then $\Xc$ solves the martingale
problem for $\Gi$. Conversely, each solution to the martingale problem
for $\Gi$ is equal in distribution to some (weak) solution of the SDE
(\ref{sde}). Thus, existence of (weak) solutions to (\ref{sde})
follows from Lemma~\ref{exislem}. Distribution uniqueness of solutions
to (\ref{sde}) follows from pathwise uniqueness, which is in turn
implied by the following comparison result.
\bl{\bf(Monotone coupling of linearly interacting diffusions)}\label{sdemon}
Let $I\sub\R$ be a closed interval, let $\sig:I\to\R$ be
H\"older-$\frac{1}{2}$-continuous, and let $b_1,b_2:I\to\R$ be Lipschitz
continuous functions such that $b_1\leq b_2$. Let $\Xc^\al$
$(\al=1,2$) be solutions, relative to the same system of Brownian
motions, of the SDE
\be
\di\Xc^\al_t(i)=\sum_ja(j,i)(\Xc^\al_t(j)-\Xc^\al_t(i))\di t
+b_\al(\Xc^\al_t(i))\di t+\sig(\Xc^\al_t(i))\di B_t(i).
\ee
$(i\in\La,\ t\geq 0,\ \al=1,2)$. Then
\be
\Xc^1_0\leq\Xc^2_0\quad\mbox{implies}\quad\Xc^1_t\leq\Xc^2_t
\quad\forall t\geq 0\quad{\rm a.s.}
\ee
\el
{\bf Proof (sketch)} Set $\De_t(i):=\Xc^1_t(i)-\Xc^2_t(i)$ and write
$x^+:=x\vee 0$. Using an appropriate smoothing of the function
$x\mapsto x^+$ in the spirit of \cite[Theorem~1]{YW71} and arguing as
in the proof of \cite[Theorem~3.2]{SS80}, one can show that
\be
E[\|\De^+_t\|_\ga]\leq (K+L)\int_0^tE[\|\De^+_s\|_\ga]\di s,
\ee
where $\|\cdot\|_\ga$ is the norm from (\ref{ganorm}), $K$ is the
constant from (\ref{gacon}), and $L$ is the Lipschitz-constant of
$b_2$. The result now follows from Gronwall's inequality.\qed

\detail{{\bf Proof} We adapt a technique due to Yamada and Watanabe
\cite{YW71} to the infinite dimensional case, very much in the spirit
of Theorem~3.2 in \cite{SS80}. Set $\De_t(i):=\Xc^1_t(i)-\Xc^2_t(i)$
and write $x^+:=x\vee 0$. We will show that $\De_0^+(i)=0$ ($i\in\La$)
implies $\De_t^+(i)=0$ for all $t\geq 0$ and $i\in\La$.

Since $\int_{0+}\frac{\di x}{x}=\infty$, it is not hard to see that we can
choose continuous functions $\rho_n:\R\to\R$ such that
\be
0\leq\rho_n(x)\leq\frac{1}{nx}1_{(0,\frac{1}{n})}(x)\quad\mbox{and}
\quad\int_0^\frac{1}{n}\!\di x\,\rho_n(x)=1.
\ee
Define twice continuously differentiable functions $\phi_n$ by
\be
\phi_n(x):=\int_0^x\!\di y\int_0^y\!\di z\,\rho_n(z)\qquad(x\in\R).
\ee
Then,
\be\ba{rll}\label{phin}
{\rm (i)}&0\leq\phi_n(x)\up x^+\qquad&\mbox{as }n\up\infty,\\[5pt]
{\rm (ii)}&0\leq\phi'_n(x)\up 1_{\{x>0\}}\qquad&\mbox{as }n\up\infty,\\[5pt]
{\rm (iii)}&0\leq|x|\phi''_n(x)\leq\ffrac{1}{n}.
\ec
Now
\be
\di\De_t(i)=\sum_ja(j,i)(\De_t(j)-\De_t(i))\di t
+\big(b_1(\Xc^1_t(i))-b_2(\Xc^2_t(i))\big)\di t
+\big(\sig(\Xc^1_t(i))-\sig(\Xc^2_t(i))\big)\di B_t(i).
\ee
It\^o's formula gives
\bc\label{dphi}
\dis\di\phi_n(\De_t(i))&=&
\dis\sum_ja(j,i)(\De_t(j)-\De_t(i))\phi'_n(\De_t(i))\di t
+\big(b_1(\Xc^1_t(i))-b_2(\Xc^2_t(i))\big)\phi'_n(\De_t(i))\di t\\[5pt]
&&\dis+\ffrac{1}{2}\big(\sig(\Xc^1_t(i))-\sig(\Xc^2_t(i))\big)^2
\phi''_n(\De_t(i))\di t\quad+\quad\mbox{martingale terms.}
\ec
Since $\sig$ is H\"older-$\frac{1}{2}$-continuous, by (\ref{phin})~(iii)
\be
\big(\sig(\Xc^1_t(i))-\sig(\Xc^2_t(i))\big)^2\phi''_n(\De_t(i))
\leq K|\De_t(i)|\phi''_n(\De_t(i))\leq\frac{1}{n}.
\ee
Taking expectations in (\ref{dphi}) and letting $n\up\infty$, using
(\ref{phin}), we get
\bc\label{Deplus}
E[\De^+_t(i)]&=&\dis\int_0^tE\big[\sum_ja(j,i)(\De_s(j)-\De_s(i))
1_{\{\De_s(i)>0\}}\di s\big]\\[5pt]
&&\dis+\int_0^tE\big[\big(b_1(\Xc^1_s(i))-b_2(\Xc^2_s(i))\big)
1_{\{\De_s(i)>0\}}\di s\big]\\[5pt]
&\leq &\dis\int_0^tE\big[\sum_ja(j,i)\De^+_s(j)\di s\big]
+L\int_0^tE\big[\De^+_s(i)\di s\big],
\ec
Here we have used that $E[\De^+_0(i)]=0$, that
$(\De_s(j)-\De_s(i))1_{\{\De_s(i)>0\}}$ can only be positive when
$\De_s(j)>\De_s(i)>0$, and that
$\big(b_1(\Xc^1_s(i))-b_2(\Xc^2_s(i))\big)1_{\{\De_s(i)>0\}}$ can only
be positive when $\Xc^1_s(i)>\Xc^2_s(i)$, in which case
\be
b_1(\Xc^1_s(i))-b_2(\Xc^2_s(i))\leq b_2(\Xc^1_s(i))-b_2(\Xc^2_s(i))
\leq L|\Xc^1_s(i)-\Xc^2_s(i)|=L\De^+_s(i),
\ee
where $L$ is the Lipschitz-constant of $b_2$. Recall the norm
$\|\cdot\|_\ga$ from (\ref{ganorm}) and let $K$ be the constant from
(\ref{gacon}). It is easy to see that (\ref{Deplus}) implies
\be
E[\|\De^+_t\|_\ga]\leq (K+L)\int_0^tE[\|\De^+_s\|_\ga]\di s.
\ee
The result now follows from Gronwall's inequality.\qed}

\bcor{\bf(Comparison of resampling-selection processes)}\label{sdecomp}
Assume that $\Xc,\ti\Xc$ are solutions to the SDE (\ref{sde}), relative
to the same collection of Brownian motions, with parameters $(a,b,c,d)$
and $(a,\ti b,c,\ti d)$ and starting in initial states $\phi,\ti\phi$,
respectively. Assume that
\be\label{mubd}
\phi\leq\ti\phi,\quad d-b\geq\ti d-\ti b,\quad d\geq\ti d.
\ee
Then
\be
\Xc_t\leq\ti\Xc_t\quad\forall t\geq 0\quad{\rm a.s.}
\ee
\ecor
{\bf Proof} Immediate from Lemma~\ref{sdemon} and the fact that by
(\ref{mubd}), $bx(1-x)-dx\leq\ti bx(1-x)-\ti dx$ for all
$x\in[0,1]$.\qed

\noi
Our next lemma shows that resampling-selection processes with finite
initial mass have finite mass at all later times. The estimate
(\ref{linex}) is not very good if $b-d<0$, but it suffices for our
purposes.
\bl{\bf\hspace{-2pt}(Summable resampling-selection processes)\hspace{-1.71pt}}\label{summab}
Let $\Xc$ be the $(a,b,c,d)$-resem-process started in
$x\in[0,1]^\La$ with $|x|<\infty$. Set $r:=(b-d)\vee 0$. Then
\be\label{linex}
E^x\big[|\Xc_t|\big]\leq|x|e^{rt}\qquad(t\geq 0),
\ee
and $|\Xc_t|<\infty$ $\forall t\geq 0$ a.s.
\el
{\bf Proof} Without loss of generality we may assume that $b\geq d$;
otherwise, using Corollary~\ref{sdecomp}, we can bound $\Xc$ from
above by a braco-process with a higher $b$. Set $r:=b-d$ and put
$\Yi_t(i):=\Xc_t(i)e^{-rt}$. By It\^o's formula,
\be\label{dYi}
\di\Yi_t(i)=\sum_ja(j,i)(\Yi_t(j)-\Yi_t(i))\,\di t
-be^{-rt}\Xc_t(i)^2\di t+e^{-rt}\sqrt{c\Xc_t(i)(1-\Xc_t(i))}\,\di B_t(i).
\ee
Set $\tau_N:=\inf\{t\geq 0:|\Xc_t|\geq N\}$. Integrate (\ref{dYi}) up
to $t\wedge\tau_N$ and sum over $i$. The motion terms yield
\be\ba{l}
\dis\int_0^{t\wedge\tau_N}\!\!\sum_{ij}a(j,i)(\Yi_s(j)-\Yi_s(i))\,\di s\\[5pt]
\dis\quad
=\int_0^{t\wedge\tau_N}\!\!\sum_j\big(\sum_ia(j,i)\big)\Yi_s(j)\,\di s
-\int_0^{t\wedge\tau_N}\!\!\sum_i\big(\sum_ja^\dgg(i,j)\big)\Yi_s(i)\,\di s=0,
\ec
where the infinite sums converge in a bounded pointwise way since
$|Y_s|\leq N$ for $s\leq\tau_N$. It follows that
\be\label{ttN}
|\Yi_{t\wedge\tau_N}|
=|x|-b\sum_i\int_0^{t\wedge\tau_N}\!\!\Xc_s(i)^2e^{-rs}\di s
+\sum_i\int_0^{t\wedge\tau_N}\!\!
\sqrt{c\Xc_s(i)(1-\Xc_s(i))}\,e^{-rs}\di B_s(i),
\ee
provided we can show that the infinite sum of stochastic integrals
converges. Indeed, for any finite $\De\sub\La$, by the It\^o isometry,
\detail{To see that the It\^o isometry also holds for integrals up to
random time $t\wedge\tau$, where $\tau$ is a stopping time, write
\be\ba{l}
\dis E\Big[\Big|\int_0^{t\wedge\tau}\!\!X_s\di B_s\Big|^2\Big]
=E\Big[\Big|\int_0^t1_{\{s\leq\tau\}}X_s\di B_s\Big|^2\Big]
=E\Big[\int_0^t|1_{\{s\leq\tau\}}X_s|^2\di s\Big]\\[5pt]
\dis\qquad=E\Big[\int_0^t1_{\{s\leq\tau\}}|X_s|^2\di s\Big]
=E\Big[\int_0^{t\wedge\tau}\!\!|X_s|^2\di s\Big],
\ec
where we used the usual It\^o isometry and the fact that
$s\mapsto 1_{\{s\leq\tau\}}X_s$ is adapted.}
\be\ba{l}
\dis\sum_{i\in\De}E\Big[\Big|\int_0^{t\wedge\tau_N}\!\!
\sqrt{c\Xc_s(i)(1-\Xc_s(i))}\,e^{-rs}\di B_s(i)\Big|^2\Big]\\[5pt]
\dis\qquad=c\sum_{i\in\De}E\Big[\int_0^{t\wedge\tau_N}\!\!
\Xc_s(i)(1-\Xc_s(i))e^{-2rs}\di s\Big]
\leq cE\Big[\int_0^{t\wedge\tau_N}\!\!|\Xc_s|\di s\Big]\leq ctN,
\ec
which shows that the stochastic integrals in (\ref{ttN}) are absolutely
summable in $L^2$-norm. It follows from (\ref{ttN}) that
\be\label{stoex}
E^x[|\Xc_{t\wedge\tau_N}|]e^{-rt}
\leq E^x[|\Xc_{t\wedge\tau_N}|e^{-r(t\wedge\tau_N)}]
=E^x[|\Yi_{t\wedge\tau_N}|]\leq|x|.
\ee
Now $NP^x[\tau_N\leq t]\leq|x|e^{rt}$ for all $t\geq 0$, which
shows that $\tau_N\up\infty$ as $N\up\infty$ a.s. Letting $N\up\infty$
in (\ref{stoex}) we arrive at (\ref{linex}).
\qed

\noi
We conclude this section with two results on the continuity of $\Xc$ in its
initial state.
\bl{\bf (Convergence in law)}\label{cola}
Assume that $\Xc^{(n)},\Xc$ are $(a,b,c,d)$-resem-processes, started in
$x^{(n)},x\in[0,1]^\La$, respectively. Then $x^{(n)}\to x$ implies
$\Xc^{(n)}\Rightarrow\Xc$.
\el
{\bf Proof} By Lemma~\ref{exislem}, the laws $\Li(\Xc^{(n)})$ are
tight and each cluster point of the $\Li(\Xc^{(n)})$ solves the
martingale problem for $\Gi$ with initial state $x$. Therefore, by
uniqueness of solutions to the martingale problem,
$\Xc^{(n)}\Rightarrow\Xc$.\qed
\bl{\bf(Monotone convergence)}\label{Xcmon}
Let $\Xc^{(n)},\Xc$ be $(a,b,c,d)$-resem-processes started in
$x^{(n)},x\in[0,1]^\La$, respectively, such that
\be
x^{(n)}\up x\qquad\mbox{as}\quad n\up\infty.
\ee
Then $\Xc^{(n)},\Xc$ may be defined on the same probablity space such that
\be
\Xc^{(n)}_t(i)\up\Xc_t(i)\qquad\forall i\in\La,\ t\geq 0\quad\mbox{as}
\quad n\up\infty\quad{\rm a.s.}
\ee
\el
{\bf Proof} Let $\Xc^{(n)},\Xc$ be solutions of the SDE (\ref{sde})
relative to the same system of Brownian motions. By
Corollary~\ref{sdecomp}, $\Xc^{(n)}\leq\Xc^{(n+1)}$ and
$\Xc^{(n)}\leq\Xc$ for all $n$. Write $\De^{(n)}_t:=\Xc_t-\Xc^{(n)}_t$
and set $\tau^{(n)}_\eps:=\inf\{t\geq 0:\De^{(n)}_t\geq\eps\}$.
A calculation as in the proof of Lemma~\ref{sdemon} shows that
\be
\di\|\De^{(n)}_t\|_\ga\leq(K+b)\|\De^{(n)}_t\|\di t\quad
+\quad\mbox{martingale terms.}
\ee
It follows that
\be
E\big[\|\De^{(n)}_{t\wedge\tau^{(n)}_\eps}\|_\ga\big]
\leq\|x-x^{(n)}\|_\ga e^{(K+b)t}.
\ee
Now $\eps P[\tau^{(n)}_\eps\leq t]\leq\|x-x^{(n)}\|_\ga e^{(K+b)t}$
from which we conclude that $\tau^{(n)}_\eps\up\infty$ as $n\up\infty$
for every $\eps>0$.\qed

\section{Dualities}\label{dusec}

\subsection{Duality and self-duality}\label{dualproof}
 
{\bf Proof of Theorem~\ref{duth}~(a)} We first prove the statement for
finite $x$. We apply Theorem~\ref{errorth}. Our duality function is
\be
\Psi(x,\phi):=(1-\phi)^x\qquad(x\in\Ni(\La),\ \phi\in[0,1]^\La).
\ee
We need to check that the right-hand side in (\ref{apdual}) is zero,
i.e., that
\be
G\Psi(\cdot,\phi)(x)=\Gi^\dgg\Psi(x,\cdot)(\phi)
\qquad(\phi\in[0,1]^\La,\ x\in\Ni(\La)),
\ee
where $G$ be the generator of the $(a,b,c,d)$-braco-process, defined
in (\ref{GdefC3}), and $\Gi^\dgg$ is the generator of the
$(a^\dgg,b,c,d)$-resem-process, defined in (\ref{GidefC3}). Note that
since $x$ is finite, $\Psi(x,\cdot)\in\Ci^2_{\rm
fin}([0,1]^\La)$. We check that
\be\ba{l}
\dis G\Psi(\cdot,\phi)(x)\\[5pt]
\dis=\sum_{ij}a(i,j)x(i)\{(1-\phi(j))-(1-\phi(i))\}(1-\phi)^{x-\de_i}
+b\sum_ix(i)\{(1-\phi(i))-1\}(1-\phi)^x\\[5pt]
\dis\hspace{5pt}+c\sum_ix(i)(x(i)-1)\{1-(1-\phi(i))\}
(1-\phi)^{x-\de_i}+d\sum_ix(i)\{1-(1-\phi(i)\}(1-\phi)^{x-\de_i}\\[5pt]
\dis=-\sum_{ij}a^\dgg(j,i)(\phi(j)-\phi(i))x(i)(1-\phi)^{x-\de_i}
-b\sum_i\phi(i)(1-\phi(i))x(i)(1-\phi)^{x-\de_i}\\[5pt]
\dis\hspace{5pt}+c\sum_i\phi(i)(1-\phi(i))x(i)(x(i)-1)
(1-\phi)^{x-2\de_i}+d\sum_i\phi(i)x(i)(1-\phi)^{x-\de_i}\\[5pt]
\quad=\Gi^\dgg\Psi(x,\cdot)(\phi)
\qquad(\phi\in[0,1]^\La,\ x\in\Ni(\La)).
\ec
Set
\be
\Phi(x,\phi):=G\Psi(\cdot,\phi)(x)=\Gi^\dgg\Psi(x,\cdot)(\phi)
\qquad(\phi\in[0,1]^\La,\ x\in\Ni(\La)).
\ee
It is not hard to see that there exists a constant $K$ such that
\be
|\Phi(x,\phi)|\leq K\Big(1+|x|^2\Big)
\qquad(\phi\in[0,1]^\La,\ x\in\Ni(\La)).
\ee
\detail{Here, one needs to use that $\sup_i\sum_ja(i,j)<\infty$.}
Therefore, condition (\ref{Phint2}) is satisfied by (\ref{kmom}). 

To generalize the statement from finite $x$ to general
$x\in\Ei_\ga(\La)$, we apply Lemma~\ref{mocolem}. Choose finite
$x^{(n)}$ such that $x^{(n)}\up x$ and couple the
$(a,b,c,d)$-braco-processes $X^{(n)},X$ with initial conditions
$x^{(n)},x$, respectively, such that $X^{(n)}\up X$. Then, for each
$t\geq 0$ and $\phi\in[0,1]^\La$,
\be
E^\phi[(1-\Xc_t)^{x^{(n)}}]\down E^\phi[(1-\Xc_t)^{x}]
\quad\mbox{as }n\up\infty,
\ee
and
\be
E[(1-\phi)^{X^{(n)}_t}]\down E[(1-\phi)^{X_t}]\quad\mbox{as }n\up\infty,
\ee
where we used the continuity of the function $x\mapsto(1-\phi)^x$
with respect to increasing sequences.\qed

\noi
{\bf Proof of Theorem~\ref{duth}~(b)} We first prove the statement under
the additional assumption that $\phi$ and $\psi$ are summable. Recall
that by Lemma~\ref{summab}, if $\Xc_0$ is summable then $\Xc_t$ is
summable for all $t\geq 0$ a.s. Let
$S:=\{\phi\in[0,1]^\La:|\phi|<\infty\}$ denote the space of
summable states. We apply Theorem~\ref{errorth}. Our duality function
is
\be
\Psi(\phi,\psi):=\ex{-\frac{b}{c}\li\phi,\psi\re}\qquad(\phi,\psi\in S).
\ee
Let $\Gi,\Gi^\dgg$ denote the generators of the
$(a,b,c,d)$-resem-process and the $(a^\dgg,b,c,d)$-resem-process, as
in (\ref{GidefC3}), respectively. We need to show that the right-hand
side in (\ref{apdual}) is zero, i.e., that
$\Gi\Psi(\cdot,\psi)(\phi)=\Gi^\dgg\Psi(\phi,\cdot)(\psi)$. It is not
hard to see that $\Psi(\cdot,\psi),\Psi(\phi,\cdot)\in\Ci_{\rm
sum}([0,1]^\La)$ for each $\psi,\phi\in S$. We calculate
\detail{Indeed, set
\be
\ga_{ij}(\phi):=\difif{\phi_i}{\phi(j)}\ex{-\frac{b}{c}\li\phi,\psi\re}
=(\ffrac{b}{c})^2\psi(i)\psi(j)\ex{-\frac{b}{c}\li\phi,\psi\re}
\qquad(\phi\in(0,1)^\La).
\ee
Then
\be
\sum_{ij}|\ga_{ij}(\phi)|
\leq(\ffrac{b}{c})^2\sum_{ij}\psi(i)\psi(j)=|\psi|^2,
\ee
which shows that $\ga\in\ell^1(\La^2)$. Moreover, if
$\phi_n\to\phi$ pointwise, then
\bc
\dis\|\ga(\phi)-\ga(\ti\phi)\|^2_2
&\leq&\dis(\ffrac{b}{c})^2\sum_{ij}\psi(i)\psi(j)
\Big|\ex{-\frac{b}{c}\li\phi_n,\psi\re}
-\ex{-\frac{b}{c}\li\phi,\psi\re}\Big|\\[5pt]
&\leq&\dis(\ffrac{b}{c})^3|\psi|^2\sum_k\psi(k)|\phi_n(k)-\phi(k)|\to 0,
\ec
which shows that $\phi\mapsto\ga(\phi)$ is continuous in the norm on
$\ell^1(\La^2)$. For the first derivatives the calculations are
similar but easier.}
\be\ba{r@{\,}c@{\,}l}\label{selfcheck}
\dis\Gi\Psi(\cdot,\psi)(\phi)
&=&\Big\{\dis\sum_{ij}a(j,i)(\phi(j)-\phi(i))(-\ffrac{b}{c})\psi(i)
+b\sum_i\phi(i)(1-\phi(i))(-\ffrac{b}{c})\psi(i)\\[5pt]
&&\dis\phantom{\Big\{}+c\sum_i\phi(i)(1-\phi(i))
(-\ffrac{b}{c})^2\psi(i)^2-d\sum_i\phi(i)(-\ffrac{b}{c})\psi(i)\Big\}
\ex{-\ffrac{b}{c}\li\phi,\psi\re}\\[5pt]
&=&\dis-\ffrac{b}{c}\Big\{\sum_{ij}a(j,i)\phi(j)\psi(i)
-\Big(\sum_ja(j,i)\Big)\sum_i\phi(i)\psi(i)\\[5pt]
&&\dis\phantom{-\ffrac{b}{c}\Big\{}
+b\sum_i\phi(i)(1-\phi(i))\psi(i)(1-\psi(i))
-d\sum_i\phi(i)\psi(i)\Big\}\ex{-\ffrac{b}{c}\li\phi,\psi\re}\\[5pt]
&=&\dis\Gi^\dgg\Psi(\phi,\cdot)(\psi).
\ec
It is not hard to see that there exists a constant $K$ such that
\be
|\Gi\Psi(\cdot,\psi)(\phi)|\leq K|\phi|\,|\psi|\qquad(\phi,\psi\in S).
\ee
\detail{In fact, one even has
\be
|\Gi\Psi(\cdot,\psi)(\phi)|\leq K|\psi|\qquad(\phi,\psi\in S).
\ee
(and analogously $\leq K|\phi|$). This follows from the facts that
\be
\sum_i\phi(i)(1-\phi(i))\psi(i)(1-\psi(i))
\leq\sum_i\phi(i)\psi(i)\leq|\psi|,
\ee
and
\be
\sum_{ij}a(i,j)\phi(i)\psi(j)\leq\sum_j\big(\sum_ia(i,j)\big)\psi(j)
\leq|a|\,|\psi|,
\ee
where $|a|:=\sup_j\sum_ia(i,j)<\infty$.}
Therefore, condition (\ref{Phint2}) is implied by Lemma~\ref{summab},
and Theorem~\ref{errorth} is applicable. To generalize the result to
general $\phi,\psi\in[0,1]^\La$, we apply Lemma~\ref{Xcmon}.\qed

\subsection{Subduality}

Fix constants $\bet\in\R$, $\ga\geq 0$. Let
$\Mi(\La):=\{\phi\in\half^\La:|\phi|<\infty\}$ be the space of finite
measures on $\La$, equipped with the topology of weak convergence, and
let $\Yi$ be the Markov process in $\Mi(\La)$ given by the unique
pathwise solutions to the SDE
\be\label{supersde}
\di\Yi_t(i)=\sum_ja(j,i)(\Yi_t(j)-\Yi_t(i))\,\di t
+\bet\Yi_t(i)\,\di t+\sqrt{2\ga\Yi_t(i)}\,\di B_t(i)
\ee
$(t\geq 0,\ i\in\La)$. Then $\Yi$ is the well-known super random walk
with underlying motion $a$, growth parameter
$\bet$ and activity $\ga$. One has \cite[Section~4.2]{Daw93}
\be\label{laplaceC3}
E^\phi\big[\ex{-\li\Yi_t,\psi\re}]=\ex{-\li\phi,\Ui_t\psi\re}
\ee
for any $\phi\in\Mi(\La)$ and bounded nonnegative $\psi:\La\to\R$, where
$u_t=\Ui_t\psi$ solves the semilinear Cauchy problem
\be\label{udif}
\dif{t}u_t(i)=\sum_ja(j,i)(u_t(j)-u_t(i))+\bet u_t(i)
-\ga u_t(i)^2\qquad(i\in\La,\ t\geq 0)
\ee
with initial condition $u_0=\psi$. The semigroup $(\Ui_t)_{t\geq 0}$
acting on bounded nonnegative functions $\psi$ on $\La$ is called the
log-Laplace semigroup of $\Yi$.

We will show that $(a,b,c,d)$-braco-process and the super random walk
with underlying motion $a^\dgg$, growth parameter
$b-d+c$ and activity $c$ are related by a duality
formula with a nonnegative error term. In analogy with words such as
subharmonic and submartingale, we call this a subduality relation.
\bp{\bf(Subduality with a branching process)}\label{subdup}
Let $X$ be the $(a,b,c,d)$-braco-process and let $\Yi$ be the super
random walk with underlying motion $a^\dgg$, growth parameter $b-d+c$
and activity $c$. Then
\be\label{subdu}
E^x\big[\ex{-\li\phi,X_t\re}]\geq E^\phi\big[\ex{-\li\Yi_t,x\re}]
\qquad(x\in\Ei_\ga(\La),\ \phi\in\Mi(\La)).
\ee
\ep
{\bf Proof} We first prove the statement
for finite $x$. We apply Theorem~\ref{errorth} to $X$ and $\Yi$
considered as processes in $\Ni(\La)$ and $\Mi(\La)$,
respectively. The process $\Yi$ solves the martingale problem for the
operator
\bc\label{Hidef}
\Hi f(\phi)&:=&\dis\sum_{ij}a^\dgg(j,i)(\phi(j)-\phi(i))\dif{\phi(i)}f(\phi)
+(b-d+c)\sum_i\phi(i)\dif{\phi(i)}f(\phi)\\[2pt]
&&\dis+c\sum_i\phi(i)\diff{\phi(i)}f(\phi)
\qquad(\phi\in[0,1]^\La),
\ec
defined for functions $\phi$ in the space $\Ci^2_{\rm
fin,b}\half^\La$ of bounded $\Ci^2$ functions on $\half^\La$
depending on finitely many coordinates. Our duality function is
$\Psi(x,\phi):=e^{-\li\phi,x\re}$. We observe that
$\Psi(x,\cdot)\in\Ci^2_{\rm fin,b}\half^\La$ for all $x\in\Ni(\La)$
and calculate
\bc
G\Psi(\cdot,\phi)(x)
&=&\dis\Big\{\sum_{ij}a(i,j)x(i)\big(e^{\phi(i)-\phi(j)}-1\big)
+b\sum_ix(i)\big(e^{-\phi(i)}-1\big)\\[5pt]
&&\dis\phantom{\Big\{}+c\sum_ix(i)(x(i)-1)\big(e^{\phi(i)}-1\big)
+d\sum_ix(i)\big(e^{\phi(i)}-1\big)\Big\}\ex{-\li\phi,x\re},
\ec
and
\bc
\Hi\Psi(x,\cdot)(\phi)
&=&\dis\Big\{\sum_{ij}a^\dgg(j,i)x(i)(\phi(i)-\phi(j))
-(b-d+c)x(i)\phi(i)\\[5pt]
&&\dis\phantom{\Big\{}+c\sum_ix(i)^2\phi(i)\Big\}
\ex{-\li\phi,x\re}
\ec
$(x\in\Ni(\La),\ \phi\in\Mi(\La))$. It is not hard to see that there
exists a constant $K$ such that
\be
|G\Psi(\cdot,\phi)(x)|\leq K|x|^2\quad\mbox{and}\quad
|\Hi\Psi(x,\cdot)(\phi)|\leq K|x|^2\,|\phi|
\qquad(x\in\Ni(\La),\ \phi\in\Mi(\La)).
\ee
\detail{For the estimate involving $G$, one should use relations such as
\be
\big|x(i)\big(e^{\phi(i)}-1\big)\ex{-\li\phi,x\re}\big|
=\big|x(i)\big(1-e^{-\phi(i)}\big)\ex{-\li x-\de_i,\phi\re}\big|\leq x(i).
\ee}
and therefore condition~(\ref{Phint2}) is implied by (\ref{kmom})
and the elementary estimate $E[|\Yi_t|]\leq e^{(b-d+c)t}|\phi|$.
One has
\be\ba{l}
\dis G\Psi(\cdot,\phi)(x)-\Hi\Psi(x,\cdot)(\phi)=\Big\{\sum_{ij}a(i,j)x(i)
\big(e^{\phi(i)-\phi(j)}-1-(\phi(i)-\phi(j))\big)\\[5pt]
\qquad\dis+b\sum_ix(i)\big(e^{-\phi(i)}-1+\phi(i)\big)
+c\sum_ix(i)(x(i)-1)\big(e^{\phi(i)}-1-\phi(i)\big)\\[5pt]
\qquad\dis+d\sum_ix(i)\big(e^{\phi(i)}-1-\phi(i)\big)\Big\}
\ex{-\li\phi,x\re}\geq 0,
\ec
and therefore, for finite $x$, (\ref{subdu}) is implied by
Theorem~\ref{errorth}. The general case follows by approximation,
using Lemma~\ref{mocolem}.\qed

\section{The maximal processes}\label{maxsec}

\subsection{The maximal branching-coalescing process}

Using Proposition~\ref{subdup} we can now prove
Theorem~\ref{maxbraco}.\medskip

\noi
{\bf Proof of Theorem~\ref{maxbraco}} Choose $x^{(n)}\in\Ei_\ga(\La)$
such that $x^{(n)}(i)\up\infty$ for all $i\in\La$. By
Lemma~\ref{mocolem}, the $(a,b,c,d)$-braco processes $X^{(n)}$ started
in $x^{(n)}$, respectively, can be coupled such that $X^{(n)}_t\leq
X^{(n+1)}_t$ for each $t\geq 0$. Define
$X^{(\infty)}=(X^{(\infty)}_t)_{t\geq 0}$ as the $\ov\N^\La$-valued
process that is the pointwise increasing limit of the $X^{(n)}$.
By Proposition~\ref{subdup} and (\ref{laplaceC3}),
\be
E\big[1-\ex{-\li\eps\de_i,X^{(n)}_t\re}\big]
\leq 1-\ex{-\li\eps\de_i,\Ui_tx^{(n)}\re}\qquad(t,\eps\geq 0,\ i\in\La).
\ee
where $(\Ui_t)_{t\geq 0}$ is the log-Laplace semigroup of the super
random walk with underlying motion $a^\dgg$, growth parameter
$r:=b-d+c$ and activity $c$. It follows that
\detail{Here we use that $e^{-\eps z}\up 1$ and therefore
$\eps^{-1}(1-e^{-\eps z})=\int_0^ze^{-\eps y}\di y\up z$ as $\eps\down 0$.}
\be\label{nbou}
E[X^{(n)}_t(i)]
=\lim_{\eps\down 0}\eps^{-1}E\big[1-\ex{-\li\eps\de_i,X^{(n)}_t\re}\big]
\leq\lim_{\eps\down 0}\eps^{-1}\big(1-\ex{-\li\eps\de_i,\Ui_tx^{(n)}\re}\big)
=\Ui_tx^{(n)}(i)
\ee
$(t\geq 0,\ i\in\La)$. Using the explicit solution of (\ref{udif}) for
constant initial conditions, it is easy to see that
$\Ui_tx^{(n)}\up\Ui_t\infty$, where
\be
\Ui_t\infty:=\left\{\ba{cl}\frac{r}{c(1-e^{-rt})}\quad&
\mbox{if }r\neq 0,\\[5pt] \frac{1}{ct}\quad&\mbox{if }r=0.\ea\right.
\ee
(See formula~(\ref{ovU}).) Letting $n\up\infty$
in (\ref{nbou}) we arrive at Theorem~\ref{maxbraco}~(b). Moreover, we
see that
\be
E\big[\|X^{(\infty)}_t(i)\|_\ga\big]\leq \Ui_t\infty\sum_i\ga_i<\infty
\qquad(t>0),
\ee
and therefore $X^{(\infty)}_t\in\Ei_\ga(\La)$ a.s.\ for each
$t>0$. Part~(a) of the theorem now follows from Lemma~\ref{mocolem}.
Using Theorem~\ref{duth}~(a) and the continuity of the function
$x\mapsto(1-\phi)^x$ with respect to increasing sequences, reasoning
as in (\ref{fikst}), we see that
\be\label{maxthin}
P[\Thin_\phi({X^{(\infty)}_t})=0]
=P^\phi[\Xc^\dgg_t=0]\qquad(\phi\in[0,1]^\La,\ t\geq 0),
\ee
where $\Xc^\dgg$ denotes the $(a^\dgg,b,c,d)$-resem-process. Since
formula (\ref{maxthin}) determines the distribution of
$X^{(\infty)}_t$ uniquely, the law of $X^{(\infty)}_t$ does not depend
on the choice of the $x^{(n)}\up\infty$ $(t\geq 0)$. This completes the
proof of part (c) of the theorem.

\detail{Since there is only one $(a,b,c,d)$-braco-process with a
given $\Ei_\ga(\La)$-valued initial condition, also the law of
$(X^{(\infty)}_t)_{t\geq 0}$ does not depend on the choice of
$x^{(n)}\up\infty$.}

To prove part (d), fix $0\leq s\leq t$. Choose
$y_n\in\Ei_\ga(\La)$, $y_n(i)\up\infty$ $\forall i\in\La$
and let $\ti X^{(n)}$ be the $(a,b,c,d)$-braco-process
started in $\ti X^{(n)}_0:=X^{(\infty)}_{t-s}\vee y_n$. Then
$\ti X^{(n)}_0\geq X^{(\infty)}_{t-s}$ and therefore, by
Lemma~\ref{complem}, $\ti X^{(n)}_s$ and $X^{(\infty)}_t$ may be
coupled such that $\ti X^{(n)}_s\geq X^{(\infty)}_t$. By part (c) of
the theorem, $\ti X^{(n)}_s$ and $X^{(\infty)}_s$ may be coupled
such that $\ti X^{(n)}_s\up X^{(\infty)}_s$ and therefore
$X^{(\infty)}_s$ and $X^{(\infty)}_t$ may be coupled such that
$X^{(\infty)}_s\geq X^{(\infty)}_t$.

It follows that $\Li(X^{(\infty)}_t)\down\ov\nu$ for some probability
measure $\ov\nu$ on $\Ei_\ga(\La)$. Set $\rho:=\Li(X^{(\infty)}_1)$
and let $(S_t)_{t\geq 0}$ denote the semigroup of the
$(a,b,c,d)$-braco-process. Recall the definition of 
$\Ci_{\rm Lip,b}(\Ei_\ga(\La))$ above (\ref{Lico}). One has
\be
\int\ov\nu(\di x)f(x)=\lim_{t\to\infty}\int\rho(\di x)S_tf(x)
\ee
for every $f\in\Ci_{\rm Lip,b}(\Ei_\ga(\La))$. Therefore, since $S_t$
maps $\Ci_{\rm Lip,\,b}(\Ei_\ga(\La))$ into itself,
\be
\int\ov\nu(\di x)S_sf(x)=\lim_{t\to\infty}\int\rho(\di x)S_tS_sf(x)
=\int\ov\nu(\di x)f(x)\qquad(s\geq 0),
\ee
for every $f\in\Ci_{\rm Lip,\,b}(\Ei_\ga(\La))$, which shows that
$\ov\nu$ is an invariant measure. If $\nu$ is another invariant
measure, then $\Li(X^{(\infty)}_t)\geq\nu$ for all $t\geq 0$. Letting
$t\to\infty$, we see that $\ov\nu\geq\nu$, proving part~(e) of the
theorem. Part~(f) has already been proved in the introduction.\qed

\subsection{The maximal resampling-selection process}

The proof of Theorem~\ref{maxresem}~(a)--(c) is similar to the proof
of Theorem~\ref{maxbraco}, but easier. Recall that
Theorem~\ref{maxresem}~(d) is proved in Section~\ref{methods}.\medskip

\noi
{\bf Proof of Theorem~\ref{maxresem}~(a)--(c)} Part~(a) can be proved
in the same way as Theorem~\ref{maxbraco}~(d), using
Lemma~\ref{Xcmon}. The proof of part~(b) goes analogue to the proof of
Theorem~\ref{maxbraco}~(e). To see why (\ref{sdemom}) holds, note that
for any $\phi\in[0,1]^\La$, by Theorem~\ref{duth}~(a),
\be
\int\ov\mu(\di\phi)(1-\phi)^x
=\lim_{t\to\infty}P^1[\Thin_{\Xc_t}(x)=0]
=\lim_{t\to\infty}P^x[\Thin_1(X^\dgg_t)=0].
\ee
To complete the proof of part~(c) we must show that $\ov\mu$ is
nontrivial if and only if the $(a^\dgg,b,c,d)$-process survives. Using
subadditivity (Lemma~\ref{subad}) it is easy to see that the
$(a^\dgg,b,c,d)$-process survives if and only if
$P^{\de_i}[X^\dgg_t\neq 0\ \forall t\geq 0]>0$ for some
$i\in\La$. Formula (\ref{sdemom}) implies that
$\int\ov\mu(\di\phi)\phi(i)=P^{\de_i}[X^\dgg_t\neq 0\ \forall t\geq
0]$, which shows that $\ov\mu=\de_0$ if and only if the
$(a^\dgg,b,c,d)$-process survives. If $\ov\mu\neq\de_0$ then the
measure $\ov\mu$ conditioned on $\{\phi:\phi\neq 0\}$ is an invariant
measure of the $(a,b,c,d)$-resem-process that is stochastically larger
than $\ov\mu$. By part~(b), this conditioned measure is $\ov\mu$
itself, thus $\ov\mu(\{0\})=0$, i.e., $\ov\mu$ is nontrivial.\qed

\detail{If you start the process in the conditioned measure then you
will of course never get extinct, since in that case $P[\Xc_t=0]$
would depend on $t$, contradicting stationarity. The conditioned
measure gives a higher expectation to every increasing function $f$,
since $f(\phi)\geq f(0)$ for every $\phi\in[0,1]^\La$.}

\section{Convergence to the upper invariant measure}\label{consec}

\subsection{Extinction versus unbounded growth}\label{growsec}

In this section we prove Lemma~\ref{exgroC3}. It has already been proved
in Section~\ref{methods} that $e^{-\frac{b}{c}|\Xc_t|}$ is a
submartingale. Therefore, if $b>0$, then $|\Xc_t|$ converges a.s.\ to
a limit in $[0,\infty]$. If $b=0$ then it is easy to see that
$|\Xc_t|$ is a nonnegative supermartingale and therefore also in this
case $|\Xc_t|$ converges a.s. Thus, all we have to do is to show that
$\lim_{t\to\infty}|\Xc_t|$ takes values in $\{0,\infty\}$ a.s.\
(Proposition~\ref{zeroinf} below), and that $\Xc$ gets extinct in
finite time if the limit is zero (Lemma~\ref{finex}). Throughout this
section, $c>0$ and $\Xc$ is the $(a,b,c,d)$-resem-process starting in
an initial state $\phi\in[0,1]^\La$ with $|\phi|<\infty$.
\bl{\bf(Finite time extinction)}\label{finex}
One has $\Xc_t=0$ for some $t\geq 0$ a.s.\ on the event
$\lim_{t\to\infty}|\Xc_t|=0$.
\el
{\bf Proof} Choose $x^{(n)}\in\Ei_\ga(\La)$ such that
$x^{(n)}(i)\up\infty$ for all $i\in\La$. Let $X^{(n)\dgg}$ denote the
$(a^\dgg,b,c,d)$-braco-process started in $x^{(n)}$ and let
$X^{(\infty)\dgg}$ denote the maximal
$(a^\dgg,b,c,d)$-braco-process. By Theorem~\ref{duth}~(a) and
Theorem~\ref{maxbraco}~(b), \be\ba{l}\label{linbo}
\dis P^\phi[\Xc_t\neq 0]
=\lim_{n\up\infty}P^\phi[\Thin_{\Xc_t}(x^{(n)})\neq 0]
=\lim_{n\up\infty}P[\Thin_{\phi}(X^{(n)\dgg}_t)\neq 0]\\
\dis\quad=P[\Thin_{\phi}(X^{(\infty)\dgg}_t)\neq 0]
\leq E\big[|\Thin_{\phi}(X^{(\infty)\dgg}_t)|\big]
=\li\phi,E[X^{(\infty)\dgg}_t]\re\leq|\phi|\Ui_t\infty,
\ec
where $\Ui_t\infty$ is the function on the right-hand side in
(\ref{explicit}). Choose $\eps>0$ and $t_0>0$ such that
$\eps\Ui_{t_0}\infty\leq\frac{1}{2}$. Let $(\Fi_t)_{t\geq 0}$ denote
the filtration generated by $\Xc$. By (\ref{linbo}),
\be\label{haeps}
\ffrac{1}{2}1_{\txt\{|\Xc_t|\leq\eps\}}\leq P[\Xc_{t+t_0}=0|\Fi_t]
\leq P[\exists s\geq 0\ \mbox{s.t. }\Xc_s=0|\Fi_t].
\ee
Now
\be\label{onext}
1_{\txt\{\lim_{s\to\infty}\Xc_s=0\}}
\leq\liminf_{t\to\infty}1_{\txt\{|\Xc_t|\leq\eps\}},
\ee
while
\be\label{right}
P[\exists s\geq 0\ \mbox{s.t. }\Xc_s=0|\Fi_t]
\to 1_{\txt\{\exists s\geq 0\ \mbox{s.t. }\Xc_s=0\}}
\quad\mbox{as }t\to\infty\quad{\rm a.s.},
\ee
by convergence of right-continuous martingales and the fact that the
left-hand side is right-continuous by a general property of strong
Markov processes described in Section~\ref{S:zeone} from
Chapter~\ref{C:ren}. Letting $t\to\infty$ in (\ref{haeps}), using
(\ref{onext}) and (\ref{right}), we find that
$\ffrac{1}{2}1_{\{\lim_{s\to\infty}\Xc_s=0\}}\leq 1_{\{\exists s\geq
0\ \mbox{s.t. }\Xc_s=0\}}$ a.s.\qed

\noi
To finish this section, we need to prove:
\bp{\bf(Convergence to zero or infinity)}\label{zeroinf}
Assume that $\La$ is infinite. Then
$\lim_{t\to\infty}|\Xc_t|\in\{0,\infty\}$ a.s.
\ep
Since the proof of Proposition~\ref{zeroinf} is rather long we break
it up into a number of steps. At each step, we will skip the proof if
it is obvious but tedious. Our first step is:
\bl{\bf(Integrable fluctuations)}\label{intfluc}
One has
\be\label{flucint}
\int_0^\infty\sum_i\Xc_t(i)(1-\Xc_t(i))\,\di t<\infty
\ee
a.s.\ on the event $\lim_{t\to\infty}|\Xc_t|\in\half$.
\el
{\bf Proof} For any $\psi\in\half^\La$ with $|\psi|<\infty$ one has
$e^{-\li\cdot,\psi\re}\in\Ci^2_{\rm sum}([0,1]^\La)$ and (compare
(\ref{selfcheck}))
\be\ba{r@{\,}c@{\,}l}
\dis\Gi\ex{-\li\cdot,\psi\re}(\phi)
&=&\dis\Big\{-\sum_i\phi(i)\sum_ja^\dgg(j,i)(\psi(j)-\psi(i))\\[5pt]
&&\dis+\sum_i\phi(i)(1-\phi(i))\big(c\psi(i)^2-b\psi(i)\big)
+d\sum_i\phi(i)\psi(i)\Big\}\ex{-\li\phi,\psi\re}.
\ec
Since $\Xc$ solves the martingale problem for $\Gi$,
\be\label{mefo}
E\Big[\int_0^t\Gi\ex{-\li\cdot,\psi\re}(\Xc_s)\di s\Big]
=E\big[\ex{-\li\Xc_t,\psi\re}\big]-\ex{-\li\phi,\psi\re}\qquad(t\geq 0).
\ee
Choose $\la>0$ such that $c\la^2-b\la=:\mu>0$ and $\psi_n\in\half^\La$
with $|\psi_n|<\infty$ such that $\psi_n\up\la$. Then the bounded
pointwise limit of the function
$i\mapsto\sum_ja^\dgg(j,i)(\psi_n(j)-\psi_n(i))$ is zero and therefore,
taking the limit in (\ref{mefo}), using Lemma~\ref{summab}, we find
that
\detail{Set $\eta_n(i):=\sum_ja^\dgg(j,i)(\psi_n(j)-\psi_n(i))$. Then
\be
E\Big[\int_0^t\sum_i\Xc_s(i)\eta_s(i)\di s\Big]
=\sum_iE\Big[\int_0^t\Xc_s(i)\di s\Big]\eta_n(i),
\ee
where by  Lemma~\ref{summab}, $i\mapsto E[\int_0^t\Xc_s(i)\di s]$
is a finite measure on $\La$.}
\be\label{onebo}
E\Big[\int_0^t\sum_i\Big\{\mu\Xc_s(i)(1-\Xc_s(i))+\la d\Xc_s(i)\Big\}
\ex{-\la|\Xc_s|}\di s\Big]=E\big[\ex{-\la|\Xc_t|}\big]-\ex{-\li\phi,\psi\re}.
\ee
Letting $t\up\infty$, using the fact that the right-hand side of
(\ref{onebo}) is bounded by one, we see that
\be
\int_0^\infty\sum_i\Big\{\mu\Xc_t(i)(1-\Xc_t(i))+\la d\Xc_t(i)\Big\}
\ex{-\la|\Xc_t|}\,\di t<\infty\quad{\rm a.s.},
\ee
which implies (\ref{flucint}).\qed

\bl{\bf(Process not started with only zeros and ones)}\label{midstart}
For every $0<\eps<\ffrac{1}{4}$ there exists a $\de,r>0$ such that
\be\label{stay}
P^\phi\big[\Xc_t(i)\in(\eps,1-\eps)\ \forall t\in[0,r]\big]\geq\de
\qquad(i\in\La,\ \phi\in[0,1]^\La,\ \phi(i)\in(2\eps,1-2\eps)).
\ee
\el
{\bf Proof} Since $\sup_i\sum_ja(i,j)<\infty$ and all the components
of the $(a,b,c,d)$-resem-process take values in $[0,1]$, the maximal
drift that the $i$-th component $\Xc_t(i)$ can experience (both in the
positive and negative direction) can be uniformly bounded. Now the
proof of (\ref{stay}) is just a standard calculation, which we skip.\qed

\detail{Set $\tau:=\inf\{t\geq 0:\Xc_t(i)\not\in(\eps,1-\eps)\}$.
Define $f\in\Ci^2_{\rm fin}([0,1]^\La)$ by $f(\phi):=\phi(i)(1-\phi(i))$.
Since $\Xc$ solves the martingale problem for $\Gi$,
\be
E[f(\Xc_{t\wedge\tau})]=f(\phi)-\int_0^{t\wedge\tau}\!\!\Gi f(\Xc_s)\di s
\ee
and therefore, since $\Gi f\in\Ci([0,1]^\La)$,
\be
E[\Xc_{t\wedge\tau}(i)(1-\Xc_{t\wedge\tau}(i))]
\geq 2\eps(1-2\eps)-\|\Gi f\|_\infty t.
\ee
Therefore, for any $\de<1$ there exists a $r$ small enough, such that
\be
P[\tau>r]=P[\Xc_{r\wedge\tau}(i)\in(\eps,1-\eps)]\geq\de.
\ee}

\bl{\bf(Uniform convergence to zero or one)}\label{unione}
Almost surely on the event that $\lim_{t\to\infty}|\Xc_t|\in\half$,
there exists a set $\De\sub\La$ such that
\be
\lim_{t\to\infty}\inf_{i\in\De}\Xc_t(i)=1\quad\mbox{and}
\quad\lim_{t\to\infty}\sup_{i\in\La\beh\De}\Xc_t(i)=0.
\ee
\el
{\bf Proof} Imagine that the statement does not hold. Then, by the
continuity of sample paths, with positive probability
$\lim_{t\to\infty}|\Xc_t|\in\half$ while there exists
$0<\eps<\ffrac{1}{4}$ such that for every $T>0$ there exists $t\geq T$
and $i\in\La$ with $\Xc_t(i)\in(2\eps,1-2\eps)$. Using
Lemma~\ref{midstart} and the strong Markov property, it is then not
hard to check that with positive probability
$\lim_{t\to\infty}|\Xc_t|\in\half$ while there exist infinitely many
disjoint time intervals $[t_k,t_k+r]$ and points $i_k\in\La$ such
that $\Xc_t(i_k)\in(\eps,1-\eps)$ for all $t\in[t_k,t_k+r]$. This
contradicts Lemma~\ref{intfluc}.\qed

\bl{\bf(Convergence to one on a finite nonempty set)}\label{Denonemp}
Almost surely on the event $\lim_{t\to\infty}|\Xc_t|\in(0,\infty)$,
the set $\De$ from Lemma~\ref{unione} is finite and nonempty.
\el
{\bf Proof} It is clear that $\De$ is finite a.s.\ on the event
$\lim_{t\to\infty}|\Xc_t|<\infty$. Now imagine that $\De$ is
empty. Then, a.s.\ on the event $\lim_{t\to\infty}|\Xc_t|>0$, there
exists a random time $T$ such that $\Xc_t(i)\leq\frac{1}{2}$ for all
$t\geq T$ and $i\in\La$. Since $z(1-z)\geq\frac{1}{2}z$ on
$[0,\frac{1}{2}]$, it follows that a.s.\ on the event
$\lim_{t\to\infty}|\Xc_t|>0$,
\be
\int_T^\infty\sum_i\Xc_t(i)(1-\Xc_t(i))\di t
\geq\frac{1}{2}\int_T^\infty|\Xc_t|=\infty.
\ee
We arrive at a contradiction with Lemma~\ref{intfluc}.\qed

\noi
{\bf Proof of Proposition~\ref{zeroinf}} Let $\De$ be the random set
from Lemma~\ref{unione}. We will show that $\De=\La$ a.s.\ on the
event $\lim_{t\to\infty}|\Xc_t|\in(0,\infty)$.  In particular, by
Lemma~\ref{Denonemp}, if $\La$ is infinite this implies that the event
$\lim_{t\to\infty}|\Xc_t|\in(0,\infty)$ has zero probability.  Assume
that with positive probability $\lim_{t\to\infty}|\Xc_t|\in(0,\infty)$
and $\De\neq\La$. By Lemma~\ref{Denonemp}, $\De$ is nonempty, and
therefore by irreducibility there exist $i\in\La\beh\De$ and $j\in\De$
such that $a(i,j)>0$ or $a(j,i)>0$. If $a(i,j)>0$ then by the fact
that the counting measure is an invariant measure for the Markov
process with jump rates $a$ and by the finiteness of $\De$, there must
also be an $i'\in\La\beh\De$ and $j'\in\De$ such that
$a(j',i')>0$. Thus, there exist $i,j\in\La$ such that $a(j,i)>0$ and
with positive probability $\lim_{t\to\infty}\Xc_t(i)=0$, and
$\lim_{t\to\infty}\Xc_t(j)=1$. It is not hard to see that this
violates the evolution in (\ref{sde}). (We skip the details.)\qed

\detail{Imagine that there exist $i,j\in\La$ such that $a(j,i)>0$ and
with positive probability $\lim_{t\to\infty}\Xc_t(i)=0$, and
$\lim_{t\to\infty}\Xc_t(j)=1$. Then, taking C\'esaro limits using
compactness, we can construct a stationary process $\Xc$ with the
property that $P[\Xc_t(i)=0,\ \Xc_t(j)=1]>0$ for each $t$. To see that
this is impossible, set $\tau:=\inf\{t\geq 0:\Xc_t(i)\neq 0\}$. Then,
using the SDE up to time $\tau$ we find that $\Xc_\tau(i)\geq
1-e^{-a(j,i)\tau}$, which shows that $\tau=0$.}

\subsection{Convergence to the upper invariant measure}\label{omnisec}

In this section we complete the proof of Theorem~\ref{homconv},
started in Section~\ref{methods}, by proving
Lemma~\ref{omnipres}. Throughout this section, $(\La,a)$ is infinite
and homogeneous and $G$ is a transitive subgroup of ${\rm
Aut}(\La,a)$.  We fix a reference point $0\in\La$. We start with two
preparatory lemmas.
\bl{\bf(Sparse thinning functions)}\label{sparse}
Assume that $\phi_n\in[0,1]^\La$, $|\phi_n|\to\infty$. Let
$\De\sub\La$ be finite with $0\in\De$. Then it is possible to choose
constants $\la_n\to\infty$, finitely supported probability
distributions $\pi_n$ on $\La$, and $\{g_i\}_{i\in{\rm supp}(\pi_n)}$
with $g_i\in G$ and $g_i(0)=i$ such that the images
$\{g_i(\De)\}_{i\in{\rm supp}(\pi_n)}$ are disjoint, and such that
$\la_n\pi_n\leq\phi_n$.
\el
{\bf Proof} Choose $(g_i)_{i\in\La}$ with $g_i\in G$ such that
$g_i(0)=i$. Let $(\xi^{\rm s}_t)_{t\geq 0}$ be the random walk on
$\La$ that jumps from $i$ to $j$ with the symmetrized jump rates
$a^{\rm s}(i,j)=a(i,j)+a^\dgg(i,j)$. By irreducibility and symmetry,
$P^i[\xi^{\rm s}_t=j]>0$ for all $t>0,\ i,j\in\La$. Put
\be
\Ga^\eps_i:=\{j\in\La:P^i[\xi^{\rm s}_1=j]\geq\eps\}\qquad(i\in\La).
\ee
We can choose $\eps>0$ small enough such that
\be\label{Gap}
j\not\in\Ga^\eps_i\quad\mbox{implies}\quad g_i(\De)\cap g_j(\De)=\emptyset
\quad(i,j\in\La).
\ee
To see this, set $\de:=\min_{k\in\De}P^0[\xi^{\rm s}_{\frac{1}{2}}=k]$
and put $\eps:=\de^2$. Imagine that $\exists k\in g_i(\De)\cap
g_j(\De)$. Then $P^i[\xi^{\rm s}_1=j]\geq P^i[\xi^{\rm
s}_{\frac{1}{2}}=k]P^k[\xi^{\rm s}_{\frac{1}{2}}=j]\geq\de^2=\eps$ by
the symmetry of the random walk and homogeneity, and therefore
$j\in\Ga^\eps_i$. Now choose inductively $i_1,i_2,\ldots\in\La$ such
that
\be
\phi_n\mbox{ assumes its maximum over }\La\beh\bigcup_{l=1}^k\Ga^\eps_{i_l}
\mbox{ in }i_{k+1}.
\ee
Then $g_{i_1}(\De),g_{i_2}(\De),\ldots$ are disjoint by (\ref{Gap}).
Since $K:=|\Ga^\eps_i|$ is finite and does not depend on $i$,
\be
\sum_{l=1}^\infty\phi_n(i_l)\geq\frac{|\phi_n|}{K},
\ee
and we can choose $k_n$ such that
\be
\la_n:=\sum_{l=1}^{k_n}\phi_n(i_l)\asto{n}\infty.
\ee
Setting
\be
\pi_n:=\frac{1}{\la_n}\phi_n1_{\{i_1,\ldots,i_{k_n}\}}
\ee
yields $\la_n$ and $\pi_n$ with the desired properties.\qed

\noi
Let $(\xi_t)_{t\geq 0}$ and $(\xi^\dgg_t)_{t\geq 0}$ denote the random
walks on $\La$ that jump from $i$ to $j$ with rates $a(i,j)$ and
$a^\dgg(i,j)$, respectively. Then, for any $\De\sub\La$, the sets
\be\label{RDe}
R\De:=\{i\in\La:P^i[\xi_t\in\De]>0\}\quad\mbox{and}
\quad R^\dgg\De:=\{i\in\La:P^i[\xi^\dgg_t\in\De]>0\}\quad(t>0)
\ee
of points from which $\xi$ and $\xi^\dgg$ can enter $\De$ do not
depend on $t>0$. Indeed
\be\label{Rdef}
R\De=\big\{i:\exists n\geq 0,\ i_0,\ldots,i_n\mbox{ s.t.\ }
i_0=i,\ i_n\in\De,\ a(i_{l-1},i_l)>0\ \forall l=1,\ldots,n\big\}
\ee
and similarly for $R^\dgg\De$. In our next lemma, for $x\in\N^\La$ and
$\De\sub\La$ we let $x|_\De:=(x_i)_{i\in\De}$ denote the restriction
of $x$ to $\De$.
\bl{\bf(Points from which 0 can be reached)}\label{reached}
If $\mu$ is a $G$-homogeneous and nontrivial probability measure on
$\N^\La$, then
\be
\mu\big(\{x:x|_{R\{0\}}=0\}\big)=0.
\ee
\el
{\bf Proof} Let $Y$ be a $\N^\La$-valued random variable with law $\mu$.
We will show that for any $\De\sub\La$,
\be\label{RR}
P\big[Y|_{R^\dgg R\De}=0\big]=P\big[Y|_{R\De}=0\big].
\ee
Assume that (\ref{RR}) does not hold. Then there exists an $i\in
R^\dgg R\De\beh R\De$ such that with positive probability $Y(i)\neq 0$
and $Y|_{R\De}=0$. Since the random walk $(\xi^\dgg_t)_{t\geq 0}$
cannot escape from $R\De$ this implies that for any $t>0$
\be
P^i\big[Y(\xi^\dgg_0)\neq 0,\ Y(\xi^\dgg_s)=0\ \forall s\geq t\big]>0,
\ee
which contradicts the fact that $(Y(\xi^\dgg_t))_{t\geq 0}$ is stationary.
This proves (\ref{RR}). Continuing this process, we see that
\be
P\big[Y|_{R\{0\}}=0\big]=P\big[Y|_{R^\dgg R\{0\}}=0\big]
=P\big[Y|_{RR^\dgg R\{0\}}=0\big]=\cdots
\ee
By irreducibility, the sets $R\{0\},R^\dgg R\{0\},RR^\dgg R\{0\},\ldots$
increase to $\La$, and therefore, since $\mu$ is nontrivial,
\be
P\big[Y|_{R\{0\}}=0\big]=P\big[Y|_\La=0\big]=0.
\ee
\qed

\noi
{\bf Proof of Lemma~\ref{omnipres}} For any finite set $\De\sub\La$,
let $X^\De$ denote the $(a,b,c,d)$-braco-process with immediate
killing outside $\De$. Thus, $X^\De_t(i):=0$ for all $i\in\La\beh\De$
and $t>0$ and $(X^\De_t(i))_{i\in\De,\ t\geq 0}$ is the Markov process
in $\N^\De$ with generator $G^\De$ given by (compare (\ref{GdefC3}))
\bc\label{GDelta}
G^\De f(x)&:=&\dis\sum_{i,j\in\De}a(i,j)x(i)\{f(x+\de_j-\de_i)-f(x)\}
\,+\!\!\sum_{i\in\De,j\in\La\beh\De}a(i,j)x(i)\{f(x-\de_i)-f(x)\}\\
&&\dis+b\sum_{i\in\De}x(i)\{f(x+\de_i)-f(x)\}
+c\sum_{i\in\De}x(i)(x(i)-1)\{f(x-\de_i)-f(x)\}\\
&&\dis+d\sum_{i\in\De}x(i)\{f(x-\de_i)-f(x)\}.
\ec
It is not hard to see that if $\De_1,\ldots,\De_n$ are disjoint finite
sets, then it is possible to couple the processes $X$ and
$X^{\De_1},\ldots,X^{\De_n}$ in such a way that
\be
X_t\leq\sum_{i=1}^n X^{\De_i}_t\qquad(t\geq 0)
\ee
and the $(X^{\De_i})_{i=1,\ldots,n}$ are independent.

Let $X$ denote the $(a,b,c,d)$-braco-process and assume that
$\phi_n\in[0,1]^\La$ satisfy $|\phi_n|\to\infty$. Fix $t>0$. Assume
that $\De\sub\La$ is a finite set such that $0\in\De$ and
\be\label{Denul}
x|_{\De}\neq 0\quad\Rightarrow\quad P^x[X^\De_t(0)>0]>0.
\ee
Choose $\la_n$, $\pi_n$, and $\{g_i\}_{i\in{\rm supp}(\pi_n)}$ as in
Lemma~\ref{sparse}. Then, for deterministic $x\in\Ei_\ga(\La)$, we can
estimate
\bc\label{prodest}
\dis P^x\big[\Thin_{\phi_n}(X_t)=0\big]
&\leq&\dis P^x\big[\Thin_{\la_n\pi_n}(X_t)=0\big]\\[.2cm]
&\leq&\dis\prod_{i\in{\rm supp}(\pi_n)}
P^x\big[\Thin_{\la_n\pi_n(i)}(X^{g_i(\De)}_t(i))=0\big]\\[.2cm]
&\leq&\dis\prod_{i\in{\rm supp}(\pi_n)}
P^{T_{g_i^{-1}}x}\big[\ex{-\la_n\pi_n(i)X^\De_t(i)}\big]\\[.2cm]
&\leq&\dis\prod_{i\in{\rm supp}(\pi_n)}
P^{T_{g_i^{-1}}x}\big[\ex{-X^\De_t(i)}\big]^{\la_n\pi_n(i)},
\ec
where the $T_{g_i^{-1}}$ are shift operators as in (\ref{shift}) and we
have used that $P[\Thin_\phi(x)=0]=E[(1-\phi)^x]
=E[e^{\li\log(1-\phi),x\re}]\leq E[e^{-\li\phi,x\re}]$ for any
$\phi\in[0,1]^\La$, $x\in\N^\La$.

If $\Li(X_0)$ is $G$-homogeneous, then by (\ref{prodest}) and
H\"older's inequality,
\bc
P\big[\Thin_{\phi_n}(X_t)=0\big]
&\leq&\dis\int P[X_0\in\di x]\prod_{i\in{\rm supp}(\pi_n)}
P^{T_{g_i^{-1}}x}\big[\ex{-X^\De_t(i)}\big]^{\la_n\pi_n(i)}\\[.2cm]
&\leq&\dis\prod_{i\in{\rm supp}(\pi_n)}\Big(\int P[X_0\in\di x]\,
P^{T_{g_i^{-1}}x}\big[\ex{-X^\De_t(i)}\big]^{\la_n}\Big)^{\pi_n(i)}\\[.2cm]
&=&\dis\int P[X_0\in\di x]\,P^x\big[\ex{-X^\De_t(0)}\big]^{\la_n},
\ec
and therefore, by (\ref{Denul}) and the fact that $\la_n\to\infty$,
\be\label{basest}
\limsup_{n\to\infty}P\big[\Thin_{\phi_n}(X_t)=0\big]
\leq P\big[X_0|_\De=0\big].
\ee
Put
\be
\De_k:=\bigcup_{n=0}^k\big\{i:\exists i_0,\ldots,i_n\mbox{ s.t.\ }
i_0=i,\ i_n=0,\ a(i_{l-1},i_l)>\ffrac{1}{k}\ \forall l=1,\ldots,n\big\}.
\ee
Then the $\De_k$ satisfy (\ref{Denul}) and $\De_k\up R\{0\}$ as
$k\up\infty$, where $R\{0\}$ is defined in (\ref{Rdef}). Therefore,
inserting $\De=\De_k$ in (\ref{basest}) and taking the limit
$k\up\infty$, using Lemma~\ref{reached}, we arrive at
(\ref{posthin}).\qed

\chapter[The contact process seen from a typical site]{The contact process seen from a typical infected site}\label{C:typ}

\renewcommand{\labelenumi}{\arabic{enumi}.}

\section{Introduction and main results}

\subsection{Contact processes on countable groups}\label{S:defi}

The aim of this chapter is to study contact processes on rather general
lattices. In particular, we are interested in the way how a certain
property of the lattice, namely subexponential growth, influences the
behavior of the process.

To keep things reasonably simple, we assume that the lattice $\La$ is
a countably infinite group with group action $(i,j)\mapsto ij$ and
unit element $0$, also referred to as the origin. Each site $i\in\La$
can be in one of two states: healthy or infected. Infected sites
become healthy with {\em recovery rate} $\de\geq 0$. An infected site
$i$ infects another site $j$ with {\em infection rate} $a(i,j)\geq
0$. We assume that the infection rates are invariant with respect to
the left action of the group, summable, and statisfy a condition that
is a bit stronger than irreducibility:
\be\ba{rl}\label{assum}
{\rm (i)}&a(i,j)=a(ki,kj)\qquad\qquad(i,j,k\in\La),\\[5pt]
{\rm (ii)}&\dis|a|:=\sum_ia(0,i)<\infty,\\[5pt]
{\rm (iii)}&\bigcup_{n\geq 0,\ m\geq 0}A^nA^{-m}
=\bigcup_{n\geq 0,\ m\geq 0}A^{-n}A^m=\La,\\[5pt]
&\mbox{where }A:=\{i\in\La:a(0,i)>0\}.
\ec
Here we adopt the convention that sums over $i,j,k$ always run over
$\La$, unless stated otherwise. For $i\in\La$ and $A,B\sub\La$ we put
$AB:=\{ij:i\in A,\ j\in B\}$, $iA:=\{i\}A$, $Ai:=A\{i\}$,
$A^{-1}:=\{i^{-1}:i\in A\}$, $A^0:=\{0\}$, $A^n:=AA^{n-1}$ $(n\geq
1)$, and $A^{-n}:=(A^{-1})^n=(A^n)^{-1}$. We let $|A|$ denote the
cardinality of $A$. Note that property~(\ref{assum})~(iii) is
equivalent to the statement that for any two sites $i,j$ there exists
a site $k$ from which both $i$ and $j$ can be infected, and a set $k'$
that can be infected both from $i$ and from $j$.

If $\La$ has a finite symmetric generating set $\De$, then the (left)
Cayley graph $\Gi=\Gi(\La,\De)$ associated with $\La$ and $\De$ is the
graph with vertex set $\Vi(\Gi):=\La$ and edges
$\Ei(\Gi):=\{\{i,j\}:i^{-1}j\in\De\}$. Examples of Cayley graphs are
the $d$-dimensional integer lattice $\Z^d$ $(d\geq 1)$ with edges
between points at distance one, or the regular tree $\Tr{d}$ $(d\geq
2)$ in which every vertex has $d+1$ neighbors. On Cayley graphs, one
often considers {\em symmetric nearest-neighbor} infection rates of
the form $a(i,j)=\lir1_{\{i^{-1}j\in\De\}}$, with $\lir>0$. In this
case, $\lir$ is simply referred to as `the' infection rate.

Let $\eta_t$ be the set of all infected sites at time $t\geq 0$. Then
$\eta=(\eta_t)_{t\geq 0}$ is a Markov process in the space
$\Pc(\La):=\{A:A\sub\La\}$ of all subsets of $\La$, called the contact
process on $\La$ with infection rates $a=(a(i,j))_{i,j\in\La}$ and
recovery rate $\de$, or shortly the {\em $(\La,a,\de)$-contact
process}. If $\de>0$, then by rescaling time we may set $\de=1$, so it
is customary so assume that $\de=1$. If $\de=0$ then $\eta$ is a
special case of first-passage percolation (see \cite{Kes86}). We equip
$\Pc(\La)\cong\{0,1\}^\La$ with the product topology and the
associated Borel-\si-field $\Bi(\Pc(\La))$, and let $\Pc_{\rm
fin}(\La):=\{A\sub\La:|A|<\infty\}$ denote the subspace of finite
subsets of $\La$.

The contact process can be constructed with the help of Harris'
\cite{Har78} {\em graphical representation}. Let $\poi=(\poi^{\rm
r},\poi^{\rm i})$ be a pair of independent, locally finite random
subsets of $\La\times\R$ and $\La\times\La\times\R$, respectively,
produced by Poisson point processes with intensity $\de$ and local
intensity $(j,k,t)\mapsto a(j,k)$, respectively. This is usually
visualized by plotting $\La\times\R$ with $\La$ horizontally and $\R$
vertically. Points $(i,s)\in\poi^{\rm r}$ and $(j,k,t)\in\poi^{\rm i}$
are marked with a recovery symbol $\ast$ at $(i,s)$ and an infection
arrow from $(j,t)$ to $(k,t)$, respectively. For $C,D\sub\La\times\R$,
say that there is a {\em path} from $C$ to $D$, denoted by $C\leadsto
D$, if there exist $n\geq 0$, $i_0,\ldots,i_n\in\La$, and
$t_0\leq\cdots\leq t_{n+1}$ with $(i_0,t_0)\in C$ and
$(i_n,t_{n+1})\in D$, such that
$\{i_k\}\times[t_k,t_{k+1}]\cap\poi^{\rm r}=\emptyset$ for all
$k=0,\ldots,n$ and $(i_{k-1},i_k,t_k)\in\poi^{\rm i}$ for all
$k=1,\ldots,n$. Thus, a path must walk upwards in time, may follow
arrows, and must avoid recoveries. For given $A\in\Pc(\La)$ and
$t_0\in\R$, put
\be
\eta^{A\times\{t_0\}}_t:=\{i\in\La:A\times\{t_0\}\leadsto(i,t_0+t)\}
\qquad(t\geq 0).
\ee
Then $\eta^{A\times\{t_0\}}=(\eta^{A\times\{t_0\}}_t)_{t\geq 0}$ is a
copy of the $(\La,a,\de)$-contact process started in
$\eta^{A\times\{t_0\}}_0=A$. For brevity, we put
$\eta^A:=\eta^{A\times\{0\}}$. The graphical representation couples
processes with different initial states in such a way that
\be\label{addit}
\eta^A_t\cup\eta^B_t=\eta^{A\cup B}_t\qquad(A,B\in\Pc(\La),\ t\geq 0).
\ee

Define {\em reversed infection rates} $a^\dgg$ by
$a^\dgg(i,j):=a(j,i)$ $(i,j\in\La)$. Say that $a$ is {\em symmetric}
if $a=a^\dgg$. For $A\in\Pc(\La)$ and $t_0\in\R$, put
\be
\eta^{\dgg\,A\times\{t_0\}}_t:=\{i\in\La:(i,t_0-t)\leadsto A\times\{t_0\}\}
\qquad(t\geq 0).
\ee
Then $\eta^{\dgg\,A\times\{t_0\}}
=(\eta^{\dgg\,A\times\{t_0\}}_t)_{t\geq 0}$ is a copy of the
$(\La,a^\dgg,\de)$-contact process started in
$\eta^{\dgg\,A\times\{t_0\}}_0=A$. For brevity, we put
$\eta^{\dgg\,A}:=\eta^{\dgg\,A\times\{0\}}$. Since for any $s\leq t$
and $A,B\in\Pc(\La)$, the event
\be
\big\{\eta^{A\times\{s\}}_{u-s}\cap\eta^{\dgg\,B\times\{t\}}_{t-u}
=\emptyset\big\}=\big\{A\times\{s\}\not\leadsto B\times\{t\}\big\}
\ee
does not depend on $u\in[s,t]$, it follows that the $(\La,a,\de)$-contact
process and the $(\La,a^\dgg,\de)$-contact process are dual in the sense that
\be\label{dual}
P[\eta^A_t\cap B=\emptyset]=P[A\cap\eta^{\dgg\,B}_t=\emptyset]
\qquad(A,B\in\Pc(\La),\ t\geq 0).
\ee
For any $C\sub\La\times\R$, say that $C\leadsto\infty$ if there is an
infinite path with times $t_k\up\infty$ starting in $C$, and define
$-\infty\leadsto D$ analogously. Instead of $\{(i,s)\}\leadsto\:$ and
$\:\leadsto\{(j,t)\}$, simply write $(i,s)\leadsto\:$ and
$\:\leadsto(j,t)$. We say that the $(\La,a,\de)$-contact process
$\eta$ {\em survives} if
\be\label{surv}
\rho(A):=P\big[\eta^A_t\neq\emptyset\ \forall t\geq 0\big]
=P[A\times\{0\}\leadsto\infty]>0
\ee
for some, and hence for all $\emptyset\neq A\in\Pc_{\rm fin}(\La)$. If
$\eta$ does not survive then we say that it {\em dies out}. Set
$\de_{\rm c}=\de_{\rm c}(\La,a):=\sup\{\de\geq 0:
\,\mbox{the $(\La,a,\de)$-contact process survives}\}$. Then the
$(\La,a,\de)$-contact process survives for $\de<\de_{\rm c}$ and dies
out for $\de>\de_{\rm c}$. One has $\de_{\rm c}\leq|a|$. If $\La$ is
finitely generated, then moreover $\de_{\rm c}>0$ (see
Section~\ref{S:surv}).

\subsection{Long-time behavior}

Since the $(\La,a,\de)$-contact process is an attractive spin system,
it has an {\em upper invariant law} $\ov\nu$, i.e., an invariant law
that is maximal with respect to the stochastic order. It may be
constructed as $\ov\nu=P[\ov\eta_0\in\cdot\,]$, where
\be\label{oveta}
\ov\eta_t:=\{i\in\La:-\infty\leadsto(i,t)\}\qquad(t\in\R).
\ee
Note that
\be\label{nuchar}
P\big[\ov\eta_0\cap A\neq\emptyset\big]=\rho^\dgg(A)
\qquad(A\in\Pc_{\rm fin}(\La)),
\ee
where $\rho^\dgg$ denotes the survival probability of the
$(\La,a^\dgg,\de)$-contact process. It is easy to see that $\ov\nu$ is
nontrivial if and only if the $(\La,a^\dgg,\de)$-contact process
survives. Here, we say that a probability law on $\Pc(\La)$ is 
{\em nontrivial} if it gives zero probability to the empty set.

We say that a probability law $\mu$ on $\Pc(\La)$ is {\em homogeneous}
if $\mu$ is shift invariant with respect to the left action of the
group, i.e., $\mu(\{iA:A\in\Ai\})=\mu(\Ai)$ for all
$\Ai\in\Bi(\Pc(\La))$. Using duality, it can be shown that
\be\label{homcon}
\int\mu(\di A)P[\eta^A_t\in\cdot\,]\Asto{t}\ov\nu
\ee
whenever the initial law $\mu$ is homogeneous and nontrivial (see
\cite{Har76}, \cite[(VI.2.1)]{Lig85}, and
\cite[(I.1.10)]{Lig99}). Here $\Rightarrow$ denotes weak convergence
of probability laws. In particular, (\ref{homcon}) shows that if
$\ov\nu$ is nontrivial, then it is the only nontrivial homogeneous
invariant law.

The long-time behavior for nonhomogeneous initial laws is more subtle
and depends on properties of the lattice $\La$ and the infection rates
$a$, such as subexponential growth.

For the symmetric nearest-neighbor contact process on $\Z^d$ started
in a finite initial state, the following picture has been rigorously
verified. Either the process dies out in finite time, or in the long
run there is a region in space with linearly growing diameter and
deterministic limiting shape, such that most of the infected sites lie
within this region and there the process is locally in the upper
invariant law \cite{BG90}. In particular, it has been shown that the
symmetric nearest-neighbor process on $\Z^d$ exhibits {\em complete
convergence}, i.e.,
\be\label{compcon}
P[\eta^A_t\in\cdot\,]\Asto{t}\rho(A)\ov\nu+(1-\rho(A))\de_0
\qquad(A\in\Pc_{\rm fin}(\La)).
\ee
Note that if complete convergence holds and $\ov\nu$ is nontrivial,
then by monotonicity, it is the unique nontrivial invariant law. For
other contact processes on $\Z^d$ the picture is supposedly similar,
provided that the infection rates are symmetric and satisfy an
appropriate tail condition. If the infection rates are not symmetric,
there is probably still a linearly growing infected region with a
limiting shape, but this region may walk out to infinity, so that
complete convergence does not hold. (For results in the
one-dimensional case, see \cite{Sch86}.)

The behavior of the symmetric nearest-neighbor process on regular
trees $\Tr{d}$ is known to be quite different. Here, there is a second
critial value $\de'_{\rm c}<\de_{\rm c}$ such that for recovery rates
$\de\in[\de'_{\rm c},\de_{\rm c})$, the process survives globally but
not locally, i.e., $\rho(A)>0$ but $P[\exists T\geq 0
\mbox{ s.t.\ }\eta^A_t\cap\{0\}=\emptyset\ \forall t\geq T]=1$
for $\emptyset\neq A\in\Pc_{\rm fin}(\La)$. In this regime, there is a
multitude of nontrivial invariant measures and complete convergence
(obviously) does not hold \cite[Section~I.4]{Lig99}.

One would like to understand which properties of the lattices $\Z^d$
and $\Tr{d}$ are responsible for the differences in the behavior of
the contact process, and which types of behavior are possible on
general lattices $\La$. The proofs for $\Z^d$ and $\Tr{d}$ use the
structure of these lattices in an essential way, and are not easily
generalized to other lattices.

It is known that (unoriented) percolation has quite different
properties on $\Z^d$ and on $\Tr{d}$. Here, the important property of
$\Z^d$, that $\Tr{d}$ lacks, is {\em amenability}. For example, the
Burton-Keane proof of the uniqueness of the infinite cluster
\cite{BK89} works on any amenable graph. Conversely, it is conjectured
that on any nonamenable graph, there exists a range of the percolation
parameter for which the infinite cluster is not unique. (See
\cite{BS01} and \cite{LP05} some partial results in this direction.)

For our main theorem, we will need to assume that the expected number
of infected sites in a contact process grows subexponentially. If
$\La$ is finitely generated, then it turns out that the
$(\La,a,\de)$-contact process grows subexponentially if $a$ satisfies
an exponential moment condition and $\La$ itself has subexponential
growth (see Proposition~\ref{P:rate}~(d) below). Here, by definition,
a finitely generated group $\La$ has {\em subexponential growth} if
\be
\lim_{n\to\infty}\frac{1}{n}\log|\De^n|=0
\ee
for some, and hence for all finite symmetric generating sets
$\De$. Observe that $\De^n=\{i:|i|\leq n\}$ where $|i|$ denotes the
distance of $i$ to the origin in the Cayley graph
$\Gi(\La,\De)$. Subexponential growth is stronger than amenability. An
example of an amenable finitely generated group that does not have
subexponential growth is the lamplighter group. (See
\cite[Section~5]{MW89} for general facts about amenability and
subexponential growth, and \cite{LPP96} or \cite[\S~6.1]{LP05} for a
nice exposition of the lamplighter group.)

\subsection{Results}

It turns out that every $(\La,a,\de)$-contact process has a
well-defined exponential growth rate.
\bp{\bf(Exponential growth rate)}\label{P:rate}

\noi
{\bf (a)} There exists a constant $r=r(\La,a,\de)\in[-\de,|a|-\de]$ such that
 the $(\La,a,\de)$-contact process satisfies
\be\label{expr}
\lim_{t\to\infty}\,\ffrac{1}{t}\log E\big[|\eta^A_t|\big]=r
\qquad(\emptyset\neq A\in\Pc_{\rm fin}(\La)).
\ee

\noi
{\bf (b)} If the $(\La,a,\de)$-contact process survives, then $r\geq 0$.\med

\noi
{\bf (c)} $r(\La,a,\de)=r(\La,a^\dgg,\de)$.\med

\noi

{\bf (d)} Assume that $\La$ is finitely generated. Let $\De$ be a
finite symmetric generating set and let $|j|$ denote the distance of
$j$ to the origin in the Cayley graph $\Gi(\La,\De)$. Assume that
$\sum_ja(0,j)e^{\eps|j|}<\infty$ for some $\eps>0$ and that $\La$ has
subexponential growth. Then $r\leq 0$.
\ep
The proof of Proposition~\ref{P:rate} will be given in
Sections~\ref{S:exp}--\ref{S:amen}. Part~(a) follows from
subadditivity, part~(b) is trivial, and part~(c) is a consequence of
duality. Part~(d) follows from some basic large deviation
estimates. The exponential moment condition on $a$ appearing in
part~(d) can perhaps be weakened, but we conjecture that it cannot be
dropped altogether. Indeed, it seems plausible that even on $\La=\Z$,
the exponential growth rate can be positive if $a$ has a sufficiently
heavy tail.

To formulate the main results of this chapter, we must describe the
contact process as seen from a `typical' infected site at a `typical'
late time. Assume that the exponential growth rate $r$ from
Proposition~\ref{P:rate} satisfies $r\leq 0$. Recall the graphical
construction of the $(\La,a,\de)$-contact process (see
Section~\ref{S:defi}). Let $(\om,\Fi,P)$ be the probability space of
the Poisson point processes used in the graphical representation.
For $\lr>r$, we define probability measures $\hat
P^A_\lr$ on $\La\times\om\times\R_+$ by
\be\label{Campdefiv}
\hat P^A_\lr(\{i\}\times\{\di\oo\}\times\{\di t\}):=
\frac{1}{\expi_\la(A)}\,1_{\txt\{i\in\eta^A_t(\oo)\}}P(\di\oo)e^{-\la t}\di t,
\ee
where
\be\label{pidef2}
\expi_\lr(A):=\int_0^\infty\! E\big[|\eta^A_t|\big]e^{-\la t}\,\di t
\qquad(A\in\Pc_{\rm fin}(\La),\ \lr>r)
\ee
is a normalizing constant. Using the fact that $\lr>r$, it is easy
to see that $0<\expi_\la(A)<\infty$, so $\hat P^A_\lr$ is
well-defined. Note that the projection of $\hat P^A_\lr$ on
$\om\times\R_+$ is given by
\be
P^A_\lr(\La\times\{\di\oo\}\times\{\di t\})
=\frac{1}{\expi_\la(A)}\,|\eta^A_t(\oo)|P(\di\oo)e^{-\la t}\di t
\ee
In other words, this projection is is obtained from the product
measure $P(\do\oo)e^{-\la t}\di t$ on $\om\times\R_+$ by size-biasing
with the number of infected sites $|\eta^A_t(\oo)|$. Let $\iota$ and
$\tau$ denote the projections on $\La$ and $\R_+$, respectively. Then,
under the law $\hat P^A_\lr$, the random variable $\eta^A_\tau$
describes a size-biased contact process as a `typical' time $\tau$,
and $\iota$ is a `typical' infected site, chosen with equal
probabilities from $\eta^A_\tau$. The law $\hat
P^A_\lr[(\iota,\eta^A_\tau)\in\cdot\,]$ is a Campbell law,
which is closely related to the more widely known Palm laws. (For the
relation between Campbell and Palm laws, see
\cite[Section~6.4]{Eth00}.) The next lemma says that as $\la$ decreases
to $r$, under the laws $P^A_\lr$, the `typical' time $\tau$ tends in
probability to $\infty$. Thus, the limit $\la\down r$
corresponds to letting time to infinity.
\bl{\bf(Typical times)}\label{L:typtime}
For each $\emptyset\neq A\in\Pc_{\rm fin}(\La)$,
\be\label{lart2}
\hat P^A_\la\big[\tau\geq t\big]\asdto{\la}{r}1\qquad(t>0).
\ee
\el
Note that $\iota^{-1}\eta_\tau$ is the process $\eta_\tau$, viewed
{f}rom the position of the typical infected site~$\iota$. The next
theorem is the main result of this chapter. Recall the definition of
$\ov\eta$ in (\ref{oveta}).
\bt{\bf(The process seen from a typical infected site)}\label{T:palm}
Assume that the upper invariant measure of the $(\La,a,\de)$-contact
process is nontrivial and that the exponential growth rate from
Proposition~\ref{P:rate} satisfies $r=0$. Let $\emptyset\neq
A\in\Pc_{\rm fin}(\La)$. Then\med

\noi
{\bf (a)} One has
\be\label{conv}
\hat P^A_\lr\big[\iota^{-1}\eta^A_\tau\in\cdot\,\big]\Asdto{\lr}{0}
P\big[\ov\eta_0\in\cdot\,\big|\,0\in\ov\eta_0\big].
\ee
{\bf (b)} Moreover,
\be\label{isup}
\hat P^A_\lr\big[\iota^{-1}\eta^A_\tau\cap\De
=\iota^{-1}\ov\eta_\tau\cap\De\big]\asdto{\lr}{0}1
\qquad(\De\in\Pc_{\rm fin}(\La)),
\ee
and the same holds with $\ov\eta_\tau$ replaced by $\eta^\La_\tau$.
\et
Note that Theorem~\ref{T:palm} holds when $\La$ is a general countable
group, but we have only verified that its assumptions are satisfied
for certain finitely generated groups (see
Proposition~\ref{P:rate}~(d)). We remark that for fixed $\lr>0$, it is
not at all obvious (and as far as we know {\em not true}) that the
distribution $\hat P^A_\lr[\iota^{-1}\ov\eta_\tau\in\cdot\,]$ should
be the same as $P[\ov\eta_0\in\cdot\,|\,0\in\ov\eta_0]$. Thus, none of
the statements (\ref{conv}) and (\ref{isup}) trivially implies the
other one.

As a result of our methods, we can also prove the following fact,
which is of some interest on its own.
\bp{\bf(Typical particles descend from every surviving site)}\label{P:totyp}
Assume that the $(\La,a,\de)$-contact process survives and that the
exponential growth rate from Proposition~\ref{P:rate} satisfies $r=0$.
Then
\be\label{totyp}
\hat P^{\{i\}}_\lr\big[(j,0)\leadsto(\iota,\tau)\,\big|
\,(j,0)\leadsto\infty\big]\asdto{\lr}{0}1\qquad(i,j\in\La).
\ee
\ep
One of the original motivations of the present chapter was to answer the
following question. Assuming survival and subexponential growth, is it
true that for any $i,j\in\La$,
\be\label{clustiv}
P\big[\exists(k,t)\mbox{ s.t.\ }(i,0)\leadsto(k,t)\leadsto\infty
\mbox{ and }(j,0)\leadsto(k,t)\,\big|\,(i,0)\leadsto\infty,
(j,0)\leadsto\infty\big]=1\quad?
\ee
This property may be interpreted as some sort of analogue of the
uniqueness of the infinite cluster in (unoriented)
percolation. Unfortunately, we do not know how to replace the
size-biased law in (\ref{totyp}) by a law conditioned on
survival. Question (\ref{clustiv}) has been answered positively for
oriented percolation on $\Z^d$ in \cite{GH02}. As a further motivation
for (\ref{clustiv}), we note that in the one-dimensional
nearest-neighbor case, a considerably stronger statement holds.
\bl{\bf(Coupling of one-dimensional processes)}\label{L:onedim}
Consider a $(\Z,a,\de)$-contact process with $a(i,j)=0$ for $|i-j|\neq 1$.
Assume that the process survives, and assume either $\de>0$ or
$a(0,1)\wedge a(1,0)>0$. Then, for any $i,j\in\Z$,
\be
P\big[\inf\{t\geq 0:\eta^{\{i\}}_t=\eta^{\{j\}}_t\}<\infty\,\big|
\,(i,0)\leadsto\infty,\ (j,0)\leadsto\infty\big]=1.
\ee
\el

\subsection{Methods}\label{S:methods}

In this section we describe the main line of our proof of
Theorem~\ref{T:palm}~(a). The first ingredient is a chararacterization
of the laws $\hat P^A_\lr[\iota^{-1}\eta^A_\tau\in\cdot\,]$ and
$P[\ov\eta_0\in\cdot\,|\,0\in\ov\eta_0]$ in terms of the dual
$(\La,a^\dgg,\de)$-contact process $\eta^\dgg$. For simplicity, we
only present the argument for $\hat P^{\{0\}}_\lr$. Let $\expi_\la(A)$
be the normalizing constant in (\ref{pidef2}).
Recall the definition of the survival probability $\rho$ in
(\ref{surv}). We write $\ov\expi_\lr$ and $\ov\rho$ for the functions
$\expi_\lr$ and $\rho$ normalised to one in the point $\{0\}$:
\be\label{ovpi}
\ov\rho(A):=\frac{\rho(A)}{\rho(\{0\})}\quad\mbox{and}
\quad\ov\expi_\lr(A):=\frac{\expi_\lr(A)}{\expi_\lr(\{0\})}.
\ee
We let $\rho^\dgg,\expi_\lr^\dgg,\ov\rho^\dgg$, and
$\ov\expi_\lr^\dgg$ denote the analogues of $\rho,\expi_\lr,\ov\rho$,
and $\ov\expi_\lr$ for the dual $(\La,a^\dgg,\de)$-contact process.
\bl{\bf(Characterization of laws seen from an infected site)}\label{L:char}\\
{\bf(a)} One has
\be\label{nutchar}
\hat P^{\{0\}}_\lr\big[A\cap\iota^{-1}\eta^{\{0\}}_\tau=\emptyset\big]
=\ov\expi^\dgg_\lr(A\cup\{0\})-\ov\expi^\dgg_\lr(A)
\qquad(A\in\Pc_{\rm fin}(\La),\ \lr>r).
\ee
{\bf(b)} Moreover,
\be\label{nu0char}
P\big[A\cap\ov\eta_0=\emptyset\,\big|\,0\in\ov\eta_0\big]
=\ov\rho^\dgg(A\cup\{0\})-\ov\rho^\dgg(A)
\qquad(A\in\Pc_{\rm fin}(\La)).
\ee
\el
It is not hard to see that the law of a $\Pc(\La)$-valued random
variable $\eta$ is uniquely characterized by all probabilities of the
form $P[A\cap\eta=\emptyset]$ with $A\in\Pc_{\rm fin}(\La)$.
Therefore, by Lemma~\ref{L:char} and the compactness of
$\Pc(\La)\cong\{0,1\}^\La$, in order to prove Theorem~\ref{T:palm}, it
suffices to prove that under the assumptions there,
$\ov\expi^\dgg_\lr\to\ov\rho^\dgg$ pointwise as
$\lr\down 0$. In order to reduce notation, we reverse the role of
$\eta$ and $\eta^\dgg$. Thus, we will prove that pointwise
$\lim_{\lr\down 0}\ov\expi_\lr=\ov\rho$, under the assumptions that
the $(\La,a,\de)$-contact process survives and its exponential growth
rate is zero. (By (\ref{nuchar}) and Proposition~\ref{P:rate}~(c),
this is equivalent to the $(\La,a^\dgg,\de)$-contact process having a
nontrivial upper invariant law and exponential growth rate zero.)

It is not hard to show (see Section~\ref{S:mart} below) that the
$(\La,a,\de)$-contact process started in a finite initial state solves
the martingale problem for the operator
\be\label{Gdef}
Gf(A):=\sum_{ij}a(i,j)1_{\{i\in A\}}\{f(A\cup\{j\})-f(A)\}
+\de\sum_i1_{\{i\in A\}}\{f(A\beh\{i\})-f(A)\},
\ee
with domain $\Di(G):=\Si(\Pc_{\rm fin}(\La))$, where
\be\label{Gdom}
\Si(\Pc_{\rm fin}(\La)):=\{f:\Pc_{\rm fin}(\La)\to\R:
|f(A)|\leq K|A|^k+M\mbox{ for some }K,M,k\geq 0\}.
\ee
It can be shown in a few lines that $\rho$ is shift invariant, monotone
(i.e., $A\sub B$ implies $\rho(A)\leq\rho(B)$),
$\rho\in\Si(\Pc_{\rm fin}(\La))$, and
\be\label{rhoharm}
G\rho=0.
\ee
Formula (\ref{rhoharm}) says that $\rho$ is a harmonic function for
the $(\La,a,\de)$-contact process. It is not hard to see that
$\expi_\lr$ shift invariant, monotone, $\expi_\lr\in\Si(\Pc_{\rm
fin}(\La))$, and
\be\label{Gip}
G\expi_\lr(A)=\lr\expi_\lr(A)-|A|
\qquad(A\in\Pc_{\rm fin}(\La),\ \lr>r).
\ee
As a consequence, one obtains:
\bl{\bf(Cluster points of the rescaled expected population
size)}\label{L:piclus}
The functions $(\ov\expi_\lr)_{\lr>r}$ are
relatively compact with respect to the product topology on
$\R^{\Pc_{\rm fin}(\La)}$. Each pointwise limit
\be\label{pinf}
\ov\expi_r(A):=\lim_{n\to\infty}\ov\expi_{\lr_n}(A)
\qquad(A\in\Pc_{\rm fin}(\La))
\ee
along a sequence $\lr_n\down r$ is shift invariant, monotone in
$A$, satisfies $\ov\expi_r\in\Si(\Pc_{\rm fin}(\La))$, and
\be\label{Gr}
G\ov\expi_r=r\ov\expi_r.
\ee
\el
In particular, if $r=0$ and the $(\La,a,\de)$-contact process
survives, it turns out that Lemma~\ref{L:piclus} gives us enough
information to determine $\ov\expi_0$ uniquely. Combined with the next
proposition, Lemma~\ref{L:piclus} shows that $\ov\expi_\lr\to\ov\rho$
pointwise as $\lr\down 0$, thereby completing the proof of
Theorem~\ref{T:palm}.
\bp{\bf(Shift\, invariant\, monotone\, harmonic\, functions)}\label{P:monhar}
Assume\, that\, the $(\La,a,\de)$-contact process survives. Assume that
$f:\Pc_{\rm fin}(\La)\to\R$ is shift invariant, monotone, $f(\emptyset)=0$,
$f\in\Si(\Pc_{\rm fin}(\La))$, and $Gf=0$. Then there exists a constant
$c\geq 0$ such that $f=c\rho$.
\ep
We note that if $\nu$ is a homogeneous invariant measure for the
$(\La,a^\dgg,\de)$-contact process, then by duality,
$f(A):=\nu(\{A:A\cap B\neq\emptyset\})$ defines a shift invariant,
monotone, bounded harmonic function $f$ for the $(\La,a,\de)$-contact
process. Therefore, in view of (\ref{nuchar}),
Proposition~\ref{P:monhar} is a strengthening of the statement that
all homogeneous invariant measures are convex combinations of $\ov\nu$
and~$\de_0$.

In order to prove Proposition~\ref{P:monhar}, we need one more lemma.
\bl{\bf(Eventual domination of finite configurations)}\label{L:evdom}
Assume that the $(\La,a,\de)$-contact process survives. Then
\be\label{evdom}
\lim_{t\to\infty}P\big[\exists i\in\La\mbox{ s.t.\ }\eta^A_t\geq iB\,\big|
\,\eta^A_t\neq\emptyset]=1
\qquad(A,B\in\Pc_{\rm fin}(\La),\ A\neq\emptyset).
\ee
\el
Formula (\ref{evdom}) says that $\eta$ exhibits a form of extinction
versus unbounded growth. More precisely, either $\eta_t$ gets extinct
or $\eta_t$ is eventually larger than a suitable shift (depending on
$\eta_t$) of any finite configuration. We remark that
Lemma~\ref{L:evdom} is no longer true if assumption
(\ref{assum})~(iii) is replaced by the weaker assumption that
$\{i\in\La:a(0,i)>0\}$ generates $\La$.\med

\noi
{\bf Proof of Proposition~\ref{P:monhar}} Since the
$(\La,a,\de)$-contact process solves the martingale problem for $G$,
and $Gf=0$, the process $f(\eta^A_t)$ is a martingale. In particular:
\be
f(A)=E[f(\eta^A_t)]\qquad(A\in\Pc_{\rm fin}(\La),\ t\geq 0).
\ee
Equip $\La$ with an arbitrary linear ordering, and for
$A,B\in\Pc_{\rm fin}(\La)$, put
\be
\hat\imath_{A,B}:=\left\{\ba{ll}\min\{i\in\La:A\geq iB\}\quad
&\mbox{if }\{i\in\La:A\geq iB\}\mbox{ is nonempty,}\\
0&\mbox{otherwise.}\ea\right.
\ee
Since $f$ is monotone and shift invariant, we have, using Lemma~\ref{L:evdom},
\bc\label{xy}
\dis f(A)&=&\dis\lim_{t\to\infty}E[f(\eta^A_t)]\\[5pt]
&\geq&\dis\limsup_{t\to\infty}E[1_{\txt\{\exists i\in\La
\mbox{ s.t.\ }\eta^A_t\geq iB\}}
f(\hat\imath_{\eta^A_t,B}B)]\\[5pt]
&=&\dis f(B)\limsup_{t\to\infty}P[\exists i\in\La
\mbox{ s.t.\ }\eta^A_t\geq iB]
\geq f(B)\rho(A)\qquad(A,B\in\Pc_{\rm fin}(\La)).
\ec
In particular, this shows that
\be
f(B)\leq\frac{f(\{0\})}{\rho(\{0\})}<\infty\qquad(B\in\Pc_{\rm fin}(\La)),
\ee
hence $f$ is bounded. Now let $A_n,B_m\in\Pc_{\rm fin}(\La)$ be
sequences such that $\rho(A_n)\to 1$ and $\rho(B_n)\to 1$. Then, by
(\ref{xy}),
\be
\liminf_{n\to\infty}f(A_n)\geq\liminf_{n\to\infty}f(B_m)\rho(A_n)
=f(B_m)\quad\forall m,
\ee
and therefore
\be\label{limsupinf}
\liminf_{n\to\infty}f(A_n)\geq\limsup_{m\to\infty}f(B_m).
\ee
This proves that the limit
\be
\lim_{\rho(A_n)\to 1}f(A_n)=:f(\infty)
\ee
exists and does not depend on the choice of the sequence $A_n$ with
$\rho(A_n)\to 1$. By the Markov property and continuity of the
conditional expectation with respect to increasing limits of
\si-fields (see Complement 10(b) from \cite[Section~29]{Loe63} or
\cite[Section~32]{Loe78}),
\be\label{loe}
\rho(\eta^A_t)=P\big[\eta^A_s\neq 0\ \forall s\geq 0\,\big|\,\eta^A_t\big]\to
1_{\txt\{\eta^A_s\neq 0\ \forall s\geq 0\}}\quad{\rm a.s.}
\quad\mbox{as }t\to\infty.
\ee
We conclude that
\be
f(A)=\lim_{t\to\infty}E[f(\eta^A_t)]=\rho(A)f(\infty)
\qquad(A\in\Pc_{\rm fin}(\La)),
\ee
which shows that $f$ is a scalar multiple of $\rho$.\qed

\subsection{Discussion and open problems}

Palm and Campbell laws are standard tools in the study of (critical)
spatial branching processes. In this context, they can be described by
Kallenberg's backward tree technique; see, for example, \cite{Kal77}
or \cite{GW91}. In the context of contact processes, it is less
obvious that they should be of any use. For example, size-biasing with
$|\eta^{\{0\}}_t\cap\{i\}|$ for fixed $i$ and $t$ is just the same as
conditioning on $(0,0)\leadsto(i,t)$. In this case there seems to be
no easy way to prove statements about $i^{-1}\eta^{\{0\}}_t$.

However, by looking at the process seen from a randomly chosen
infected site rather than a fixed site, i.e., by looking at Campbell
laws rather than Palm laws, we can make a connection with the growth
of $E[|\eta^{\{0\}}_t|]$ as $t\to\infty$, and in this way obtain a
result. A disadvantage of this approach is that one ends up with
statements about size-biased laws, where one would probably be more
interested in laws conditioned on survival. Nevertheless, it seems
that the statements in Theorem~\ref{T:palm} do catch a phenomenon that
depends in a crucial way on a property of the underlying lattice, in
this case, subexponential growth.

We next state some open problems and questions, and then comment on them.
\begin{enumerate}
\item {\bf Problem} Replace the random time $\tau$
in by a deterministic time $t$ and prove the analogue
of Theorem~\ref{T:palm} for $t\to\infty$.\label{1}
\item {\bf Problem} Study the contact process seen from a typical
infected site in case the exponential growth rate is
positive.\label{2}
\item {\bf Problem} Study the contact proces seen from a typical
infected site chosen from a process conditioned to survive, instead of
size-biased on the number of infected sites.\label{3}
\item {\bf Problem} Prove (\ref{clustiv}) assuming survival and
subexponential growth.\label{3a}
\item {\bf Problem} Assuming survival and subexponential growth, prove
that conditional on $(i,0)\leadsto\infty$ and $(j,0)\leadsto\infty$,
eventually most sites in $\eta^{\{i\}}_t$ are also in
$\eta^{\{j\}}_t$.\label{3b}
\item {\bf Question} With the same set-up as in the previous problem,
is it even true that $\eta^{\{i\}}_t$ and $\eta^{\{j\}}_t$ are
eventually equal? (Compare Lemma~\ref{L:onedim}.)\label{3c}
\item {\bf Problem} Prove that $\de_{\rm c}>0$ for a contact process
on a group $\La$ that is not finitely generated, for example on the
hierarchical group.\label{4}
\item {\bf Problem} Give an example of a contact process on $\Z$ for
which the exponential growth rate is positive.\label{5}
\item {\bf Question} Assuming that $\La$ has exponential growth, is it
true that the $(\La,a,\de)$-contact process survives if and only if
$r(\La,a,d)>0$?\label{6}
\item {\bf Question} Does survival of the $(\La,a,\de)$-contact
process imply survival of the $(\La,a^\dgg,\de)$-contact
process?\label{8}
\end{enumerate}
If one tries to solve Problem~\ref{1} in a naive way, by mimicking the
techniques in this chapter, it seems one would have to strengthen
Proposition~\ref{P:rate}~(a) in the sense that
\be\label{difrate}
\lim_{t\to\infty}\dif{t}\log E\big[|\eta^A_t|\big]
=r\qquad(\emptyset\neq A\in\Pc_{\rm fin}(\La)).
\ee
Then it would follow that each cluster point $\ov\expi_\infty$ of the
functions $\ov\expi_t(A):=
E\big[|\eta^A_t|\big]/E\big[|\eta^{\{0\}}_t|\big]$ satisfies
$G\ov\expi_\infty=0$. However, (\ref{difrate}) does not simply follow
from subadditivity and seems hard to establish in general. Even random
times $\tau$ that are uniformly distributed on intervals $[0,T]$ seem
difficult to treat, since they would require that
$\lim_{T\to\infty}\dif{T}\log\int_0^T E\big[|\eta^A_t|\big]\di t=r$.

In order to solve Problem~\ref{2}, generalizing
Proposition~\ref{P:monhar}, one would like to show that the equation
$G\ov\pi_r=r\ov\pi_r$ has a unique shift invariant, monotone
solution $\ov\pi_r$ with $\ov\pi_r(\emptyset)=0$ and
$\ov\pi_r(\{0\})=1$ (perhaps also using that $\ov\pi_r$ is
subadditive).

Problems~\ref{3}--\ref{3b} and Question~\ref{3c} have been discussed
before. The difficulty is to replace size-biased laws by laws
conditioned on survival in statements like
Proposition~\ref{P:totyp}. Although size-biasing and conditioning are
asymptotically equivalent in a `local' sense (see
Proposition~\ref{P:sizco} below), this does not seem easy. Note that
if (\ref{clustiv}) holds for the $(\La,a^\dgg,\de)$-contact process,
then the limit law in Theorem~\ref{T:palm}~(a) may also be written as
$P[\hat\eta_0\in\cdot\,|\,-\infty\leadsto(0,0)]$, where
$\hat\eta_t:=\{i\in\La:\exists(j,s)\mbox{ s.t.\
}-\infty\leadsto(j,s)\leadsto(0,0)\mbox{ and }(j,s)\leadsto(i,t)\}$
$(t\in\R)$. This construction is similar to Kallenberg's
backward tree technique, and also somewhat reminiscent of the
construction of the the second lowest extremal invariant measure of
the contact process in \cite{SS97,SS99}.

Problem~\ref{4} seems interesting, since the hierarchical group has
found applications in population biology, and the usual comparison
with one-dimensional oriented percolation cannot work here.

Problem~\ref{5} and Question~\ref{6} are naturally motivated by
Proposition~\ref{P:rate}~(d). Related to Question~\ref{6} is the more
general question: what does the behavior of $E[|\eta_t|]$ for
$t\to\infty$ tell us about survival? Especially for critical
processes, it seems conceivable that
$\lim_{t\to\infty}E[|\eta_t|]=\infty$ while the process dies
out.

Related to this is Question~\ref{8}, which has been asked before for
branching-coalescing particle systems in \cite{AS05}. For symmetric
processes or for processes on abelian groups, the answer is obviously
positive, but in general $(\La,a)$ and $(\La,a^\dgg)$ need not be
isomorphic. However, in formula (\ref{Esym}) below, it is shown that
$E[|\eta^{\{0\}}_t|]=E[|\eta^{\dgg\,\{0\}}_t|]$ for all $t\geq 0$. (On
the other hand, dropping the assumption that $\La$ is a group, by
considering contact processes on transitive graphs that are not
unimodular, it is easy to construct examples where
$E[|\eta^{\{0\}}_t|]\neq E[|\eta^{\dgg\,\{0\}}_t|]$ and where $\eta$
survives but $\eta^\dgg$ dies out.) An example of a model on $\Z^2$
where nontriviality of the upper invariant law and survival are not
equivalent is the NEC model due to A.~Toom \cite{BG85,DLSS91}.

Related to Question~\ref{8} (compare also Question~\ref{3c}) is the
following question: is it always true that $\inf\{t\geq 0:
\eta^{\{0\}}_t\sub\ov\eta_t\}$ is a.s.\ finite? Note that if the
answer is positive, then extinction of the $(\La,a^\dgg,\de)$-contact
process implies extinction of the $(\La,a,\de)$-contact process, since
in this case $\ov\eta\equiv 0$.

\subsection{Outline}

Section~\ref{S:parta} is devoted to the proof of
Theorem~\ref{T:palm}~(a). In Section~\ref{S:mart} we prove that
contact processes started in finite initial states solve the
martingale problem for the operator $G$ in (\ref{Gdef}). We establish
Proposition~\ref{P:rate}~(a)--(c) in Section~\ref{S:exp}, and part~(d)
in Section~\ref{S:amen}. In Section~\ref{S:char}, we establish
Lemmas~\ref{L:typtime} and \ref{L:char}. In Section~\ref{S:harm}, we
prove basic facts about the functions $\rho$ and $\expi_\lr$; in
particular, formulas (\ref{rhoharm}) and (\ref{Gip}), and
Lemma~\ref{L:piclus}. In Section~\ref{S:dom}, we prove
Lemma~\ref{L:evdom}, thereby completing the proof of
Theorem~\ref{T:palm} in the case $A=\{0\}$. In Section~\ref{S:A} we
show how the arguments may be generalized to arbitrary $\emptyset\neq
A\in\Pc_{\rm fin}(\La)$.

Section~\ref{S:rest} contains proofs of all results that are not
directly needed for Theorem~\ref{T:palm}~(a). In
Section~\ref{S:sizco}, we prove that size-biasing and conditioning on
survival are equivalent in a `local' sense. Section~\ref{S:maxcoup}
contains the proofs of Theorem~\ref{T:palm}~(b) and
Proposition~\ref{P:totyp}. Section~\ref{S:coup} contains the proof of
Lemma~\ref{L:onedim}. For completeness, we prove in
Section~\ref{S:surv} the fact mentioned in the text that
$\de_{\rm c}>0$ whenever $\La$ is finitely generated.\med

\noi
{\bf\large Acknowledgements} The author thanks Geoffrey Grimmmett,
Olle H\"aggstr\"om, Russel Lyons, and Roberto Schonmann for useful
email conversations about the contact process, oriented percolation,
and amenability.

\section{The law seen from a typical particle}\label{S:parta}

\subsection{A martingale problem}\label{S:mart}

In this section we prove that the $(\La,a,\de)$-contact process
started in finite initial states solves the martingale problem for the
operator $G$ in (\ref{Gdef})--(\ref{Gdom}).
\bp{\bf\hspace{-1.8pt}(Martingale problem and moment estimate)\hspace{-1.7pt}}\label{L:mart}
For each $f\in\Si(\Pc_{\rm fin}(\La))$ and $A\in\Pc_{\rm fin}(\La)$,
the process
\be\label{MP}
M_t:=f(\eta^A_t)-\int_0^tGf(\eta^A_s)\di s\qquad(t\geq 0)
\ee
is a martingale with respect to the filtration generated by $\eta^A$.
Moreover, setting $z^{\li k\re}:=\prod_{i=0}^{k-1}(z+i)$, one has
\be\label{momest}
E\big[|\eta^A_t|^{\li k\re}\big]\leq|A|^{\li k\re}e^{k(|a|-\de)t}
\qquad(A\in\Pc_{\rm fin}(\La),\ k\geq 1,\ t\geq 0).
\ee
\ep
{\bf Proof} The proof of \cite[Proposition~8]{AS05} can in a straightforward
way be adapted to the present set-up. Set $f_k(A):=|A|^{\li k\re}$. Then
\bc\label{Gmom}
\dis Gf_k(A)&=&\dis\sum_{ij}a(i,j)1_{\{i\in A\}}1_{\{j\not\in A\}}
\{(|A|+1)^{\li k\re}-|A|^{\li k\re}\}+\de\sum_i1_{\{i\in A\}}
\{(|A|-1)^{\li k\re}-|A|^{\li k\re}\},\\[5pt]
&\leq&\dis(|a|-\de)|A|\{(|A|+1)^{\li k\re}-|A|^{\li k\re}\}
=k(|a|-\de)|A|^{\li k\re}.
\ec
Define stopping times $\tau_N:=\inf\{t\geq 0:|\eta^A_t|\geq N\}$. The
stopped process $(\eta^A_{t\wedge\tau_N})_{t\geq 0}$ has bounded jump
rates, and therefore standard theory tells us that for each $N\geq 1$
and $f\in\Si(\Pc_{\rm fin}(\La))$, the process
\be\label{MNP}
M^N_t:=f(\eta^A_{t\wedge\tau_N})-\int_0^{t\wedge\tau_N}Gf(\eta^A_s)\di s
\qquad(t\geq 0)
\ee
is a martingale. Moreover, it easily follows from (\ref{Gmom}) that
\be\label{stopest}
E\big[|\eta^A_{t\wedge\tau_N}|^{\li k\re}\big]\leq|A|^{\li k\re}e^{k(|a|-\de)t}\qquad(k\geq 1,\ t\geq 0).
\ee
It is easy to see that $f\in\Si(\Pc_{\rm fin}(\La))$ implies
$Gf\in\Si(\Pc_{\rm fin}(\La))$. Using this fact and (\ref{stopest})
for some sufficiently high $k$ (depending on $f$), one can show that
for fixed $t\geq 0$, the random variables $(M^N_t)_{N\geq 1}$ are
uniformly integrable. Therefore, letting $N\to\infty$ in (\ref{MNP}),
one finds that the process in (\ref{MP}) is a martingale. Letting
$N\to\infty$ in (\ref{stopest}) yields (\ref{momest}).\qed

\subsection{The exponential growth rate}\label{S:exp}

In this section we prove Proposition~\ref{P:rate}~(a)--(c).\med

\noi
{\bf Proof of Proposition~\ref{P:rate}~(a)} By a slight abuse of notation,
let us write (compare (\ref{pidef2}))
\be\label{pitdef}
\expi_t(A):=E\big[|\eta^A_t|\big]\qquad(A\in\Pc_{\rm fin}(\La),\ t\geq 0).
\ee
We start by showing that
\be\label{pimult}
\expi_{s+t}(\{0\})\leq\expi_s(\{0\})\expi_t(\{0\})\qquad(s,t\geq 0).
\ee
By (\ref{addit}),
\be\label{pixpi}
E\big[|\eta^A_t|\big]=E\Big[\big|\bigcup_{i\in A}\eta^{\{i\}}_t\big|\Big]
\leq\sum_{i\in A}E\big[|\eta^{\{i\}}_t|\big]=|A|E\big[|\eta^{\{0\}}_t|\big],
\ee
where in the last step we have used shift invariance. As a consequence,
\be
\expi_{s+t}(\{0\})=\int P[\eta^{\{0\}}_s\in\di A]E\big[|\eta^A_t|\big]
\leq\int P[\eta^{\{0\}}_s\in\di A]|A|E\big[|\eta^{\{0\}}_t|\big]
=\expi_s(\{0\})\expi_t(\{0\}).
\ee
This proves (\ref{pimult}). It follows that $t\mapsto\log\expi_t(\{0\})$
is subadditive and therefore, by \cite[Theorem~B.22]{Lig99}, the limit
\be\label{rdef}
\lim_{t\to\infty}\ffrac{1}{t}\log\expi_t(\{0\})=:r\in[-\infty,\infty]
\ee
exists. By monotonicity and (\ref{pixpi}),
\be
\expi_t(\{0\})\leq\expi_t(A)\leq|A|\expi_t(\{0\})\qquad(A\in\Pc_{\rm fin}(\La)).
\ee
Taking logarithms, dividing by $t$, and letting $t\to\infty$ we arrive
at (\ref{expr}). Since $\eta$ can be bounded from below by a simple
death process and from above by a branching process (see (\ref{brest})
below), one has
\be
e^{-\de t}\leq E\big[|\eta^{\{0\}}_t|\big]\leq e^{(|a|-\de)t}\qquad(t\geq 0),
\ee
which implies that $-\de\leq r\leq|a|-\de$.\qed

\noi
{\bf Proof of Proposition~\ref{P:rate}~(b)} If the $(\La,a,\de)$-contact
process survives, then
\be
\expi_t(\{0\})\geq P[\eta^{\{0\}}_t\neq 0]
\asto{t}P[\eta^{\{0\}}_s\neq 0\ \forall s\geq 0]>0,
\ee
which implies that $r\geq 0$.\qed

\noi
{\bf Proof of Proposition~\ref{P:rate}~(c)} By duality (formula (\ref{dual}))
and shift invariance,
\bc\label{Esym}
\dis E\big[|\eta^{\{0\}}_t|\big]
&=&\dis\sum_iP\big[\eta^{\{0\}}_t\cap\{i\}\neq\emptyset\big]
=\sum_iP\big[\{0\}\cap\eta^{\dgg\,\{i\}}_t\neq\emptyset\big]\\[5pt]
&=&\dis\sum_iP\big[\{i^{-1}\}\cap\eta^{\dgg\,\{0\}}_t\neq\emptyset\big]
=E\big[|\eta^{\dgg\,\{0\}}_t|\big],
\ec
which implies that $r(\La,a,\de)=r(\La,a^\dgg,\de)$.\qed

\subsection{Subexponential growth}\label{S:amen}

{\bf Proof of Proposition~\ref{P:rate}~(d)} Consider a branching
process on $\La$, started with one particle in the origin, where a
particle at $i$ produces a new particle at $j$ with rate $a(i,j)$, and
each particle dies with rate $\de$. Let $B_t(i)$ denote the number of
particles at site $i\in\La$ and time $t\geq 0$. It is not hard to see
that $\eta^{\{0\}}$ and $B$ may be coupled such that
\be\label{brest}
1_{\eta^{\{0\}}_t}\leq B_t\qquad(t\geq 0).
\ee
Let $(\xi_t)_{t\geq 0}$ be a random walk on $\La$ that jumps from $i$
to $j$ with rate $a(i,j)$, started in $\xi_0=0$.  Then it is not hard
to see that (compare \cite[Proposition~I.1.21]{Lig99})
\be\label{walk}
E[B_t(i)]=P[\xi_t=i]e^{(|a|-\de)t}\qquad(i\in\La,\ t\geq 0).
\ee
\detail{Indeed, if $b_t(i)=E[B_t(i)]$, then
\[
\dif{t}b_t(i)=\sum_ja(j,i)b_t(j)-b_t(i)
=(|a|-\de)b_t(i)+\sum_ja(j,i)b_t(j)-|a|b_t(i).
\]}
Let $\lc>0$ be a constant, to determined later. It follows from
(\ref{brest}) and (\ref{walk}) that
\bc\label{split}
\dis E\big[|\eta^{\{0\}}_t|\big]&\leq&
\dis\sum_i\big(1\wedge P[\xi_t=i]e^{(|a|-\de)t}\big)\\[5pt]
&=&|\{i\in\La:|i|\leq\lc t\}|+P[|\xi_t|>\lc t]e^{(|a|-\de)t}\qquad(t\geq 0).
\ec
Let $(Y_i)_{i\geq 1}$ be i.i.d.\ $\N$-valued random variables with
$P[Y_i=k]=\frac{1}{|a|}\sum_{j:\ |j|=k}a(0,j)$ $(k\geq 0)$, let $N$ be
a Poisson-distributed random variable with mean $|a|$, independent of
the $(Y_i)_{i\geq 1}$, and let $(X_m)_{m\geq 1}$ be i.i.d.\ random
variables with law
$P[X_m\in\cdot\,]=P[\sum_{i=1}^NY_i\in\cdot\,]$. Since the random walk
$\xi$ makes jumps whose sizes are distributed in the same way as the
$Y_i$, and the number of jumps per unit of time is Poisson distributed
with mean $|a|$, it follows that
\be\label{linsplit}
P[|\xi_t|>\lc t]\leq P\Big[\frac{1}{\lceil t\rceil}
\sum_{m=1}^{\lceil t\rceil}X_m>\lc\frac{t}{\lceil t\rceil}\Big]\qquad(t>0),
\ee
where $\lceil t\rceil$ denotes $t$ rounded up to the next integer.
By our assumptions,
\be
E\big[\ex{\eps X_m}\big]=E\big[\ex{\eps\sum_{i=1}^NY_k}\big]
=e^{-|a|}\sum_{n=0}^\infty\frac{|a|^n}{n!}E\big[\ex{\eps Y_1}\big]^n
=\ex{-|a|(1-E[e^{\eps Y_1}])}<\infty,
\ee
for some $\eps>0$. Therefore, by \cite[Theorem~2.2.3 and Lemma~2.2.20]{DZ98},
for each $R>0$ there exists a $\lc>0$ and $K<\infty$ such that
\be\label{Rex}
P\Big[\frac{1}{n}\sum_{m=1}^{n}X_m>\lc\Big]\leq K\ex{-nR}\qquad(n\geq 1).
\ee
Choosing $\lc$ such that (\ref{Rex}) holds for some $R>|a|-\de$ yields,
by (\ref{linsplit})
\be
\lim_{t\to\infty}P\big[|\xi_t|>\lc t\big]e^{(|a|-\de)t}=0.
\ee
Inserting this into (\ref{split}) we find that the exponential growth
rate $r=r(\La,a,\de)$ satisfies
\be
r\leq\limsup_{t\to\infty}\frac{1}{t}\log|\{i\in\La:|i|\leq\lc t\}|=0,
\ee
where we have used that $\La$ is subexponential.\qed

\subsection{Duality and Campbell laws}\label{S:char}

\detail{{\bf Proof that $\hat P^A_\lr$ is well-defined} It suffices to
show that $r\leq 0$ implies $E_\lr[|\eta^A_\tau|]<\infty$ for all
$\lr>r$ and $\emptyset\neq A\in\Pc_{\rm fin}(\La)$. To prove this,
observe that by (\ref{expr}), for each $\eps>0$ there exists a
$K<\infty$ such that for all $t\geq T$,
\be
E\big[|\eta^A_t|\big]\leq K\ex{(r+\eps)t}\qquad(t\geq 0).
\ee
In particular, choosing $\eps<1/\lr$ and using that $r\leq 0$, we obtain that
\be\label{welldef}
E_\lr\big[|\eta^A_\tau|\big]
=\int_0^\infty\! E\big[|\eta^A_\tau|\big]\ffrac{1}{\lr}e^{-t/\lr}\,\di t
\leq K\ffrac{1}{\lr}\int_0^\infty\!\ex{(\eps-1/\lr)t}\,\di t<\infty,
\ee
as desired.\qed}

\noi
{\bf Proof of Lemma~\ref{L:char}~(a)} This follows by writing
\be\ba{l}\label{charp}
\dis\hat P^{\{0\}}_\lr\big[A\cap\iota^{-1}\eta^{\{0\}}_\tau=\emptyset\big]
\stackrel{(1)}{=}\dis\expi_\lr(\{0\})^{-1}
\sum_i\int_0^\infty\! P\big[i\in\eta^{\{0\}}_t,
\ A\cap i^{-1}\eta^{\{0\}}_t=\emptyset\big]e^{-\la t}\,\di t\\[5pt]
\qquad\stackrel{(2)}{=}\dis\expi_\lr(\{0\})^{-1}
\sum_i\int_0^\infty\! P\big[0\in\eta^{\{i^{-1}\}}_t,
\ A\cap\eta^{\{i^{-1}\}}_t=\emptyset\big]e^{-\la t}\,\di t\\[5pt]
\qquad\stackrel{(3)}{=}\dis\expi_\lr(\{0\})^{-1}
\sum_j\int_0^\infty\!\Big\{P\big[(A\cup\{0\})\cap\eta^{\{j\}}_t\neq\emptyset\big]
-P\big[A\cap\eta^{\{j\}}_t\neq\emptyset\big]\Big\}e^{-\la t}\,\di t\\[5pt]
\qquad\stackrel{(4)}{=}\dis\expi^\dgg_\lr(\{0\})^{-1}
\sum_j\int_0^\infty\!\Big\{P\big[\eta^{\dgg\,A\cup\{0\}}_t\cap\{j\}\neq\emptyset\big]
-P\big[\eta^{\dgg\,A}_t\cap\{j\}\neq\emptyset\big]\Big\}e^{-\la t}\,\di t\\[5pt]
\qquad\stackrel{(5)}{=}\dis\expi^\dgg_\lr(\{0\})^{-1}
\int_0^\infty\!\Big\{E\big[|\eta^{\dgg\,A\cup\{0\}}_t|\big]
-E\big[|\eta^{\dgg\,A}_t|\big]\Big\}e^{-\la t}\,\di t\\[5pt]
\qquad\stackrel{(6)}{=}\dis\expi^\dgg_\lr(\{0\})^{-1}
\big\{\expi^\dgg_\lr(A\cup\{0\})-\expi^\dgg_\lr(A)\big\}
\stackrel{(7)}{=}\ov\expi^\dgg_\lr(A\cup\{0\})-\ov\expi^\dgg_\lr(A).
\ec
Here, in step (2) we have used shift invariance, in step (3) we have
changed the summation order and used that $\{0\in\eta^{\{j\}}_\tau,
\ A\cap\eta^{\{j\}}_\tau=\emptyset\}
=\{(A\cup\{0\})\cap\eta^{\{j\}}_\tau\neq\emptyset\}
\beh\{A\cap\eta^{\{j\}}_\tau\neq\emptyset\}$,
and in step (4) we have used duality (formula (\ref{dual})) and
formula (\ref{Esym}).\qed

\noi
{\bf Proof of Lemma~\ref{L:char}~(b)} We have
\bc
\dis P\big[A\cap\ov\eta_0=\emptyset\,\big|\,0\in\ov\eta_0\big]
&\stackrel{(1)}{=}&\dis P\big[0\in\ov\eta_0\big]^{-1}
P\big[0\in\ov\eta_0,\ A\cap\ov\eta_0=\emptyset\big]\\[5pt]
&\stackrel{(2)}{=}&\dis P\big[\{0\}\cap\ov\eta_0\neq\emptyset\big]^{-1}
\big\{P\big[(A\cup\{0\})\cap\ov\eta_0\neq\emptyset\big]
-P\big[A\cap\ov\eta_0\neq\emptyset\big]\big\}\\[5pt]
&\stackrel{(3)}{=}&\dis\rho^\dgg(\{0\})\big\{\rho^\dgg(A\cup\{0\})
-\rho^\dgg(A)\big\}\stackrel{(4)}{=}\ov\rho^\dgg(A\cup\{0\})-\ov\rho^\dgg(A),
\ec
where in step (3) we have used (\ref{nuchar}).\qed

\noi
As a preparation for the proof of Lemma~\ref{L:typtime}, we prove:
\bl{\bf(Expected population size)}\label{L:size}
One has $\lim_{\la\down r}\pi_\la(A)=\infty$ for all $\emptyset\neq A\in\Pc_{\rm fin}(\La)$.
\el
{\bf Proof} We start with the case $A=\{0\}$. Recall that Proposition~\ref{P:rate}~(a) is a consequence of the subadditivity of the function $t\mapsto\log E[|\eta^{\{0\}}_t|]$. In fact, subadditivity gives us a little more. By \cite[Theorem~B.22]{Lig99},
\be\label{expr2}
\lim_{t\to\infty}\,\ffrac{1}{t}\log E\big[|\eta^{\{0\}}_t|\big]=\inf_{t>0}\ffrac{1}{t}\log E\big[|\eta^{\{0\}}_t|\big]=r,
\ee
where $r=r(\La,a,\de)\in[-\de,|a|-\de]$ is the exponential growth
rate. Formula (\ref{expr2}) says that $E[|\eta^{\{0\}}_t|]=e^{r_tt}$
where $\lim_{t\to\infty}r_t=\inf_{t>0}r_t=r$. Thus, for every
$\eps>0$, there exists a $T_\eps<\infty$ such that
\be\label{twosid}
e^{rt}\leq E\big[|\eta^{\{0\}}_t|\big]\leq e^{(r+\eps)t}\qquad(t\geq T_\eps).
\ee
It follows from the lower bound in (\ref{twosid}) and monotone convergence that
\be
\lim_{\lr\down r}\expi_\lr(\{0\})=\int_0^\infty\!\! E\big[|\eta^{\{0\}}_t|\big]\,e^{-rt}\,\di t=\infty.
\ee
The generalization to arbitrary $\emptyset\neq A\in\Pc_{\rm fin}(\La)$ is immediate, since $\expi_\la$ is monotone.\qed

\noi
{\bf Proof of Lemma~\ref{L:typtime}} By Lemma~\ref{L:size},
\be
\hat P^{A}_\lr\big[\tau<t\big]=\frac{\int_0^t\!\! E\big[|\eta^{A}_s|\big]\,e^{-\lr s}\,\di s}{\int_0^\infty\!\! E\big[|\eta^{A}_s|\big]\,e^{-\lr s}\,\di s}\leq\frac{\int_0^t\!\! E\big[|\eta^{A}_s|\big]\,e^{-rs}\,\di s}{\expi_\la(A)}\asdto{\lr}{r}0.
\ee
for any $t>0$.\qed

\subsection{Harmonic functions}\label{S:harm}

In this section we prove formulas (\ref{rhoharm}) and (\ref{Gip}), and
Lemma~\ref{L:piclus}.\med

\noi
{\bf Proof of (\ref{rhoharm})} The shift invariance and monotonicity
of $\rho$ follow from the corresponding properties of the contact
process. Since $\rho$ is bounded, obviously $\rho\in\Si(\Pc_{\rm
fin}(\La))$. Since $\eta^A$ solves the martingale problem for $G$, for
any $f\in\Si(\Pc_{\rm fin}(\La))$, one has
\be
\int_0^tE[Gf(\eta^A_s)]\di s=E[f(\eta^A_t)]-f(A)\qquad(A\in\Pc_{\rm fin}(\La)),
\ee
and therefore
\be\label{Gf}
Gf(A)=\lim_{t\to 0}t^{-1}\big\{E[f(\eta^A_t)]-f(A)\big\}
\qquad(A\in\Pc_{\rm fin}(\La)).
\ee
By the Markov property,
\be
\rho(\eta^A_t)=E\big[\eta^A_s\neq 0\ \forall s\geq 0\,\big|\,\eta^A_t\big]
=E\big[\eta^A_s\neq 0\ \forall s\geq 0\,\big|\,\Fi^A_t\big],
\ee
where $(\Fi^A_t)_{t\geq 0}$ denotes the filtration generated by $\eta^A$.
It follows that $\rho(\eta^A_t)$ is a martingale, and therefore, by
(\ref{Gf}), $G\rho=0$.\qed

\noi
{\bf Proof of (\ref{Gip})} The shift invariance and monotonicity of
$\expi_\lr$ follow from the corresponding properties of the contact
process. It follows from (\ref{addit}) that
$\expi_\lr(A)\leq\expi_\lr(\{0\})|A|$, which shows that
$\expi_\lr\in\Si(\Pc_{\rm fin}(\La))$. Moreover,
\detail{Recall (\ref{welldef}), which shows that all integrals below are finite.}
\be\ba{l}
\dis t^{-1}\big\{E[\expi_\lr(\eta^A_t)]-\expi_\lr(A)\big\}\\[5pt]
\dis\quad= t^{-1}\int_0^\infty\!\Big\{E\big[|\eta^A_{t+s}|\big]
-E\big[|\eta^A_s|\big]\Big\}e^{-\la s}\,\di s\\[5pt]
\dis\quad= t^{-1}\Big\{\int_t^\infty\! E\big[|\eta^A_s|\big]
e^{-\lr(s-t)}\,\di s-\int_0^\infty\! E\big[|\eta^A_s|\big]
e^{-\lr s}\,\di s\Big\}\\[5pt]
\dis\quad= t^{-1}(e^{\lr t}-1)\int_0^\infty\! E\big[|\eta^A_s|\big]
e^{-\lr s}\,\di s
-e^{\lr t}t^{-1}\int_0^t\! E\big[|\eta^A_s|\big]
e^{-\lr s}\,\di s.
\ec
Letting $t\to 0$, using (\ref{Gf}), it follows that
\be
G\expi_\lr(A)=\lr\expi_\lr(A)-|A|
\qquad(A\in\Pc_{\rm fin}(\La),\ \lr>r),
\ee
as desired.\qed

\noi
{\bf Proof of Lemma~\ref{L:piclus}} It follows from (\ref{nutchar})
that $\ov\expi_\lr(A)\leq|A|$, which shows that the functions
$(\ov\expi_\lr)_{\lr>r}$ are relatively compact, and each pointwise
limit $\ov\expi_\infty$ along a sequence $\lr_n\down r$ satisfies
$\ov\expi_\infty\in\Si(\Pc_{\rm fin}(\La))$. Since each
$\ov\expi_{\lr_n}$ is shift invariant an monotone, the same is true
for $\ov\expi_\infty$. If $f_n,f\in\Si(\Pc_{\rm fin}(\La))$, $f_n\to
f$ pointwise, and the $f_n$ are uniformly bounded on sets of the form
$\{A\in\Pc_{\rm fin}(\La):|A|\leq K\}$, then it is not hard to see
that pointwise
\be
\lim_{n\to\infty}Gf_n=Gf.
\ee
Applying this to the functions $\ov\expi_{\lr_n}$, which satisfy the
uniform bound $\ov\expi_{\lr_n}(A)\leq|A|$, using (\ref{Gip}) and
Lemma~\ref{L:size}, we find that
\be
G\ov\expi_\infty(A)=\lim_{n\to\infty}
\frac{\lr_n\expi_{\lr_n}(A)-|A|}{\expi_{\lr_n}(\{0\})}
=\lim_{n\to\infty}\lr_n\ov\expi_{\lr_n}(A)
-\frac{|A|}{\expi_{\lr_n}(\{0\})}=r\ov\expi_r(A)
\qquad(A\in\Pc_{\rm fin}(\La)),
\ee
as required.\qed

\subsection{Eventual domination of finite configurations}\label{S:dom}

In this section we prove Lemma~\ref{L:evdom}. We start with two
preparatory lemmas.
\bl{\bf(Local creation of finite configurations)}\label{L:loc}
For each $B\in\Pc_{\rm fin}(\La)$ and $t>0$, there exists a finite
$\De\sub\La$ and $j\in\La$ such that
\be\label{inDe}
\eps:=P\big[\eta^{\{0\}}_t\supset jB\mbox{ and }
\eta^{\{0\}}_s\sub\De\ \forall 0\leq s\leq t\big]>0.
\ee
\el
{\bf Proof} It follows from assumption~(\ref{assum})~(iii) that there
exists a site $j^{-1}\in\La$ with $P\big[\eta^{\{j^{-1}\}}_t\supset B]>0$,
and therefore $P\big[\eta^{\{0\}}_t\supset jB]>0$. Since
$\bigcup_{0\leq s\leq t}\eta^{\{0\}}_s$ is a.s.\ finite, we can choose
a finite but large enough $\De$ such that (\ref{inDe}) holds.\qed

\bl{\bf(Domination of finite configurations)}\label{L:dom}
For each $B\in\Pc_{\rm fin}(\La)$, $t>0$, and $A_n\in\Pc_{\rm fin}(\La)$
satisfying $\lim_{n\to\infty}|A_n|=\infty$, one has
\be
\lim_{n\to\infty}P[\exists i\in\La\mbox{ s.t.\ }\eta^{A_n}_t\geq iB]=1.
\ee
\el
{\bf Proof} Let $\De$, $j$, and $\eps$ be as in Lemma~\ref{L:loc}. We can
find $\ti A_n\sub A_n$ such that $|\ti A_n|\to\infty$ as $n\to\infty$, and
for fixed $n$, the sets $(k\De)_{k\in\ti A_n}$ are disjoint. It follows that
\be\ba{l}
P[\exists i\in\La\mbox{ s.t.\ }\eta^{A_n}_t\geq iB]\\[5pt]
\dis\qquad\geq1-\prod_{k\in\ti A_n}\big(1-P\big[\eta^{\{k\}}_t\supset kjB
\mbox{ and }\eta^{\{k\}}_s\sub k\De\ \forall 0\leq s\leq t\big]\big)\\[5pt]
\dis\qquad=1-(1-\eps)^{|\ti A_n|}\asto{n}1,
\ec
where we have used (\ref{inDe}) and the fact that events concerning the
graphical representation in disjoint parts of space are independent.\qed

\noi
{\bf Proof of Lemma~\ref{L:evdom}} If $\de=0$, then obviously
$\lim_{t\to\infty}|\eta^A_t|=\infty$ a.s. If $\de>0$, then it is easy
to see that $\inf\{\rho(A):|A|\leq M\}<1$ for all $M<\infty$. Therefore,
by (\ref{loe}),
\be\label{exgroiv}
\eta^A_t=\emptyset\mbox{ for some }t\geq 0\quad\mbox{or}\quad|\eta^A_t|
\asto{t}\infty\qquad{\rm a.s.}
\ee
Fix $\emptyset\neq B\in\Pc_{\rm fin}(\La)$ and set
$\psi_t(A):=P[\exists i\in\La\mbox{ s.t.\ }\eta^A_t\geq iB]
\quad(A\in\Pc_{\rm fin}(\La),\ t\geq 0)$. Then, for each $t>0$,
\be
\lim_{T\to\infty}P[\exists i\in\La\mbox{ s.t.\ }\eta^A_T\supset iB]
=\lim_{T\to\infty}E[\psi_t(\eta^A_{T-t})]=\rho(A),
\ee
where we have used Lemma~\ref{L:dom} and (\ref{exgroiv}).\qed

\subsection{Generalization to arbitrary initial states}\label{S:A}

In this section, we show how the proof of Theorem~\ref{T:palm}~(a) must
be adapted to cover general initial states $\emptyset\neq A\in\Pc(\La)$.\med

\noi
{\bf Proof of Theorem~\ref{T:palm}~(a) for general initial states} For
$A,B\in\Pc_{\rm fin}(\La)$ with $A\neq\emptyset$, we observe that
$i\in BA^{-1}\desd B\cap iA\neq\emptyset$,
\detail{Indeed, $\exists a\in A,\ b\in B$ s.t.\ $i=ba^{-1}
\desd\exists a\in A,\ b\in B$ s.t.\ $ia=b\desd$.}
 and therefore
\be\label{BAinv}
\dis|BA^{-1}|=\sum_i1_{\txt\{B\cap iA\neq\emptyset\}}.
\ee
We define
\be\label{piAB}
\expi_{A,\lr}(B):=\int_0^\infty\! E\big[|\eta_tA^{-1}|\big]e^{-\lr t}\,\di t\quad\mbox{and}
\quad\ov\expi_{A,\lr}(B):=\frac{\expi_{A,\lr}(B)}{\expi_{A,\lr}(\{0\})},
\ee
and let $\expi^\dgg_{A,\lr}$ and $\ov\expi^\dgg_{A,\lr}$ denote the
analogues of $\expi_{A,\lr}$ and $\ov\expi_{A,\lr}$ for the
$(\La,a^\dgg,\de)$-contact process. Generalizing the proof of
Lemma~\ref{L:char}~(a), we find that
\be\label{charp2}
\hat P^A_\lr\big[B\cap\iota^{-1}\eta^A_\tau=\emptyset\big]
=\ov\expi^\dgg_{A,\lr}(B\cup\{0\})-\ov\expi^\dgg_{A,\lr}(B).
\ee
\detail{Indeed,
\be\ba{l}
\dis\hat P^A_\lr\big[B\cap\iota^{-1}\eta^A_\tau=\emptyset\big]
=\dis E_\lr\big[|\eta^A_\tau|\big]^{-1}
\sum_i P_\lr\big[i\in\eta^A_\tau,\ B\cap i^{-1}\eta^A_\tau=\emptyset\big]\\[5pt]
\qquad=\dis E_\lr\big[|\eta^A_\tau|\big]^{-1}
\sum_i P_\lr\big[0\in\eta^{i^{-1}A}_\tau,\ B\cap\eta^{i^{-1}A}_\tau=\emptyset\big]\\[5pt]
\qquad=\dis E_\lr\big[|\eta^A_\tau|\big]^{-1}
\sum_j\Big\{P_\lr\big[(B\cup\{0\})\cap\eta^{jA}_\tau\neq\emptyset\big]
-P_\lr\big[B\cap\eta^{jA}_\tau\neq\emptyset\big]\Big\}\\[5pt]
\qquad=\dis E_\lr\big[|\eta^A_\tau|\big]^{-1}
\sum_j\Big\{P_\lr\big[\eta^{\dgg\,B\cup\{0\}}_\tau\cap jA\neq\emptyset\big]
-P_\lr\big[\eta^{\dgg\,B}_\tau\cap jA\neq\emptyset\big]\Big\}\\[5pt]
\qquad=\dis\expi^\dgg_{A,\lr}(\{0\})^{-1}
\big\{\expi^\dgg_{A,\lr}(B\cup\{0\})-\expi^\dgg_{A,\lr}(B)\big\}
=\ov\expi^\dgg_{A,\lr}(B\cup\{0\})-\ov\expi^\dgg_{A,\lr}(B),
\ec
where we have used (\ref{BAinv}) and the fact that, in analogy with
(\ref{Esym}), for any $t\geq 0$,
\bc\label{Esym2}
\dis E\big[|\eta^{\dgg\,\{0\}}_tA^{-1}|\big]
&=&\dis\sum_jP\big[\eta^{\dgg\,\{0\}}_t\cap jA\neq\emptyset\big]
=\sum_jP\big[\{0\}\cap\eta^{jA}_t\neq\emptyset\big]\\[5pt]
&=&\dis\sum_jP\big[\{j^{-1}\}\cap\eta^{A}_t\neq\emptyset\big]
=E\big[|\eta^{A}_t|\big].
\ec}
Since $|B|\leq|BA^{-1}|\leq|A|\,|B|$ for any $A,B\in\Pc_{\rm fin}(\La)$
with $A\neq\emptyset$, it follows that
\be
\lim_{t\to\infty}\ffrac{1}{t}\log E\big[|\eta^B_tA^{-1}|\big]=r
\qquad(\emptyset\neq B\in\Pc_{\rm fin}(\La)),
\ee
where $r$ is the exponential growth rate from
Proposition~\ref{P:rate}. The proofs of (\ref{Gip}) and
Lemma~\ref{L:piclus} now carry over to the functions
$(\ov\expi_{A,\lr})_{\lr>r}$ without a change, and therefore the
arguments in Section~\ref{S:methods} show that
Theorem~\ref{T:palm}~(a) holds for general initial states
$\emptyset\neq A\in\Pc(\La)$.\qed

\section{Proofs of further results}\label{S:rest}

Recall that $\poi=(\poi^{\rm r},\poi^{\rm i})$ is the pair of Poisson
point processes used in the graphical representation. We construct
$\poi$ on the canonical probability space $\Poi:=\Pc_{\rm
loc}(\La\times\R)\times\Pc_{\rm loc}(\La\times\La\times\R)$, where
$\Pc_{\rm loc}(\La\times\R)$ and $\Pc_{\rm loc}(\La\times\La\times\R)$
denote the spaces of locally finite subsets of $\La\times\R$ and
$\La\times\La\times\R$, respectively. These spaces can in a natural
way be identified with subspaces of the spaces of locally finite
counting measures on $\La\times\R$ and $\La\times\La\times\R$,
respectively. Using this identification, we equip $\Pc_{\rm loc}
(\La\times\R)$ and $\Pc_{\rm loc}(\La\times\La\times\R)$ with the
vague topology. We equip $\Poi$ with the product topology and the
associated Borel-\si-field $\Fi$, and let $P$ be the probability
measure on $(\om,\Fi)$ such that under $P$, the coordinate functions
$\poi^{\rm r},\poi^{\rm i}$ are Poisson point processes as described
in the introduction.

We equip $\La\times\R$ and $\La\times\La\times\R$ with a group
structure by putting $(i,s)(j,t):=(ij,s+t)$ and
$(i,j,s)(k,l,t):=(ik,jl,s+t)$, respectively. In line with our earlier
notation, for any subset $\al\sub\La\times\R$, we write
$(i,s)\al:=\{(ij,s+t):(j,t)\in\al\}$. For
$\bet\sub\La\times\La\times\R$, we define $(i,j,s)\bet$
analogously. We define shift operators $\tet_{i,t}:\Poi\to\Poi$ by
\be
\tet_{i,t}(\al,\bet):=((i,t)\al,(i,i,t)\bet)
\ee
$(i\in\La,\ t\in\R,\ (\al,\bet)\in\Poi)$. Thus, $\tet_{i,t}$ shifts a
graphical representation by left-multiplication with $i$ and increasing
all times by $t$.

\subsection{Conditioning and size-biasing}\label{S:sizco}

In this section, we prove that size-biasing and conditioning on survival
are asymptotically equivalent in a `local' sense. Let
\be
\poi_t:=(\poi^{\rm r}\cap\La\times(-\infty,t]\,,
\,\poi^{\rm i}\cap\La\times\La\times(-\infty,t]\big)
\ee
denote the restriction of the Poisson point processes used in the
graphical representation to the time interval $(-\infty,t]$.
\bp{\bf(Conditioning and size-biasing)}\label{P:sizco}
Assume that the $(\La,a,\de)$-contact process survives and that the
exponential growth rate satisfies $r(\La,a,\de)=0$. Then, for any
$\emptyset\neq A\in\Pc_{\rm fin}(\La)$,
\be\label{aseq}
\hat P^A_\lr\big[\poi_t\in\,\cdot\,\big]
\asdto{\lr}{0}P\big[\poi_t\in\,\cdot\,\,\big|\,A\times\{0\}\leadsto\infty\big]
\qquad(t\in\R).
\ee
\ep
{\bf Proof} It suffices to prove the claims for $t>0$. For any $\Ai\in\Fi$, write
\be\label{decomp}
\hat P^A_\lr\big[\poi_t\in\Ai\big]
=\hat P^A_\lr\big[\poi_t\in\Ai\,\big|\,\tau\geq t]
\,\hat P^A_\lr\big[\tau\geq t\big]+\hat P^A_\lr\big[\poi_t\in\Ai,\ \tau<t\big],
\ee
and observe that
\be\ba{l}\label{pir}
\dis\hat P^A_\lr\big[\poi_t\in\Ai\,\big|\,\tau\geq t]
=\frac{\int_0^\infty\!E[|\eta^A_{t+s}|1_{\{\poi_t\in\Ai\}}]e^{-\lr s}\,\di s}{\int_0^\infty\!E[|\eta^A_{t+s}|]e^{-\lr s}\,\di s}
=\frac{E[\int_0^\infty\!E[|\eta^A_{t+s}|\,|\,\poi_t]e^{-\lr s}\,\di s1_{\{\poi_t\in\Ai\}}]}
{E[\int_0^\infty\!E[|\eta^A_{t+s}|\,|\,\poi_t]e^{-\lr s}\,\di s]}\\[12pt]
\dis\quad=\frac{E[\expi_\lr(\eta^A_t)1_{\{\poi_t\in\Ai\}}]}
{E[\expi_\lr(\eta^A_t)]}
=\frac{E[\ov\expi_\lr(\eta^A_t)1_{\{\poi_t\in\Ai\}}]}
{E[\ov\expi_\lr(\eta^A_t)]}\asdto{\lr}{0}
\frac{E[\ov\rho(\eta^A_t)1_{\{\poi_t\in\Ai\}}]}{E[\ov\rho(\eta^A_t)]},
\ec
where we have used that $\ov\expi_\lr\to\ov\rho$ pointwise as $\lr\down0$
by Lemma~\ref{L:piclus} and Proposition~\ref{P:monhar}, and bounded
convergence, using the uniform bound $\expi_\lr\leq|\,\cdot\,|$. Since
\be\ba{l}\label{ident}
\dis\frac{E[\ov\rho(\eta^A_t)1_{\{\poi_t\in\Ai\}}]}{E[\ov\rho(\eta^A_t)]}
=\frac{E[\rho(\eta^A_t)1_{\{\poi_t\in\Ai\}}]}{E[\rho(\eta^A_t)]}\\[12pt]
\dis\qquad=\frac{E[P[A\times\{0\}\leadsto\infty\,|\,\poi_t]
1_{\{\poi_t\in\Ai\}}]}{E[P[A\times\{0\}\leadsto\infty\,|\,\poi_t]]}
=P\big[\poi_t\in\Ai\,\big|\,A\times\{0\}\leadsto\infty\big],
\ec
formula (\ref{aseq}) follows from Lemma~\ref{L:typtime},
(\ref{decomp}), and (\ref{pir}).\qed

\subsection{Coupling to the maximal process}\label{S:maxcoup}

In this section we prove Theorem~\ref{T:palm}~(b) and Proposition~\ref{P:totyp}. In analogy with (\ref{Campdefiv}), we put
\be\label{Campdef2}
\hat P^{\dgg\,A}_\lr(\{i\}\times\{\di\oo\}\times\{\di t\}):=
\expi^\dgg_\la(A)^{-1}\,1_{\txt\{i\in\eta^{\dgg\,A}_t(\oo)\}}P(\di\oo)e^{-\la t}\di t,
\ee
which is well-defined for any $\emptyset\neq A\in\Pc_{\rm fin}(\La)$
and $\lr>r$. Recall that
$\eta^{\dgg\,A}_\tau=\{i\in\La:(i,-\tau)\leadsto A\times\{0\}\}$. We
can view $\eta^{\dgg\,A}_t$ as the set of all `ancestors' at time $-t$
of the set $A$ at time $0$. As before, let $\iota$ and $\tau$ denote
the projections on $\La$ and $\R_+$, respectively. Then, under the law
$\hat P^{\dgg\,A}_\lr$, the random variables $\iota$ and $\tau$
describe a `typical' ancestor of $A$ and a `typical' time $-\tau$.

In the next lemma, we shift the graphical representation $\poi$ in
such a way that the `typical' infected site and time $(\iota,\tau)$,
chosen with respect to $\hat P^{\{0\}}_\lr$, are mapped to the point
$(0,0)$. Note that under such a shift, the origin is mapped to
$\iota^{-1}$. Thus, the next lemma can be described by saying that if
we start the contact process with only the origin infected, then seen
{f}rom a typical infected site, the origin is a typical ancestor.
\bl{\bf(Origin seen from a typical infected site)}\label{L:orse}
Assume that $r(\La,a,\de)\leq 0$. Then
\be\label{typanc}
\hat P^{\{0\}}_\lr\big[(\iota^{-1},\tet_{\iota^{-1},-\tau}\poi,\tau)
\in\cdot\,\big]=\hat P^{\dgg\,\{0\}}_\lr\big[(\iota,\poi,\tau)\in\cdot\,\big].
\ee
\el
{\bf Proof} Let us write $(i,s)\stackrel{\poi}{\leadsto}(j,t)$ when
$(i,s)$ can be connected to $(j,t)$ along a path in the graphical
representation $\poi$. Then
\be\ba{l}
\dis\hat P^{\{0\}}_\lr\big[\iota^{-1}=j,
\ \tet_{\iota^{-1},-\tau}\poi\in\Ai,\ \tau\in(a,b)\big]
=\hat P^{\{0\}}_\lr\big[\iota=j^{-1},
\ \tet_{\iota^{-1},-\tau}\poi\in\Ai,\ \tau\in(a,b)\big]\\[5pt]
\dis\quad=\expi_\lr(\{0\})^{-1}
\int_a^b\!P\big[j^{-1}\in\eta^{\{0\}}_t,
\ \tet_{j,-t}\poi\in\Ai\big]e^{-\lr t}\,\di t\\[5pt]
\dis\quad=\expi_\lr(\{0\})^{-1}
\int_a^b\!P\big[(0,0)\stackrel{\poi}{\leadsto}(j^{-1},t),
\ \tet_{j,-t}\poi\in\Ai\big]e^{-\lr t}\,\di t\\[5pt]
\dis\quad=\expi_\lr(\{0\})^{-1}
\int_a^b\!P\big[(j,-t)\stackrel{\tet_{j,-t}\poi}{\leadsto}(0,0),
\ \tet_{j,-t}\poi\in\Ai\big]e^{-\lr t}\,\di t\\[5pt]
\dis\quad=\expi_\lr(\{0\})^{-1}
\int_a^b\!P\big[(j,-t)\stackrel{\poi}{\leadsto}(0,0),
\ \poi\in\Ai\big]e^{-\lr t}\,\di t\\[5pt]
\dis\quad=\expi^\dgg_\lr(\{0\})^{-1}
\int_a^b\!P\big[j\in\eta^{\dgg\,\{0\}}_t,\ \poi\in\Ai\big]e^{-\lr t}\,\di t
=\hat P^{\dgg\,\{0\}}_\lr\big[\iota=j,\ \poi\in\Ai,\ \tau\in(a,b)\big],
\ec
where we have used (\ref{Esym}).\qed

\noi
In order to prove Theorem~\ref{T:palm}~(b), we need two more lemmas.

\bl{\bf(Large populations)}\label{L:lart}
Assume that the $(\La,a,\de)$-contact process survives and that the
exponential growth rate satisfies $r(\La,a,\de)\leq 0$. Then, for any
$\emptyset\neq A\in\Pc_{\rm fin}(\La)$,
\be\label{toinf}
\hat P^A_\lr\big[|\eta^A_\tau|\geq K\big]\asdto{\lr}{0}1\qquad(K<\infty).
\ee
\el
{\bf Proof} Let $\tau_\lr$ be an exponentially distributed reandom variable with mean $1/\lr$, independent of the Poisson processes used in the graphical representation. Then
\bc
\dis\hat P^A_\lr\big[|\eta^A_\tau|\geq K\big]
&=&\dis\frac{E\big[|\eta^A_{\tau_\lr}|1_{\txt\{|\eta^A_{\tau_\lr}|\geq K\}}\big]}
{E\big[|\eta^A_{\tau_\lr}|\big]}\\[5pt]
&=&\dis\frac{E\big[|\eta^A_{\tau_\lr}|1_{\txt\{|\eta^A_{\tau_\lr}|\geq K\}}\,\big|
\,\eta^A_{\tau_\lr}\neq\emptyset\big]}{E\big[|\eta^A_{\tau_\lr}|\,\big|
\,\eta^A_{\tau_\lr}\neq\emptyset\big]}
\geq E\big[1_{\txt\{|\eta^A_{\tau_\lr}|\geq K\}}\,\big|
\,\eta^A_{\tau_\lr}\neq\emptyset\big]\asdto{\lr}{0}1,
\ec
where we have used (\ref{exgroiv}), and the fact that $|\eta^A_{\tau_\lr}|$ and
$1_{\{\eta^A_{\tau_\lr}\geq K\}}$ are positively correlated since the functions
$z\mapsto z$ and $z\mapsto 1_{\{z\geq K\}}$ are nondecreasing.\qed

\detail{Indeed, if  $f,g:\N\to\R$  are nondecreasing, then they are positively correlated with respect to any probability measure on $\N$.}

\noi
Recall that in the proof (in Section~\ref{S:methods}) of
Proposition~\ref{P:monhar}, sequences $A_n\in\Pc_{\rm fin}(\La)$ such
that $\rho(A_n)\to 1$ played an important role. Although we did not
need this fact there, the next lemma implies that for $\de>0$, actually
$\rho(A_n)\to 1$ if and only if $|A_n|\to\infty$.
\bl{\bf(High survival probabilities)}\label{L:rhone}
Assume that the $(\La,a,\de)$-contact process survives, and
$A_n\in\Pc_{\rm fin}(\La)$. Then $|A_n|\to\infty$ implies $\rho(A_n)\to 1$. 
\el
{\bf Proof} By (\ref{loe}) there exist $B_m\in\Pc_{\rm fin}(\La)$
with $\rho(B_m)\to 1$. Now if $A_n\in\Pc_{\rm fin}(\La)$ satisfy
$|A_n|\to\infty$, then by Lemma~\ref{L:dom},
\be\ba{l}
\dis\liminf_{n\to\infty}\rho(A_n)\\[5pt]
\dis\quad\geq\liminf_{n\to\infty}P\big[\eta^{A_n}_s\neq\emptyset
\ \forall s\geq t\,\big|\,\exists i\in\La\mbox{ s.t.\ }
\eta^{A_n}_t\geq iB_m\big]P\big[\exists i\in\La\mbox{ s.t.\ }
\eta^{A_n}_t\geq iB_m\big]\\[5pt]
\dis\qquad\geq\rho(B_m),
\ec
for each $t>0$ and $m$. Letting $m\to\infty$ yields the claim.\qed

\noi
We now first prove Theorem~\ref{T:palm}~(b) in the case $A=\{0\}$, and
then indicate how the arguments may be generalised to $\emptyset\neq
A\in\Pc_{\rm fin}(\La)$. We will obtain Proposition~\ref{P:totyp} as a
corollary to our proofs in the case $A=\{0\}$.\med

\noi
{\bf\boldmath Proof of Theorem~\ref{T:palm}~(b) in the case $A=\{0\}$}
By Lemma~\ref{L:orse}, we must show that for fixed $\De\in\Pc_{\rm fin}$,
the sets $\{j\in\De:(\iota,-\tau)\leadsto(j,0)\}$,
$\{j\in\De:-\infty\leadsto(j,0)\}$, and
$\{j\in\De:\La\times\{-\tau\}\leadsto(j,0)\}$ are asymptotically equal
under the laws $\hat P^{\dgg\,\{0\}}_\lr$ as $\lr\down 0$. It
suffices to show that for any $j\in\La$,
\be\label{step1}
\hat P^{\dgg\,\{0\}}_\lr\big[(\iota,-\tau)\leadsto(j,0)\,\big|
\,\La\times\{-\tau\}\leadsto(j,0)\big]\asdto{\lr}{0}1.
\ee
and
\be\label{step2}
\hat P^{\dgg\,\{0\}}_\lr\big[-\infty\leadsto(j,0)\,\big|
\,\La\times\{-\tau\}\leadsto(j,0)\big]\asdto{\lr}{0}1.
\ee
Reversing the direction of time and interchanging the roles of $\eta$
and $\eta^\dgg$, this then yields Proposition~\ref{P:totyp} as a corollary.

For any $t>0$, by Proposition~\ref{P:sizco},
\be\ba{l}
\dis\hat P^{\dgg\,\{0\}}_\lr\big[\La\times\{-\tau\}\leadsto(j,0)\big]
=\hat P^{\dgg\,\{0\}}_\lr\big[\eta^{\dgg\,\{j\}}_\tau\neq\emptyset\big]\\[5pt]
\dis\quad\leq\hat P^{\dgg\,\{0\}}_\lr\big[\eta^{\dgg\,\{j\}}_t
\neq\emptyset\big]+\hat P^{\dgg\,\{0\}}_\lr\big[\tau<t\big]
\asdto{\lr}{0}P\big[\eta^{\dgg\,\{j\}}_t\neq\emptyset\,\big|
\,-\infty\leadsto(0,0)\big].
\ec
Letting $t\to\infty$ yields
\be\label{step1a}
\limsup_{\lr\down 0}\hat P^{\dgg\,\{0\}}_\lr\big[\La\times\{-\tau\}
\leadsto(j,0)\big]\leq P\big[-\infty\leadsto(j,0)\,\big|
\,-\infty\leadsto(0,0)\big]=:\phi(j).
\ee
By Lemma~\ref{L:orse} and Theorem~\ref{T:palm}~(a),
\be\label{step1b}
\lim_{\lr\down 0}\hat P^{\dgg\,\{0\}}_\lr\big[(\iota,-\tau)
\leadsto(j,0)\big]=\lim_{\lr\down 0}
\hat P^{\{0\}}_\lr\big[j\in\iota^{-1}\eta^{\{0\}}_\tau\big]
=P\big[j\in\ov\eta_0\,\big|\,0\in\ov\eta_0\big]=\phi(j).
\ee
Combining (\ref{step1a}) and (\ref{step1b}) we arrive at (\ref{step1}).

Since conditional on $\eta^{\dgg\,\{0\}}_\tau$, the typical site $\iota$
is chosen with equal probabilities from the sites in
$\eta^{\dgg\,\{0\}}_\tau$,
\be
\hat P^{\dgg\,\{0\}}_\lr\big[(\iota,-\tau)\leadsto(j,0)\,\big|
\,\La\times\{-\tau\}\leadsto(j,0)\big]
=\hat E^{\dgg\,\{0\}}_\lr\Big[
\frac{|\eta^{\dgg\,\{j\}}_\tau\cap\eta^{\dgg\,\{0\}}_\tau|}
{|\eta^{\dgg\,\{0\}}_\tau|}\,\Big|
\,\La\times\{-\tau\}\leadsto(j,0)\Big].
\ee
Therefore, (\ref{step1}) and Lemma~\ref{L:lart} imply that
\be
\lim_{\lr\down 0}\hat P^{\dgg\,\{0\}}_\lr\big[|\eta^{\dgg\,\{j\}}_\tau|
\geq K\,\big|\,\La\times\{-\tau\}\leadsto(j,0)\big]=1\qquad(K<\infty),
\ee
which by Lemma~\ref{L:rhone} implies (\ref{step2}).\qed

\noi
{\bf Generalization to arbitrary initial states} In analogy with
(\ref{Campdef2}), we define, for any $\emptyset\neq A,B\in\Pc_{\rm fin}(\La)$,
\be\label{Campdef3}
\hat P^{\dgg\,B}_{A,\lr}(\{i\}\times\{\di\oo\}\times\{\di t\}):=
\expi^\dgg_{A,\la}(B)^{-1}\,1_{\txt\{\eta^{\dgg\,B}_t(\oo)\cap iA\neq\emptyset\}}P(\di\oo)e^{-\la t}\di t,
\ee
where $\expi^\dgg_{A,\la}(B)^{-1}$ is defined below (\ref{piAB}). Note
that this is a probability measure by (\ref{BAinv}). As before, let
$\iota$ denote the projection on $\La$. Then, under the law $\hat
P^{\dgg\,B}_{A,\lr}$, the random variable $\iota$ describes a
`typical' site such that $\iota A\times\{-\tau\}\leadsto
B\times\{0\}$. By an obvious analogue of Lemma~\ref{L:orse}, we must
prove the following generalisations of (\ref{step1}) and
(\ref{step2}):
\be\ba{rc}\label{step12}
{\rm (i)}&\dis\hat P^{\dgg\,\{0\}}_{A,\lr}\big[\iota A\times\{-\tau\}
\leadsto(j,0)\,\big|\,\La\times\{-\tau\}\leadsto(j,0)\big]\asdto{\lr}{0}1,\\[5pt]
{\rm (ii)}&\dis\hat P^{\dgg\,\{0\}}_{A,\lr}\big[-\infty
\leadsto(j,0)\,\big|\,\La\times\{-\tau\}\leadsto(j,0)\big]\asdto{\lr}{0}1.
\ec
Define a measure $\ti P^{\dgg\,\{0\}}_{A,\lr}$ on
$\La\times\La\times\om\times\R_+$ by
\be\label{tiPdef}
\ti P^{\dgg\,\{0\}}_{A,\lr}(\{k\}\times\{i\}\times\{\di\oo\}\times\{\di t\}):=
\hat E^{\dgg\,\{0\}}_{A,\lr}\Big[|\eta^{\dgg\,\{0\}}_\tau\cap iA|^{-1}
1_{\txt\{k\in\eta^{\dgg\,\{0\}}_\tau\cap iA\}}1_{\txt\{i\}\times\{\di\oo\}\times\{\di t\}}\Big].
\ee
Let $\kappa,\iota:\La\times\La\times\om\times\R_+\to\La$ denote the
projections on the first and second coordinate, respectively. Then,
under the law $\ti P^{\dgg\,\{0\}}_{A,\lr}$, the random variable
$\kappa$ describes a site chosen with equal probabilities from
$\eta^{\dgg\,\{0\}}_\tau\cap\iota A$. Therefore, in order to prove
(\ref{step12}), it suffices to prove:
\be\ba{rc}\label{step123}
{\rm (i)}&\dis\ti P^{\dgg\,\{0\}}_{A,\lr}\big[(\kappa,-\tau)
\leadsto(j,0)\,\big|\,\La\times\{-\tau\}\leadsto(j,0)\big]\asdto{\lr}{0}1,\\[5pt]
{\rm (ii)}&\dis\ti P^{\dgg\,\{0\}}_{A,\lr}\big[-\infty\leadsto(j,0)\,\big|
\,\La\times\{-\tau\}\leadsto(j,0)\big]\asdto{\lr}{0}1.
\ec
We claim that $\ti P^{\dgg\,\{0\}}_{A,\lr}\big[(\kappa,\poi,\tau)\in\cdot\,]$
has a density with respect to
$\hat P^{\dgg\,\{0\}}_\lr\big[(\iota,\poi,\tau)\in\cdot\,]$ that is
uniformly bounded away from $0$ and $\infty$, and therefore (\ref{step123})
follows from (\ref{step1}) and (\ref{step2}). Indeed, by (\ref{Campdef3})
and (\ref{tiPdef}),
\be\ba{l}
\dis\ti P^{\dgg\,\{0\}}_{A,\lr}(\{k\}\times\La\times\{\di\oo\}\times\{\di t\})\\[5pt]
\dis\quad=\expi^\dgg_{A,\lr}(\{0\})^{-1}
\sum_iE\big[|\eta^{\dgg\,\{0\}}_\tau\cap iA|^{-1}
1_{\txt\{k\in\eta^{\dgg\,\{0\}}_\tau\cap iA\}}
1_{\txt\{\eta^{\dgg\,\{0\}}_\tau\cap iA\neq\emptyset\}}1_{\txt\{\di\oo\}}\big]e^{-\lr t}\,\di t\\[5pt]
\dis\quad=Z\expi^\dgg_\lr(\{0\})^{-1}
E\big[F1_{\txt\{k\in\eta^{\dgg\,\{0\}}_\tau\}}1_{\txt\{\di\oo\}}\big]e^{-\lr t}\,\di t
=\hat E^{\dgg\,\{0\}}_\lr\big[ZF(k)1_{\txt\{k\}\times\{\di\oo\}\times\{\di t\}}\big],
\ec
where $Z:=\expi^\dgg_\lr(\{0\})/\expi^\dgg_{A,\lr}(\{0\})$ satisfies $|A|^{-1}\leq Z\leq 1$ and
\be
F(k):=\sum_i|\eta^{\dgg\,\{0\}}_\tau\cap iA|^{-1}1_{\txt\{k\in iA\}}
=\sum_{i\in kA^{-1}}|\eta^{\dgg\,\{0\}}_\tau\cap iA|^{-1}
\ee
satisfies $1\leq F(k)\leq |A|$.\qed

\subsection{Coupling of one-dimensional processes}\label{S:coup}

{\bf Proof of Lemma~\ref{L:onedim}} For any point $(i,s)$ such that
$(i,s)\leadsto\infty$, set
\be
r_{s,t}(i):=\max\{j\in\Z:(i,s)\leadsto(j,t)\leadsto\infty\}\qquad(t\geq s).
\ee
Then $(r_{s,t}(i))_{t\geq s}$ is the right-most path to infinity
starting at $(i,s)$. By symmetry and the nearest-neighbor property, it
suffices to show that for any $(i,s)$ and $(j,s)$ such that
$(i,s)\leadsto\infty$ and $(j,s)\leadsto\infty$, there exists a $t\geq
s$ such that $r_{s,t}(i)=r_{s,t}(j)$. Imagine that this is not the
case. Then, for any $i\in\Z$ and $s\leq t$, the maximum
\be
R_{s,t}(i):=\max\{j\in\Z:r_{t,u}(j)=r_{s,u}(i)\mbox{ for some }u\geq t\}
\ee
exists. Set
\be
\chi_s:=\{i\in\Z:R_{s,s}(i)=i\}\qquad(s\in\R).
\ee
It is not hard to see that $R_{s,t}$ maps $\Z$ into $\chi_t$ and that
$R_{s,t}:\chi_s\to\chi_t$ is one-to-one. We claim that
$R_{s,t}:\chi_s\to\chi_t$ is with positive probability not surjective
if $s<t$. Indeed, since we are assuming that $\de>0$ or $a(0,1)\wedge
a(1,0)>0$, it is easy to see that with positive probability there
exist $i,j,k\in\chi_t$ with $i<j<k$ such that
\be
\max\{i'\in\Z:(0,s)\leadsto(i',t)\}=i\quad\mbox{and}
\quad\max\{k'\in\Z:(1,s)\leadsto(k',t)\}=k.
\ee
It follows that $R_{s,t}(0)=i$ and $R_{s,t}(1)=k$, and therefore, since
$R_{s,t}$ is monotone, there is no $n\in\Z$ with $R_{s,t}(n)=j$. 

This `obviously' violates stationarity. More formally, fix $s<t$ and
define $f:\Z\times\Z\to\R$ by
\be
f(i,j):=P[i\in\chi_s, j\in\chi_t,\ j=R_{s,t}(i)].
\ee
Then
\be\ba{l}
\dis\sum_jf(0,j)=P[0\in\chi_s, \exists j\in\chi_t\mbox{ s.t.\ }j
=R_{s,t}(0)]\\[5pt]
\dis\quad=P[0\in\chi_s]>P[\exists i\in\chi_s\mbox{ s.t.\ }0\in\chi_t,\ 0
=R_{s,t}(i)]=\sum_if(i,0).
\ec
Since $\sum_jf(0,j)=\sum_jf(-j,0)=\sum_if(i,0)$ (this equality is a
special case of the mass transport principle; see \cite{Hag97},
\cite[Section~3]{BLPS99}, or \cite[Chapter~7]{LP05}), we arrive
at a contradiction.\qed

\subsection{Survival on finitely generated groups}\label{S:surv}

In this section we prove:
\bl{\bf(Survival for low recovery rates)}
If $\La$ is finitely generated, then $\de_{\rm c}>0$.
\el
{\bf Proof} Let $\De$ be a finite generating set for $\La$. Since
$\{i:a(0,i)>0\}$ generates $\La$, there exists a finite subset
$A\sub\{i:a(0,i)>0\}$ that generates $\De$, and thereby all of
$\La$. Therefore, we can find $i_0,i_1,\ldots\in\La$, all different,
such that $\inf_{k\geq 0}a(i_k,i_{k+1})>0$. We will use comparison to
oriented site percolation to show that $P((i_0,0)\leadsto\infty)>0$ if
$\de$ is sufficiently small. Fix $T>0$. Call a point $(n,m)$ with
$n,m\in\N^2$ good if in the grapical representation, in the time
interval $[Tm,T(m+1))$, there is an arrow from $i_{n}$ to $i_{n+1}$
and there are no recoveries in $i_{n}$ and $i_{n+1}$. By choosing $T$
large enough and $\de$ small enough, the probability that a point is
good can be made arbitrarily high, uniformly in $n$. If this
probability is larger than the critical parameter for independent
2-dimensional oriented site percolation, then with positive
probability there is an upward path along good points, and therefore
the contact process survives.\qed


\begin{thebibliography}{BCGH95}

\bibitem[Apo69]{Apo69}
T.M.~Apostol.
{\em Calculus, Vol.~II}.
Wiley, 1969.

\bibitem[AS05]{AS05}
S.R.~Athreya and J.M.~Swart.
Branching-coalescing particle systems.
{\em Prob.\ Theory Relat.\ Fields.}~131(3), 376--414, 2005.

\bibitem[BCGH95]{BCGH95}
J.-B.~Baillon, Ph.~Cl{\'e}ment, A.~Greven, and F.~den Hollander.
On the attracting orbit of a non-linear transformation arising from
  renormalization of hierarchically interacting diffusions.\ I.\ The compact
  case.
{\em Canad.\ J.\ Math.}~47(1):\ 3--27, 1995.

\bibitem[BCGH97]{BCGH97}
J.-B.~Baillon, Ph.~Cl{\'e}ment, A.~Greven, and F.~den Hollander.
On the attracting orbit of a non-linear transformation arising from
  renormalization of hierarchically interacting diffusions. II.\ The
  non-compact case.
{\em J.\ Funct.\ Anal.}~146:\ 236--298, 1997.

\bibitem[BEM03]{BEM03}
J.~Blath, A.M.~Etheridge, and M.E.~Meredith.
Coexistence in locally regulated competing populations.
University of Oxford, preprint, 2003.

\bibitem[BES04]{BES02}
N.H.~Barton, A.M.~Etheridge, and A.K.~Sturm.
Coalescence in a random background.
{\em Ann. Appl. Probab.}~14(2), 754--785, 2004.

\bibitem[BG85]{BG85}
C.~Bennett and G.~Grinstein.
Role of irreversibility in stabilizing complex and nonergodic behavior in local interacting discrete systems.
{\em Phys.\ Rev.\ Lett.}~55, 657--660, 1985.
 
\bibitem[BG90]{BG90}
C.~Bezuidenhout and G.~Grimmett.
The critical contact process dies out.
{\em Ann.\ Probab.}~18(4), 1462--1482, 1990.

\bibitem[BGN91]{BGN91}
D.J.~Barsky, G.R.~Grimmett, and C.M.~Newman.
Percolation in half-spaces: Equality of critical densities and continuity of the percolation probability.
{\em Probab.\ Theory Relat.\ Fields}~90(1), 111--148, 1991.

\bibitem[BK89]{BK89}
R.M.~Burton and M.~Keane.
Density and uniqueness in percolation.
{\em Commun.\ Math.\ Phys.}~121, 501--505, 1989.

\bibitem[BLPS99]{BLPS99}
I.~Benjamini, R.~Lyons, Y.~Peres, O.~Schramm.
Group-invariant percolation on graphs.
{\em Geom.\ Funct.\ Anal.}~9(1), 29--66, 1999.

\bibitem[BS01]{BS01}
I.~Benjamini and O.~Schramm.
Percolation in the hyperbolic plane.
{\em J.~Am.~Math.\ Soc.}~14, 487--507, 2001.

\bibitem[CDG04]{CDG04}
J.T.~Cox, D.A.~Dawson, and A.~Greven.
Mutually catalytic super branching random walks: Large finite systems
and renormalization analysis.
{\em Mem.\ Am.\ Math.\ Soc.}~809, 2004. 

\bibitem[CFG96]{CFG96}
J.T.~Cox, K.~Fleischmann, and A.~Greven.
Comparison of interacting diffusions and an application to their ergodic theory.
{\em Probab.\ Theory Relat.\ Fields}~105, 513--528, 1996.

\bibitem[CG86]{CG86}
J.T.~Cox and D.~Griffeath.
Diffusive clustering in the two dimensional voter model.
{\em Ann.\ Probab.}~14(2), 347--370, 1986.

\bibitem[CG94]{CG94}
J.T.~Cox and A.~Greven.
Ergodic theorems for infinite systems of locally interacting diffusions.
{\em Ann. Probab.}~22(2), 833--853, 1994.

\bibitem[Che87]{Che87}
M.F.~Chen.
Existence theorems for interacting particle systems with non-compact
state space.
{\em Sci.\ China Ser.\ A}\ 30, 148--156, 1987.

\bibitem[Daw77]{Daw77}
D.A.~Dawson.
The critical measure diffusion process.
{\em Z.\ Wahrscheinlichkeitstheor.\ Verw.\ Geb.}~40, 125--145, 1977.

\bibitem[Daw93]{Daw93}
D.A.~Dawson.
Measure-valued Markov processes.
{\em Ecole d'Eté de probabilités de Saint-Flour XXI. Lect. Notes Math.}~1541, 1-260, Springer, Berlin, 1993.

\bibitem[DDL90]{DDL90}
W.~Ding, R.~Durrett, and T.M.~Liggett.
Ergodicity of reversible reaction diffusion processes.
{\em Probab.\ Theory Relat.\ Fields} 85(1), 13--26, 1990.

\bibitem[DE68]{DE68}
D.A.~Darling and P.~Erd\H os.
On the recurrence of a certain chain.
{\em Proc.\ Am.\ Math.\ Soc.}~19(1):\ 336-338, 1968.

\bibitem[DEFMPX02a]{DEFMPX02a}
D.A.~Dawson, A.~Etheridge, K.~Fleischmann, L.~Mytnik, E.A.~Perkins, and J.~Xiong.
Mutually catalytic branching in the plane: Finite measure states.
{\em Ann. Probab.}~30(4), 1681--1762, 2002. 

\bibitem[DEFMPX02b]{DEFMPX02b}
D.A.~Dawson, A.~Etheridge, K.~Fleischmann, L.~Mytnik, E.A.~Perkins, and J.~Xiong.
Mutually catalytic branching in the plane: Infinite measure states.
{\em Electron. J. Probab.}~7, paper no.~15, 61 pp., 2002.

\bibitem[Deu89]{Deu89}
J.-D.~Deuschel.
Invariance principle and empirical mean large deviations of the critical Ornstein-Uhlenbeck process.
{\em Ann. Probab.}~17(1), 74--90, 1989.

\bibitem[DF97a]{DF97a}
D.A.~Dawson and K.~Fleischmann.
A continuous super-Brownian motion in a super-Brownian medium.
{\em J.\ Theoret.\ Probab.}~10(1), 213--276, 1997.

\bibitem[DF97b]{DF97b}
D.A.~Dawson and K.~Fleischmann.
Longtime behavior of a branching process controlled by branching catalysts.
{\em Stoch.\ Process.\ Appl.}~71(2), 241--257, 1997.

\bibitem[DF02]{DF02}
D.A.~Dawson and K.~Fleischmann.
Catalytic and mutually catalytic super-Brownian motions
In R.C.~Dalang (ed.) {\em Seminar on stochastic analysis, random fields and applications III. Proceedings of the 3rd seminar, Ascona, Switzerland, September 20-24, 1999.} Prog.\ Probab.~52, 89--110, Birkh\"auser, Basel, 2002.

\bibitem[DFMPX03]{DFMPX03}
D.A.~Dawson, K.~Fleischmann, L.~Mytnik, E.A.~Perkins, and J.~Xiong.
Mutually catalytic branching in the plane: Uniqueness.
{\em Ann.\ Inst.\ Henri Poincaré, Probab.\ Stat.}~39(1), 135--191, 2003.

\bibitem[DG93a]{DG93a}
D.A.~Dawson and A.~Greven.
Hierarchical models of interacting diffusions: Multiple time scale
  phenomena, phase transition and pattern of cluster-formation.
{\em Probab.\ Theory Related Fields}~96(4):\ 435--473, 1993.

\bibitem[DG93b]{DG93b}
D.A.~Dawson and A.~Greven.
Multiple time scale analysis of interacting diffusions.
{\em Probab.\ Theory Related Fields}~95(4):\ 467--508, 1993.

\bibitem[DG96]{DG96}
D.A.~Dawson and A.~Greven.
Multiple space-time scale analysis for interacting branching models.
{\em Electron.\ J.\ Probab.}~1, paper no.\ 14, 84 pp., 1996.

\bibitem[DG99]{DG99}
D.A.~Dawson and A.~Greven.
Hierarchically interacting Fleming-Viot processes with selection
and mutation: Multiple space time scale analysis and quasi-equilibria.
{\em Electron.\ J.\ Probab.}\ 4, Paper no.~4, 81 p., 1999.

\bibitem[DGV95]{DGV95}
D.A.~Dawson, A.~Greven and J.~Vaillancourt.
Equilibria and quasi-equilibria for infinite collections of
interacting Fleming-Viot processes.
{\em Trans.\ Amer.\ Math.\ Soc.}~347(7): 2277--2360, 1995.

\bibitem[DK96]{DK96}
P.~Donnelly and T.G.~Kurtz.
A countable representation of the Fleming-Viot measure-valued diffusion.
{\em Ann. Probab.}~24(2), 698--742, 1996.

\bibitem[DK99]{DK99}
P.~Donelly and T.G.~Kurtz.
Genealogical processes for Fleming-Viot models with selection and
recombination.
{\em Ann.\ Appl.\ Probab.}~9, 1091--1148, 1999.

\bibitem[DLSS91]{DLSS91}
B.~Derrida, J.L.~Lebowitz, E.R.~Speer, and H.~Spohn.
Dynamics of an anchored Toom interface.
{\em J.\ Phys.\ A: Math.\ Gen.}~24, 4805--4834, 1991.

\bibitem[DP98]{DP98}
D.A.~Dawson and E.A.~Perkins.
Long-time behavior and coexistence in a mutually catalytic branching model.
{\em Ann. Probab.}~26(3), 1088-1138, 1998.

\bibitem[DS95]{DS95}
R.~Durrett and R.~Schinazi.
Intermediate phase for the contact process on a tree.
{\em Ann.\ Probab.}~23(2), 668--673, 1995.

\bibitem[Dyn02]{Dyn02}
E.~Dynkin.
{\em Diffusions, superdiffusions and partial differential equations.}
Vol. 50 of American Mathematical Society Colloquium Publications.
AMS, Providence, RI, 2002.

\bibitem[DZ98]{DZ98}
A.~Dembo and O.~Zeitouni.
{\em Large Deviations Techniques and Applications. 2.\ ed.}
Springer, New York, 1998.

\bibitem[EF98]{EF98}
A.~Etheridge and K.~Fleischmann.
Persistence of a two-dimensional super-Brownian motion in a catalytic medium.
{\em Probab.\ Theory Relat.\ Fields}~110(1), 1--12, 1998.

\bibitem[EK86]{EK}
S.N.~Ethier and T.G.~Kurtz.
{\em Markov Processes; Characterization and Convergence.}
John Wiley \& Sons, New York, 1986.

\bibitem[EK04]{EK02}
J.~Engl\"ander and A.E.~Kyprianou.
Local extinction versus local exponential growth for spatial branching processes.
{\em Ann. Probab.}~32(1A), 78--99, 2004.

\bibitem[EP99]{EP99}
J.~Engl\"ander and R.G.~Pinsky.
On the construction and support properties of measure-valued diffusions on $D\subseteq\R^d$ with spatially dependent branching.
{\em Ann.\ Probab.}, 27(1): 684--730, 1999.

\bibitem[ER91]{ER91}
N. El Karoui and S.~Roelly.
Propri\'et\'es de martingales, explosion et repr\'esentation
de L\'evy- Khintchine d'une classe de processus de branchement \`a
valeurs mesures. 
{\em Stoch.\ Proc.\ Appl.}~38(2):\ 239--266, 1991.

\bibitem[ET02]{ET02}
J.~Engl\"ander and  D.~Turaev.
A scaling limit for a class of superdiffussions.
{\em Ann. Probab.}~30(2), 683--722, 2002.

\bibitem[Eth00]{Eth00}
A.~Etheridge.
{\em An Introduction to Superprocesses.}
University Lecture Series~20,
AMS, Providence, 2000.

\bibitem[EU48]{EU48}
C.J.~Everett and S.~Ulam.
Multiplicative systems in several variables, I.
Los Alamos Scientific Laboratory, LA-683, 1948.

\bibitem[Ewe04]{Ewe04}
W.J.~Ewens.
{\em Mathematical Population Genetics. I: Theoretical Introduction. 2nd ed.}
Interdisciplinary Mathematics~27.\ Springer, New York, 2004.

\bibitem[FG94]{FG94}
K.~Fleischmann and A.~Greven.
Diffusive clustering in an infinite system of hierarchically interacting diffusions.
{\em Probab.\ Theory Relat.\ Fields}~98(4), 517--566, 1994.

\bibitem[Fit88]{Fit88}
P.J.~Fitzsimmons.
Construction and regularity of measure-valued branching processes.
{\em Isr.~J.~Math.}~64(3):\ 337--361, 1988.

\bibitem[Fit91]{Fit91}
P.J.~Fitzsimmons.
Correction to ``Construction and regularity of measure-valued branching processes''.
{\em Isr.~J.~Math.}, 73(1): 127, 1991.

\bibitem[Fit92]{Fit92}
P.J.~Fitzsimmons.
On the martingale problem for measure-valued Markov branching processes.
pp~39--51 in: {\em Seminar on Stochastic Processes, 1991.}
Progress in Probability~29, Birkh\"auser, Boston, 1992.

\bibitem[FK99]{FK99}
K.~Fleischmann and A.~Klenke.
Smooth density field of catalytic super-Brownian motion.
{\em Ann.\ Appl.\ Probab.}~9, 298--318, 1999.

\bibitem[FS03]{FSsup}
K.~Fleischmann and J.M. Swart.
Extinction versus exponential growth in a supercritical
  super-Wright-Fischer diffusion.
{\em Stoch. Proc. Appl.}~106(1):\ 141--165, 2003.

\bibitem[FS04]{FStrim}
K.~Fleischmann and J.M. Swart.
Trimmed trees and embedded particle systems.
{\em Ann. Probab.}~32(3a):\ 2179--2221, 2004.

\bibitem[GH02]{GH02}
G.~Grimmett and P.~Hiemer.
Directed percolation and random walk.
Pages~273--297 in: V.~Sidoravicius (ed.), {\em In and Out of Equilibrium.}
Prog.\ Probab.~51.
Birkh\"auser, Boston, 2002

\bibitem[GKW99]{GKW99}
A.~Greven, A.~Klenke, and A.~Wakolbinger.
The longtime behavior of branching random walk in a catalytic medium.
{\em Electron.\ J.\ Probab.}~4, paper no.\ 12, 80 pp., 1999.

\bibitem[GKW01]{GKW01}
A.~Greven, A.~Klenke, and A.~Wakolbinger.
Interacting Fisher-Wright diffusions in a catalytic medium.
{\em Probab.\ Theory Related Fields}~120(1):\ 85--117, 2001.

\bibitem[GLW05]{GLW05}
A.~Greven, V.~Limic, and A.~Winter.
Representation theorems for interacting Moran models and interacting Fisher-Wright models, and applications.
{\em Electron. J. Probab.}~10, paper no.~39, 1286--1358, 2005.

\bibitem[GM90]{GM90}
G.R.~Grimmett and J.M.~Marstrand.
The supercritical phase of percolation is well behaved.
{\em Proc.\ R.\ Soc.\ Lond., Ser.~A}~430, no.~1879, 439--457, 1990.

\bibitem[GW91]{GW91}
L.G.~Gorostiza and A.~Wakolbinger.
Persistence criteria for a class of critical branching particle systems 
in continuous time.
{\em Ann.\ Probab.}~19(1), 266--288, 1991.

\bibitem[Hag97]{Hag97}
O.~H\"aggstr\"om.
Infinite clusters in dependent automorphism invariant percolation on trees.
{\em Ann.\ Probab.}~25(3), 1423--1436, 1997.

\bibitem[Har63]{Har63}
T.E.~Harris.
{\em The theory of branching processes.}
Springer, Berlin, 1963.

\bibitem[Har76]{Har76}
T.E.~Harris.
On a class of set-valued Markov processes.
{\em Ann.\ Probab.}~4, 175--194, 1976.

\bibitem[Har78]{Har78}
T.E.~Harris.
Additive set-valued Markov processes and graphical methods.
{\em Ann.\ Probab.}~6, 355--378, 1978.

\bibitem[HS98]{HS98}
F.~den Hollander and J.M.~Swart.
Renormalization of hierarchically interacting isotropic diffusions.
{\em J.\ Stat.\ Phys.}~93:\ 243--291, 1998.

\bibitem[Jir64]{Jir64}
M.~Ji\v rina.
Branching processes with measure-valued states.
In {\em Trans. Third Prague Conf. Information Theory, Statist.
          Decision Functions, Random Processes (Liblice, 1962)},
          pages 333--357, Czech.\ Acad.\ Sci., Prague, 1964.

\bibitem[Kal76]{Kal76}
O.~Kallenberg.
{\em Random Measures.}
Akademie-Verlag, Berlin, 1976.

\bibitem[Kal77]{Kal77}
O.~Kallenberg.
Stability of critical cluster fields.
{\em Math.\ Nachr.}~77, 7--43, 1977.

\bibitem[Kal83]{Kal83}
O.~Kallenberg.
{\em Random measures}, 3rd rev.\ and enl.\ ed.
Akademie-Verlag, Berlin, 1983.

\bibitem[Kes86]{Kes86} 
H.~Kesten.
Aspects of first passage percolation. Pages 125--264 in: 
P.L.~Hennequin (ed.),
{\em \'Ecole d'\'et\'e de probabilit\'es de Saint-Flour XIV - 1984.}
Lect.\ Notes Math.~1180, Springer, Berlin, 1986.

\bibitem[Kle96]{Kle96}
A.~Klenke.
Different clustering regimes in systems of hierarchically interacting
  diffusions.
{\em Ann.\ Probab.}~24(2):\ 660--697, 1996.

\bibitem[KN97]{KN97}
S.M.~Krone and C.~Neuhauser.
Ancestral processes with selection.
{\em Theor.\ Popul.\ Biol.}~51(3), 210--237, 1997.

\bibitem[Kol33]{Kol33}
A.~Kolmogorov.
{\em Grundbegriffe der Wahrscheinlichkeitsrechnung.}
Ergebnisse der Mathematik, 1933.

\bibitem[Law05]{Law05}
G.F.~Lawler.
{\em Conformally Invariant Processes in the Plane.}
AMS, 2005.

\bibitem[Lie81]{Lie81}
A.~Liemant.
Kritische Verzweigungsprozesse mit allgemeinem Phasenraum.~IV.
{\em Math.\ Nachr.}~102:\ 235--254, 1981.

\bibitem[Lig85]{Lig85}
T.M.~Liggett.
{\em Interacting Particle Systems}.
Springer, New York, 1985.

\bibitem[Lig96]{Lig96}
T.M.~Liggett.
Branching random walks and contact processes on homogeneous trees.
{\em Probab.\ Theory Relat.\ Fields}~106(4), 495--519, 1996.

\bibitem[Lig99]{Lig99}
T.M.~Liggett.
{\em Stochastic Interacting Systems: Contact, Voter and Exclusion Process.}
Springer, Berlin, 1999.

\bibitem[Loe63]{Loe63}
M.~Lo\`eve.
{\em Probability Theory} 3rd ed.
Van Nostrand, Princeton, 1963.

\bibitem[Loe78]{Loe78}
M.~Lo\`eve.
{\em Probability Theory II} 4th ed.
Graduate Texts in Mathematics~46.
Springer, New York, 1978.

\bibitem[LP05]{LP05}
R.~Lyons and Y.~Peres.
{\em Probability on trees and networks.}
Draft available from 
http://mypage.iu.edu/$\ti\ $rdlyons/prbtree/prbtree.html, 2005.

\bibitem[LPP96]{LPP96}
R.~Lyons, R.~Pemantle, and Y.~Peres.
Random walks on the lamplighter group.
{\em Ann.\ Probab.}~24(4), 1993--2006, 1996.

\bibitem[LS81]{LS81}
T.M.~Liggett and F.~Spitzer.
Ergodic theorems for coupled random walks and other systems with
locally interacting components.
{\em Z. Wahrsch. verw. Gebiete} 56, 443--468, 1981.

\bibitem[Mou92]{Mou92}
T.S.~Mountford.
The ergodicity of a class of reversible reaction-diffusion processes.
{\em Probab.\ Theory Relat.\ Fields} 92(2), 259--274, 1992.

\bibitem[MT95]{MT95}
C.~M\"uller and R.~Tribe.
Stochastic p.d.e.'s arising from the long range contact and
long range voter processes.
{\em Probab.\ Theory Relat.\ Fields} 102(4), 519--545, 1995.

\bibitem[MW89]{MW89}
B.~Mohar and W.~Woess.
A survey on spectra of infinite graphs.
{\em Bull.\ Lond.\ Math.\ Soc.}~21(3), 209--234, 1989. 

\bibitem[Neu90]{Neu90}
C.~Neuhauser.
An ergodic theorem for Schl\"ogl models with small migration.
{\em Probab.\ Theory Relat.\ Fields} 85(1), 27--32, 1990.

\bibitem[NS80]{NS80}
M.~Notohara and T.~Shiga.
Convergence to genetically uniform state in stepping stone models of population genetics.
{\em J.\ Math.\ Biology}~10, 281--294, 1980.

\bibitem[Pat88]{Pat88}
A.L.T.~Paterson.
{\em Amenability.}
AMS, Providence, 1988.

\bibitem[Paz83]{Paz83}
A.~Pazy.
{\em  Semigroups of Linear Operators and Applications
to Partial Differential Equations.}
Springer, New York, 1983.

\bibitem[Pem92]{Pem92}
R.~Pemantle.
The contact process on trees.
{\em Ann.\ Probab.}~20(4), 2089--2116, 1992.

\bibitem[Pen04]{Pen04}
C.~Pen\ss el.
{\em Interacting Catalytic Feller Diffusions: Finite System Scheme and Renormalisation.}
Logos, Berlin, 2004.

\bibitem[Rev84]{Rev84}
D.~Revuz.
{\em Markov chains. 2nd rev. ed.}
North-Holland Mathematical Library, Vol~ 11.
North-Holland, Amsterdam, 1984.

\bibitem[RW87]{RW87}
L.C.G.~Rogers and D.~Williams.
{\em Diffusions, Markov Processes, and Martingales, Volume 2: Ito
  Calculus.}
Wiley, Chichester, 1987.

\bibitem[Sch72]{Sch72}
F.~Schl\"ogl.
Chemical reaction models and non-equilibrium phase transitions.
{\em Z.~Phys.}\ 253, 147--161, 1972.

\bibitem[Sch86]{Sch86} R.H.~Schonmann.
The asymmetric contact process.
{\em J.\ Stat.\ Phys.}~44, 505--534, 1986.

\bibitem[Sch98]{Sch98}
F.~Schiller.
{\em Application of the Multiple Space-Time Scale Analysis on a
  System of R-valued, Hierarchically Interacting, Stochastic Differential
  Equations.}
Master thesis, Universtity Erlangen-N\"urnberg, 1998.

\bibitem[SF83]{SF83}
S.~Sawyer and J.~Felsenstein.
Isolation by distance in a hierarchically clustered population.
{\em J.\ Appl.\ Probab.}~20:\ 1--10, 1983.

\bibitem[Sha88]{Sha88}
M.~Sharpe.
{\em General Theory of Markov Processes}.
Academic Press, Boston, 1988.

\bibitem[Shi80a]{Shi80}
T.~Shiga.
An interacting system in population genetics.
{\em J.\ Math.\ Kyoto Univ.} 20(2), 213--242, 1980.

\bibitem[Shi80b]{Shi80b}
T.~Shiga.
An interacting system in population genetics, II.
{\em J.\ Math.\ Kyoto Univ.}~20(4), 723--733, 1980.

\bibitem[Shi81]{Shi81}
T.~Shiga.
Diffusion processes in population genetics.
{\em J.~Math.\ Kyoto Univ.}\ 21, 133--151, 1981.

\bibitem[Shi92]{Shi92}
T.~Shiga.
Ergodic theorems and exponential decay of sample paths for certain interacting diffusion systems.
{\em Osaka J.\ Math.}~29, 789--807, 1992.

\bibitem[Smo83]{Smo83}
J.~Smoller.
{\em Shock Waves and Reaction-diffusion Equations.}
Vol.~258 of Grund\-lehren Math.\ Wiss., Springer, New York, 1983.

\bibitem[SS80]{SS80}
T.~Shiga and A.~Shimizu.
Infinite dimensional stochastic differential equations and their applications.
{\em J.~Math.\ Kyoto Univ.}\ 20, 395--416, 1980.

\bibitem[SS97]{SS97}
M.~Salzano and R.H.~Schonmann.
The second lowest extremal invariant measure of the contact process.
{\em Ann.\ Probab.}~25(4), 1846--1871, 1997.

\bibitem[SS99]{SS99}
M.~Salzano and R.H.~Schonmann.
The second lowest extremal invariant measure of the contact process. II.
{\em Ann.\ Probab.}~27(2), 845--875, 1999.

\bibitem[Sta96]{Sta96}
A.M.~Stacey.
The existence of an intermediate phase for the contact process on trees.
{\em Ann.\ Probab.}~24(4), 1711--1726, 1996.

\bibitem[SU86]{SU86}
T.~Shiga and K.~Uchiyama.
Stationary states and their stability of the stepping stone model
involving mutation and selection.
{\em Probab.\ Theory Relat.\ Fields} 73, 87--117, 1986.

\bibitem[Swa99]{Swa99}
J.M.~Swart.
{\em Large Space-Time Scale Behavior of Linearly Interacting
  Diffusions}.
PhD thesis, Katholieke Universiteit Nijmegen, 1999.
http://helikon.ubn.kun.nl/ mono/s/swart\un{\ }/largspscb.pdf.

\bibitem[Swa00]{Swa00}
J.M.~Swart.
Clustering of linearly interacting diffusions and universality of
  their long-time limit distribution.
{\em Prob.\ Theory Related Fields}~118:\ 574--594, 2000.

\bibitem[WG74]{WG74}
H.W.~Watson and F.~Galton.
On the probability of the extinction of families.
{\em J.\ Anthropol.\ Inst.\ Great Britain and Ireland}~4, 138--144, 1874.

\bibitem[WK74]{WK74}
K.G.~Wilson and J.~Kogut.
The renormalization group and the $\eps$-expansion.
{\em Phys.\ Rep.}~12C, 75--200, 1974.

\bibitem[YW71]{YW71}
T.~Yamada and S.~ Watanabe.
On the uniqueness of solutions of stochastic differential equations.
{\em J.\ Math.\ Kyoto Univ.}~11:\ 155--167, 1971.

\end{thebibliography}
\end{document}